\documentclass{report}

\usepackage{amsfonts,amsmath,graphicx}

\begin{document}
\pagenumbering{roman}
\begin{center}
{\Large \bf Mathematical Models in Biology}
\vskip 0.2in
\centerline{By}
\vskip 0.1in
\centerline{\large  BARBARA CATHRINE MAZZAG}
\centerline{\large B.A. (University of California, Santa Cruz), 1995} 
\centerline{\large M.S. (University of California, Davis), 2000 }
\vskip 0.3in
\centerline{\large  DISSERTATION}
\vskip 0.1in
\centerline{\large Submitted in partial satisfaction of the requirements
for the degree of}
\vskip 0.1in
\centerline{\large  DOCTOR OF PHILOSOPHY} 
\vskip 0.1in
\centerline{\large in}  
\vskip 0.1in
\centerline{\large APPLIED MATHEMATICS}
\vskip 0.2in
\centerline{\large in the}
\vskip 0.1in  
\centerline{\large  OFFICE OF GRADUATE STUDIES}
\vskip 0.1in
\centerline{\large of the}
\vskip 0.1in
\centerline{\large  UNIVERSITY OF CALIFORNIA}
\vskip 0.1in
\centerline{\large  DAVIS} 
\end{center}
\vskip 0.1in
{\large Approved:}
\begin{center}
\centerline{\underbar{\hskip 2.5in}}   
\vskip 0.15in
\centerline{\underbar{\hskip 2.5in}}
\vskip 0.15in
\centerline{\underbar{\hskip 2.5in}}   
\vskip 0.2in
\centerline{\large Committee in Charge}
\vskip 0.2in  
\centerline{\large 2002}
\end{center}
\newpage
\large
\tableofcontents
\listoffigures
\listoftables


\newpage
\large
{\Large \bf ACKNOWLEDGMENTS} \\

My most grateful thanks go to my advisor, Prof. Alex Mogilner for all of
his help and understanding.  I feel extremely fortunate for the
opportunity to work with him and learn from him.  He has introduced me to
mathematical biology, has encouraged me to participate in conferences and
an internship, and in general has shown me how exciting and lively work
in this area can be.  I would like to thank him for providing me with a
lot of practical help in research, for sharing his insights and intuition
with me, and for the financial and emotional support he has given me.  
This dissertation is based on collaborations with him but, naturally, the
errors in this work are mine alone. This document would not exist without
his help. 

Prof. Goodhill has generously provided a place for me while I was
spending a six-month-long internship in his laboratory.  I cannot thank
him enough for teaching me to use many valuable research tools, for
stimulating discussions, for a lot of personal help while I was at
Georgetown, and for his continued interest and detailed comments on our
collaborative project.  Some day I hope to live up to his excellence in
writing.

Prof. Barakat also has given me many very useful comments and criticisms 
on our joint work as well as on my Masters' thesis.  I also appreciate 
the excellent background he has given me on fluid dynamics in his course.    

I would like to thank Prof. Angela Cheer and Prof. Albert Fannjiang for 
reading this work and for providing me with their comments and 
suggestions.  

I am also very grateful to my personal system's administrator and support 
team, Tyler Evans.  Writing this dissertation would have been much easier 
without having to follow him to Eureka, but it would have been impossible 
to complete it without him.


\newpage
\begin{center}
\underline{\bf \Large Abstract}
\end{center}

\medskip

\large
\noindent

Aerotaxis is the particular form of chemotaxis in which oxygen plays the 
role of both the attractant and the repellent.  Aerotaxis occurs without 
methylation adaptation, and it leads to fast and complete aggregation 
toward the most favorable oxygen concentration.  Biochemical pathways of 
aerotaxis remain largely elusive, however, aerotactic pattern formation is 
well documented.  This allows mathematical modeling to test plausible 
hypotheses about the biochemical mechanisms.  Our model demonstrates that 
assuming fast, non-methylation adaptation produces theoretical results 
that are consistent with experimental observations.  We obtain analytical 
estimates for parameter values that are difficult to obtain 
experimentally.  

Chemotaxis in growth cones differs from gradient sensing in other animal 
cells, because growth cones can change their attractive or repulsive 
response to the same chemical gradient based on their internal calcium or 
cAMP levels.  We create two models describing different aspects of growth 
cone guidance.  One model describes the internal switch that determines 
the direction of movement.  However, this model allows chemotaxis under 
certain conditions only, so a second model is created to propose a 
mechanism that allows growth cone guidance in any environment.  

Endothelial cells go through extensive morphological changes when exposed 
to shear stress due to blood flow.  These morphological changes are 
thought to be at least partially the result of mechanical signals, such as 
deformations, transmitted to the cell structures.  Our model describes an 
endothelial cell as a network of viscoelastic Kelvin bodies with 
experimentally obtained parameters.  Qualitative predictions of the model 
agree with experiments.


\newpage
\pagestyle{myheadings}
\pagenumbering{arabic}
\markright{  \rm \normalsize CHAPTER 1. \hspace{0.5cm}
 MATHEMATICAL MODELS IN BIOLOGY }
\large
\chapter{Introduction}
\thispagestyle{myheadings}

Biology has gone through an extraordinary change in the past century,
partially due to increasingly advanced methods of being able to collect
data, and partially because of the sophistication in the quantitative
analysis of this data.  These changes are particularly striking in
molecular and cellular biology, where incredibly complex interactions are
revealed to be at the basis of all cell functions such as sensing,
movement or reproduction.  It is precisely the complexity of experimental
observations that necessitates a more accurate and in-depth analysis.  
The need for a quantitative understanding of biological phenomena has lead
to different modeling approaches.  Many times highly advanced numerical
simulations are created, and some research is aimed at highly realistic
computer models of entire signal transduction pathways, or even entire
cells.  A very different, but equally valid approach, is to simplify
possibly very complicated interactions to a smaller set of key components
which lends itself easier to analytical models.

This dissertation is concerned with models of the latter type, namely,
with arriving at biologically meaningful results from mathematical models
based on a simplification of experimental observations.  The hope of such
models is that if they do indeed capture the key principles of the
underlying the phenomenon, then a quantitative understanding of these
principles leads to new information which has not been uncovered by the
experiments.  In successful models the mathematical analysis leads to
insights which are unattainable (or very difficult to attain)  
experimentally.

The first two chapters of this work are centered around a common theme:  
gradient sensing.  Chapter \ref{aerotaxis} discusses a model of pattern
formation due to bacteria searching for optimal oxygen concentrations.  
This work is based on experiments conducted by Zhulin et al. \cite{Z} on
pattern formation of such bacteria.  Our mathematical model presented in
this chapter confirms that the experiments are produced by a novel form of
gradient sensing, and in addition, it offers some experimentally testable
predictions.

Chapter \ref{axon} presents two models of how signal transduction events
lead to the orientation of a neuron toward the appropriate target.  This
is an inherently difficult problem because of the limitations on
experimental data available.  In this dissertation two mathematical models
of neuronal gradient sensing are developed, each aiming at understanding a
different aspect of this question.  Their advantages and disadvantages are
discussed in detail in this chapter, and it is concluded that further work
is necessary in this area.  This chapter of my dissertation is a
collaboration with Prof. Geoffrey J.  Goodhill from the Neuroscience
Department of Georgetown University Medical Center who introduced me to
the biological background of growth cone guidance.

The third main topic of this work, presented in Chapter \ref{endothelial},
also describes sensing on the cellular level, but in this case the signal
is not biochemical, but mechanical.  The mathematical model describes
morphological changes in the cells lining the blood vessels when they are
exposed to different types of shear stresses (induced by different types
of flow over the cell surface).  The predictions of the model agree
qualitatively with the experiments, however, improvements of the model
described in this chapter are necessary in order to make quantitative
predictions.  This work was done in collaboration with John S. Tamaresis
from the Graduate Group in Applied Mathematics and Prof. Abdul Barakat
from the Mechanical and Aeronautical Engineering Department of the
University of California, Davis.

Each chapter contains the biological terminology and data relevant to the
topic, as well as a discussion and conclusion of the mathematical model
presented. The three projects are quite distinct both biologically and
mathematically, therefore separate conclusions appeared to be most
appropriate.


\newpage
\pagestyle{myheadings}
\chapter{Aerotaxis} \label{aerotaxis}
\thispagestyle{myheadings}
\markright{  \rm \normalsize CHAPTER 2. \hspace{0.5cm}
  MATHEMATICAL MODELS IN BIOLOGY}

\section{Introduction}

The study of cell motility is a broad subject with numerous applications
ranging from understanding how nerve cells find their place in the
developing brain to understanding wound healing.  Of all different forms
of motile cell behavior, bacterial chemotaxis is the best understood.
Bacterial aerotaxis differs from conventional chemotaxis in a number of
interesting ways that were highlighted by experiments conducted at the
Loma Linda Medical School by Zhulin et al. \cite{Z}.

Some differences between conventional chemotaxis and aerotaxis are already
known, but the biochemical signal transduction pathways involved in
aerotaxis are not.  The purpose of this work is to demonstrate that the
unusual patterns found in Zhulin's aerotaxis experiments are consistent
with a novel form of taxis without slow adaptation.  The main result is a
mathematical model based on fast adaptation that characterizes aerotactic
behavior and uses experimentally obtained parameters. Analytical and
numerical results of the model are compared to experimental data.

The biological background is introduced in Section \ref{conv_chem} which
explains the terms chemotaxis and aerotaxis in detail and describes the
differences between the two. Current knowledge of the biochemical pathways
involved in chemotaxis is also explained.

A detailed description of the Zhulin aerotaxis experiments and explanation
of the observed pattern follows, as well as questions that arise from
these experiments.  We conclude with descriptions of various chemotaxis
models and their limitations when applied to aerotaxis experiments.  In
particular, we introduce the Keller-Segel model, and discuss the modeling
assumptions.  This is followed by a summary of Gr\"unbaum's work on
approximating a general class of equations describing random walk
behaviors.  We used his analysis to show why most conventional chemotaxis
models cannot be applied to the Zhulin aerotaxis experiments.  Some other
mathematical models are discussed briefly (Tranquillo \cite{TL}, Barkai \&
Leibler \cite{BL} ), and we argue that no existing models in the
literature provide an appropriate framework for aerotaxis.

Section \ref{model1} contains our mathematical model of aerotaxis.  The
terms of the simple advection-reaction equations are explained, and the
main question is the determination of the turning rates (reaction terms).  
A phenomenological justification is given for the choice of the particular
terms.  Appendix \ref{app_aerotax} provides the biological reasoning
behind choosing such form for the turning frequencies and provides a
simple mathematical model for a receptor which could produce such turning
rates.  Section \ref{model1} also explains how non-dimensonalization and
scaling were obtained for the model.

Section \ref{num_sim1} contains the results, both analytic and
numerical.  We show the numerical simulations of the aerotactic band
formation.  The interpretation of numerical results emphasizes how
the model matches and predicts the band formation.  Analytical
solutions are given for the steady state of the system.  Various
parameter values which are experimentally not easily measurable are
estimated.

In Section \ref{conclusion1} we summarize our findings and talks
about potential future projects related to this topic.

\clearpage
\section{Background} \label{conv_chem}

\subsection{Conventional chemotaxis} 

\begin{figure}[h]
\centerline{\includegraphics[width=0.8\textwidth]{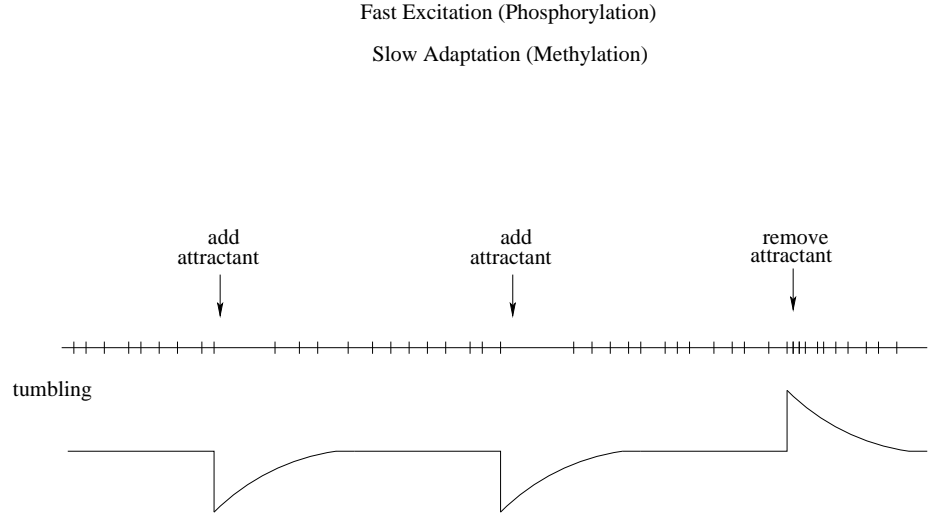}}
\caption[Turning frequencies]{Changes in turning frequencies as a
response to attractants
and repellents.  Fast excitation and slow adaptation are demonstrated.
(Figure based on Bray, \cite{B})}
\label{turn_freq.eps}
\end{figure}

Bacterial chemotaxis is a term used for motility in the direction of
higher nutrient concentrations (such as sugars and amino acids) and
away from repellents.  In a neutral environment (i.e. one with
uniform chemical concentrations) bacteria swim smoothly in a given
direction for a period of time, then go through a period of abrupt
changes of direction, called tumbling.  The sequence of 'runs' and
'tumbles' results in a random walk. This random walk becomes biased
when attractant (or repellent) is added (or removed).  Additional
attractants suppress the frequency with which tumbling occurs;
therefore, the straight runs lengthen in the direction of the
highest attractant concentration.  This allows bacteria to move up
the gradient.  Removing repellents has the same effect.  On the
other hand, adding repellents or removing attractants both increase
the frequency of tumbling and facilitate the movement of bacteria
down the gradient \cite{B}.  Figure \ref{turn_freq.eps} shows how
turning frequencies change in response to attractants or repellents.
Keeping attractant (and repellent) concentrations constant results
in adaptation of turning rates.

The swimming mechanism of cells can vary greatly, resulting in different
types of swimming not discussed in this thesis.  However, the underlying
principle of a biased random walk is prevalent in all forms of bacterial
motility.

Chemotactic movement is necessary for bacterial cells because
concentration gradients are detected by temporal comparison.  This
means that rather than being able to measure concentrations at the
different ends of the cell, bacteria must move through their
environment to be able to detect concentration changes.  In order to
have a temporal sensing mechanism, it is crucial that cells retain
some information about the previous environment; in other words, it
is necessary for the cell to have some sort of memory.

The mechanism for temporal sensing is also dependent on adaptation
to the current level of attractants and nutrients. This allows cells
to remain sensitive to concentration changes in a wide range of
chemical environments.  Adaptation and memory are related concepts,
since memory is a consequence of a slow adaptation mechanism which
allows cells to retain information about the previous environment
for a period of time.  During adaptation, the cell remains in a
chemical state determined by the attractant (repellent)
concentration before, and this time lag between the current and the
past states serves as the memory of the cell.  Any model of
conventional chemotaxis must address the issues of sensitivity and
memory.

The signal transduction pathways in bacteria such as {\it
Escherichia coli} have been widely studied \cite{BL}. There are two
important chemical reactions, phosphorylation and methylation, which
are responsible for controlling the tumbling frequency, and which
act on very different time scales.  Phosphorylation is a very fast
reaction (on the order of milliseconds \cite{B}) that causes a fast
response (on the order of 0.1 sec) to changes in the attractant or
repellent concentration, while methylation acts on the same time
scale as a single run (1-3 sec), and it is responsible for the
adaptation to current chemical concentrations.

\begin{figure}[h]
\centerline{\includegraphics[width=0.8\textwidth]{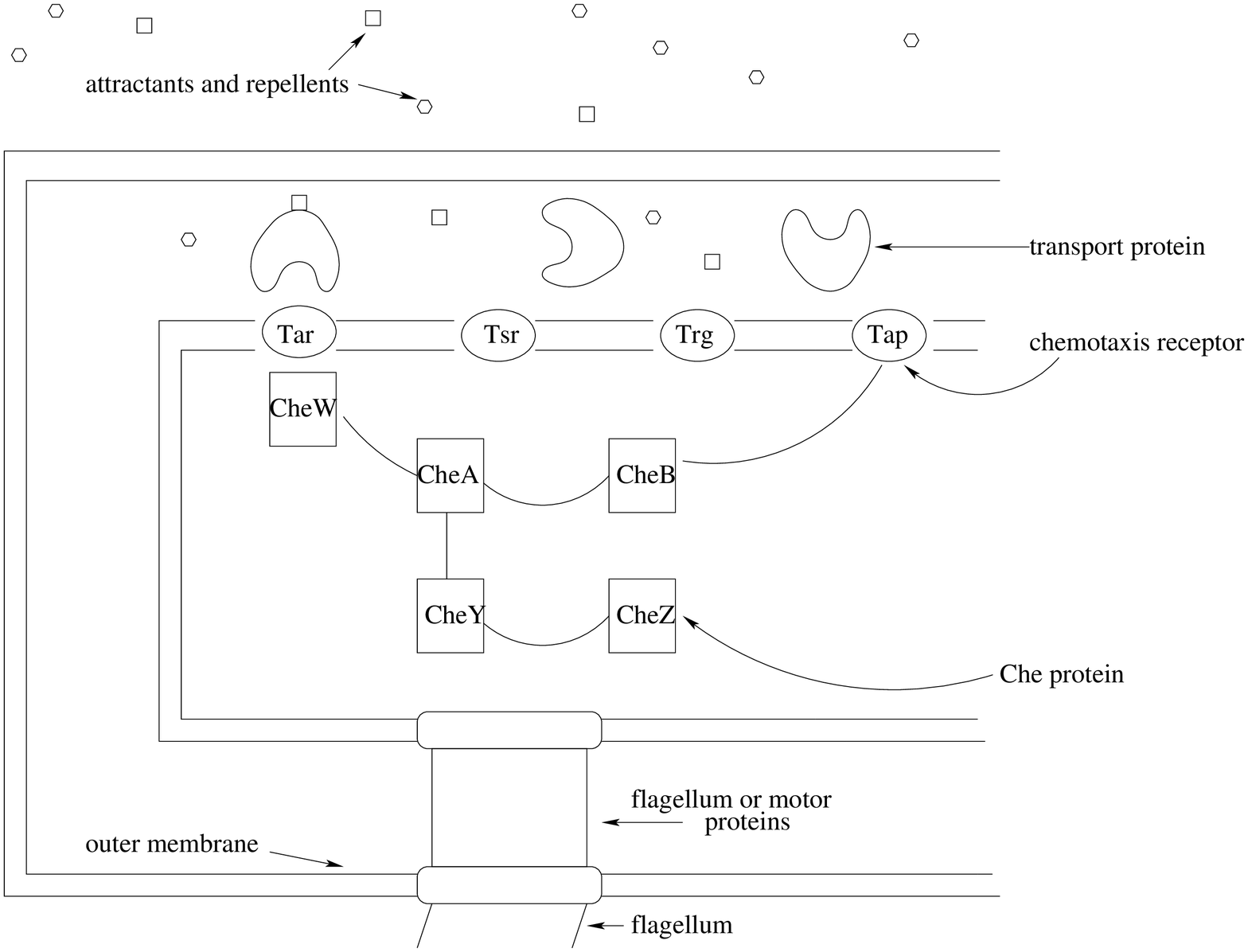}}
\caption[Signal transduction pathway.]{Signal transduction pathway in
bacteria. (Figure based on Bray, \cite{B}.)}
\label{sig_trans.eps}
\end{figure}

A simple model of the signal transduction pathway shown in Figure
\ref{sig_trans.eps} involves a number of proteins, CheA, CheW, CheY,
etc.  The signaling works as follows: ligands (attractants and
repellents) bind to transmembrane chemoreceptors of the bacteria.  
The cytosolic sides of the receptors are attached to kinases CheA
and CheW.  CheA phosphorylates both itself and CheY, which is a
protein that acts as a messenger inside the cell.  The role of CheY
is to increase the tumbling frequency by reacting with the motor.  
CheY eventually dephosphorylates with the help of CheZ.  Attractant
binding slows down CheY phosphorylation, which in turn leads to
suppressed tumbling and longer runs in the direction of attractants.  
Methylation, on the other hand, speeds up CheY phosphorylation which
leads to an eventual return to the base tumbling rates, i.e.
adaptation. Methylation occurs through a pair of enzymes, CheR and
CheB, which add and remove methyl groups.  As the rate of
methylation of the receptors slowly increases, the tumbling
frequency increases as well, and the cell returns to the original
turning rates. There is also a coupling between the phosphorylation
and methylation pathways, and it is the phosphorylation of CheB by
CheA.  This results in an increase in the demethylation activity of
CheB, so this also contributes to the suppression of tumbling rates
\cite{BL}.

\subsection{Aerotaxis}

Aerotaxis is a specific type of taxis (directed movement) in which
both the attractant and the repellent are particular concentrations
of oxygen.  Very low and very high concentrations both act as
repellents whereas some intermediate concentrations attract
bacteria.  The particular range of desirable concentration depends
on the species.  Aerotaxis was first discovered in 1676 by van
Leeuwenhoek who noticed aggregation of cells underneath the surface
of a solution.  Later, in 1881, Englemann observed that bacteria
aggregate at the edges of coverglass and around the air bubbles
trapped underneath \cite{T}. Initially, aerotaxis was considered a
chemotactic response to oxygen, but further studies, summarized in
Taylor's 1983 paper \cite{T}, reveal several aspects in which
aerotaxis and conventional chemotaxis differ drastically.

First of all, we list these differences between conventional
chemotaxis and aerotaxis, then explain their meaning and
significance below.  One of the most notable differences is that
aerotaxis is metabolism-dependent.  This means that while bacteria
do consume the oxygen, chemotactic bacteria can be attracted to
nutrients which they are unable to metabolize \cite{A}. There is
also quite a bit of evidence \cite{R, TZ, Z}, that aerotaxis is an
example of so called "energy taxis" in which cells monitor their
internal energy balance and react to optimize it.  (In contrast,
chemotaxis is based on monitoring and optimizing the nutrient
availability in the external environment.) We discuss the notion of
energy taxis further below.  The third factor distinguishing
aerotaxis and chemotaxis is the signal transduction pathway.  This
includes differences in the receptors utilized, as well as the fact
that aerotactic movement does not have a methylation-dependent
adaptation \cite{T}.  The precise mechanism of adaptation in
aerotactic bacteria is currently unknown.

In order to explain energy taxis, we must introduce some new terms.  
Proton motive force refers to the electrochemical potential difference
across the membrane \cite{T}, and it is produced by linking electron
transport due to respiration to translocation of protons \cite{T}.  The
coupling between the proton motive force and the electron transport system
is tight, and currently it is not known which of the two acts as a signal
for the cell's behavior \cite{R}. However, it is believed \cite{R, TZ}
that aerotactic bacteria respond to internal changes in the proton motive
force or the elector transport system and not directly to the
extracellular oxygen levels.  This is supported by a series of experiments
\cite{T,TZ} in which other signals that changed electron transport also
resulted in behavioral responses.  For example, Taylor and Zhulin
\cite{TZ} cite cases in which metabolized substrates elicit tactic
response, as do chemicals which are able to donate electrons to or accept
electrons from the electron transport system.

\begin{figure}[h]
\centerline{\includegraphics[width=0.6\textwidth]{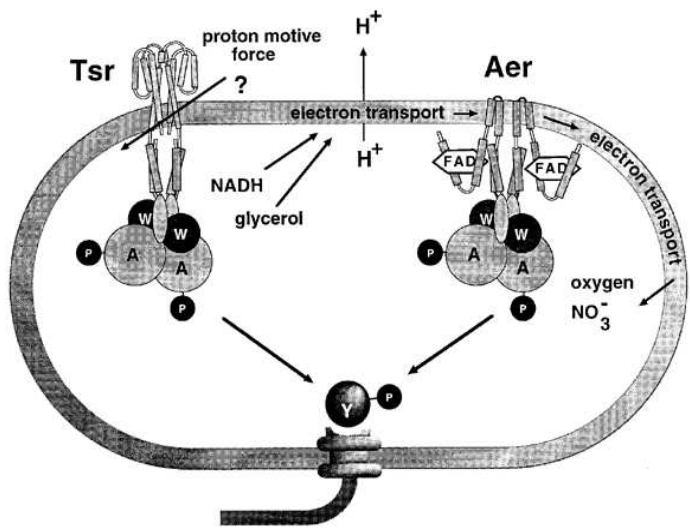}}
\caption[Aer and Tsr receptors.]{Aer and Tsr receptors.  (Figure from
Taylor and Zhulin, \cite{TZ})}
\label{tay_zhu2.eps}
\end{figure}

There are two receptors that have been proven to act as signal
transducers \cite{R}.  The two receptors, Aer and Tsr are shown in
Figure \ref{tay_zhu2.eps}.  Aer is a novel receptor which plays no
role in conventional chemotaxis.  Rebbapragada et al. \cite{R}
demonstrate that when mutants with a deactivated aer gene are
exposed to an oxygen gradient, they find optimal oxygen
concentrations much slower and less efficiently than wild type
(non-mutant)  bacteria.  Meanwhile, the same mutants still exhibit
normal chemotaxis indicating that the signal transduction pathways
for chemotaxis and aerotaxis differ. When expression of Aer is
restored, the aerotactic behavior returns.

Tsr, the other aerotaxis receptor, also works as a receptor for
conventional chemotaxis, sensing external environments as well as
monitoring internal pH.  The aer tsr double mutants were not capable
of aerotactic sensing; however, upon restoring one or both of the
Tsr and Aer receptors, aerotactic response returns \cite{R}. The
signal transduction pathway of aerotaxis is not well understood,
although it is believed to converge with the phosphorylation pathway
of conventional chemotaxis \cite{R2}. CheA, CheW and CheY are part
of the signal transduction pathway for aerotaxis \cite{R2}, but
there is evidence that adaptation is methylation-independent
\cite{T}), and it is much faster than the adaptation response in
conventional chemotaxis \cite{T}.

Aerotaxis is thought to be beneficial, because finding the
appropriate concentrations of oxygen is essential for the metabolism
of some species and can, in fact, be a more immediate need than
finding the appropriate nutrient levels \cite{T}. According to
Taylor and Zhulin \cite{TZ}, the aerotactic response might prevent
bacteria from getting trapped in anaerobic, growth-limiting
environments.  In their hypothesis, bacteria living in conditions
which support growth would mostly rely on their chemotactic behavior
and would use aerotaxis when the maintenance of optimal internal
energy levels becomes impossible. The mechanism for aerotaxis is
expected to operate on simpler principles than that of chemotaxis,
since it is "designed" to find a well-defined range of a single
chemical, oxygen.  On the other hand, chemotaxis allows bacteria to
choose between different types of nutrients of possibly different
concentrations and quality, as well as allowing adaptation to a wide
range of concentrations \cite{R, TZ, Z}.

\subsection{Aerotaxis experiments with {\it Azospirillum 
brasilense}}

\begin{figure}[h]
\centerline{\includegraphics[width=0.6\textwidth]{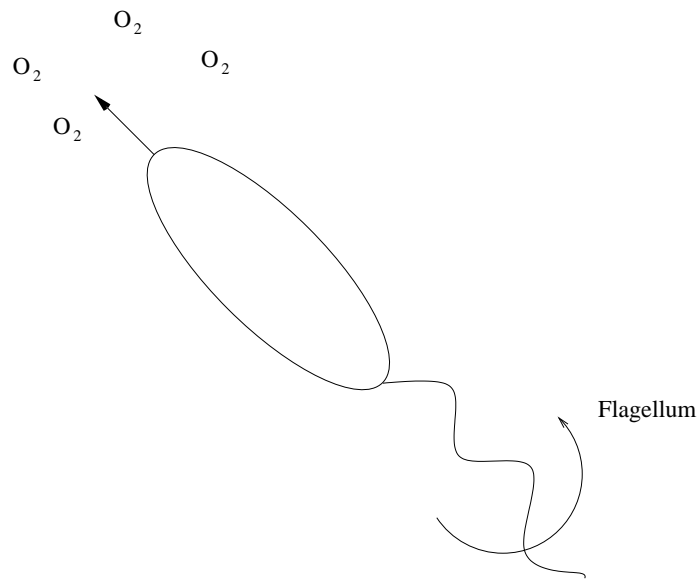}}
\caption[Azospirillum brasilense.]{Schematic figure of Azospirillum
brasilense. One flagellum is attached to the ellipsoidal body.}
\label{azospirillum.ps}
\end{figure}

{\it Azospirillum brasilense} is a 1-2 $\mu m$ nitrogen-fixing
plant-associated bacterium \cite{Z} with an ellipsoidal body to
which one flagellum is attached.  (See Figure \ref{azospirillum.ps}
for a schematic diagram of {\it A. brasilense}.)  Counterclockwise
rotation of the flagellum produces forward motion, while clockwise
rotation reverses the direction.  The essentially one dimensional
movement of {\it A. brasilense} makes it a simple organism to model.  
It is generally accepted that its positive aerotaxis (attraction) is
a response to changes in the proton motive force, and in Zhulin's
hypothesis \cite{Z} this is also the signal for negative aerotaxis.

{\it A. brasilense} is aerobic, but it prefers very low
concentrations of dissolved oxygen.  Zhulin's experiments
demonstrated by spatial and temporal assays that oxygen indeed acts
both as repellent and attractant, and Zhulin has provided evidence
\cite{Z} that both are linked to monitoring the proton motive force
inside the cell.  As before, one can reason that aerotaxis provides
an additional advantage for {\it A. brasilense} by guiding it to an
optimal range of oxygen concentration.  Nitrogen fixing can be
accomplished only in environments where the oxygen concentration is
below 1 percent \cite{Z}, and in these environments {\it A.
brasilense} remains capable of maintaining aerobic metabolism
\cite{Z}.

\begin{figure}[h]
\includegraphics[width=0.45\textwidth]{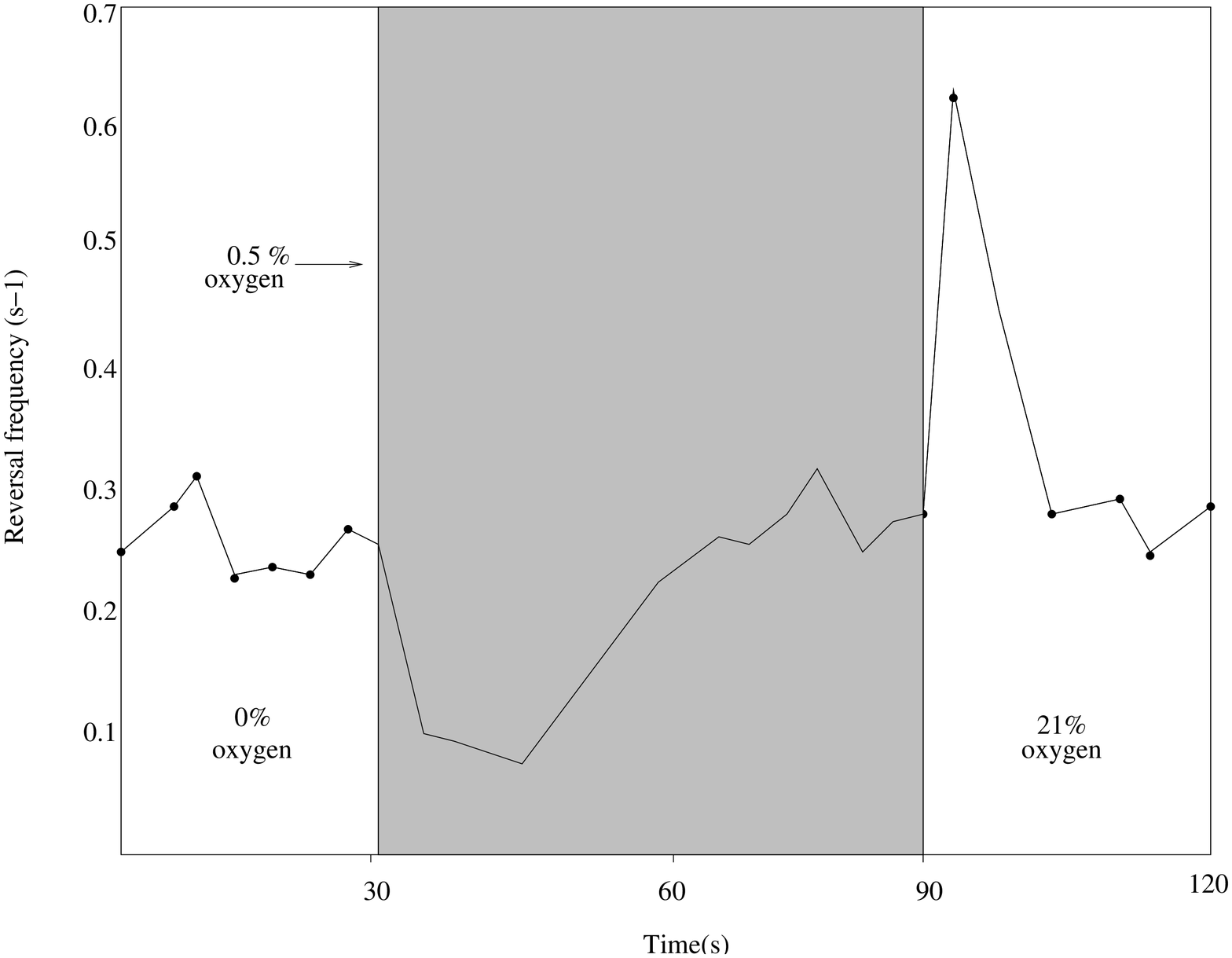} \qquad
\includegraphics[width=0.45\textwidth]{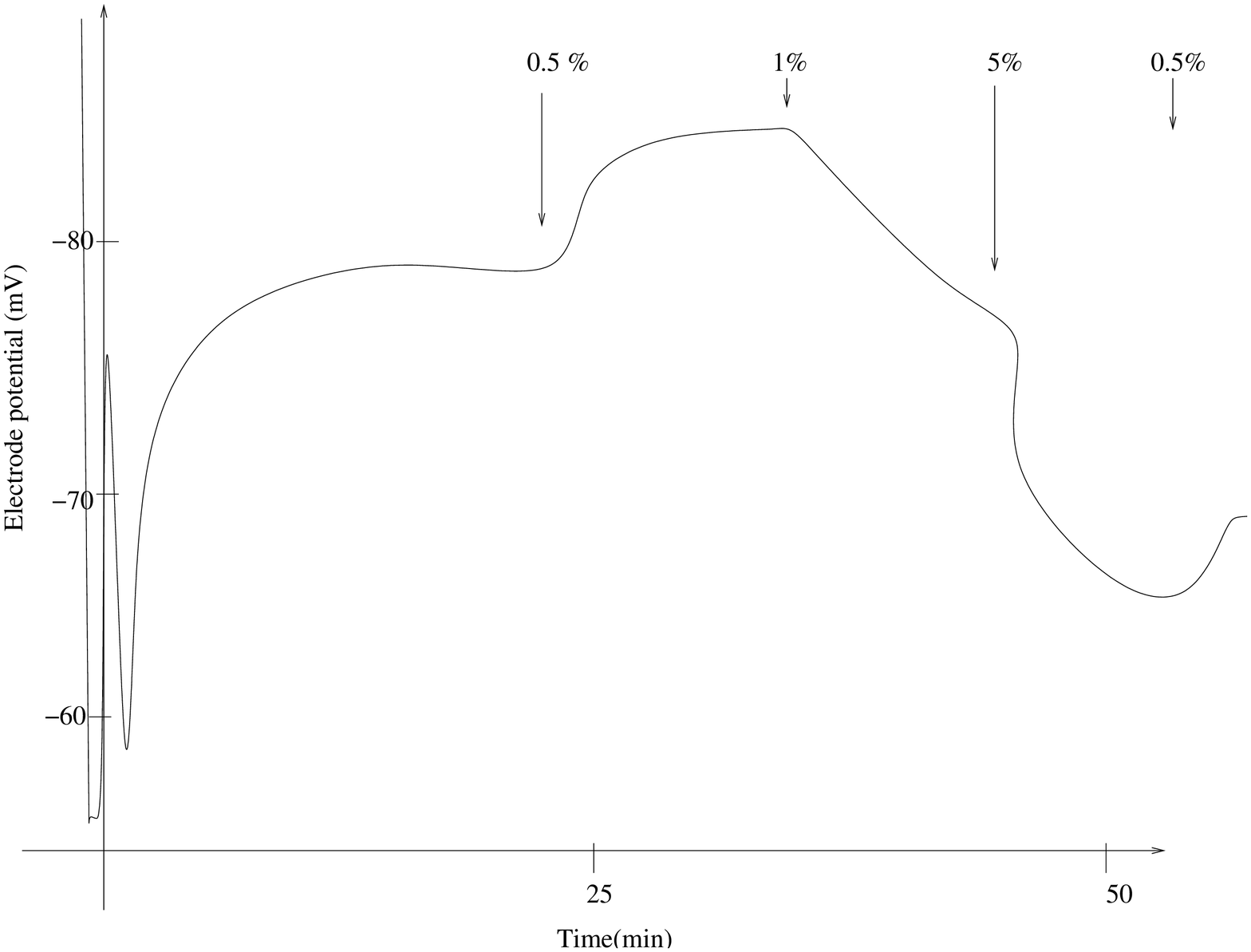}
\caption[Temporal assays.]{The graph on the left shows the turning
frequency as a function  of time.  Significant changes occur as the
oxygen concentration jumps to  0.5 percent from no oxygen, and as the
oxygen concentration jumps to 21 percent.  The figure on the right
shows the proton motive force as a function of time. (Figure based 
on Zhulin,\cite{Z}) } 
\label{tem_assays}
\end{figure}

In the temporal assays for aerotaxis, a small droplet of bacterial
suspension was spread on a slide \cite{Z} then exposed to different oxygen
concentrations.  Figure \ref{tem_assays} illustrates the results of these
experiments.  It shows that the bacterial turning frequency remained
constant at approximately 0.28 $s^{-1}$ for most changes in oxygen
concentration.  There were two instances where significant changes were
observed in the turning frequency. The reversal frequency dropped to 0.09
$s^{-1}$ when cells were ventilated with 0.5 percent oxygen after being
equillibrated to nitrogen (no oxygen).  The other noticeable change occurred
when cells adapted to 0.5 percent oxygen were exposed to 5 percent oxygen.  
In this case, the reversal frequency jumped to 0.49 $s^{-1}$.  Table
\ref{turn_freq} summarizes the results of these experiments.

\begin{table}[h]
\centering
\begin{tabular}{|c c c|}
\hline
Change in oxygen concentration (percentage) & {} & Reversal frequency
($s^{-1}$) \\ 
\hline
21 to 100 & ........... & 0.32 $\pm$ 0.04 \\
0 to 21   & ........... & 0.27 $\pm$ 0.04 \\
100 to 21 & ........... & 0. 24 $\pm$ 0.06 \\
21 to 0  &  ........... & 0. 22 $\pm$ 0.05 \\
0 to 0.5 &  ........... & 0.09 $\pm$  0.03 \\
0.5 to 5  & ........... & 0.49 $\pm$ 0.09 \\
\hline
\end{tabular}
\caption[Temporal assay.]{Turning frequencies as the oxygen
concentration changes.  (Table from Zhulin, \cite{Z})}
\label{turn_freq}
\end{table}

These findings indicate that cells are attracted to 0.5 percent
oxygen concentration, since both increases and decreases in the
oxygen concentration caused negative aerotactic response.  Zhulin
also calculated the proton motive force for various oxygen
concentrations and found that the highest values of proton motive
force corresponded to the 0.3-0.5 percent oxygen concentration
range.

In the spatial assay for aerotaxis, a flat 50 by 2 by 0.1 mm
capillary tube was filled with a solution containing {\it A.
brasilense} distributed uniformly, with no dissolved oxygen in the
solution.  At the open end of the capillary, the oxygen
concentration was maintained at various fixed levels, and this
induced a band formation inside the capillary within 50 seconds to 3
minutes, depending on the oxygen concentration. The oxygen diffused
into the capillary and was consumed at a certain rate by the
bacteria.  The bacteria aggregated to the favorable oxygen
concentration, and the region where their density was high was seen
as the band inside the capillary.  If the oxygen concentration was
100 percent, the band moved further away from the meniscus. If
nitrogen replaced oxygen, then the band moved closer to the open
end, and eventually disappeared.  If air was introduced again, then
the band reappeared. An unusual feature of these spatial assays was
the steepness of the gradients produced.  Oxygen concentrations
change from 20 percent to zero oxygen in about $ 1.6$ mm.

\begin{figure}[h]
\centerline{\includegraphics[width=0.6\textwidth]{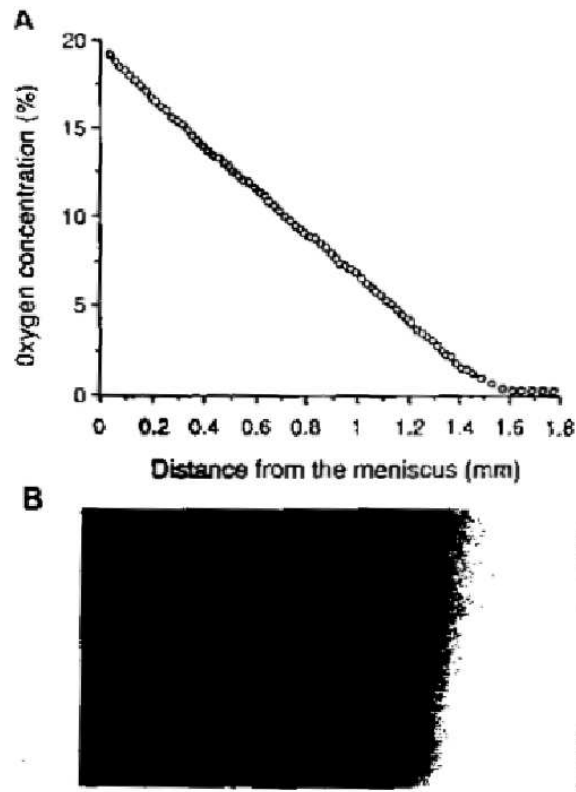}}
\caption[Band formation.]{Figure A shows the oxygen concentration as
a function of space.  Figure B shows the aerotactic band on the same
scale.  The band forms about 1.6 mm form the meniscus where the
oxygen concentration is 0.3-0.5 percent. (Figure based on Zhulin,
\cite{Z})} 
\label{spat_assay2.eps}
\end{figure}

In the spatial assay, the swimming speed of bacteria was observed.
Inside the band, the cells swam with an average speed of $ 49 \frac
{\mu m} {s} $ whereas outside the band the average speed was $14 -
22 \frac {\mu m} {s} $, depending on the oxygen concentration. The
cells swimming in either direction would swim straight through the
band then reverse direction immediately once they passed the band.  
The bacteria did not reverse their direction of swimming while
inside the band \cite{Z}.  The bacterial density inside the band was
nearly a hundred times that in front of the band (between the band
and the meniscus).  The bacterial density behind the band remained
approximately constant.  Figure \ref{spat_assay2.eps} shows the
bacterial band which evolves in the spatial assay, and, on the same
scale, the oxygen concentration corresponding to the various
bacterial densities.

Our mathematical model will have to answer the following questions
that arise based on the experiments: Can we create a model which
leads to the evolution an aerotactic band with similar
characteristics?  We would like to test whether a model based on the
simple turning rules which are suggested by the experiments would be
able to produce the observed pattern.  In particular, the hope is
that of the model would demonstrate that the measurable quantities,
such as the bacterial density inside the band and outside the band,
the width of the band, and its distance from the meniscus, can be
explained by assuming our simple turning rule.
   
Some deeper questions to pose beyond this model would be: What kind
of signal transduction mechanism do aerotactic bacteria have?  In
particular, what sort of receptors make it possible for oxygen to
act both as an attractant and a repellent?  What sort of properties
of this signal transduction mechanism allow the cells to react to
such steep concentration gradients?  How can the turning frequency
of the cells change so abruptly at the boundaries of the aerotactic
band?  What is the underlying mechanism that allows the turning
frequency to respond immediately?  These questions are likely to be
answered by further experiments.

\subsection{Mathematical models for conventional chemotaxis} 
\label{chemotaxis_models}

This section provides a brief overview of the mathematical tools
used to describe conventional chemotaxis and briefly discusses why
these existing models cannot be employed to study aerotaxis.

The first mathematical models for bacterial chemotaxis were created
in the early 1970's by Keller and Segel \cite{KS,S}.  In their work
\cite{KS}, they show how local gradient detection by individual
cells produces fluctuations in their path, and how the average over
many cells corresponds to a macroscopic flux.  They derive the
macroscopic flux equations based on individual bacteria.  Although
in this paper they make incorrect assumptions about the way bacteria
are able to detect gradients (namely, they assume that bacteria can
compare concentrations at their head and tail), they arrive at a
very useful expression for the macroscopic flux: \begin{eqnarray*}
J(x) = - \mu (\frac {db} {dx} ) + \chi b \frac {dL} {dx}
\end{eqnarray*} In this equation $\mu(L)$ is the "diffusion", or
motility coefficient, and $\chi (L)$ is the chemotactic coefficient.
$b(x,t)$ is the bacterial density, and $L(x,t)$ is the concentration
of the ligands.  The rate of change of the bacterial density
$b(x,t)$ can be given by $$ \frac {\partial b} {\partial t} = -
\nabla J $$ with the appropriate boundary condition.  In one
dimension, this becomes the well known Keller-Segel chemotaxis
equation: \begin{eqnarray} \frac {\partial b} {\partial t} = - \frac
{\partial} {\partial x} ( - \mu \frac {\partial b} {\partial x} +
\chi b \frac {\partial L} {\partial x} ) \label{KSeqn}
\end{eqnarray} Again, boundary conditions must be imposed.  It is
clear in the equation (\ref{KSeqn}) that the first term on the
right-hand side is diffusion due to random motility, and the second
term is due to chemotactic flux.  Usually it is also assumed that
the spatial gradient of the concentration is also small, and
therefore $\frac {\partial L} {\partial x}$ can be approximated by a
constant, and absorbed in $\chi$, the chemotactic coefficient.

In a later paper \cite{S}, Segel derives the same chemotactic
equation (\ref{KSeqn}) based on changes in receptor configuration as
a result of attractant (or repellent) binding.  In this model, there
is no need to use the incorrect hypothesis that bacteria are able to
compare concentrations at different parts of their body.  Instead,
this model is developed by writing down a separate equation for
left- and right-moving bacteria in different receptor
configurations.

In order to to arrive at (\ref{KSeqn}) from this system of
equations, Segel must use the assumption that spatial gradients are
small.  Almost the same analysis is summarized more lucidly by
Gr\"unbaum; therefore, in order to understand where the assumption
of small spatial gradients is needed, it is appropriate to look at
Gr\"unbaum's article ``Advection-diffusion equations for internal
state-meditated random walks"  \cite{G}.

Gr\"unbaum gives a simple argument regarding the validity of
parabolic advection-diffusion approximations to hyperbolic advection
equations.  He starts with a system of one-dimensional advection
equations describing left- and right-moving individuals, similar to
the Segel system in \cite{S}.  \begin{eqnarray} \frac {\partial b^+}
{\partial t} + v \frac {\partial b^+} {\partial x} = \sigma^- b^- -
\sigma^+ b^+ \nonumber \\ \frac {\partial b^-} {\partial t} - v
\frac {\partial b^-} {\partial x} = \sigma^+ b^+ - \sigma^- b^-
\label{Gadv} \end{eqnarray}

Boundary conditions are omitted again for the time being.  The
turning rates are denoted by $\sigma^{+}$, $\sigma^{-}$ and velocity
by $v$.  The velocity is chosen to be positive for the left-moving
bacteria, and negative for the right-moving bacteria.  Total
bacterial density is given by $b(x,t) = b^+ (x,t) + b^-(x,t)$; in
other words, the total bacterial density is the sum of the right-
and left-moving terms.  Gr\"unbaum also defines the density flux,
$J(x,t)$ as $J(x,t) = v (b^+ (x,t) - b^-(x,t))$, and introduces two
new variables, $\sigma_0 = \frac {1} {2} (\sigma^+ + \sigma^-)$ and
$\Delta \sigma = \frac {1} {2} (\sigma^+ - \sigma^-)$.  By taking
the sum and difference of the two equations in (\ref{Gadv}) and
making the appropriate substitutions, he arrives at the following
two equations: \begin{eqnarray} \frac {\partial b} {\partial t} = -
\frac {\partial J} {\partial x} \label{cons2} \\ \frac {1} {2
\sigma_0} \frac {\partial J} {\partial t} + J = \frac {1} {2
\sigma_0} (-v^2 \frac {\partial b} {\partial x} + 2 \Delta \sigma v
b) \label{Gdiff} \end{eqnarray} In equation (\ref{Gdiff}), the
objective is to estimate the term $J_t$ using the conservation
equation (\ref{cons2}).

Gr\"unbaum does this the following way.  Assume that spatial
derivatives are small, then a small $J_x$ implies that $b_t$ is
small, or the bacterial density varies on a slow time scale.  The
solution to (\ref{Gdiff}) is the linear combination of the
homogeneous and inhomogeneous solutions.  The inhomogeneous terms
all contain $b(x,t)$ or a spatial derivative of $b$; therefore, they
vary on a slow time scale. However, the solution to the homogeneous
equation \begin{eqnarray*} \frac {1} {2 \sigma_0} \frac {\partial J}
{\partial t} + J = 0 \end{eqnarray*} will be an exponential term
acting on a faster time scale.  This implies that after a short
initial period, $J$ approaches a quasi-equilibrium.  This gives the
following equation: \begin{eqnarray} J \approx \frac {1} {2
\sigma_0} (-v^2 \frac {\partial b} {\partial x} + 2 \Delta \sigma v
b) \end{eqnarray}

By substituting this approximation for the flux into (\ref{cons2}),
one can recover the Keller-Segel equation: \begin{eqnarray} \frac
{\partial b} {\partial t} \approx \frac {\partial} {\partial x} (\mu
\frac {\partial b} {\partial x} - \chi b) \label{Geqn}
\end{eqnarray} with \begin{eqnarray} \mu = \frac {v^2} {2 \sigma_0}
\nonumber \\ \chi = \frac {v \Delta \sigma} {\sigma_0}
\label{Gconst} \end{eqnarray} The advantage of this formulation over
the Keller-Segel equation is clear, since (\ref{Gconst}) shows that
the two important coefficients, the motility coefficient and the
chemotactic coefficient, can be expressed in terms of empirically
observable quantities, namely the cell velocity and the turning
rates. It is also clear from Gr\"unbaum's analysis that the
assumption of a small spatial gradient is necessary in order to
arrive at the Keller-Segel chemotaxis equation, (\ref{Geqn}).

In addition to the assumption of a small spatial gradient, another
difficulty with this model is that an exact measure of the
difference of turning frequencies, $\Delta \sigma$, is very
difficult to obtain empirically, or to approximate analytically. An
analytical expression for $\chi$, the taxis coefficient, can be
derived in terms of the characteristic time scale of a run, the
bacterial speed and the characteristic attractant concentration.  
This analytical result does not rely on having to measure the
turning frequency, but is rather obtained by a perturbation method
which separates time scales of the reaction.  In this perturbation
analysis, the characteristic run time is obtained based on the
assumption that adaptation is slow.  This is not a correct
assumption in the aerotaxis experiments.

Another significant development in mathematical modeling of
chemotaxis came in the late 1980's and early 1990's.  A large effort
was made to unite knowledge about receptor dynamics and signal
transduction pathways (many times referred to as 'internal state
dynamics') to random motility \cite{DT, MT}.  In many of these
models (e.g. \cite{MT}), not only bacterial chemotaxis but
chemotaxis of animal cells is considered, making the internal state
dynamics far more complex than in the bacterial case. There were
several models, for example those by Barkai and Leibler \cite{BL}
and by Tranquillo et al. \cite{DT, MT} which included very detailed
models for the biochemical mechanism.

One such detailed model of chemotaxis in an animal cell is due to
Tranquillo \cite{MT}. The general approach of Tranquillo and his
collaborators is as follows. Identify a simplified (but still
realistic) scheme to describe all possible receptor states and
intracellular chemicals that govern the turning behavior.  All rate
constants and parameters are approximated based on empirical data.  
Each of the receptor states are interdependent stochastic variables,
and their time evolution is calculated using a multivariate
probability density function.

The jump processes between the various states can be approximated by
continuous stochastic differential equations by making certain
assumptions about the system.  The probability function must be
linearized in order to be represented by an analogous system of
stochastic differential equations.  However, the linearization
procedure implicitly assumes small fluctuations in the stochastic
variables.  This, again, is an assumption that would be violated in
the aerotaxis experiments.

In order to analyze the stochastic differential equations,
deterministic equations are derived either using averages of the
stochastic variables or transformations from the stochastic
variables to new, deterministic variables.  With this method, the
Fokker-Planck equations are derived and used to estimate cell
movement on longer time scales \cite{MT}.

These models (\cite{DT, MT}) are highly sophisticated biologically and
mathematically alike, although lack the satisfying simplicity of the
earlier Keller-Segel equations.  As mentioned above, they also assume
certain features of the internal state dynamics that would be violated in
steep attractant gradient; therefore, they cannot be used for our
purposes.

Barkai \& Leibler also create a model which is very closely built on
experimentally obtained data.  In ``Robustness in simple biochemical
networks" \cite{BL} they examine the question of receptor adaptation in
bacteria, and they propose a quantitative model in which a wide range of
biochemical parameters are admissible.  The model Barkai \& Leibler
propose is a system of ordinary differential equations that is based on
the accurate description of the possible receptor states.  The key to the
model is a small network that has two states: an active and an inactive
state.  In the active state the external signal leads to a fast response,
whereas in the inactive state there is no response.  The shift between the
active and inactive states is a slower process, which corresponds to the
receptor methylation.  The system exhibits perfect adaptation if the
active state is independent of the magnitude of the outside stimulus.  A
model based on the same principles is presented in the chapter on animal
cell chemotaxis, Section \ref{perfect}.  The contribution of Barkai \&
Leibler is significant for bacterial chemotaxis, because their work shows
that fine-tuning the model parameters is unnecessary in order to achieve
adaptation.  However, in our model for aerotaxis there is no methylation
adaptation, and we argue that if there is adaptation, it must occur on the
fast time scale.

There is also some literature that focuses directly on aerotactic
behavior and models for aerotaxis \cite{H, HPK}.  In \cite{HPK},
certain bioconvection patterns are described quantitatively by
developing a model of aerotaxis. Kessler's experiments described in
this paper involve an initially well-stirred suspension in which
cells swim upwards toward oxygen, then, after the top layer becomes
sufficiently denser than the bottom layer, an instability occurs.
The overturning instability evolves into the observed patterns.  
Some of the phenomenon is similar to the Zhulin aerotaxis
experiments, namely, here bacteria consuming oxygen create the
gradient while the oxygen level is fixed at the meniscus.  However,
in the Kessler experiments the bacterial convection stirs the
solution, and the bacteria carry oxygen into deeper layers of the
solution.

The mathematical model describing these experiments consists of a
conservation equation for the cells and a reaction-diffusion
equation for the oxygen concentration.  The Navier-Stokes equations
are not required, since no bulk fluid flow is assumed.  The authors
arrive at the same equation as the Keller-Segel equation for the
bacteria.  The whole system is: \begin{eqnarray} \frac {\partial b}
{\partial t} = - \nabla \cdot [b( u + v) - D \cdot \nabla b] \\
\frac {\partial L} {\partial t} = - \nabla \cdot ( L u - D_L \nabla
L)-k b \label{HPKeqn1} \end{eqnarray}

In this model $L(x,y,z,t)$ is the concentration of oxygen,
$b(x,y,z,t)$ is the density of bacteria, as above, $u$ is the fluid
velocity taken to be zero, $v(\Theta)$ is the average cell velocity
which is a function of a dimensionless measure of $L$ denoted
$\Theta$, $D$ and $D_L$ are the diffusion coefficients of the
bacteria and of oxygen, respectively, and, finally, $k(\Theta)$ is
the coefficient of nutrient consumption by the bacteria.  The
boundary conditions are applied at $z=0$ and at $z=-h$, the bottom
of the suspension.  They are given by:  \begin{eqnarray} L=L_0 \
\hbox{\rm at $ z=0$,} \ \frac {\partial L} {\partial z} = 0 \
\hbox{\rm at $z=-h$} \\ v_z b - D \frac {\partial b} {\partial z} =0
\ \hbox{\rm at $z=0$ and $-h$} \end{eqnarray}

These boundary conditions mean that there is no flux of oxygen or
bacteria at any of the boundaries.  The initial condition is given
for a well-stirred solution and uniform suspension: \begin{eqnarray}
L(z, 0)= L_0 \\ b(z, 0) = b_0 \end{eqnarray}

By non-dimensionalizing the equations, the authors arrive at a
perturbation problem whose analysis leads to a good quantitative
description of the experiments.  There are several reasons why this
model would not be valid for the Zhulin aerotaxis experiments.  
First, just like the Keller-Segel equation, Hillesdon et al. make
use of the small gradient assumption implicitly.  Also, the Kessler
experiments are examples of kinesis rather than taxis.  In kinesis
the movement is determined by local concentration of attractant, not
by concentration gradients, as it is in taxis.

In his review article \cite{H}, Hill summarizes other models of
various types of chemotaxis related to pattern formation (gyrotaxis,
geotaxis, phototaxis).  These equations make very similar
assumptions to the models discussed above, namely, conservation
equations for the bacteria involving gradients of the flux due to
swimming and random motility. When, in addition, randomness in drift
speed and direction of motility are introduced, the Fokker-Planck
equation is used to describe tactic behavior.

In all the above mentioned models, there are two important aspect of
chemotaxis that must be taken into account: the gradient of
attractant concentration and the adaptation time.  For small
gradients and fast adaptation, it is simple to find approximations
to the hyperbolic advection equations describing bacterial motion.  
Gr\"unbaum also shows that for slow adaptation in shallow gradients,
one is able to simplify the equations to a parabolic system.  
However, if the attractant gradients are large, then these
approximations do not work for slow adaptation any more.  In this
case, the only possible approach is to try Monte Carlo simulations.
The model of aerotaxis presented in the next chapter uses the fact
that while gradients are large, the adaptation time in our case is
fast.  This allows us to solve the original hyperbolic system.
 
\begin{figure}[h]
\centerline{\includegraphics[width=0.6\textwidth]{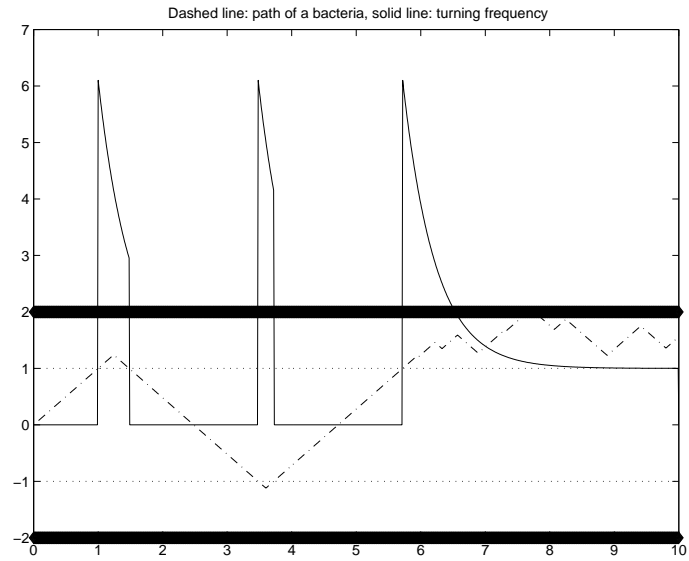}}
\caption[Monte-Carlo simulation.]{Monte-Carlo simulation.}
\label{mc2.ps}
\end{figure}

It is also possible to give a heuristic argument for why
conventional chemotaxis models that rely on slow adaptation cannot
be used to model the Zhulin aerotaxis experiments.  The following
Monte-Carlo simulation,Figure \ref{mc2.ps} illustrates this 
reasoning.

The wide solid lines (at -2 and 2) represent the sides of the
capillary tube (so the bacteria are confined to this region), and
the dotted lines (at -1 and 1) represent the favorable oxygen
concentrations.  The following rules govern the movement of each
cell:  \begin{itemize} \item a cell moves straight to the left or to
the right with a constant velocity, $v$; \item turning frequency
inside the favorable region (between -1 and 1) is $\sigma = 0$;  
\item the cell leaves the band at some random time, $\tau$, and at
this time the turning frequency jumps form $\sigma = 0$ to $\sigma =
c$; \item characteristic time of adaptation is $t_a$; \item
adaptation to the baseline turning frequency of $\sigma=0$ is
exponential, and it is given by $\sigma = c e^{- \frac {t-\tau }
{t_a}}$. \end{itemize}

The adaptation in these simulations is assumed to be slow, which
means that $\frac {1} {c}$, the time of a straight run, is of the
same order of magnitude as $t_a$, the characteristic time of
adaptation.  Later, in our model of aerotaxis we assume that $t_a <<
\frac {1} {c}$.

In these simulations, the particle's velocity and turning frequency
are given deterministically.  When the cell is outside the band, the
time of turning, $\tau$, is determined based on the difference of a
uniformly generated random number between 0 and 1 and the turning
frequency.

In the figure, we can see a typical run of the simulation.  Most of
the time the bacterium stays inside the optimal oxygen
concentration, because upon leaving the band, its turning frequency
jumps from 0 to c, and it is likely to turn back into the band.  
However, outside the band the turning frequency is large for a
period of time ($t_a$) due to slow adaptation, and it frequently
causes the cell to keep tumbling and getting trapped outside the
optimal environment.  Running the simulation up to 10,000 times, the
bacterial density inside the band was only three times the density
outside the band.  This is clearly very different from the 100:1
ratio observed experimentally.

One must conclude that there are no existing models of bacterial
chemotaxis (other than Monte Carlo simulations)  that can describe
the behavior in steep gradients.  Most existing chemotaxis models
also assume a slow adaptation of the turning rates.  Exact
mathematical descriptions of turning rates based on slow adaptation
are very difficult to analyze, and in order to create tractable
equations, one must make approximations.  The approximations involve
assumptions of small spatial gradients, since this allows continuity
of internal state variables.  Assuming slow adaptation and a steep
spatial gradient, no approximations are possible leading to simple
mathematical expression. This suggests that it would be futile to
attempt to model the Zhulin experiments with already existing
chemotaxis equations.  However, since aerotaxis is known to have
fast adaptation, mathematical expression of the turning rates is
much simpler; thus, we can develop a different model in which one
need not rely on approximations based on shallow gradients.

\clearpage
\section{Model} \label{model1}

\subsection{Mathematical model for aerotaxis}

Now we can present the mathematical model for aerotaxis which
describes the pattern formation as a result of steep concentration
gradients.  A brief derivation of the equations is given below.  
The bacteria's movement is governed by advection-reaction equations.  
The advection term describes the directed swimming of bacteria,
while the reaction terms denote the turning of bacteria in response
to the oxygen gradients.  The evolution of the ligand (oxygen)
concentration is given by a reaction-diffusion equation which is
coupled to the advection equations through the term describing
oxygen consumption by the bacteria.

We assume that bacterial movement is one-dimensional.  Although some
turning is possible, this assumption is based on empirical evidence
of {\it A. brasilense} swimming.  The cells are observed to swim
either forward, or, upon changing their direction, backward. The
hypothesis that the turning rates are dependent on oxygen
concentration and the oxygen gradient is the main assumption of the
model, and it is discussed at length below.  The boundary conditions
imposed simply mean that all left-moving cells turn to the right at
the left boundary, and similarly, all right-moving cells turn to the
left at the right boundary.  The assumption of the conservation of
the number of bacteria is used later, and it justifies the lack of
birth and death terms in the equations.  Since the band forms on the
order of minutes, this is a reasonable omission.

\begin{eqnarray}
\frac {\partial r} {\partial t} =  \frac{\partial (-v r)} {\partial
x} - f_{rl} r + f_{lr} l \label{main1} \\
\frac {\partial l} {\partial t} = \frac {\partial (vl )} {\partial x} +
f_{rl} r - f_{lr} l \label{main2} \\
\frac {\partial L} {\partial t} = D \frac {\partial^2 L} {\partial
x^2} - k (r + l) \label{main3}
\end{eqnarray}
\begin{eqnarray*} r(0)=l(0)  \\
r(a) = l(a)  \end{eqnarray*}

$r(x,t)$  - right moving bacteria \\
$l(x,t)$  - left moving bacteria \\
$L(x,t)$ - ligand (oxygen) concentration \\
$f_{rl}(L)$ - rate of turning from
right-moving  to left-moving cell \\
$f_{lr}(L)$ - rate of turning from
left-moving to right-moving cell \\
$v(L(x))$ - bacterial speed \\
$D$ - diffusion coefficient for oxygen \\
$k$ - rate of oxygen consumption by bacteria \\
$a$ - length of the capillary \\

The initial condition for the left- and right-moving bacterial
densities is the same constant for all positions, $x$, such that the
sum of the two populations is the total bacterial density.  The
initial condition for the oxygen concentration is $L_0$ at the left
boundary, and zero everywhere else.

\begin{figure}[h]
\centerline{\includegraphics[width=0.6\textwidth]{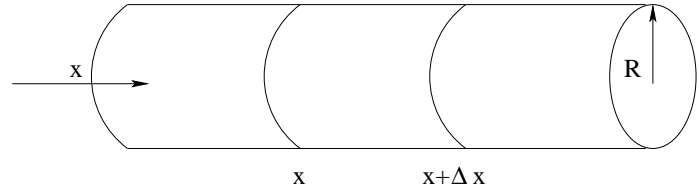}}
\caption[Derivation of the advection equation.]{A volume element that
bacteria swim through.}
\label{volume.eps}
\end{figure}

A brief justification (similar to the reasoning of Segel in
\cite{S}) of using the advection equation (\ref{main1}) is as
follows. Consider the cells swimming in a capillary of a fixed
cross-sectional area, $R$.  Let us look at a short section of the
capillary, from $x$ to $x+\Delta x$.  The density of right-moving
cells in this section is given by $r R \Delta x$. Then the rate of
change of the $r$, the right-moving bacterial density, with time
respect to time is equal to (i) change due to reversal of direction
to become left-moving, and, similarly, left-moving bacteria turning
to become right-moving; (ii) change due to cells swimming into and
out of the slice.

The change due to turning is quite straight-forward.  It is the turning
rate of cells times the bacteria in the given volume, or $f_{lr} l R
\Delta x - f_{rl} r R \Delta x$.  The term due to cell swimming is the
rate at which cells flow into the cell, and the rate at which they flow
out, $$v(x)  r(x) R  - v(x+ \Delta x) r(x + \Delta x) R. $$ This leads
to $$ \frac {\partial (r R \Delta x)} {\partial t} = f_{lr} l R \Delta x
- f_{rl} r R \Delta x - R [ r(x+ \Delta x) v(x + \Delta x) - r(x) v(x)].  
$$ Dividing by $ R \Delta x$ and letting $\Delta x$ approach zero, we
obtain $$\frac {\partial r} {\partial t} = - \frac{\partial (v r)}
{\partial x} - f_{rl} r + f_{lr} l $$ Equation (\ref{main2}) can be
obtained in a similar fashion.

The most important question regarding the model is the determination of
the the turning frequencies.  There are two questions to be answered: what
function of $L$ are the turning frequencies?  What is the biological
evidence in support of this mathematical expression?  Both of these
questions are answered below.

\subsection{Expression for the turning frequencies}

\begin{figure}[h]
\centerline{\includegraphics[width=0.8\textwidth]{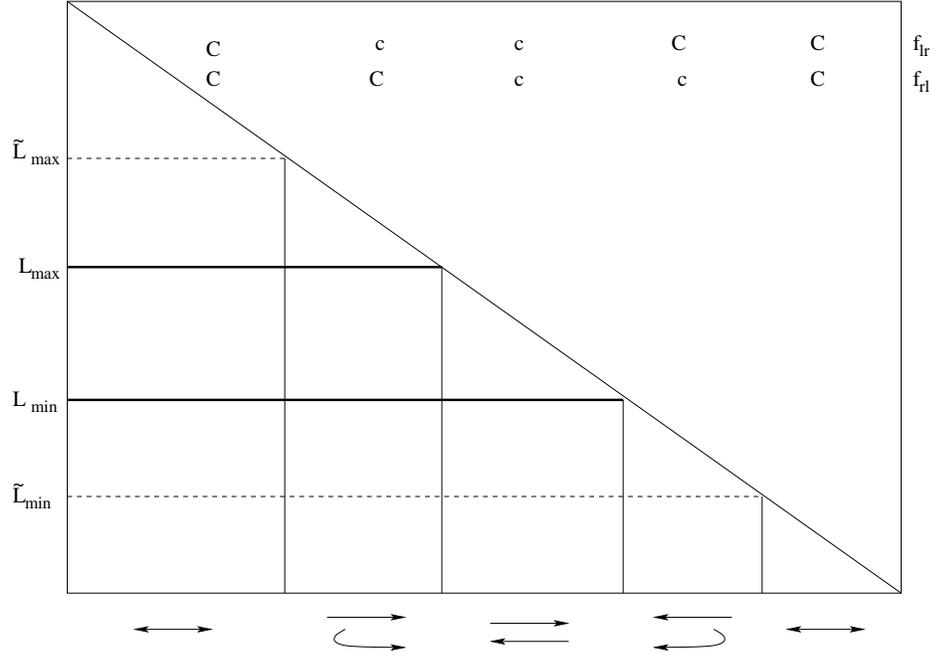}}
\caption[Turning rates.]{Turning rates of the bacteria.}
\label{new_turn_freq.eps}
\end{figure}

We must consider the mathematical expression for the turning rates.  
It is known from the experiments that there is a range of
concentrations where bacteria do not turn around.  We can call this
range $L_{min}$ to $L_{max}$.  This is the range where the proton
motive force (PMF) is the highest; therefore, it is the preferred
concentration range where the bacterial band develops.  Outside this
range, there must be some threshold values, $\tilde{L}_{min}$ and
$\tilde{L}_{max}$ such that if a cell is between $\tilde{L}_{min}$
and $L_{min}$, or between $L_{max}$ and $\tilde{L}_{max}$, it turns
back inside the band.  The bacterial frequencies are shown in Figure
\ref{new_turn_freq.eps}.

However, the turning rates cannot depend on the values of oxygen
concentration alone, since if that were the case, a cell getting to
$\tilde{L}_{min}$ from outside the band would have to turn, and
would therefore never enter the favorable oxygen range.  This
suggests that the bacteria must be able to retain some additional
information about the environment, for example the gradient of the
ligand.  This would allow a cell arriving to $\tilde{L}_{min}$ to
determine whether to keep on swimming (if it came from outside of
the band) or to turn around (if it came from within the band).

\begin{figure}[h]
\centerline{\includegraphics[width=0.8\textwidth]{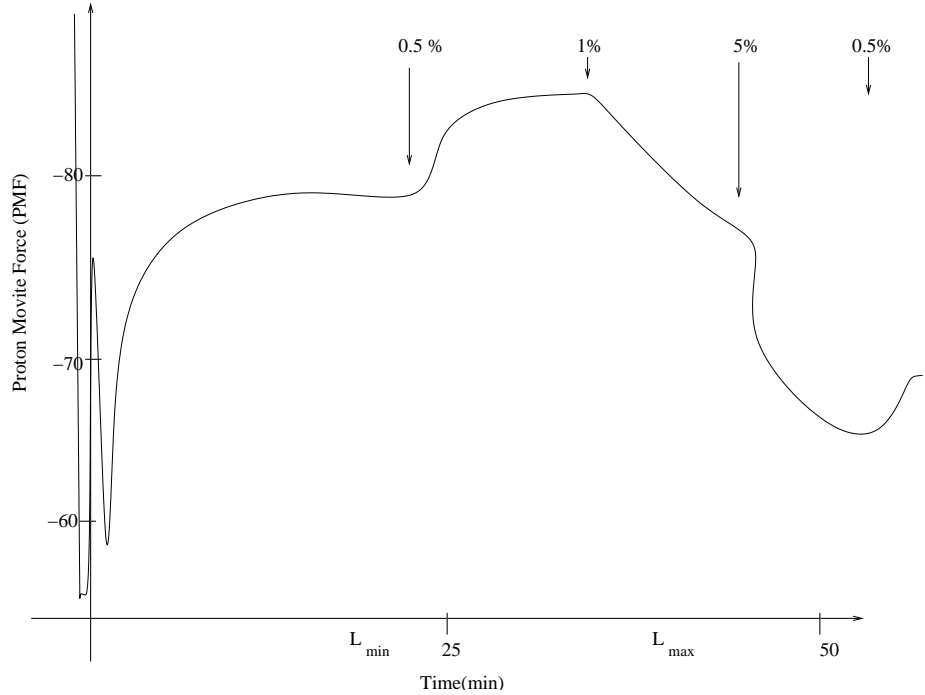}}
\caption[PMF vs. oxygen.]{The figure shows the proton motive force
(PMF) versus the oxygen concentration. The highest proton motive force
is observed between $L_{min}$ and $L_{max}$.  PMF is increasing
between $\tilde{L}_{min}$ and $L_{min}$ and decreases between
$L_{max}$ and $\tilde{L}_{max}$. (Figure based on Zhulin, \cite{Z}.)}
\label{pmf.eps}
\end{figure}

In our hypothesis, this additional information comes from monitoring the
proton motive force inside the cell.  As it is shown in Figure
\ref{pmf.eps}, the PMF has its largest value, and it is a constant, when
the cell is between $L_{min}$ and $L_{max}$.  Outside the band, the PMF
has a low value, and it is also a constant.  The turning signal for a cell
is the negative temporal gradient of the PMF.  Positive temporal gradients
and the high constant value of PMF suppresses tumbling. Cells are able to
detect temporal gradients of PMF by swimming through a spatial gradient of
oxygen (which is linearly proportional to time).  The biological
justification for such turning rates is given in Appendix
\ref{app_aerotax}. Further details on a model of a receptor producing such
turning rates are also included.

The easiest way to construct such turning rates is by defining piecewise
linear functions as is shown in the above figure (Fig.  
\ref{new_turn_freq.eps}).  We express the turning rates of the
left-turning bacteria and the turning rate of the right-turning bacteria
separately, so now the two rates only depend explicitly on the ligand
concentration, and their dependence on the oxygen gradient is implicit.  
(If a general turning rate was given for all bacteria, this rate would
explicitly depend on the gradient of oxygen.) This simple choice for the
turning frequencies makes the system of equations almost linear, with the
only non-linearity resulting from the dependence of turning rates on the
ligand gradient.  Now we can explicitly write down the turning rates.

\begin{equation}
        f_{rl} = \left\{ {\begin{array} {lllll}
        C, &  L < \tilde{L}_{min}  \\
        c, &  \tilde{L}_{min} < L < L_{min} \\
        c, &  L_{min} < L < L_{max} \\
        C, &  L_{max} < L < \tilde{L}_{max}  \\
        C, &  \tilde{L}_{max} < L
\end {array}} \label{f1}
\right.
\end{equation}

\begin{equation}
f_{lr} = \left\{ {\begin{array} {lllll}
        C, &  L < \tilde{L}_{min} \\
        C, &  \tilde{L}_{min} < L < L_{min} \\
        c, &  L_{min} < L < L_{max} \\
        c, &  L_{max} < L < \tilde{L}_{max}  \\
        C, &  \tilde{L}_{max} < L 
\end {array}} \label{f2}
\right.
\end{equation}

In all the above equations $C$ is some larger turning rate than $c$.  The
exact value of $C$ and $c$ will be specified later.  Now with equations
(\ref{main1}), (\ref{main2}), (\ref{main3}), (\ref{f1}) and (\ref{f2}),
we have the complete system of equations describing bacterial swimming.

\subsection{Non-dimensionalization and scaling}

We have the full system of equations for the aerotaxis experiments
including turning rates and boundary conditions.
\begin{eqnarray*}
\frac {\partial r} {\partial t} =  \frac{\partial (-v r)} {\partial
x} - f_{rl} r + f_{lr} l  \\
\frac {\partial l} {\partial t} = \frac {\partial (vl )} {\partial x} + 
f_{rl} r - f_{lr} l  \\
\frac {\partial L} {\partial t} = D \frac {\partial^2 L} {\partial x^2} -
k (r + l)
\end{eqnarray*}

\begin{equation} \nonumber
        f_{rl} = \left\{ {\begin{array} {lllll}
        C, &  L < \tilde{L}_{min} \\
        c, &  \tilde{L}_{min} < L < L_{min} \\
        c, &  L_{min} < L < L_{max} \\
        C, &  L_{max} < L < \tilde{L}_{max}  \\   
        C, &  \tilde{L}_{max} < L 
\end {array}}
\right.
\end{equation}
\begin{equation} \nonumber
f_{lr} = \left\{ {\begin{array} {lllll}
        C, &  L < \tilde{L}_{min} \\
        C, &  \tilde{L}_{min} < L < L_{min} \\
        c, &  L_{min} < L < L_{max} \\
        c, &  L_{max} < L < \tilde{L}_{max}  \\
        C, &  \tilde{L}_{max} < L 
\end {array}}
\right.
\end{equation}
\begin{eqnarray*} r(0)=l(0)  \\
r(a) = l(a)  \end{eqnarray*}

Before the system can be analyzed, it is important to non-dimensionalize
all variables.  The following parameters are important: \\

Size of capillary tube: 50x2x0.1 mm \\
Preferred $[O_2]$: 0.3-0.5 percent \\   
Band width: 0.2 mm \\
Distance of band from capillary tube end: 1.6 mm \\
Time of band formation: 50 sec - 3 min \\
$[O_2]$ outside the capillary tube: 21 percent \\
Speed: $ 40 \frac {\mu m} {sec} $\\
Diffusion coefficient: $ 2 \cdot 10^{-9} \frac {m^2} {sec} $\\
Turning frequency: $ 1 sec^{-1}$ \\
Rate of oxygen consumption: $ 3 \cdot 10^{-11} \frac {\mu M} {(cell)
(sec)} $ \\

The appropriate spatial scale can be found by estimating the
distance which would allow the bacteria to outrun the invasion of
the oxygen.  The diffusion of oxygen implies that $x \sim \sqrt{D t}
$, whereas the distance for the escape of bacteria is $x \sim v t$.  
From these two expressions we can approximate the time for the cells
to outrun diffusion, it is $t_0 \sim \frac {v} {D} \approx $ 5 sec,
which means that $x_0 \sim v t_0 \approx $ 100 $\mu m$.  This
suggests that the time scale should be on the order of about 10
seconds (which agrees with the experiments where the band develops
between 50 seconds and 3 minutes).  For the spatial scale, we choose
2 mm, since this is the length of the region where the bacteria are
found, but from the scaling argument it is clear that the resolution
must be smaller than 100 $\mu m$.  Based on this we arrive at the
characteristic scales summarized in Table \ref{char_scales}.

\begin{table}[h!]
\centering
\begin{tabular}{|c c c|}
\hline
Measurement &   {}  &  One unit \\
\hline
Length &        ...................... & 2 mm \\
Time &          ...................... &  10 sec  \\
Oxygen concentration &  ......................  &1 $\frac {\mu M} {ml}$
\\
Bacterial concentration &......................  & $ 2 \cdot 10^7 \frac
{cells} {ml}$ \\
\hline
\end{tabular}
\caption[Characteristic scales.]{Table gives the units of
measurement for length, time, oxygen concentration and bacterial
concentration.}
\label{char_scales}
\end{table}

\begin{table}[h] 
\centering
\begin{tabular}{|c|c|c|}
\hline
Measurement & Dimensional quantity & Non-dimensional value \\
\hline
Speed &  $ 40 \frac {\mu m} {sec}$  &  0.2 \\
Diffusion coefficient &  $ 2 \cdot 10^{-9} \frac {m^2} {sec}$ & 0.01 \\
Turning frequency &  $ 1 sec^{-1}$ &  10  \\
Rate of oxygen consumption & $ 3 \cdot 10^{-11} \frac {\mu M} {(cell)
(sec)}$ &  $4 \cdot 10^{-3} $ \\
\hline
\end{tabular}
\caption[Non-dimensional parameter values.]{Non-dimensional parameter
values.}
\label{nondim_param}
\end{table}

These scales give us the non-dimensionalized values for our parameters,
shown in Table \ref{nondim_param}.  

Values for $\tilde{L}_{min}$ and $\tilde{L}_{max}$ are not readily
measured experimentally, so several different parameter values were
used in the simulations.  The non-dimensionalization gives the
following differential equations:
$$\frac {\partial r} {\partial t} =  \frac{\partial (- r)} {\partial
x} - f_{rl} r + f_{lr} l  $$
$$\frac {\partial l} {\partial t} = \frac {\partial (l )} {\partial x} +
f_{rl} r - f_{lr} l $$
$$ \frac {\partial L} {\partial t} =  \frac {\partial^2 L} {\partial x^2} 
- \kappa (r + l)$$
\begin{equation} \nonumber
        f_{rl} = \left\{ {\begin{array} {lllll}
        C', &  L < \tilde{L}_{min} \\
        c', &  \tilde{L}_{min} < L < L_{min} \\
        c', &  L_{min} < L < L_{max} \\
        C', &  L_{max} < L < \tilde{L}_{max}  \\
        c', &  \tilde{L}_{max} < L
\end {array}}
\right.
\end{equation}
\begin{equation} \nonumber
f_{lr} = \left\{ {\begin{array} {lllll}
        C', &  L < \tilde{L}_{min} \\
        c', &  \tilde{L}_{min} < L < L_{min} \\
        c', &  L_{min} < L < L_{max} \\
        c', &  L_{max} < L < \tilde{L}_{max}  \\
        C', &  \tilde{L}_{max} < L
\end {array}}
\right.
\end{equation}
\begin{eqnarray*} r(0)=l(0)  \\
r(1) = l(1)  \end{eqnarray*}

Here, $\kappa = \frac {k t_0 b_0} {L_0}$ where $k$ is the rate of
oxygen consumption by bacteria, $t_0$ is the time scale, $b_0$ and
$L_0$ are the scales of the original bacterial and oxygen
concentrations, respectively. The non-dimensional values for the
turning rates, $C'$ and $c'$ are given by $C' = C t_0$, and $c'=c
t_0$.

\clearpage
\section{Results} \label{num_sim1}

\subsection{Numerical simulations}

This chapter describes the main results of the model.  The most
important result is the series of computer simulations showing the
band development as it is observed in the actual experiments.
  
Analytical solutions for the steady state are possible, but
numerical results are given to see the time evolution of the
solutions.  In order to have numerical stability, the equations
describing the left-moving bacteria were solved with a
forward-differencing scheme, the equations for the right-moving
bacteria with a backward-differencing scheme.  The diffusion
equation for the oxygen was discretized with a forward-time,
centered space (FTCS) scheme.  The boundary condition at the right
is that all right-moving bacteria become left-moving, and at the
left boundary all the left-moving bacteria become right-moving.  
The domain was discretized by 40 grid points, representing the full
length of the capillary, 2 mm.  The number of grid points was chosen
to be 40 so the developing aerotactic band would have sufficient
resolution.  The size of the time step, $ \Delta t = 0.01$ was
chosen such that the solution to the diffusion equation would be
stable.  The initial conditions in all simulations were that the
bacterial density (both left- and right-moving) is a constant,
scaled to 1.  The initial oxygen concentration is also scaled.

\subsection{Numerical results}

In the following figures, one can follow the development of the
aerotactic band.

\begin{figure}[p]
\centerline{\includegraphics[width=0.6\textwidth]{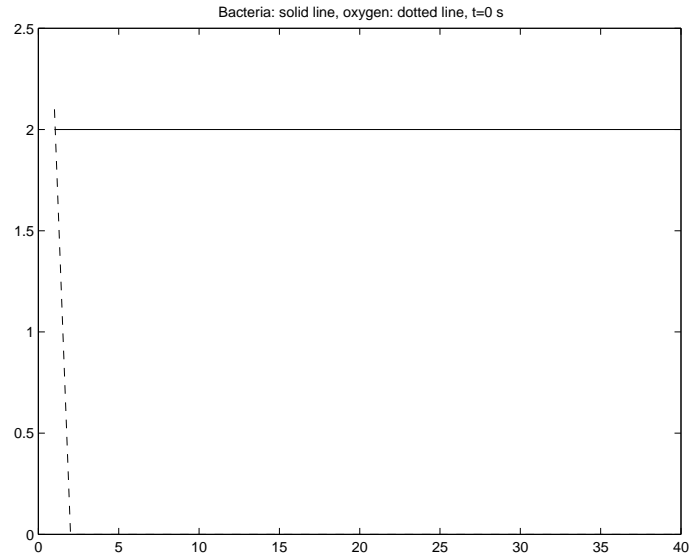}}
\caption[Aerotactic band formation 1.]{Initial condition.  Before the
aerotactic band forms, the bacterial density is uniform. (Bacterial
density is given by solid line.) Oxygen is zero everywhere, except at the 
left boundary.  (Oxygen concentration is given by dashed line.) }
\label{aerotax11.ps}
\end{figure}

\begin{figure}[p]
\centerline{\includegraphics[width=0.6\textwidth]{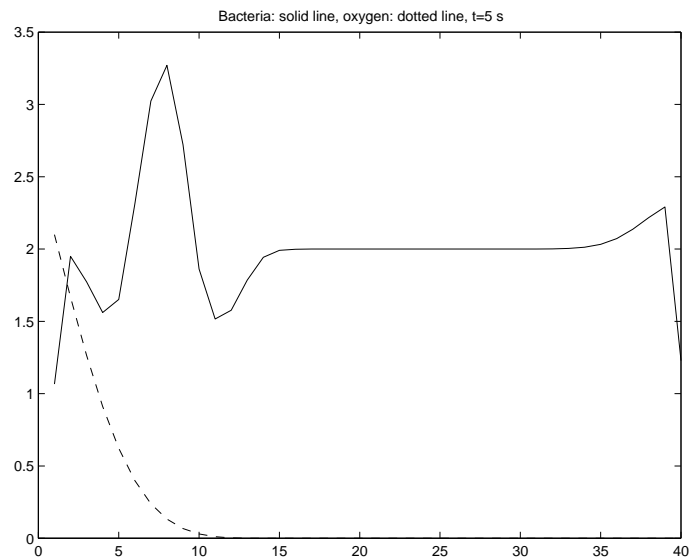}}
\caption[Aerotactic band formation 2.]{Bacteria beginning to aggregate
at the favorable oxygen concentration. Sharp oxygen gradient is
developing. }
\label{aerotax12.ps}
\end{figure}

\begin{figure}[p]
\centerline{\includegraphics[width=0.6\textwidth]{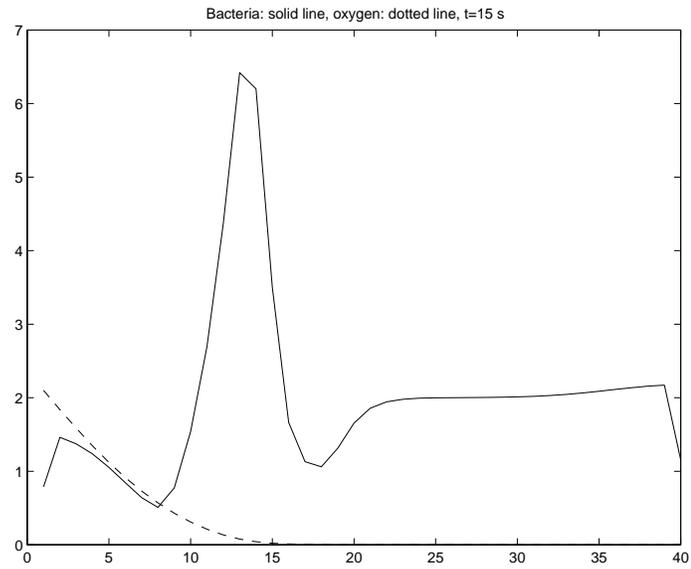}}
\caption[Aerotactic band formation 3.]{After 15 seconds the band is
clearly visible. Bacterial density  in front of the band still
changing.}
\label{aerotax13.ps}
\end{figure}

\begin{figure}[p]
\centerline{\includegraphics[width=0.6\textwidth]{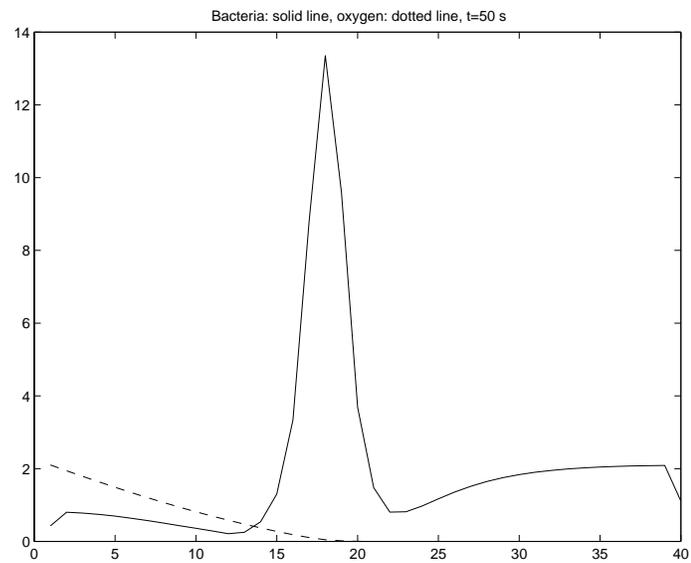}}
\caption[Aerotactic band formation 4.]{Aerotactic band after 50
seconds.}
\label{aerotax14.ps}
\end{figure}

\begin{figure}[h]
\centerline{\includegraphics[width=0.6\textwidth]{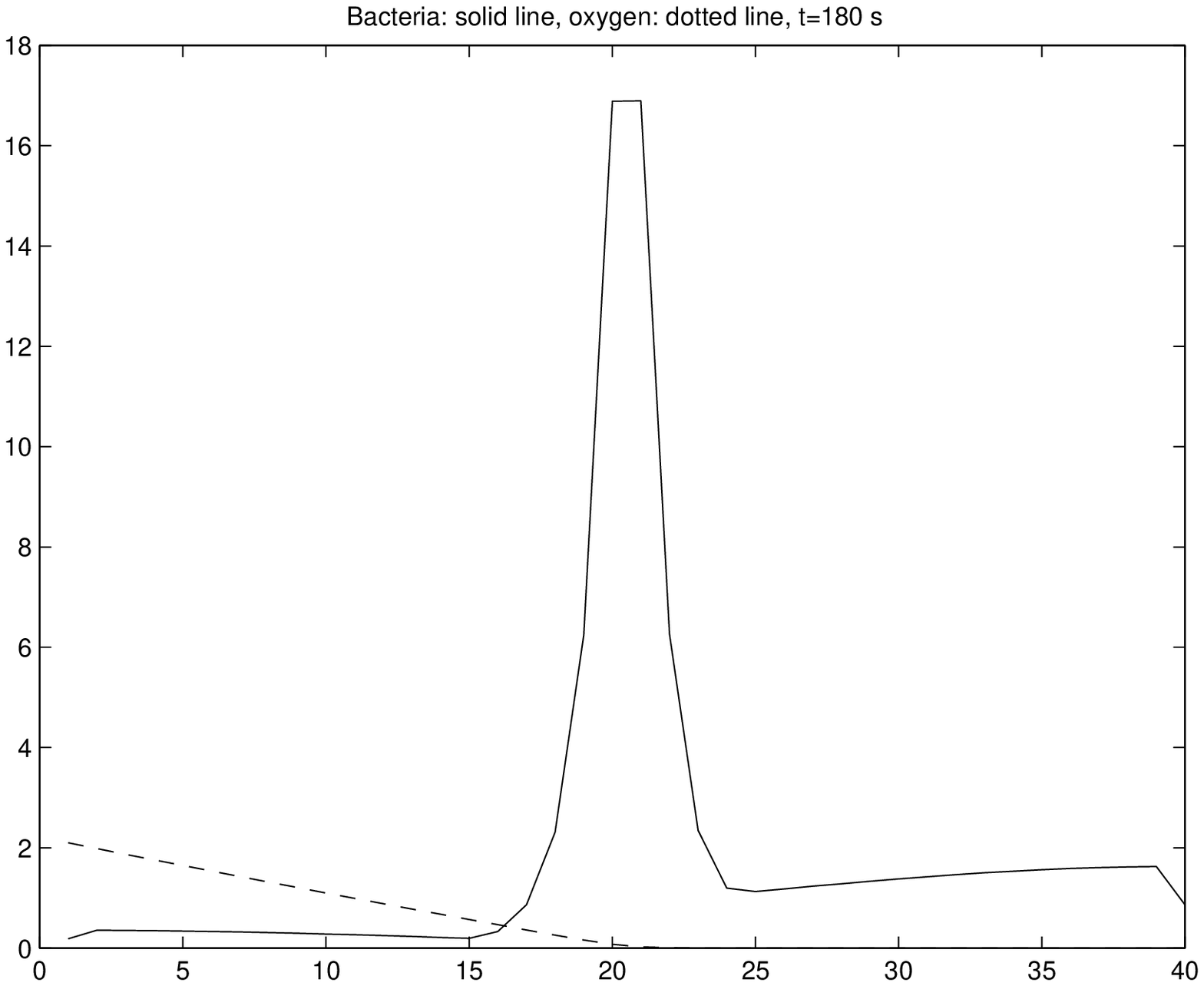}}
\caption[Aerotactic band formation 5.]{Quasi steady state of bacterial
band. Ratio of bacterial densities in front of the band and inside the
band agree with experimental measurements.}
\label{aerotax15.ps} 
\end{figure}

In the figures one can follow the development of the aerotactic
band. Initially, (Figure \ref{aerotax11.ps}) the bacterial density
is uniform everywhere.  The oxygen concentration is a constant at
the first grid point, and zero on all other gridpoints.  The
apparent gradient of oxygen is due to the software package, Matlab
connecting the first and second gridpoints.  After 5 seconds, in
Figure \ref{aerotax12.ps}, the cells at the open end of the
capillary that are exposed to high levels of oxygen concentration
start swimming toward the lower, optimal oxygen concentrations.  
Meanwhile, some of the bacteria at the back that are close enough to
the oxygen are also able to detect the optimal oxygen range, and the
aggregation begins from both sides.  Oxygen also begins to diffuse
through the solution, and a sharp gradient evolves.

In the figure which shows the aerotactic band after 15 seconds
(Figure \ref{aerotax13.ps}), it is clear that most cells from the
open end of the capillary have already aggregated to the band.  
There is low bacterial density on either side of the band, because
all cells at these positions are able to swim into the optimal
range.  Cells at the back of the band never sense the oxygen, and
the bacterial density here remains practically unchanged.  Once the
band reaches its steady state (after 3 minutes), it is clear that
almost all bacteria from the front have aggregated to the band,
(Figure \ref{aerotax15.ps}). The ratio of the cell density inside
the band to the front of the band is more than a 100, and the
density inside the band to behind the band is about 10.  Both of
these values agree with the experimental data.  The time scale of
the band formation, the width of the band and the distance of the
band from the meniscus also give good agreement with the
experimental measurements. The total bacterial density is conserved
in all the simulations.  $\tilde{L}_{min} = 0.2$ and
$\tilde{L}_{max} = 0.7$ were used in the numerical computations.

\begin{figure}[h]
\centerline{\includegraphics[width=0.5\textwidth]{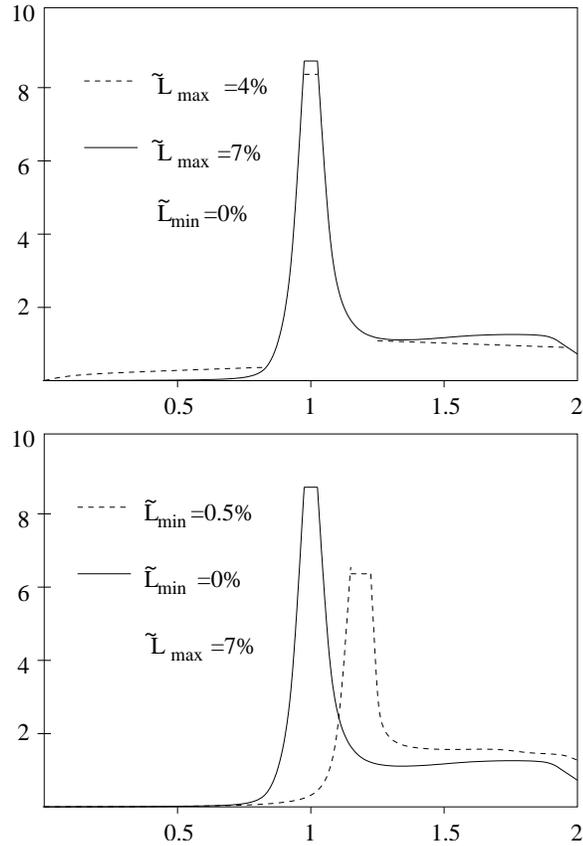}}  
\caption[Numerical experiments with $\tilde{L}_{max}$ and
$\tilde{L}_{min}$]{The top figure shows how the final bacterial band 
changes as $\tilde{L}_{max}$ is changed.  The bottom figure shows
the numerical simulations with $\tilde{L}_{max}$ fixed, and
$\tilde{L}_{min}$ changing.}
\label{l_max_min.eps} 
\end{figure}

The simulations can also help to determine the optimal values of
$\tilde{L}_{max}$ and $\tilde{L}_{min}$.  (Figure
\ref{l_max_min.eps})  This is useful, because these parameters are
difficult to measure experimentally.  By increasing the value of
$\tilde{L}_{max}$, the bacterial density inside the band increases
too, and the the aggregation from the front of the band is more
complete.  This happens because a larger $\tilde{L}_{max}$
effectively increases the region in which cells are in favorable
conditions.  The cells do not have to swim very far from the
meniscus to be able to sense the favorable oxygen concentration,
thus more of them get trapped in the optimal zone.

Lowering the value of $\tilde{L}_{min}$ has a similar effect.  The
bacterial density is higher for low values of $\tilde{L}_{min}$, and
the aggregation from the back of the band is more complete in this
case.  As before, this happens because cells now get inside the band
by detecting lower oxygen concentrations.  For a higher value of
$\tilde{L}_{min}$, sensing the same oxygen concentration would act
as a repellent, where in this case (for the low value) it becomes an
attractant.  This change in $\tilde{L}_{min}$ also results in the
main bacterial density shifting toward lower oxygen concentration.

\subsection{Analytical results}

Here we obtain analytically explicit expressions for asymptotically
stable stationary distributions of bacteria and oxygen. These
distributions are possible to find because the nonlinear equations
(\ref{main1}, \ref{main2}, \ref{main3}, \ref{f1}, \ref{f2}) of the
energy taxis model are piece-wise linear.  In the steady state, the
system of differential equations can be solved on intervals and
replaced with algebraic equations. We solve the equations in two
different cases.

In the first case, we look at the fully developed band, after all
the bacteria have aggregated here.  (An estimate below shows
approximately how much time must pass for this case to be valid.)  
In this case we can look at three regions, and solve the system of
ordinary differential equations in these regions, using continuity
and conservation arguments to match the solutions at the boundaries.  
In the general case, we assume that $ L_{max} < L_0$ for the outside
oxygen concentration, $L_0$, and $L_{min} = 0$. In the next section,
(``Special steady state solutions") we discuss situations in which
the above conditions do not hold.  We provide solutions for $L_{min}
< L_0 < L_{max}$ and $L_{min} < L_0$.

In the second case, we look at what happens if we only wait 3
minutes (in other words, wait just enough time for the bacteria to
develop the band.) In this case, it is possible to show that the
bacteria behind the band have not had sufficient time to feel the
effect of the oxygen gradient; therefore, we can assume that the
bacterial concentration there does not change at all.  This allows
us to solve the equations in two regions only, because the region
behind the band does not need to be considered.  The second case is
discussed in the section ``Quasi steady state solutions".

\subsubsection{\large Steady state solutions}

\begin{figure}[h]
\centerline{\includegraphics[width=0.6\textwidth]{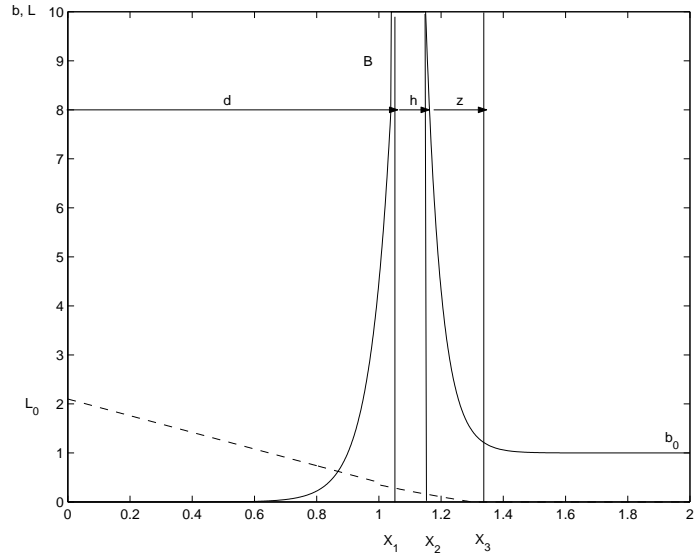}}
\caption[Steady state distribution of oxygen and bacteria.]{Steady
state distribution of the oxygen concentration and the bacterial
density for the case when $L_0 > L_{max} > L_{min}$.} 
\label{sstate2.ps}
\end{figure}

Let us now first consider the steady state solution when $ L_{max} <
L_0$ and $L_{min} = 0$. There are three regions in this case: region
I, in front of the band, from $x_0=0$ to $x_1=d$; region II, inside
the band, from  $x_1=d$ to $x_2=d+h$; and region III, behind the band,
from $x_2=d+h$ to $x_3=d+h+z$.  Figure \ref{sstate2.ps} illustrates 
this scenario.

Let $b(x)$ and $L(x)$ be stationary linear densities of cells and
oxygen, respectively, where the $x$-axis has its origin at the
meniscus and is directed inside the capillary. Let the oxygen
concentration reach values $L_{min}$ and $L_{max}$ at distances
$x_1$ and $x_2$ from the meniscus, respectively. Let $s =
v/(f_{rl}-f_{lr}) \simeq 100 \mu m$, $k = \kappa/D$. We use $s$ as
our length scale.

In region I the bacterial density does not change with respect to
time and is therefore given by:
$$-vr'-f_{rl}r+f_{lr}l = 0, \ \ \ vl'+f_{rl}r-f_{lr}l = 0.$$ 
By noting that the steady state densities of left- and right-moving
bacteria are the same and adding the two equations, we obtain:
$$r = l \sim e^{x/s}.$$
Let $B$ be the density of bacteria inside the band and $d$ be the
distance of
the band from the open end of the capillary.  Then,
$$b = B e^{(x-d)/s}$$
The equation for oxygen in this region is:
$$L'' = kb = kB e^{(x-d)/s}$$ 

Here we can use the fact that the oxygen concentration outside the
capillary is $L_0$, we integrate to find the oxygen concentration in
this region, $L_I$:
$$L_I(x) = L_0 - c_1 x + kBs^2 [e^{(x-d)/s} - e^{-d/s}].$$

In region II we can write the same expression for the bacterial
density as the expression in region I:
$$-vr'-f_{rl}r+f_{lr}l = 0, \ \ \ vl'+f_{rl}r-f_{lr}l = 0.$$
We note that the bacterial density in this region is a constant;
therefore,
we have: 
$$ r = l= B/2; \ \ \  b = B.$$
The equation for oxygen in this region is given by:
$$L'' = kB$$ 
We integrate this expression and use the information that at one edge of
the band the oxygen concentration is $L_{max}$ and the fact that the band
is $d$ distance away from the meniscus to find:  
$$L_{II}(x) = L_{max} - c_2 (x-d) + {1 \over 2}kB(x-d)^2.$$

In region III: 
$$-vr'-f_{rl}r+f_{lr}l = 0, \ \ \ vl'+f_{rl}r-f_{lr}l = 0; \ \ \ r = l
\sim e^{-x/s}.$$
Then, we use the continuity of the bacterial density at the boundary
of the band, $x_2$ to write down the equation for the density in
this region:
$$b = B e^{(x_2-x)/s}.$$
For oxygen we substitute this expression into the bacterial density:
$$L'' = kb = kB e^{(x_2-x)/s}$$
To find the equation for oxygen we integrate, and use the fact that the
oxygen concentration at $x_2=d+h$ is $L_{min}$: 
$$L_{III}(x) = L_{min} - c_3 (x-d-h) + kBs^2 [e^{(d+h-x)/s} - 1].$$

Here we have 7 unknown parameters: $B, c_1, c_2, c_3, d, h, z$. To
find them we use the following boundary conditions: \\

1) Continuity of bacterial density at $x = x_3$ : $B e^{((d+h) -
(d+h+z))/s} = e^{(-z/s)}= b_0.$\\ 
2) No oxygen behind $x = x_3$ (which provides constant bacterial
density there): $L(d+h+z) = 0$.\\ 
3) No flux of oxygen to $x > x_3$: $L'(x_3) = 0$.\\ 
4-5) Continuity of oxygen at the boundaries between the
zones: $L_I(x_1) = L_{max}$ and $ L_{II}(x_2) = L_{min}$.\\ 
6-7) Continuity of flux of oxygen at the boundaries between the zones:
$L_I'(x_1) = L_{II}'(x_1)$ and $ L_{II}'(x_2) = L_{III}'(x_2)$. \\

These conditions lead to the following system of algebraic equations   
using the notation $\lambda = e^{z/s}$: (Note: Throughout the following
calculations, we will use $ k=0.003, b_0=2, L_{min}=0.003,L_{max}=0.005,
s=1$.)

\begin{eqnarray}
B  = b_0 \lambda \label{1} \\
L_{min} - c_3 z + kBs^2 [\lambda^{-1} - 1] = 0 \label{2}\\
-c_3+kBs \lambda^{-1} = 0 \label{3} \\
L_0 - c_1 d + kBs^2 [1 - e^{-d/s}] = L_{max} \label{4} \\
L_{max} - c_2 h + {1 \over 2}kBh^2 = L_{min} \label{5}\\
-c_1+kBs = -c_2 \label{6} \\
-c_2+kBh = -c_3+kBs \label{7}
\end{eqnarray}

From (\ref{3}) and using (\ref{1}): $c_3 = kb_0s$. From (\ref{6},
\ref{7}): $$c_1 = kb_0(s+h\lambda), \ \ \ c_2 =
kb_0(s+h\lambda-s\lambda).$$

We can find analytic approximations for z, d and h.  We start by
approximating z.  Substituting $c_1, c_2, c_3$ as functions of $z,
h$ into (\ref{2}), we find: $$ z/s = (L_{min} +
kb_0s^2(1-e^{z/s}))/(kb_0s^2) $$ Let $\alpha = \frac {L_{min}}
{kb_0s^2} + 1$, and $ \zeta = z/s $.  Then we can rewrite the
expression for $ z/s$ as follows: $$ \alpha - \zeta = e^{-\zeta} $$
Approximating this expression we find that $$ z \simeq \alpha s $$
(The approximation can be obtained as follows: We assume that $
e^{-\zeta} << 1$. Then we can express $\zeta$ in the form $\zeta =
\alpha - \epsilon $, where $\epsilon$ is a small parameter.  This
parameter can be found from:  $ \epsilon = e^ {-\alpha + \epsilon}$.  
In the zeroth approximation $\epsilon \simeq e^{-\alpha}$. Thus, $ z
\simeq \alpha s$.) Numerically, $ z \simeq 1.5 $ and $ e^{-\zeta}
\simeq 0.2231 $ which verifies the original assumption of $
e^{-\zeta} << 1$.

Next, to find h we use (\ref{5}) to obtain the expression:  $$
(h/s)^2 - 2 \frac {\lambda -1} {\lambda} (h/s) - \frac {2 (L_{max} -
L_{min})} {kb_0 s^2 \lambda} =0$$ From above, $\lambda = e^{z/s}
\simeq 4.4817$, so we have $\frac {\lambda -1} {\lambda} \sim 1$.  
Let $ y = (h/s)$, $u=\frac {2 (L_{max} - L_{min})} {kb_0 s^2
\lambda} \simeq 0.1488 $.  Then, from the equation $$ y^2 -2 y - u
=0$$ using the smallness of $u$ in comparison with 1 we arrive at
two roots, $ y_1 \simeq {u \over 2} , y_2 \simeq 2 - {u \over 2}$.  
Only the second root corresponds to a biologically meaningful
situation, since the width of the bacterial band must be of the same
order as the length of the straight runs.  Choosing the first
solution would predict the band width to be much narrower;
therefore, the bacteria would always have to swim across the band
without turning.  If we substitute $y_2$ into the expression for
$h$, then we get, in terms of the original variables, $ h= s(2 -
\frac { (L_{max} - L_{min})} {kb_0 s^2 \lambda}) \simeq 1.8512 $

Finally, we approximate $ d/s $.  From (\ref{4}):  $$ e^{-d/s} +
(d/s) \frac {c_1} {kb_0 \lambda s} = \frac{L_0 - L_{max}} {kb_0
\lambda s^2} + 1 $$ Since the value of d is expected to be large in
comparison with s, we can approximate the above expression by
letting $e^{-d/s} \simeq 0$. We obtain: $$ d/s = \frac { \frac {L_0
- L_{max} } {kb_0s^2} + \lambda } {1+ \lambda h/s} \simeq \frac {
\frac {L_0} {kb_0s^2} + \lambda } {\lambda h/s} = \frac { \frac
{L_0} {kb_0s^2 \lambda} + 1 } {h/s}$$ We have $ \frac {L_0} {kb_0s^2
\lambda} + 1 \simeq 8.4377 >> 1$ Then, using $ h \simeq 2s , d=
4.2188 $.  This also verifies the assumption that $e^{-d/s} \simeq
0$.

In summary, we get the following values:
\begin{itemize}
\item $ h \simeq 2 s \simeq 185 \mu m $ 
\item $ z \simeq (1 + \frac {L_{min}} {kb_0s^2})s \simeq 150 \mu m $ 
\item $ d \simeq \frac {L_0} {2 kb_0 s \lambda} s \simeq 370 - 1860
\mu m$ using the values $L_0 = 0.2 $ and $L_0 = 1.00$, respectively
\end{itemize}

\begin{figure}[h]
\centerline{\includegraphics[width=0.3\textwidth]{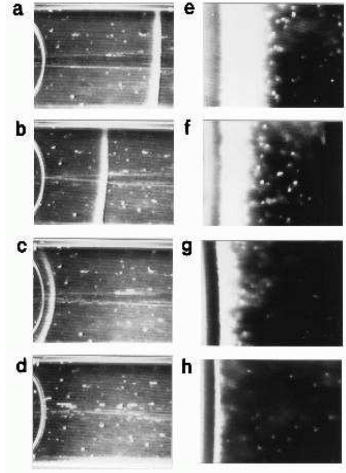}}
\caption[Spatial assays for different oxygen
concentrations.]{Spatial assays for different oxygen concentrations. 
A, 100 percent. B, 21 percent. The width of the band is independent
of $L_0$. The distance between the meniscus and the band increases
with increasing $L_0$. (Figure from Zhulin, \cite{Z})} 
\label{spat_assay.eps}
\end{figure}

Conclusion:\\ (1) Because of the absence of detailed measurements of
bacterial density profile, we do not have data on $z$. \\ (2) The
band's width, $h$, is close to $ 200 \mu m$, exactly as predicted by
our theory if the favorable value for velocity is chosen. \\ (3) The
band was observed to form 400 to 900 $ \mu m$ from the meniscus
corresponding to 20 percent to 100 percent oxygen concentrations,
which is in agreement with our rough estimate of 370 to 1860 $\mu
m$\\ (4) According to the estimate, $d$ is proportional to $L_0$.  
In the experiment, d was observed to grow at greater values of
$L_0$. Thus, we capture the qualitative dependence of d on $L_0$.  
\\ (5) h does not depend on $L_0$. \\ Both conclusions (4) and (5)
are in agreement with Figure \ref{spat_assay.eps} from Zhulins's 
paper, \cite{Z}.

\subsubsection{\large Special steady state solutions}

As discussed above, in this section we look at special cases of the
steady state solution where the assumption that $L_0 > L_{max} >
L_{min} $ no longer holds.  This section extends
the analytical results to two other cases, namely to the case were
$L_{min} < L_0 < L_{max}$ and when $ L_0 < L_{min}$.

\begin{figure}[h]
\centerline{\includegraphics[width=0.6\textwidth]{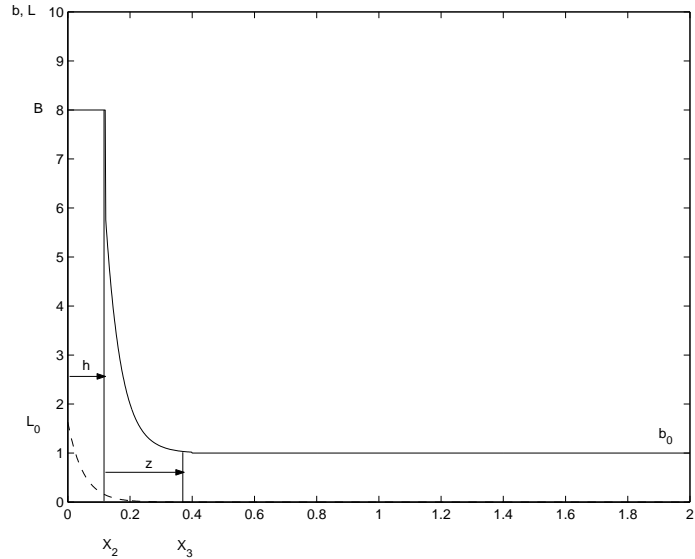}}
\caption[Special steady state distribution of oxygen and
bacteria 1.]{Special steady
state distribution of the oxygen concentration and the bacterial
density for the case when $L_0 < L_{max}$. } 
\label{sstate3.ps}
\end{figure}  

First case: $L_0$ falls between $L_{min}$ and $L_{max}$, as shown in
Figure \ref{sstate3.ps}.  In this case we expect to find two regions
only.  In region I, the band is from $x_1=0$ to $x_2=h$, and region
II, from $x_2=h$ to $x_3=h+z$, is behind the band.

In region I the equation for the bacterial density is a constant; 
therefore,
it is given by:
$$-vr'-f_{rl}r+f_{lr}l = 0, \ \ \ vl'+f_{rl}r-f_{lr}l = 0; \ \ \ r = l =
B/2$$
Here $B$ is density of bacteria inside the band.
The equation for oxygen in this region, using $L(0) = L_0$, and
substituting the above expression for the bacterial density:
$$L'' = kb = kB/2, \ \ \
L_I(x) = L_0 + c_1 x + 1/2 k B x^2$$ 
The other boundary condition,  $L(h) = L_{min}$  helps to determine
$c_1$:
$$ L_I = L_0 + ( \frac {L_{min}- L_0} {h} - {1 \over 2} kBh) x + {1
\over 2} kBx^2$$

In region II we have the following equations.  For the bacterial
density we have exponential decay:
$$-vr'-f_{rl}r+f_{lr}l = 0, \ \ \ vl'+f_{rl}r-f_{lr}l = 0; \ \ \ r = l
\sim e^{-x/s}.$$
Then, using the continuity of the bacterial density at $x_2=h$ we get:
$$b = B e^{(h-x)/s}$$
For oxygen we have the following equations:
$$L'' = kb = kB e^{(h-x)/s}.$$
By integrating the equation twice, we have:
$$L_{II}(x) = kBs^2e^{(h-x)/s}+c_2(x-h)+c_3$$  Using
$ L(h) = L_{min} , L(h+z) = 0$ we get
$$ L_{II}(x) = kBs^2[e^{(h-x)/s} + \frac {kBs^2(1- e^{-z/s)}-L_{min}} {z}
(x-h) + L_{min} - kBs^2$$

To find B,h, z we need the following conditions:\\

(1) Continuity of bacterial density at $x=x_3$.  \\
(2-3) Continuity of flux at the boundaries of the regions.\\
$L'_I(x_2) = L'_{II}(x_2), L'_{II}(x_3)=0$

\begin{eqnarray}
B  = b_0 e^{z/s} \label{8} \\
\frac {L_{min}-L_0} {h} + {1 \over 2} kBh = -kBs + \frac
{kBs^2(1-e^{-z/s)} - L_{min}} {z} \label{9} \\
- k s B e^{-z/s} + \frac {kBs^2(1-e^{-z/s)}-L_{min}} {z} = 0 \label{10}
\end{eqnarray}

Substituting the expression for $B$ into (\ref{10}) we get:

$$ - kb_0s + \frac {kb_0s^2e^{z/s}-kb_0s^2 - L_{min}}{z} = 0$$ 

Simplifying this equation we arrive at: $$ e^{z/s} - {z \over s} - 1
- \frac {L_{min}}{kb_0s^2} = 0 $$ Recalling that $\alpha = 1 + \frac
{L_{min}} {k b_0 s^2}$ and $\zeta = z/s$ we can simplify the above
expression to: $$ e^{\zeta} - \zeta - \alpha = 0$$ An analytical
approximation of $\zeta$ is difficult to achieve, since all the
terms of the equation are expected to be of the same order of
magnitude.  Based on a numerical estimate, $\zeta = 0.85 $. This,
indeed, verifies that all terms of the above equation are of order
1. In terms of the original variables we get $ z \simeq 0.85 s
\simeq 85 \mu m$

We can also substitute the expression for $B$ into (\ref{9}), and we
get:

$$ \frac {L_{min}- L_0} {h} + {1 \over 2} k b_0 e^{z/s} h = - kb_0
e^{z/s} s + \frac {kb_0s^2 e^{z/s} - kb_0s^2 -L_{min}} {z} $$ This
gives us a quadratic equation in $h$ which can be written in a
simplified form by introducing $\beta = kb_0 e^{z/s}$.  $$ h^2+ 2 (s
+ \frac { L_{min}+kb_0} {z \beta} - \frac {s^2 } {z}) h + \frac
{L_{min}-L_0} {\beta} = 0 $$ We know the numerical values of all
terms in the equation, so we can estimate $h$.

Conclusions:\\
(1) As before, experimental data does not exist for $z$.\\
(2) Values of $h$ can be estimated for various values of $L_0$.  (For
example, for $L_0=0.35$ percent, $h \approx 35 \mu m$.)  In our
estimates there is a dependence of $h$ on $L_{min} - L_0$ which could
be tested by running experiments with various species that have
different preferences for $L_{min}$, the minimum tolerable oxygen
concentration.\\  

Second case: $L_0$ falls under $L_{min}$.  We expect that there is no
bacterial band developing in this scenario since the outside oxygen
concentration, $L_0$, is below the minimum preferred concentration.
There is only one region in this case, region I, from $x=0$ to $x=z$. 

In region I the equation for the bacterial density is:
$$-vr'-f_{rl}r+f_{lr}l = 0, \ \ \ vl'+f_{rl}r-f_{lr}l = 0; \ \ \ r = l =
\sim e^{-x/s}$$
$B$ is density of bacteria at the boundary.
$ b = B e^{-x/s}.$
The equation for oxygen is given by:
$$L'' = kb = k B e^{-x/s} \ \ \
L_I(x) = k B s^2 e^{-x/s} + c_1 x + c_2$$  We use that the oxygen
concentration outside the capillary is $L_0$: $L(0) = L_0$.
$$L_I(x) = k B s^2 e^{-x/s} + c_1 x + L_0 - kBs^2$$

To find B, $c_1$ and z we need the following conditions:\\
(1) Continuity of bacterial density at z: $ B e^{-z/s} = b_0$ \\
(2) Continuity of flux of oxygen at $x=z$:
$$ k B s e^{-z/s} + (k B s^2)/z (e^{-z/s} - 1) + L_0/z = 0$$
(3) Boundary condition at $x= z$: $L_I(z) = 0$

By making the appropriate substitutions we arrive at: $$ kb_0sz+
kb_0s^2(1 - e^{z/s}) = L_0$$ Again, recalling $\zeta = z/s$, a
simple form of the expression is $$ e^{\zeta} - \zeta = 1 - \frac
{L_0} {kb_0s^2} $$ Here we must use a numerical estimate as well,
since we again expect the terms to be the same order of magnitude.

Conclusion:\\
(1) A numerical estimate for z is possible, but, as mentioned above,
it cannot be compared to experimental values. \\

This concludes the steady state solutions.  Now we can examine how the
solution behaves on the time scale over which the band develops.   
 
\subsubsection{\large Quasi steady state solutions}

\begin{figure}[h]
\centerline{\includegraphics[width=0.5\textwidth]{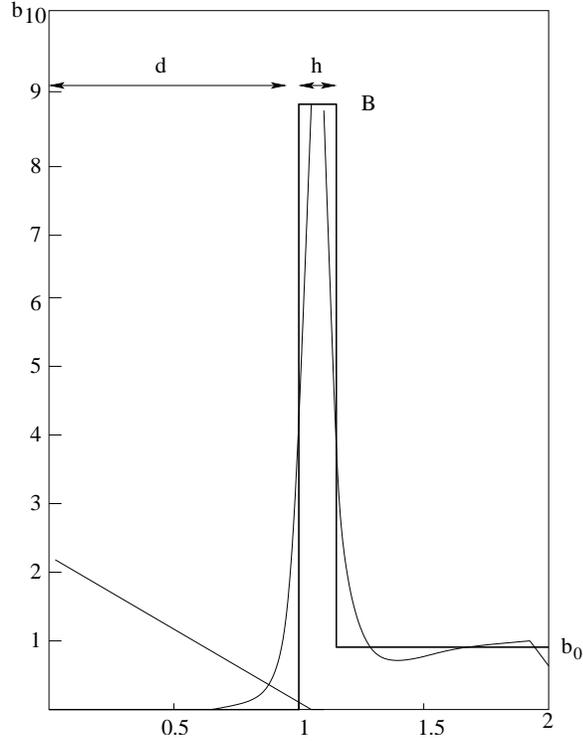}}
\caption[Quasi steady state solution.]{Quasi steady state solution of
the bacterial density.}
\label{quas_ss.eps}
\end{figure}

In the previous section we provided an analytical solution for the
steady state in three different regions: in front of, inside, and
behind the band.  In this section we obtain an analytical solution
for the time interval 50 seconds to 3 minutes, during which the
bacterial band forms, shown on Figure \ref{quas_ss.eps}. We cannot
obtain transient solutions for our system, but we know that the
steady state solutions will be valid only after a very long time
period (hours).  On the other hand, the experimental observations
are for a much shorter time scale, on the order of minutes, so we
would like to find the quasi steady state solutions for this time
period.  We will show below that during this time period only the
bacteria in front of the band have time to move inside the band, so
in fact, we only have two regions: one in front of the band and one
inside the band.  We check our heuristic analytical results against
the experimental findings.

In order to demonstrate that only the bacteria in front of the band
are involved in the band formation, we must estimate the time for the
cells from behind the band to move into the band by pure diffusion in
the absence of oxygen.  The diffusion coefficient is given by: $ D
\sim \frac {v^2} {f} \simeq \frac {(20 \mu m)^2} {0.5 s^{-1} } \simeq
8 \cdot 10^2 \frac {\mu m^2} {s} \simeq 8 \cdot 10^{-4} \frac {m m^2}
{s}$.  Now we can estimate the time of diffusion, $ t \sim \frac {L^2}
{D} \simeq \frac {(1.5 mm)^2} {8 \cdot 10^{-4} m m^2 / s} \simeq 2.8
\cdot 10^3 s \simeq 47 $ min. This confirms that if we only want to
look at the time interval up to 3 minutes, we only need to consider
two regions, region I, from $x_0=0$ to $x_1=d$, which is in front of
the band, and region II, from $x_1=d$ to $x_2=d+h$, inside the band.

Region I:\\

The equation for oxygen in this region is:
$$ D \frac {d^2 L} {d x^2} = 0 $$
with boundary conditions: 
\begin{eqnarray}
L(0) = L_0  \label{b1}\\
L(d) = L_{max} \label{b2} 
\end{eqnarray} 
Using (\ref{b1}) we arrive at the equation for the oxygen in region I:
$$L(x) =L_0 - c_1 x$$

We assume that there are no bacteria in this region because they have
all aggregated to the band.\\

Region II:

The equation for oxygen in this region is:
$$ D \frac {d^2 L} {d x^2} = \chi B $$
With boundary conditions:
\begin{eqnarray} 
L(d) = L_{max} \label{b3} \\
L(d+h) = 0 \label{b4}
\end{eqnarray}
Using (\ref{b3}) the equation for oxygen in this region is
$$ L = L_{max} - c_1 (x-d) + \frac {\chi B} {2 D} (x-d)^2 $$

The number of bacteria in the band is constant, $b = B$.    

We are trying to find four parameters, $ c_1, d, h$ and $B$.\\
\begin{eqnarray}
c_1 = \frac {L_0- L_{max}} {d}  \label{b6} \\
L_{max} - c_1 h + \frac {\chi B} {2 D} h^2 = 0 \label{b7} \\
L'_I(d)=L'_{II}(d) \label{b8} \\
L'_{II}(d+h) = 0 \label{b9} \\
Bh = b_0 (d+ h) \label{b10}  
\end{eqnarray}

We obtained (\ref{b6}) using (\ref{b2}), and (\ref{b7}) from
(\ref{b4}). In (\ref{b8}) the conservation of flux at the left edge of
the band is a redundant condition, $ - c_ 1 = -c_1 + (d-d)$.  The
expression (\ref{b9}) gives $-c_1 + \frac {\chi B} {D} h = 0$.  
Equation (\ref{b10}) follows from the conservation of bacteria.  
Solving (\ref{b6}) and (\ref{b9}) for $c_1$ and setting them equal
gives $ \frac {L_0-L_{max}} {d} = \frac {\chi B} {D} h $. We can solve
this for $d$, $$ d = \frac {D(L_0 - L_{max} )} {\chi B h } $$ From
(\ref{b7}), again using (\ref{b9}) we arrive at $ L_{max} - \frac
{\chi B} {2 D} h^2 = 0 $ which we can use to obtain $$ h = \sqrt{\frac
{2 D L_{max}} {\chi B }} $$ From (\ref{b10}) we obtain $ B = b_0 \frac
{d+h} {h} $.

Now we introduce a new variable, $k = \frac {\chi} {D} $. Using the new
notation we can rewrite as $d = \frac { (L_0 - L_{max})} {k B h} \simeq
\frac {L_0} { k b_0 (d+h) }$.  We impose the condition that $d$, the
distance in front of the band must be much larger than the width of the
band, $d >> h$ which implies that $ d \simeq \frac {L_0} {k b_0 d}$. This
gives $$ d \simeq \sqrt{ \frac {L_0} {k b_0} } \simeq 8 $$ if $L_0 =
0.2$.  If $L_0 = 1$, then $d \simeq 17$. 

The expression $d \simeq \frac {L_0} {k b_0 d}$ also can be obtained
in a more intuitive way.  According to our assumption, the oxygen flux
is equal to the oxygen consumed, or: $ D \frac {L_0} {d} = \chi B h $.  
This implies $ \frac {L_0} {d} = \frac {\chi B} {D} h = k B h = k b_0
\frac {d} {h} h$ (by using the expression for the conservation of
bacterial cells). From this we arrive at $ d \simeq \sqrt{ \frac {L_0}
{k b_0} }$

Our assumption on $d$ allows us to rewrite $ B \simeq b_0 \frac {d}
{h}$. Then we can substitute this expression for $B$ into the
expression for $h$ we had earlier, $ h \simeq \sqrt{\frac {2 D
L_{max}} {\chi B} } \simeq \sqrt{\frac {2 L_{max}} {k b_0 d/h} }$.
Solving this quadratic expression for $h$ we obtain $h \simeq \frac
{2 D L_{max} } {\chi b_0 d } \simeq \frac {2 L_{max}} {k b_0 d } =
\frac {2 L_{max} } {kb_0} \sqrt{\frac {kb_0} {L_0}} $.  For $ L_0 =
0.2$ this gives $ h =0.4$, and for $L_0 = 1$, $ h= 0.2$.

Since we have an expression for both $d$ and $h$ in terms of the
original variables, we can find $c_1$ from (\ref{b6}) and $ B$ from
$ B \simeq b_0 \frac {d} {h} $

Conclusions: \\
(1) The distance of the band from the meniscus, $d$, is given as a
function of the outside oxygen concentration, $L_0$: $$d \simeq \frac
{L_0} {k b_0 d}$$
This formula provides estimates $d \simeq 0.8$ mm for $L_0 = 0.2$ and
$d \simeq 1.7$ mm for $L_0=1$.  These values are in agreement with the
figures from the Zhulin experiment.  (See Figure \ref{sstate2.ps})\\
(2) The band width, $h$ is given as a function of $L_0$ and $L_{max}$.
More experimental data would be necessary to determine whether the
dependence of $h$ on $L_{max}$ is correct.  The two values for the
band width, $h = 0.4$ mm and $h=0.2$ mm for $L_0 = 0.2$ and $L_0= 1$,
respectively, are also in agreement with empirical findings.   

\clearpage
\section{Conclusions} \label{conclusion1}

Our main accomplishment is the development of a mathematical model for
aerotaxis in a steep attractant gradient for fast-adapting bacterial
turning frequencies.  As discussed in Section \ref{chemotaxis_models},
most of the literature is focused on shallow attractant gradients
because this allows approximations and the use of the Keller-Segel-type
chemotaxis equations (\ref{KSeqn}).  As we note, the only possible
approach to steep gradients and slow adaptation is Monte-Carlo
simulations.

The significance of our model for biologists is its prediction of the
appropriate pattern formation by assuming novel receptor mechanics and
fast adaptation, while the same patterns are impossible to obtain by
models of conventional chemotaxis. This indicates that aerotaxis (and
probably other forms of energy taxis) uses a distinct biochemical
mechanism to achieve motility. The separation of biochemical pathways
supports the notion \cite{TZ} that energy taxis is a ``flight response"
which can overrule chemotactic behavior  thus giving organisms with
aerotaxis a selective advantage over purely chemotactic bacteria.

An immediate benefit of our model are its estimates for the maximal and
minimal tolerable oxygen concentrations, $\tilde{L}_{max}$ and
$\tilde{L}_{min}$, respectively.  The analytical solutions also show
dependence of the band width and the distance of the band from the
meniscus on other parameters which makes our predictions experimentally
testable.

A possible limitation of the model is the assumption that adaptation has
to be faster than the length of the run.  A further question to
investigate is the dynamic of the model if the run length and adaptation
occur on the same time scale.  The option that adaptation is longer than
the time scale of the run is ruled out by our Monte-Carlo simulations.  
Determining whether there is an optimal characteristic adaptation time for
given attractant gradients would also be scientifically valuable.


\newpage
\pagestyle{myheadings}
\markright{  \rm \normalsize CHAPTER 3. \hspace{0.5cm}
  MATHEMATICAL MODELS IN BIOLOGY}
\chapter{Growth cone guidance} \label{axon}
\thispagestyle{myheadings}

\section{Introduction}

Sensory input from the environment and responses to such inputs are
characteristic of all animal cells.  Chemotaxis, the movement along
effector gradients toward sources of attractants (e.g. nutrients or certain
guidance molecules) or away from chemical repellents, provides an example
of such a process.  In fact, even bacterial cells are capable of chemotaxis
(see Chapter \ref{aerotaxis}), although the mechanism underlying chemotaxis
of eukaryotic cells and of bacteria are very different.  Bacterial cells,
because of their extremely small size, are unable to detect concentration
differences across their body.  They swim through gradients, and compare
concentration differences in time.  Animal cells, on the other hand, are
thought to genuinely measure concentration differences at different parts
of the cell body.  Yet another difference is that the signal transduction
pathways involved in bacterial chemotaxis of some species (e.g {\it E.  
coli}) are very well known, and accurate mathematical models of bacterial
chemotaxis exist as well \cite{BL}.  Signal transduction pathways involved
in chemotaxis of animal cells have not been so completely uncovered, and it
is possible that there are significant differences between the pathways as
well as in the mechanisms of gradient sensing between neutrophils, growth
cones or slime molds ({\it Dictyostelium discoideum}).

Understanding gradient sensing is an interesting theoretical question
in its own right, even without considering potential biological
applications of successful theoretical models.  Many different
hypotheses can be formulated regarding how concentration differences
are compared in a cell.  Dallon \& Othmer \cite{Da} discuss several
distinct mechanisms, such as comparing receptor occupancy at different
sides of the cell or setting up an internal gradient in response to
external gradients.  Theoretical models help to distinguish between
these, and other possibilities.  Many authors \cite{M,Na,PD} also
believe that the amplification of the external signal is of key
importance in understanding chemotaxis.  It is not well understood how
a 2-3 \% concentration change over its body in a noisy background is
sufficient to unambiguously orient a cell in a gradient.  While
adaptation is known to be an important feature of bacterial
chemotaxis, it is unclear whether adaptation occurs in all animal
cells.

This work focuses on chemotaxis of growth cones, partially because of the
biological significance of the question, and partially because of exciting
experimental findings in growth cone chemotaxis discussed in Section
\ref{bio_backg}.  Growth cones are the highly sensitive and extremely
dynamic tip of axons.  They are finger-like protrusions that are able to
detect ligand gradients, and initiate the axon movement in them.  Growth
cone guidance is most important during development when axons can cover
enormous distances to their final place in the nervous system, and also
during axon regeneration following injuries.  A successful theoretical
description of growth cone chemotaxis could lead to very important medical
applications.  Growth cones are also unique among animal cells, because
they can respond to gradients of the same chemical differently, depending
on the receptors expressed in the cell, and also, depending on the
internal chemical state of the cell.  How the internal chemical state
might influence growth cone guidance is explained in the experiments with
{\it Xenopus} spinal neurons, described below.

The next section provides the necessary biological background, and
explains experimental data relevant to the models formulated later.  The
second part of the Background, Section \ref{math_backg} describes previous
theoretical models of chemotaxis of eukaryotic cells.  Two papers, one by
Meinhardt \cite{M}, and the other one by Levchenko \& Iglesias are in the
focus.  Next, we formulate two mathematical models of chemotaxis, and
provide analysis and numerical results of the models.  The results are
interpreted, and the we discuss the merits and shortcomings of both
models.  New experimental data on adaptation of growth cones \cite{SP1}
and the limitations of the two models suggests some directions for further
work in this area.

\clearpage
\section{Background}

\subsection{Biological background}
\label{bio_backg}

Our project aims to create a mathematical model of early events in growth
cone response to extracellular effectors, because these events have a
clear significance in axon guidance.  During development, axons must find
their appropriate targets in the nervous system, and they accomplish this
with the help of their growth cones.  Growth cones are able to detect
spatial ligand gradients, and turn toward the source of the attractant.
The axon body follows the guidance of growth cones in this direction.
How growth cones are able to 'measure' gradients, even in cases where the
concentration differences are very small, is not well understood, and is
currently an important research topic.

In spite of extensive research, there are many obstacles in
understanding chemotaxis.  Gradient sensing involves an enormous
number of signaling pathways, and some of the pathways are unique to
certain organisms.  Many years of research has been devoted to
studying {\it Dictyostelium} cells, for example, while growth cone
signaling pathways are a relatively new topic of investigation.  The
hope is that the underlying principles are the same in all eukaryotic
cells, however, important exceptions exist.  It is generally accepted
\cite{PD} that chemotactic cells adapt to persistent stimulus, i.e.
certain signaling events occur transiently shortly after the stimulus
is changed.  This appears necessary in order to explain how cells can
orient in a wide range of external stimulus.  However, it is not clear
whether adaptation occurs in growth cones, as {\it in vivo} they
navigate in attractant concentration gradients that may not span many
orders of magnitude.  Although this question is not resolved, experts
in the field believe growth cones adapt too (Mu-ming Poo, personal
communications).

We will focus our attention on growth cone responses in netrin-1
gradients, because netrin-1 is the best characterized molecule which has
proven to exert a tropic effect {\it in vivo}.  Part of the problem of
gradient detection is the determination of the signal transduction events
immediately following ligand binding, and the reconstruction of the signal
transduction pathway.  We focus on the early events, because it is well
established that in animal cell chemotaxis sensing and motility are
independent processes \cite{PD}, and we make the assumption that the
decision to turn toward or away from a substance is made at an early stage
of the pathway.  There is also strong evidence suggesting that the
decision is made by the time a common second messenger, cAMP is produced.  
We do not consider events happening downstream from this point, although
we recognize that many other downstream events must take place for the
actual turning response to be executed.

Song and Poo \cite{SP1,SP2} note that there are two categories of
guidance cues, type I and type II.  Members of the first group are
characterized by the termination of the turning response when
extracellular $Ca^{2+}$ is depleted, and by the co-activation of two
pathways: the PI-3 kinase and PLC-$\gamma$ pathways.  Whether the
turning response is attractive or repulsive, depends on the level of
cAMP activity.  In contrast, type II.  guidance cues are independent
of extracellular calcium and PI-3 kinase.  The level of cAMP activity
does not alter the turning response in guidance molecules of this
group.  The experiments of the Song, Poo, and Tessier-Lavigne labs
\cite{M,SP,SP2, So}, demonstrate that turning response of growth cones
toward a netrin-1 source strongly depends on internal cAMP levels.
Several results \cite{HN,Zh} indicate the dependence of the turning
response on internal calcium dynamics as well, therefore netrin-1 is a
type I guidance molecule.  Whether the turning response to netrin-1 is
attractive or repulsive also depends on the receptor expressed in the
cell.  DCC is a receptor that leads to attractive turning while UNC-5
leads to repulsion.  However, in the same cell the turning response
can be altered by changing the cytosolic cAMP concentration.  For
example, in a neuron expressing DCC receptors attractive turning can
be changed to repulsion by lowering cAMP amounts.  Some of the
experiments investigating the calcium and cAMP dependence of the
response are described in further detail, because we hypothesize that
understanding calcium and cAMP dynamics leads to understanding how
attractive or repulsive turning decision is made in growth cones.  
Additional information is provided about signal transduction
components, in particular, about the interaction of the calcium and
cAMP pathways.

\begin{figure}[h!]
\centerline{\includegraphics[width=0.8\textwidth]{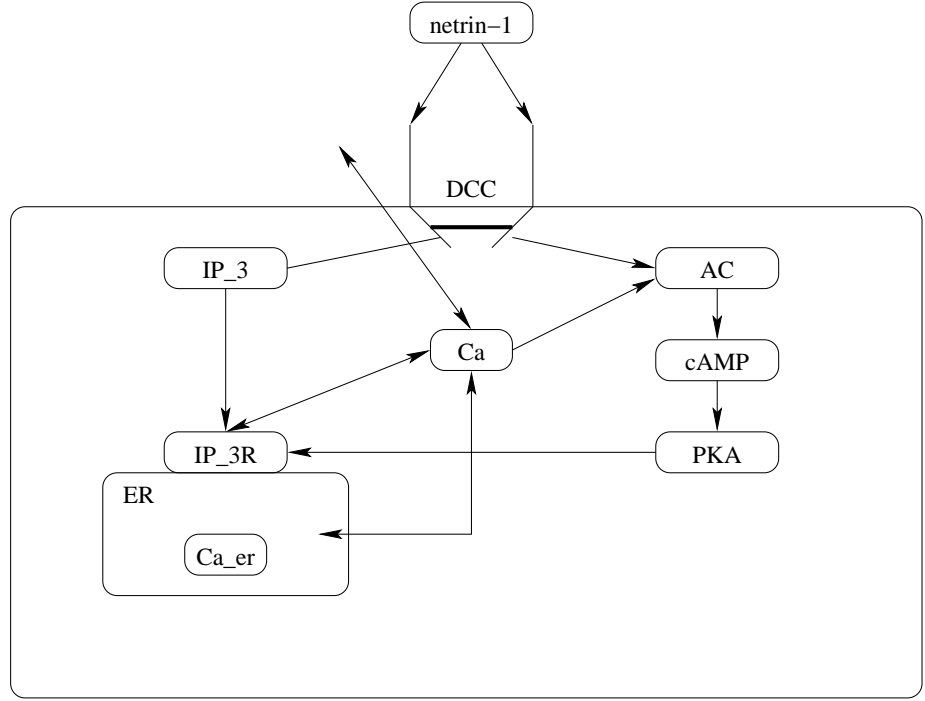}}
\caption[Part of the signal transduction pathway in growth cones.]{Part
of the signal transduction pathway in growth cones.  Abbreviations:
{\bf AC}: adenylate cyclase,{\bf Ca}: cytosolic calcium ion,
$Ca^{2+}$,{\bf $Ca_{er}$}:  calcium ion in the endoplasmic reticulum,
{\bf cAMP}: cyclic AMP,{\bf DCC}: Deleted in Colorectal Cancer,{\bf
ER}: endoplasmic reticulum,{\bf $IP_3$}: inositol 1,4,5
triphosphate,{\bf L}: ligand, netrin-1,{\bf PKA}: protein kinase A. 
References: \cite{Al,BH,Ber,HN,M,MST,SP,So}} 
\label{small3.eps}
\end{figure}

A number of experiments were performed whose results provide constraints
for a model of how ligand binding leads to a turning response. The basic
idea, that high concentrations of cytosolic cAMP correspond to attractive
turning, and low concentrations to repulsive turning, has been widely
accepted, and is reviewed by Song and Poo, \cite{SP2}.  This relationship
is clearly demonstrated by the experiments of Ming et al. on {\it Xenopus}
spinal neurons \cite{M} which show that attraction to a netrin-1 gradient
can be abolished by bath addition of Rp-cAMP, an antagonist of cAMP which
causes lowered cAMP levels in the cell. Figure \ref{small3.eps}. provides
information on the relevant parts of the signal transduction pathway.

Song and Poo also note \cite{SP2} that for each guidance cue, there is a
characteristic range of cAMP which determines whether attractive or
repulsive turning is exhibited.  When cAMP activity falls below the
critical range, repulsive turning is observed and when cAMP activity is
above the critical range, attractive turning is induced.  The critical 
range of cAMP must be low for guidance cues that induce attraction, but 
further reduction of the cAMP converts attractive to repulsive turning.  

The justification to talk about a "critical range" or "threshold" comes
from experiments done on {\it Xenopus} spinal neurons by Ming et al.  
\cite{M}.  They want to determine whether the cAMP concentration changes
turning behavior gradually, or if at a certain concentration attractive
response changes to repulsive response (and vice versa).  They test this
by administering an increasing amount of Rp-cAMP to the bath in addition
to the same netrin-1 gradient.  What they find is that there is a clear
transition from attraction to repulsion consistent with the idea of the
"switch".  On the other hand, enhancing the effectiveness of cAMP by
adding Sp-cAMP to the bath (a chemical which enhances the effect of cAMP)
has no effect on the angle of turning.  The turning angle is statistically
the same during attractive and repulsive turning \cite{M, So}.

Ming et al. also examine the effects of low netrin-1 concentrations, and
they find that for sufficiently low concentrations there is no turning
response. These experiments give us the following important pieces of
information regarding growth cone turning in a netrin-1 gradient: (Note:  
a concentration of 5-10 $\frac {\mu g} {mL}$ at the pipette which results
in about a thousand times the concentration at the cell leads to clear
turning response while 0.5 and 1.5 $\frac {\mu g} {mL}$ at the pipette
leads to no response.)

\begin{itemize} \item There is a certain low concentration of netrin-1
which is insufficient to elicit any sort of response from the growth cone.
\item Under normal circumstances, netrin-1 binding to the receptors induces
turning toward the higher netrin-1 concentration, but lowering cytosolic
cAMP levels can change the attractive response to repulsive turning. \item
Experiments with other molecules (rMAG) which induce repulsion, show that
this turning response can also be reversed and changed into attractive
turning by increasing the cytosolic cAMP levels. \item The turning response
is a 'switch-like' behavior as demonstrated by the fact that the angle of
turning cannot increase by increasing the internal cAMP levels, and cannot
be decreased by decreasing the cAMP levels, rather, decreasing cAMP levels
changes attractive turning into repulsive turning abruptly. \end{itemize}

The abrupt, switch-like change in the turning response suggests that
the internal cAMP levels must reach a value that we can call the
'threshold', above which level the growth cone is able to turn into
the gradient, and below which the gradient is repulsive.  There is no
empirical evidence on how cAMP is distributed inside the cell during
attraction or repulsion. The reversal of behavior is achieved by bath
additions of substances which are known to lower or increase cAMP
levels.  If we assume that cAMP levels are uniform in the cell, and
the attraction or repulsion is only dependent on whether the cAMP
levels are above or below threshold, then quite clearly, we need
another mechanism to provide the spatial information about which
direction the growth cone must turn.  On the other hand, if we assume
that cAMP is elevated to different levels in the cell, and this is
the only mechanism that drives the turning, we arrive at a
contradiction again, because raising the internal cAMP levels by bath
addition of Sp-cAMP should (depending on the amount of Sp-cAMP)
increase the cAMP level everywhere in the cell above the threshold
level which would terminate any type of response.  Based on this
heuristic argument, it is unlikely that cAMP activation alone could
be responsible for the observed behavior.  Let us examine what is
known about the $Ca^{2+}$ pathway.

First of all, many experimentalists note that depleting the external
calcium completely, any sort of turning response is abolished \cite{HN,
SP}.  Hong et al. \cite{HN} examine how different calcium levels regulate
the turning response in {\it Xenopus} neurons, by blocking plasma membrane
channels, channels of the internal storage, and by triggering calcium
efflux from the internal storage by bath addition of ryanodine.  They
observe that high level internal $Ca^{2+}$ signals (resulting from normal
extracellular calcium concentrations) in a netrin-1 gradient lead to
attractive turning, while low level $Ca^{2+}$ signals (resulting from a
depleted extracellular calcium level) lead to repulsive turning.  In these
experiments the authors cannot indicate whether the $Ca^{2+}$ level
elevation is local in all growth cones, although in the larger ones they
are able to observe a transient calcium gradient on the side facing
netrin-1.

A natural question, given the role of calcium as an intracellular
messenger, is the spatial distribution of calcium during growth cone
turning.  Both Hong et al. \cite{HN} and Zheng \cite{Zh} conduct
experiments on {\it Xenopus} neurons to show that spatially restricted
calcium elevation on one side of the growth cone is sufficient to
trigger a turning response.  Zheng is able to quantify the amount of
calcium present by his technique of releasing caged calcium with a
laser beam, and he is able to calculate the amount of caged calcium by
controlling the frequency of the laser.  He shows that the base level
of cytosolic calcium plays an important role, because the same amount
of increase induces attraction in a cell which has 'normal' cytosolic
calcium levels, and repulsive turning in a cell which has lowered
cytosolic calcium levels.  Hong et al. manipulate calcium induced
calcium release (CICR) by creating a ryanodine gradient across the
growth cones. Ryanodine is known to activate the receptors of the
internal calcium storage, the endoplasmic reticulum, ER.  The calcium
release from the stores is autocatalytic.  Hong et al. do not quantify
the amount of released calcium, only note that repulsive turning
corresponds to shallower gradients and a lowered cytosolic calcium
level.

The interaction between the cAMP and calcium pathways is examined by Hong
et al. They induce attractive turning by bath addition of ryanodine,
which leads to increased calcium release, and subsequently override the  
attractive response by blocking a downstream component, PKA, of the cAMP
pathway.  Similarly, they can overturn a calcium induced repulsion by bath
addition of Sp-cAMP.

There is also a steady calcium elevation in some growth cones 10 -15
minutes after the netrin-1 gradient is set up, which only terminates when
the netrin-1 gradient is no longer present.  This was not observed in all
growth cones, and it is possible that this is a result of other processes
(for example actin dynamics) that are unrelated to the initial decision
making regarding the direction of turning \cite{HN}.

We can summarize the relevant facts about the calcium pathway as
follows:

\begin{itemize}
\item Influx from the extracellular medium is important, because depleting
the calcium here abolishes the turning response.
\item Local elevation of calcium mediates a response, but whether it is   
attractive or repulsive, depends on the absolute amount of calcium, not
just the gradient of calcium across the growth cone.  High levels of
calcium lead to attractive turning in the side of the cell where calcium
levels are elevated, but elevation of calcium to a only a low level leads
to repulsive turning from the side with the elevated calcium.
\item cAMP can change the turning response initiated by a calcium      
gradient.  If the growth cone is responding by repulsion, elevating the
cAMP levels will change the response to attraction, and vice versa,
attraction induced by high calcium levels on one side can turn into
repulsion, if cAMP is lowered.
\end{itemize}

Calcium and cAMP are common second messengers and are known to interact
with each other several ways.  cAMP is produced by adenylate cyclase,
and different adenylate cyclase isoforms can be activated or inhibited
by calcium.  Adenylate cyclase in {\it Xenopus} spinal neurons is
calcium-activated (Y. Gorbunova, personal communication).  Calcium
activates three types of AC, namely AC1, AC8 and AC3, but currently it
is not know which of these AC isoforms is found in {\it Xenopus} spinal
neurons.  Each of these isoforms responds to different calcium levels.
AC3 only appears activated by un-physiologically high concentrations of
calcium.  AC8 also responds to calcium in concentrations 5 to 10 times
higher than the concentration necessary to activate AC1, and is
considered a "pure calcium detector" \cite{De}.  AC1, on the
other hand responds to simultaneous activation by G-proteins and
calcium concentrations in the physiological 0.1-1 $\mu M$ range \cite{MD}, 
and is considered "a coincidence detector".  This is
consistent with experimental data on {\it Xenopus} neurons, therefore 
we consider interactions of AC1 and calcium in our model.  Also in
favor of AC1 is the fact that this particular isoform is abundant in
the central nervous system, mainly in the brain, particularly during
development when growth cone guidance is especially relevant \cite{De}. 

Another known effect connecting the two pathways is the phosphorylation of
IP$_3$ receptors by protein kinase A or PKA (downstream product of cAMP)
which leads to changes in calcium flux from the endoplasmic reticulum, ER
into the cytosol.  In some systems it increases calcium flux, in some
others the phosphorylation leads to a decrease.  We hypothesize that in our
case the flux increases, and also, that phosphorylation due to PKA to be
proportional to adenylate cycles.  This concludes all the biological
information necessary to understand the development of our theoretical
model of axon guidance.

\subsection{Theoretical models of gradient sensing}
\label{math_backg}

Theoretical descriptions of how cells can detect sometimes only a few
percent change over their body and orient themselves toward attractants
have a long history.  Previous models (Moghe \& Tranquillo \cite{MT1},
Tranquillo \& Lauffenburger \cite{TL}) focus on a phenomenological
description of chemotactic movement.  Although these models do include one
intracellular messenger, their goal is not the accurate description of the
signal transduction underlying gradient sensing, rather, they attempt a
complete description of leukocyte chemotaxis from sensing to directional
movement.  Although these models made important contributions to the
understanding of chemotactic movement, their approach to model chemotaxis
is significantly different from later models which separate sensing from
motility, and build on detailed information about signal transduction
events.  Since these models, separation of sensing and motility has emerged
as an important principle in understanding chemotaxis \cite{PD}, discussed
below.   Although  we do not describe the Moghe \& Tranquillo and
Tranquillo \& Lauffenburger models in detail, they are similar to models 
discussed in Section \ref{chemotaxis_models}.

Recently, several sophisticated types of theoretical models of the
signal transduction events underlying bacterial and leukocyte
chemotaxis have been developed \cite{BL,Da,LI,Me}. However, chemotaxis
in eukaryotic cells is still not well characterized for a number of
reasons.  One problem is the sheer number of connections and pathways
that exist in these organisms, which makes the design of a tractable
model of all signaling molecules involved in gradient detection and
motility impossible with traditional methods.  Another problem is the
difficulty in deciding which experimental data is widely applicable to
all chemotactic cells, which data may be true for {\it Dictyostelium},
for instance, but not growth cones.  An example of such a case is the
question whether growth cones adapt.  Developing and maintaining
polarity may not be important for growth cone chemotaxis either.  
Difficulty can also arise from the correct interpretation of the
experiments.  For example, Meinhardt \cite{Me} bases his model
partially on observations about the dynamic nature of the membrane
protrusions of {\it Dictyostelium}, however, later these protrusions
were proved to be unnecessary for sensing.

So far no theoretical model of chemotaxis in growth cones has been
developed, in spite of the intriguing recent data on the signal
transduction mechanism presented by Ming et al. \cite{M} and Song \& Poo
\cite{SP2}, among others.

Parent \& Devreotes \cite{PD} established a widely accepted
characterization of chemotaxis.  Their work provides important criteria for
every model of chemotactic sensing and movement must meet.  
\begin{itemize} \item Extreme sensitivity and the ability to detect a
concentration difference of as little as 2 \% between front and back
of the cell in a range of absolute concentrations; \item Polarity: when the
cell is exposed for a period of time to the same gradient the rear becomes
less sensitive;  \item Directional sensing is not essential for movement;  
\item Movement is not necessary for sensing (i.e. it is not like bacterial
chemotaxis)  \item Adaptation: transient response is observed in response
to uniform changes in the attractant concentration, while responses at the
leading edge are persistent; eventhough uniform changes lead to transient
response, immobile cells are still able to sense an unchanging gradient.
\end{itemize} Although Parent and Devreotes do not emphasize this, the
amplification of the signal (which is necessary to explain the enormous
sensitivity)  is also a common goal of theoretical descriptions of
chemotaxis.

We review two fundamentally different theoretical models of growth 
cone sensing, one by Meinhardt \cite{Me}, and the other one by 
Levchenko \& Iglesias \cite{LI}.  We discuss how these models address 
the criteria set by Parent \& Devreotes, and some implications and 
limitations of these models.   Three other models exploring some aspect 
of chemotactic sensing are also mentioned.  

Perhaps one of the most widely quoted and most widely criticized model of
recent years is Meinhardt's chemotaxis model \cite{Me}, which attempts to
give a general framework of chemotactic sensing.  His model is based on two
main observations.  Firstly, that chemotactic cells are extremely sensitive
and are able to detect only a few percent change in the attractant
concentration over the cell body, regardless of the absolute concentration
of attractants.  Secondly, that sensing is a dynamic process involving
quickly changing protrusions of the cell membrane, called pseudopods. Even
when no external stimulus is present, pseudopods of {\it Dictyostelium} can
travel around the cell circumference or, in other cases, happen in
synchrony on opposite sides of the the cell.  Meinhardt's model seeks to
reproduce these characteristic patterns of pseudopod extension. Meinhardt,
echoing Parent \& Devreotes, assumes that sensing and motility are
independent processes. Based on the observations, he sets four criteria for
the model: high sensitivity; sensitivity in a wide range of attractant
concentrations; polarization of the cell adapts to changes in the
orientation of the external signal; intracellular pattern formation
continues even in the absence of an external signal.

Meinhardt's model describes what he calls an abstract "intrinsic pattern
forming system" that is responsible for orienting the cell.  Meinhardt
first addresses the question of sensitivity.  He proposes a Turing-like
mechanism, in which a global inhibitor and a local activator amplify a
small change in the attractant concentration.  The activator enhances its
own production as well as the production of the inhibitor.  Although the
mathematical details are not given, one can assume that this is a standard
reaction-diffusion model, in which the inhibitor diffuses faster than the
activator, therefore the range of inhibitor is larger than that of the
activator.  Such mechanisms produce a stable pattern which does not
respond to later fluctuations in the attractant concentration.  The size
of the response is independent of the initial attractant gradient.  The
model, therefore, explains how sensitivity can be independent of the
absolute ligand concentration.  However, an important shortcoming of this
mechanism is that once the cell orients itself, it is unable to respond to
new stimulus, because of the stable pattern of the internal signaling
system.  Meinhardt explores several ways a stable pattern can be re-set,
and checks the predictions of each method against experimental findings.

First, if the half life of the inhibitor is longer than that of the
activator, oscillations occur.  First, the activator level peaks,
subsequently the inhibitor accumulates and ends the activation, and this
gives one full cycle of oscillation.  This model implies that sensing is
not continuous, rather, it is possible in certain time intervals
corresponding to the phase when the activator levels are low.  However,
there is experimental evidence to the contrary at least in {\it
Dictyostelium}.  Another way to destabilize a pattern is by reducing the
range of the inhibitor to a region smaller than the whole cell surface.  
This allows more protrusions to appear, but often this orients the cell in
the wrong direction.  

Patterns produced by reaction-diffusion systems are also adjustable when
an upper bound is imposed on the production rate of the activator.  If, in
addition, the activator is slowly diffusing, then the activated region can
move around the cell surface, but it also becomes broad, unlike the
appearance of the protrusions observed.  If the activator is almost
non-diffusible, then the protrusions have some very desirable
characteristics.  Namely, a steeper attractant gradient leads to more
preferential protrusions in the cell, and without the external gradient
protrusions appear randomly distributed over the cell surface.  However,
in this case the appearing peaks cannot be shifted.  In order to explain
the random appearance of pseudopods over the cell surface in the absence
of attractants, Meinhardt keeps the assumption that the production of the
activator has an upper bound and that the activator diffuses very slowly.  
This model still does not adjust to new attractant gradients, however.  
To achieve this, Meinhardt assumes the existence of a second inhibitor
that diffuses slowly and that acts on a slow time scale.  The second
inhibitor accumulates over time where the peak of the activator is
located, and it destroys the activator peak.  This process readjusts the
cell, and it allows the formation of new protrusions.

The model is analyzed exclusively numerically.  The equations given by
Meinhardt are ordinary differential equations describing the evolution
of the activator and two inhibitors at the surface of the cell broken
into n sections.  $a$ is the activator, $b$ is the global inhibitor with
fast diffusion and $c$ is the local inhibitor that acts on a slow time
scale.  The constants of the equation are:  $s$ - ligand
concentration; several different functions are used here to reflect
random fluctuations or external asymmetry in different sections of the
cell surface; $b_a$ - basic production rate of the activator; $s_c$ -
Michaelis-Menten constant for the local inhibitor; $s_a$ - saturation
constant of the activator; $ r_a$ - decay rate of the activator; $ r_b$
- decay rate of the global inhibitor;  $b_c$ - production rate of the
local inhibitor;  $r_c$ - decay rate of the local inhibitor.
\begin{eqnarray}
\frac {da_i} {dt} = \frac {s_i(a_i^2/b + b_a)} {(s_c + c_i) (1 + s_a    
a_i^2)} - r_a a_i \nonumber  \\
\frac {db} {dt} = r_b \sum_{i=1}^n a_i /n - r_b b \nonumber  \\  
\frac {dc} {dt} = b_c a_i - r_c c_i
\label{Meinhardteqn}
\end{eqnarray}

These equations contain the implicit assumption that the activator, $a$
and the local inhibitor, $c$ do not diffuse while the global inhibitor,
$b$ diffuses so rapidly that it is given by the same function in each
part of the cell membrane.  The term $$ \frac {s_i(a_i^2/b + b_a)} {(s_c
+ c_i) (1 + s_a a_i^2)} = \frac {a_i^2 + b_a b} {b (1+ s_a a_i^2)} \frac
{s_i} {s_c + c_i} $$ shows that $a_i$, the activator is autocatalytic,
and as its level increases, it saturates.  Both $b$ and $c$ inhibit the
production of $a$.  The activator decays at the rate $r_a$.  $a$
promotes the production of both $b$ and $c$.  The production of $b$
depends on the average level of the activator in the cell, while the
production of the local activator only depends on the local activator
level, $a_i$.  The mathematical treatment of the model in this paper is 
superficial, as the author draws on his extensive experience of 
reaction-diffusion systems.   

Several aspects of the model have been attacked. It is unclear what
biochemical mechanism the proposed pattern forming system \cite {Me}
represents, and which components of the signal transduction mechanism play
the role of the activator and the inhibitors.  Meinhardt does suggest a
molecular interpretation, but he calls it "tentative" as well.  Levchenko
\& Iglesias note that although the existence of a second inhibitor is a
theoretical solution, it makes the model less likely to be biologically
accurate \cite{LI}.  Another problem is that with the Meinhardt model
persistent activation at a given small region of the cell membrane is
impossible.  Instead, activated regions move around the circumference of
the cell.  This is consistent with the dynamic fluctuations of the
pseudopods that are at the basis of the Meinhardt model, however, these
fluctuations have since been shown to be unnecessary both for the
adaptability and persistence of signaling in {\it Dictyostelium} cells
\cite{LI}.  In spite of these and other criticism, the model is widely
known, because it does provide an appealing framework for the processes
that results in chemotactic sensing.  Specifically, the model provides a
mechanism for amplification and the readjustment of the signaling
pathways, so the cell can reorient itself in a changing external stimulus.  
Meinhardt's model is also attractive, because it is minimal, in the sense
that it explains chemotactic sensing with assuming the fewest possible
signaling components.  This makes the model easy to grasp and test against
experimental data.  There have been few theoretical alternatives proposed
which can so consistently and clearly explain the most important features
of chemotactic orientation.

Such a theoretical alternative is offered by Levchenko \& Iglesias
\cite{LI}.  Their interpretation of the experimental data, therefore,
their modeling goals, are different from those of Meinhardt.  They
believe that chemotactic signaling pathways must be able to adapt to
spatially uniform increases in the external attractant concentration and
they must also be able to signal persistently when graded stimulus is
presented to the cell.  The hypothesis that chemotactic cells adapt to
uniform increases in the stimulus directly contradicts Meinhardt's basic
assumption regarding the necessity of dynamic pseudopods for successful
sensing.  This hypothesis is based on experiments with {\it
Dictyostelium} and neutrophils where certain signaling components called
phosphoinositides are shown to adapt perfectly.  There are no
differences in the other theoretical goals.  Levchenko and Iglesias also
want to account for the huge internal amplification of the external
signal, the reorientation of cells when a different stimulus is
presented, sensitivity in a wide range of concentrations, and the
independence of motility and sensing.

\begin{figure}[h!]
\centering
\centerline{\includegraphics[width=0.18\textwidth]{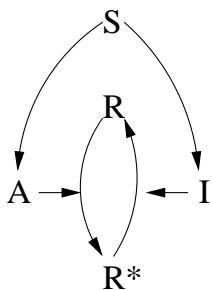}}
\caption[Perfect adaptation according to Levchenko and Iglesias.] 
{Illustration of the mechanism proposed by Levchenko and
Iglesias.  This is the mechanism for perfect adaptation.  The signal, S 
activates both the activator, A and the inhibitor, I.  The output is the 
activated form of the response element, $R^*$.}
\label{lev_igl1.eps}
\end{figure}

Levchenko and Iglesias believe that the most important feature of 
chemotactic sensing is adaptation to uniform changes in the stimulus, 
because this is what allows cells to orient in a wide range of external 
attractant concentrations.  Thus, they begin by building a model for 
adaptation. The pathway they assume is shown in figure 
\ref{lev_igl1.eps}.  Both the activator, $A$ and the inhibitor, $I$ are 
activated by the signal, $S$.  

Based on this signaling scheme, one can write down the following 
equations.  

\begin{eqnarray}
\frac {dR^*} {dt} = -k_{-r} I R^* + k_r A R \nonumber \\
\frac {dA} {dt} = -k_{-a} A + k'_a S(A_{tot} - A) \nonumber \\
\frac {dI} {dt} = -k_{-i} I + k'_i S (I_{tot} - I)
\label{lev_igl1}
\end{eqnarray}

The quantities, $A_{tot}$ and $I_{tot}$ refer to the total amount of 
activator and inhibitor available. Activation of $A$ and $I$ depend on 
the signal, while the inactivation is constitutive.  $R$ could be found 
by the conservation $R_{tot} = R + R^*$.  By assuming that the available 
substrate for $S$ is always much larger than $A$ and $I$, i.e, that 
$A_{tot} \gg A$ and $I_{tot} \gg I$, and by non-dimensionalizing the 
equations, Levchenko and Iglesias arrive at the equations
\begin{eqnarray}
\frac {dr} {dt} = - \beta ir + a(1-r) \nonumber \\
\frac {dA} {dt} = -(a -s) \nonumber \\
\frac {dI} {dt} = - \alpha (i-s)   
\label{lev_igl1'} 
\end{eqnarray}

The two new constants are given as follows.  $$\alpha = k_{-i}/k_{-a}$$
and $$\beta = \frac {(k_{-r}/k_r) (k_{-a}/k_a)} {k_{-i}/k_i}$$  It is
easy to see that the steady state of the active response element,
$r_{ss}$ is given by the ratio of the activator and inhibitor
concentration, and it is independent of the signal: $$ r_{ss} = \frac
{a/i} {a/i + \beta} $$ The authors show numerical simulations for the
adaptation of the response element.

\begin{figure}[h!]
\centering
\centerline{\includegraphics[width=0.28\textwidth]{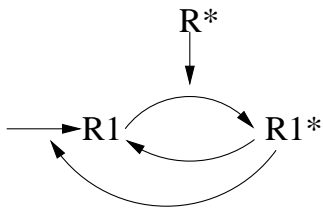}} 
\caption[Amplification according to Levchenko and Iglesias.] {Illustration 
of the mechanism proposed by 
Levchenko and
Iglesias.  This is the mechanism for amplification.  $R^*$ promotes the 
production of $R1^*$, the signaling component which is amplified in this
scheme.  $R1^*$ is autocatalytic, because it controls the substrate for
its own precursor, $R1$.}
\label{lev_igl2.eps}
\end{figure}

Next, the authors consider amplification of the signal.  This is
illustrated in Figure \ref{lev_igl2.eps}.  The goal is to amplify the
production of the output, $R1^*$, They note that there are several ways to
achieve signal amplification: for example, increasing the amount of enzyme
or the substrate.  Increasing the enzyme concentration would correspond to
increasing the concentration of $R^*$, while increasing the substrate
concentration corresponds to increasing the amount of $R1$.  The figure
shows that the authors choose amplification by increasing the substrate.  
By letting $R1^*$ control the production of $R1$, they create what they
call "substrate-supply positive feedback".  The significance of this
mechanism is that amplification only occurs when $R^*$ is turned on.  The
implication of this is better understood when the entire signal
transduction pathway is considered.

\begin{figure}[h!]
\centering
\centerline{\includegraphics[width=0.28\textwidth]{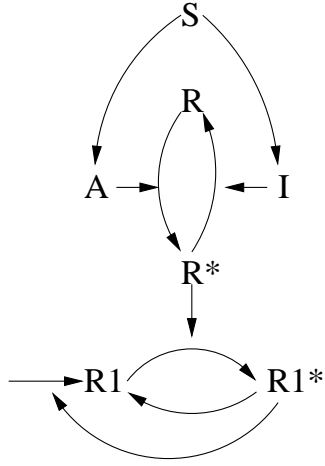}}
\caption[Signal transduction in Levchenko and Iglesias.] {Illustration of 
the mechanism proposed by Levchenko and 
Iglesias.  This is the mechanism for perfect adaptation.   This figure 
shows the whole signaling pathway.  Signaling is possible only when $S$ 
is present.  $R^*$ adapts to persisting stimulus, and $R1^*$ is the 
final, amplified output.}
\label{lev_igl3.eps}
\end{figure}

Adaptation and amplification occur at different levels of this signal
transduction pathway.  Let us examine how the level of $R1^*$ changes
according to this scheme.  As we have seen in the first scheme, $R^*$
can only be activated when a signal is present, and $R^*$ can adapt,
or return to a base level.  Only when there is a signal, $S$, can
there be an amplified response, $R1^*$.  Therefore, the autocatalytic
activation is dependent on the existence external signal, and it
cannot grow unboundedly.  More importantly, there is no need for an
inhibitor to end the signal and force the signaling pathways to be
readjustable, like in Meinhardt's model.

By assuming Michaelis-Menten kinetics for the production of $R1^*$, and
assuming that the production of $R1$ is mediated by another enzyme, $E$,
the authors give the following equations for $R1^*$, $R1$ and $E$.
\begin{eqnarray}
\frac {dR1^*} {dt} = - k_{-2} R1 + \frac {k_2 R1 R^*} {k_M + R1} 
\nonumber \\
\frac {dR1} {dt} = -k_{-1} R1 + (k_1 + k_2 E) R1^* \nonumber \\
\frac {dE} {dt} = -k_{-e} E + k_e (E_{tot}-E) R1^* 
\label{lev_igl2}  
\end{eqnarray}

However, so far the model only describes temporal dynamics, and does not
show how an internal spatial gradient of might develop.  If signaling
takes place entirely locally, then adaptation means that all components of
the signaling pathway return to the same base level everywhere, and no
internal gradient develops.  To answer this question, the authors assume
that all reactions described in equations \ref{lev_igl1} and
\ref{lev_igl2} (with the appropriate initial conditions) happen in n
separate compartments along the cell membrane and, in addition, the
inhibitor is allowed to diffuse.  They add the appropriate terms, $k_D
(I_{j+1} + I_{j_1} - 2 I_j) $ to the equation describing the evolution of
the inhibitor in the jth compartment.  (They also assume that there is not
flux at the two endpoints.)  Now, because the greatest level of inhibitor
is produced where the external signal is the strongest, each compartment
experiences a slightly different ratio $ a/i$, so the steady state of the
activated response element will be different in each compartment.  This
allows the development of a gradient in the internal signaling system.  
Numerical simulations illustrate how the model works, and verify that both
perfect adaptation and an amplified graded response are possible.  The
authors also enter a lengthy discussion of how this theoretical scheme is
mapped onto the actual biochemical pathways of amoebae and neutrophils.  
At the heart of their argument are G-protein activated phosoinositide
kinases and phosphatases.

The Levchenko \& Iglesias model shares a lot in common with Meinhardt's
model.  They both use an activator-inhibitor system where the activator
is assumed non-diffusive and the inhibitor diffuses quickly.  However,
in Levchenko \& Iglesias the production of the inhibitor is linked to
the signal rather than the activator, and in this model the
amplification is a "substrate-supply feedback" mechanism that is
switched on only in the presence of the external signal as well.  
Another significant difference between the models is that in Levchenko
\& Iglesias adaptation and amplification happen at separate levels of
signal transduction. These features solve some problems of the Meinhardt
model.  Levchenko and Iglesias offer a plausible theoretical framework
for the understanding of chemotaxis, and the authors make some
experimentally verifiable predictions regarding the nature of the
inhibitor and activator. 

The two mathematical models discussed so far offer the two most 
comprehensive conceptual approaches to chemotactic sensing.  Some other 
models must also be mentioned because their results illuminate 
particular aspects of chemotaxis.  

Dallon and Othmer \cite{Da} developed a mathematical model which
carefully analyzes chemotactic signals between {\it Dictyostelium
discoideum} (slime mold) cells.  Their novel approach focuses on the
signal characteristics at the boundary of the cells.  The authors use
these characteristics to make predictions regarding which, out of four
potential mechanisms, is the most feasible to act as a signal to
initiate chemotactic orientation.  The four mechanisms to orient cells
in a ligand gradient which had been proposed in the literature are as
follows.  Spatial sensing: the cell measures the concentration
difference or the difference in the number of occupied receptors in
the front and the back of the cell.  Differential force mechanism: the
cell adhesion to the substratum and to other cells depends on the
level of chemoattractant.  Pseudo-spatial mechanism:  the cell extends
pseudopods to convert the spatial gradient in the attractant into a
temporal gradient.  Spatio-temporal sensing: the external attractant
gradient sets up an intracellular gradient.  Adaptation is possible
with this mechanism, by allowing the internal gradient to decay if the
external gradient is unchanging.  \cite{Da} Dallon and Othmer focus on
distinguishing between the spatial, pseudo-spatial and spatio-temporal
mechanisms.

When {\it Dictyostelium} cells are starved, they start secreting cAMP
which acts as a chemotactic signal to initiate aggregation of cells.  
cAMP diffuses, and it is broken down by two chemicals: mPDE inside the
cell and ePDE outside the cell .  Externally, these are the only
reactions included in the model.  The signal transduction inside the cell
is based on a previous model which postulates two pathways: an excitable
one and an inhibitory one.  These two pathways regulate the production of
cAMP by the cell.  The model assumes two cylindrical cells which are
homogeneous in the vertical direction, so the problem can be solved in
the plane.  After non-dimensionalization and the use of the conservation
for mPDE and ePDE, the model consists of the following equations.  
Outside the cells cAMP diffuses and it is degraded by ePDE.  At the
external boundary of the cells, the outward flux of cAMP is equal to the
degradation due to mPDE and secretion.  Inside the cell there is another
reaction-diffusion equation for cAMP accounting for the cAMP diffusion
and its degradation due to mPDE.  At the internal boundary, the inward
flux for cAMP is equal to the basal production, stimulated production
minus the secretion.  The relevant components of the signal transduction
mechanism are membrane-bound, and their evolution is given by three
ordinary differential equations.

A first set of numerical simulations examines the early phase, when
only one cell is signaling and the receiver cell is inactive.  The
difference in the cAMP concentrations at the front and back of the
receiving cell are shown for various distances between the cells,
and for various activity levels of mPDE.  Increasing the activity of
mPDE, which corresponds to an increase in the attractant levels
everywhere, results in a decreased difference between the front and
back cAMP concentrations.  This implies that basing the orientation
on the difference in attractant concentrations at the front and back
of the cell is not advantageous, so spatial sensing is unlikely.  
In order to use the pseudo-spatial mechanism the cell must detect
the time rate of change in the cAMP concentration.  However, the
front-to-back ratio of this is essentially a constant during the
early phase of the signaling, and it only peaks when the rates of
change in cAMP are negative both in the front and the back.  This
suggests that pseudo-spatial sensing would also work poorly.  
However, the front-to-back ratio of cAMP increases with increasing
activity levels of mPDE, and this aspect of the signal is also not
hampered by cAMP levels increasing everywhere around the receiving
cell.  Based on the signal characteristics of the first set of
simulations, spatio-temporal sensing is the most feasible
alternative.

Next, the authors also include the internal signal transduction.  
During the initial time frame when cells must orient, the
qualitative cAMP profiles are the same, although the peaks shift.  
The overall time evolution of cAMP also changes.  A second peak of
cAMP appears that corresponds to the production of cAMP by the
receiving cell.  Similarly to the first set of simulations, spatial
sensing still does not work well, because the difference between the
front and back concentrations is too small.  However, under certain
assumptions regarding the effectiveness of mPDE, the pseudo-spatial
mechanism might be useful.  In general, the front-to-back difference
in the time derivatives is too small to orient cells, but the ratio
of the rates does have a peak at a later time which might be a
signal to initiate cell movement.  The authors conclude that a
pseudo-spatial mechanism cannot be excluded, but it may control
initiation of movement rather than orientation.  As before, the
front-to-back ratio of cAMP concentrations gives a clear signal
again, which further supports the notion of the spatio-temporal
sensing.  Although the internal cAMP gradients are weak, because 
cAMP diffuses quickly, a stable gradient is established which may 
be amplified by some other part of the signal transduction 
mechanism.   

The most significant argument of the article is that purely spatial
mechanism is ineffective for organisms which must orient in a wide range
of concentrations, and this would restrict cells to navigate only in
specific concentration ranges.  However, this may be the case for nerve
growth cones.  It is also important to note that measuring the
front-to-back ratio of the attractant concentrations is an effective
sensing mechanism regardless of the absolute concentration.  Many aspects
of the problem are particular to chemotactic sensing in {\it
Dictyostelium}: the degradation of the external cAMP gradient, and the
fact that cells themselves can change the cAMP profile by secreting cAMP.  
Although the particular mechanism might be very different for growth
cones, it is likely that here also a spatio-temporal sensing based on a
steady internal gradient is at work.

Finally, two articles ought to be mentioned that were published very
recently, in 2001.  They reflect the renewed interest in describing
chemotactic sensing.  Narang et al. \cite{Na} attempt to formulate a
model of chemotactic sensing based on an accurate description of what
they believe to be the most relevant part of the signal transduction
pathway.  The model addresses two questions, namely, the sensitivity of
the chemotactic cells, and why the cellular response is only dependent
on the attractant gradient while it is independent of the absolute
concentration of the attractant.  The biochemical mechanism the authors
focus on is similar to the one examined by Levchenko \& Iglesias.  
Unlike Meinhardt or Levchenko \& Iglesias who create a theoretical
model first, Narang et al. focus on identifying the part of the signal
transduction pathway responsible for a polarized response in an
attractant gradient.  They pinpoint certain membrane phosphoinositides
that respond to uniform attractant changes transiently and spatially
uniformly inside the cell, while they maintain a polarized distribution
in a gradient.

The model has four variables: active receptors, $R_{10}$;  membrane
phosphoinositides, $P$; cytosolic inositides and phosophates, $I$; and
stored phosphoinositides, $P_s$.  All four species are allowed to
diffuse.  The cell is assumed to be two dimensional and disk shaped,
but the diffusion inside the cytosol is assumed to be fast, so radial
gradients are ignored.  This simplifies the model to one spatial
variable, $\theta$, the angle between the leading edge and a given
point of the membrane.  At the heart of the model is a Meinhardt-type
activator-inhibitor system in which membrane phosphoinositides play the
role of activators, and the cytosolic inositides are the inhibitors.  
As in Meinhardt's model, the activator is autocatalytic and it diffuses
slower than the inhibitor.  The receptor dynamics are assumed to follow
perfect adaptation, based on the work of Barkai \& Leibler \cite{BL}.  
This assumption has been disproved by recent data showing that perfect
adaptation of receptors is not characteristic of eukaryotic cells. The
stored phosphoinositides do not play a crucial role in the dynamics of
the model.

Numerical simulations show that the pathway responds to spatially uniform
increases transiently.  This is a result of the assumed perfect adaptation
of the receptors.  Simulations in a graded stimulus demonstrate a stable
and amplified response of $P$, the membrane inositides.  This response is
also consistent with the behavior of chemotactic cells, as it is expected
based on Meinhardt's analysis of activator-inhibitor systems.  Further
numerical experiments also verify that size of the response, i.e the size
of the membrane inositide peak, is independent of the mean external
gradient concentration.  The authors of the article do obtain
amplification of the signal; chemotactic response in a wide range of
external attractant concentration; and adaptation to spatially uniform
stimulus.  Narang et al., however, do not resolve the question of how this
amplified response can be readjusted.  The article proposes that a calcium
surge could result in a destruction of the phosphoinositide peak, and the
authors mention that they have obtained promising preliminary data
supporting this hypothesis.  Although it is impossible to decide based on
the brief description of it in the article, the calcium surge may
correspond to the action of the second inhibitor Meinhardt proposes.

Postma and Van Haastert \cite{PVH} investigate the limitations on the
localization and amplification of intracellular responses by analyzing
the diffusion of second messenger molecules.  During chemotactic
sensing, second messengers must transmit signals from the cell membrane
to the cytoskeleton and various locations inside the cytosol.  The 
speed of the signal transmission depends on the diffusion speed of the 
second messenger, however, fast diffusion leads to the loss of spatial 
information.  This is the first dilemma addressed in the paper.  A 
related question is the amplification of the signal.  Linear signal 
transduction always produces shallower second messenger gradients than 
the original stimulus, therefore a strong local amplification is needed.  
Postma and Van Haastert propose a mechanism that enhances second 
messenger gradients.  

They consider two models for second messenger production.  In the first
scheme the cell is considered to be cylindrical, and second messengers
are produced at one end of the cylinder.  The molecules are allowed to
diffuse and decay.  In the second scheme, a spherical cell is
considered in a linear gradient of the external chemoattractant.  
Production, diffusion and degradation of the second messenger occur at
the cell membrane.  The authors find that the dispersion range of the
second messenger is given by the expression $\lambda =
\sqrt{D_m/k_{-1}}$ where $D_m$ is the diffusion coefficient of the
molecule and $k_{-1}$ is the rate of its degradation.  This expression
implies that fast diffusing second messengers are only able to localize
if their half life is short.  

Next, Postma and Van Haastert propose a model for signal amplification.  
They introduce nonlinearity by assuming that a component of the signal
transduction pathway translocates between the cytosol and the membrane.  
The active receptors stimulate the already membrane-bound effector
molecules which begin the production of second messengers.  Then, as
the second messenger concentration increases locally, more effector
molecules are recruited from the cytosol to the membrane, resulting in
the amplification of the original signal.  The effector translocation
can be considered a positive feedback, or local activation.  By
depleting the cytosol of effector molecules global inhibition is
introduced, therefore the system is another example of an
activator-inhibitor model.  The authors do not model adaptation of the
pathway and the readjustment of the cell.  Numerical simulations
demonstrate that the diffusion-translocation model is able to amplify
an external gradient about tenfold.  However, this amplification is
smaller than experimentally observed values.  The magnitude of
amplification of the model also depends on the gradient in receptor
activity, as stronger external gradients are enhanced more.  In shallow
attractant gradients the model needs to be improved, and in these
situations the authors assume an additional mechanism:  translocation
of another molecule from the cytosol to the membrane which activates
the production of second messengers.  The  model offers two important new
concepts: the analysis of the dispersion characteristics of second
messengers, and that translocation of certain components of the signal
transduction pathway can also act as a positive feedback mechanism.

\clearpage
\section{Mathematical models and results}

Two mathematical models are presented in this section.  Each model has 
focused on a different aspect of chemotactic sensing, and both have 
limitations in explaining growth cone guidance.  The goals and 
shortcomings of each model are discussed in two separate subsections.  
Numerical and analytical results are presented.  

\subsection{cAMP-adenylate cyclase switch}

Based on the experimental findings, the biochemically accurate description
of the cAMP switch seemed like the most important goal because of its key
role in determining the turning response.  The main feature of the model
must be the switch-like, "all or none" response.  As the turning angle
remains the same regardless of the external netrin-1 concentration, we
concluded that the size of the cAMP response should be independent of the
size of the external stimulus.

Experimental evidence \cite{M, SP, So} also suggests that in growth
cones, similarly to other animals cells, we can decouple the
cytoskeletal reorganization and other downstream parts of the signal
transduction pathway from chemotactic sensing.  Turning occurs over
the period of minutes whereas the local elevation of calcium and the
increase in cAMP in response to stimulus happen much quicker, on the
order of seconds.  This allows us to consider the simplified pathway
as shown in Figure \ref{small3.eps}.

We also believed that because growth cones only encounter effector
concentrations over one order of magnitude, they would only need to sense
a gradient within a given range of concentrations.  This implied that
adaptation to a wide range of attractant concentrations would not need to
be considered.

Related to size of the cAMP response is the question of how a possibly
very small spatial gradient of an attractant (or repellent) in the
extracellular medium is amplified into a large internal gradient.  The
amplification must exist, because decisive response can only be expected
if the internal signal is clear.  This implies that even if the receptor
occupancy on different sides of the growth cones is similar, there must be
large enough differences in downstream parts of the signal transduction
pathway for the growth cone to turn in the appropriate direction.  We
assumed that the amplification occurs at the level of cAMP, if it is the
unambiguous biochemical signal for turning. There are several questions
related to gradient amplification, such as the ability to clearly
distinguish noise from signal, and physical limitations on receptor
activation in very low and very high ligand concentrations, but we did not
plan to address these questions.

This model is based on the dynamics of cAMP and $Ca^{2+}$. Namely, high
levels of cAMP correspond to attractive turning for some guidance
molecules, and for these substances, lowering the cAMP concentration means
switching to repulsive turning; a gradient of cytosolic calcium leads to
attraction, while the same internal calcium gradient at a lower overall
calcium level induces repulsion.  We want to give a plausible explanation
of how cAMP concentrations are lowered and increased in a growth cone, and
include realistic cAMP and $Ca^{2+}$ interactions.

To summarize, we wanted the model in which \begin{itemize}
\item sensing is modeled independently of motility
\item the size of cAMP response is independent of the stimulus
\item there is large internal amplification of the external stimulus
\item there is no adaptation
\item cAMP and $Ca^{2+}$ dynamics are realistic. \end{itemize}

A direct approach would be including the spatial and temporal dynamics of
all parts of the signal transduction pathway represented in Figure
\ref{small3.eps}.  This signal transduction pathway could be described by
a system of partial differential equations with the appropriate rate
constants and diffusions coefficients.  These constants, in general, are
difficult to find in the literature, and such a system would have a very
large number of unknowns.  Therefore, the first step in our work is the
reduction of the system.

We want to find the simplest mechanism that can account for a sharp switch
in cAMP levels.  The idea of a smooth external gradient of netrin-1
inducing a sharp, discontinuous response in the internal cAMP
concentration is very closely related to the idea of a smooth gradient
giving rise to thresholds during development, as discussed by Lewis, Slack
and Wolpert \cite{LSW}.

In order to model this mechanism, we made further reductions in the
pathway to be considered, and only focused on the interaction between
cytosolic calcium and adenylate cyclase (AC), the enzyme that produces
cAMP.  There is ample evidence \cite{C,De,MC} of a wide range of such
interactions.  By focusing on these, we hypothesize that these are the
most important nonlinear interactions contributing to the switch.  Such a
simplification is based on all other processes happening on a faster time
scale.  This assertion should be checked again if other signaling
molecules would be added to the model.

As discussed in the ``Biological background" section, we assume that the
adenylate cyclase isoform found in growth cones is calcium-activated, and
it can simultaneously be activated by the receptor, via G-proteins.  We
also consider a positive feedback loop on AC by assuming that the protein
produced by cAMP, called PKA enhances the calcium flux from the cytosolic
calcium stores.  This assumption is based on phosphorylation of the
receptors on the calcium stores by PKA, also mentioned in the ``Biological
background".  We assume that the increase in the calcium flux is
proportional to the concentration of AC.  This implies that there is only
amplification between AC and PKA, and also, that production of cAMP, then
production of PKA happens on a faster time scale then the calcium-AC or
calcium-PKA interactions.  We formulate the model for calcium and
adenylate cyclase, so it is important to comment on how the AC level is
related to the internal cAMP concentration.  A large amount of cAMP is
produced by active adenylate cyclase, and consequently cAMP is degraded.  
No other processes regulate cAMP, therefore, the concentration of cyclic
AMP is proportional to the concentration of active AC.  The amplification
of cAMP due to its production makes the actual switch mechanism is even
more dramatic than the results of our model show.  Our system of
differential equations for calcium, denoted by $C$ and the activated form
of AC, denoted by $A$ is:
\begin{eqnarray}
\frac {dC} {dt} = \overbrace {k_0 \frac {L} {k_{n1}+L}}^{1} -
\overbrace{k1 \frac
{C^2} {K_p^2 + C^2}}^{2} + \overbrace{k_2 (C_b - C)}^{3} \nonumber \\ 
+\overbrace{(k_f + k_3 A) \cdot \frac {L C (C_{er}
- C) } {C  + k_a L}}^{4} \nonumber \\ 
\frac {dA} {dt} = \underbrace{ k_4 \frac {L} {k_{n2} + L} }_{5} \cdot
\underbrace{ \frac {C_m C^4} {K_r^5 + C_m C^4} }_{6} \cdot \underbrace{
(A_t - A)}_{7} - \underbrace{k_5 A}_{8}  \nonumber \\
C(0) = C_b \nonumber \\
A(0) = 0
\label{ca_ac_eqns}
\end{eqnarray}

The equation for calcium dynamics draws heavily on previous models of
calcium dynamics \cite{GDB,KD, MS, Ta}.  The first equation describes the
time evolution of cytosolic calcium concentration, $[Ca^{2+}]$.  As ligand
binds, there is a calcium flux from outside the cell which saturates with
increasing ligand concentrations (1).  The cytosolic calcium is
continuously pumped into the endoplasmic reticulum (ER), as shown in (2).  
The pump is believed to transport two calcium ions per cycle, hence the
second order form \cite{GDB,KD,MS}. There are a number of mechanisms that
maintain the cytosolic calcium level near the resting value, $Ca_b$. These
include passive leak between the cytosol and ER, the cytosol and the
extracellular medium and calcium buffering.  These mechanisms are
summarized in (3).

The last term, (4) describes the calcium flux from the ER into the cytosol
due to the activation of the IP$_3$ channels.  This term is similar to the
analogous term in Tang et al \cite{Ta} who show that flux from the
internal storage, the endoplasmic reticulum (ER) in all models based
receptor-kinetics can be written in the same form.  It is experimentally
shown that calcium has a dual role in the dynamics of the IP$_3$
receptors: the initial fast increase of calcium leads to the opening of
the channels, but consequently calcium also contributes to the slow
closing of the channel.  In our model there is a flux purely due to the
direct activation of the IP$_3$ channels which is proportional to the
calcium concentration difference between the ER and the cytosol,
$C_{er}-C$.  This flux reaches maximal value at a calcium concentration
dependent on the ER calcium level and $k_a$, and is small for both low and
high calcium concentrations.  IP$_3$ is taken to be proportional to the
ligand concentration.  We also include additional calcium flux from the ER
due to the phosphorylation IP$_3$ channels by PKA.  PKA is assumed to be
proportional to the activated form of adenylate cycles.

\begin{table}[h!] 
\centering
\begin{tabular}{|c|c|c|c|}
\hline
Constant  & Value  & Value in lit. & Reference \\
\hline
\hline
$k_0$ & 7 $\frac {1} {s}$ & 6.61  $\frac {1} {s}$ & \cite{MS} \\
\hline
$k_{n1}$ & 1 $\mu M$  &  0.1 $\mu M$ & \cite{MS} \\
\hline
$k_1$ & 5 $\frac {\mu M} {s}$ &  5 $\frac {\mu M} {s}$ & \cite{MS} \\
\hline
$K_p$ & 0.15 $\mu M$  & 0.15 $\mu M$  & \cite{MS} \\
\hline
$k_2$ & 10 $\frac {1} {2}$ &  2 $\frac {1} {2}$ & \cite{JK} \\
\hline
$C_b$ & 0.1 $\mu M$ & 0.08 $\mu M$ & \cite{BI} \\
\hline
$k_f$ & 10 $\frac {1} {\mu M s}$&  &  \\
\hline
$k_3$ & 1 $\frac {1} {\mu M s}$ & &  \\
\hline
$C_{ER}$ & 7 $\mu M$ & 6.3 $\mu M$ & \cite{BI} \\
\hline
$k_a$ & 1 & 1 & \\
\hline
$k_4$ & 2 $\frac {1} {s}$ &  & \\
\hline
$k_{n2}$ & 1 $\mu M$ &  & \\
\hline
$C_m$ & 20 $\mu M$ & 20 $\mu M$ & \cite{BI} \\
\hline
$K_r$ & 1 $\mu M$ & & used \cite{C} \\
\hline
$A_t$ & 20 $\mu M$ & 0.02 $\mu M$ & \cite{BI} \\
\hline
$k_5$ & 1 $\frac {1} {s}$ &  & \\
\hline
\end{tabular}
\caption{Parameter values for the calcium-adenylate cyclase switch model}
\label{Parameters}
\end{table}

The total amount of adenylate cyclase in the cell is fixed on a short time
scale.  The inactivated part of the total given in term (7) becomes 
activated if simultaneously stimulated by netrin-1 and calcium. 
Activation by the ligand is assumed to have first-order kinetics (5).
The calcium activation is mediated by calmodulin, and (8) is fitted to
experimental data \cite{C}.  Its fourth-order form is likely to reflect   
that four calcium ions are necessary to form $Ca_4 CaM$, the
calcium-calmodulin complex which is responsible for the activation of AC.
Activated adenylate cyclase, AC decays linearly with rate $k_5$ (8).

The parameter values are shown in Table \ref{Parameters}.  Wherever it was
possible, we used values that have been established in the literature,
with the exception of the total adenylate cyclase concentration, $A_t$.
The value for $K_r$ was obtained by fitting experimental data in \cite{C}.
We tried several values, ranging over two orders of magnitude, for the
parameters which we could not obtain from the literature.  In all cases,
the qualitative behavior of the system remained the same.  The equations
\ref{ca_ac_eqns} were not non-dimensionalized.  Comments on the parameter 
value for $A_t$ and the reason for leaving the model in its dimensional 
form are given at the end of the section.  

It is important to mention how the model might change, if other isoforms
turn out more important in growth cone guidance.  Because AC3 responds to
very high calcium levels, it is unlikely to play any role, however, it is
possible that AC8 is present.  AC8 can also activated by G-proteins as
well as somewhat higher calcium concentration than what is required for
AC1.  This would change the activation term in the differential equation
for adenylate cyclase, for example to term (5) + term (6)  (instead of the
current activation term which is (5)(6)), but based on numerical
experiments it is still possible to find a parameter range for which we
see a similar behavior as the one in our current model. This implies that
the main mechanism of the cAMP switch does not crucially depend on our
current hypothesis about the AC isoform in growth cones.

\begin{figure}[h!]
\centering
\centerline{\includegraphics[width=0.8\textwidth]{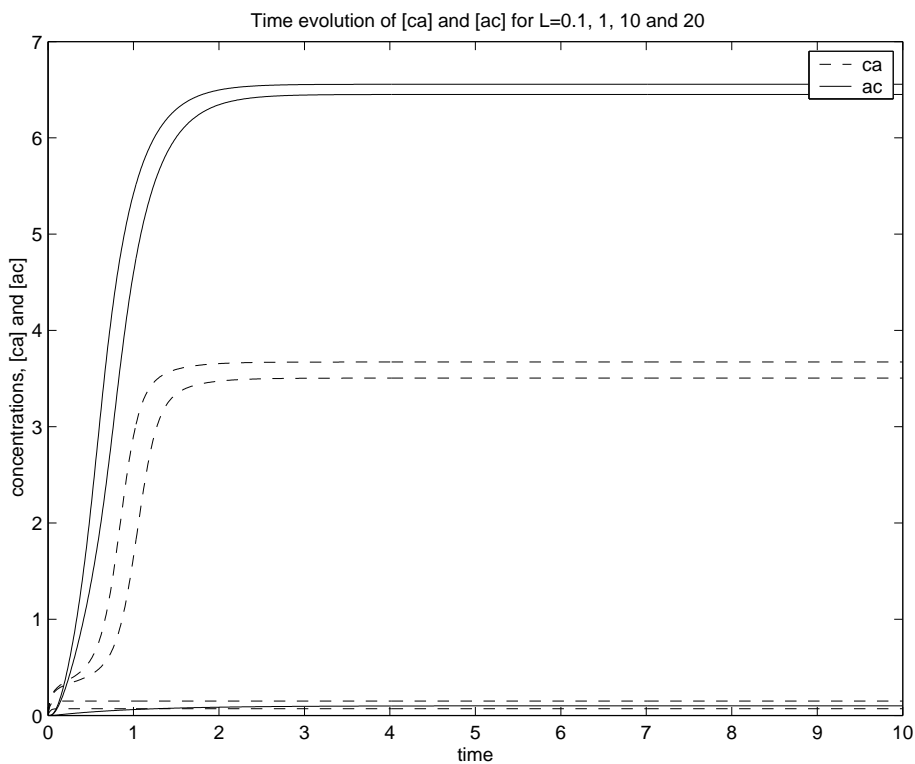}}
\caption[Time evolution of Ca and AC.] {Time evolution of calcium and 
adenylate cyclase for four
different netrin-1 concentrations. Time units: seconds}
\label{ca_ac10.ps}
\end{figure}

We solved the equations numerically using Matlab.  The system is in
dimensional form, and we used the parameter values as they are given in
table \ref{Parameters} without their units.  The numerical simulation is
run for 10 seconds.  Solutions for various ligand concentrations are shown
in Figure \ref{ca_ac10.ps}.  As some parameter values are disputable, the
following discussion on the behavior of the model is limited to the
qualitative behavior.  Numerical values are included in the discussion of
the results only in order to make the discussion easier, but these values
are not claimed to provide realistic information.  Increasing ligand
concentrations corresponds to increases in the cytosolic concentration of
calcium and adenylate cyclase.  It is clear that changing the ligand
concentration from 0.1 $\mu M$ to 1 $\mu M$ does not change the steady
state value of adenylate cyclase and $Ca^{2+}$ significantly, just as
changing the ligand concentration from 10 $\mu M$ to 20 $\mu M$ does not.  
However, there is a significant jump between the steady state values as
the ligand concentration increases from 1 to 10 $\mu M$.  This suggests
that the ligand concentration is the bifurcation parameter in the
differential equations for adenylate cyclase and $Ca^{2+}$.

\begin{figure}[h!]
\includegraphics[width=0.45\textwidth]{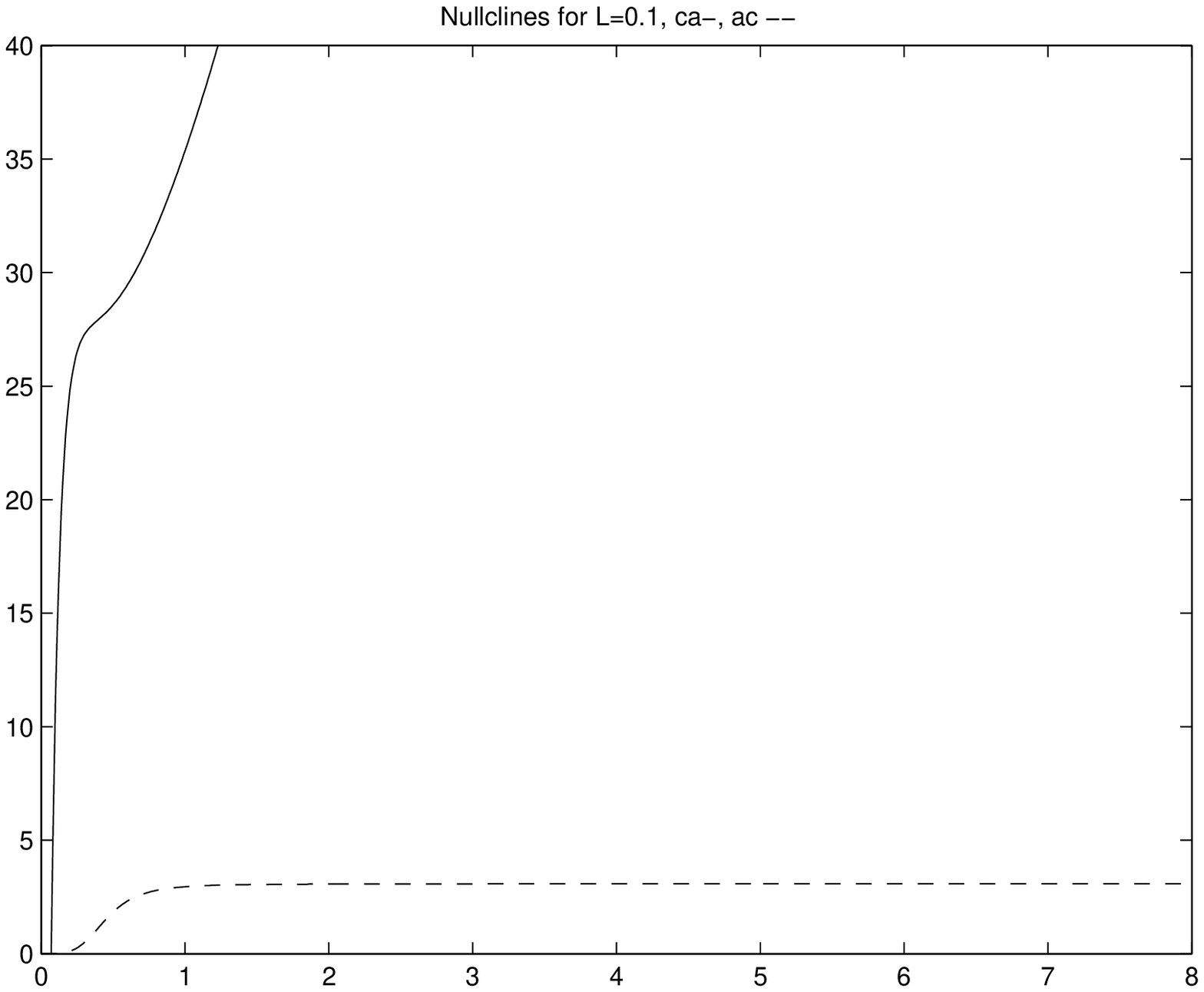} \qquad
\includegraphics[width=0.45\textwidth]{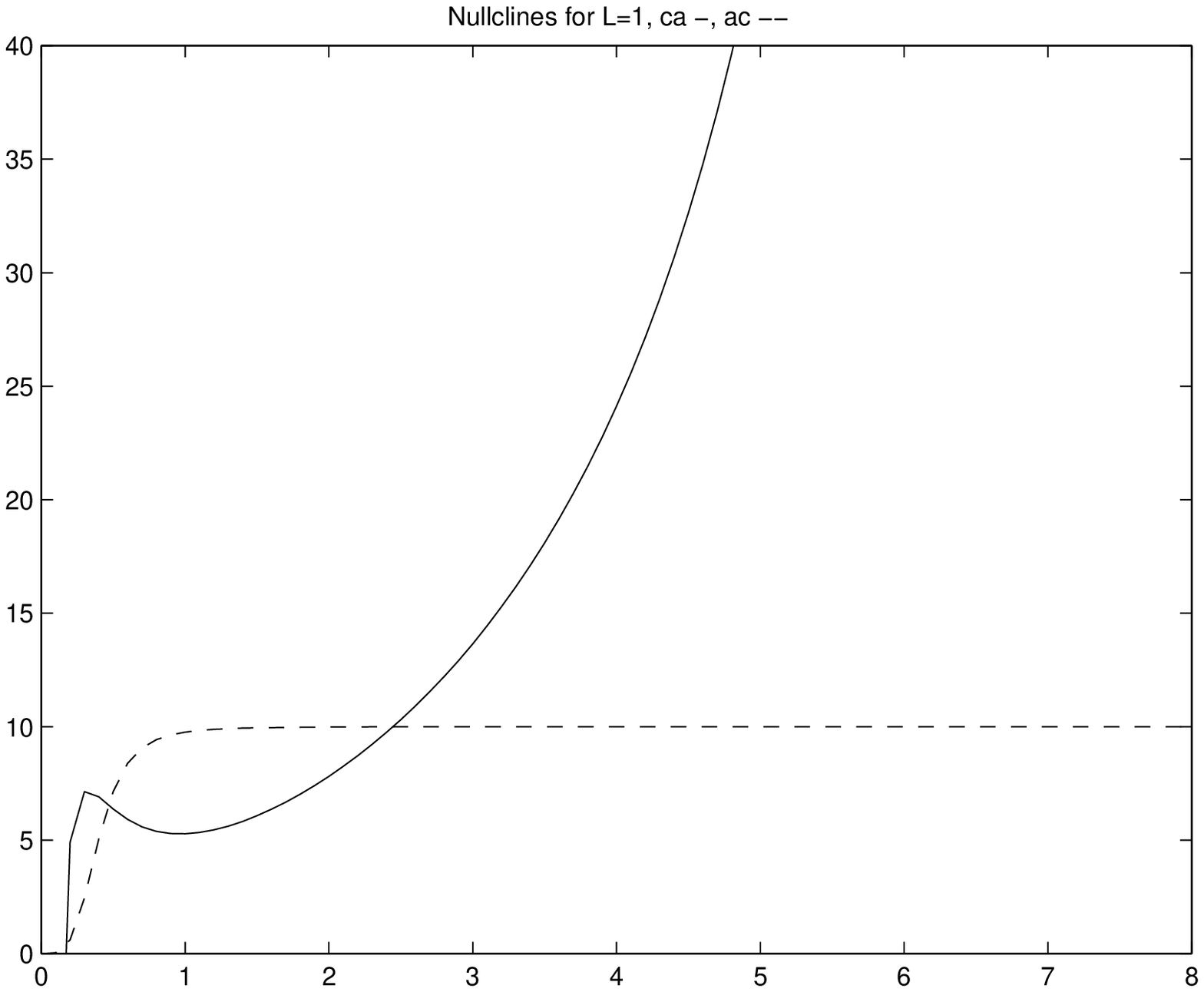} \\
\centerline{\includegraphics[width=0.45\textwidth]{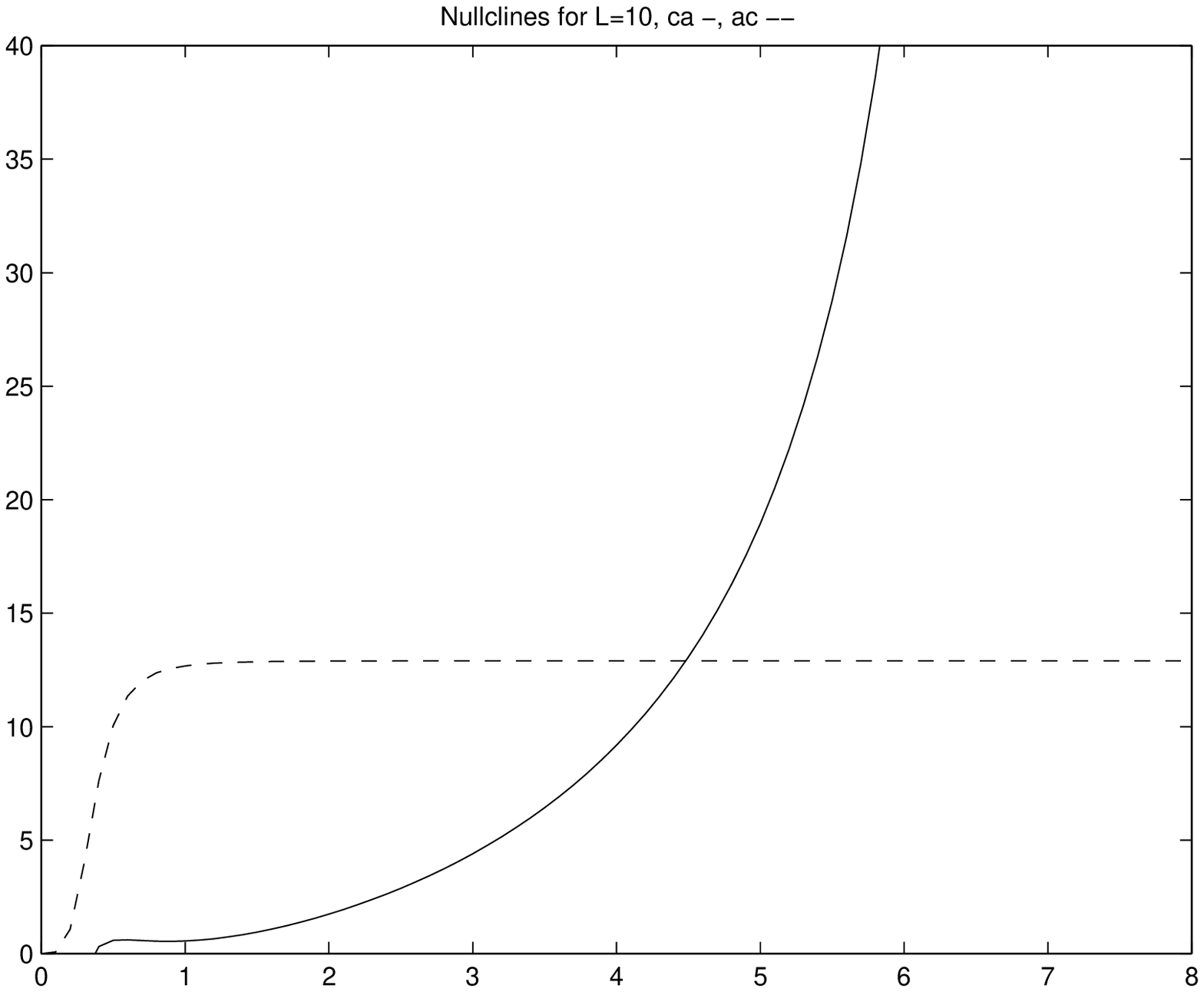}}
\caption[Nullclines.]{Nullclines for values of the netrin-1 concentration, 
L=0.1, 1
and 10 from left to right. The solid curve gives the nullcline for
calcium and the dashed curve is the nullcline for adenylate cyclase.}
\label{Nullclines}
\end{figure}

This is verified when the nullclines of the system are plotted in Figure
\ref{Nullclines}.  For small values, there is only one small stable steady
state.  Regardless of the initial conditions, the cytosolic AC (and
calcium) concentration will approach the same small value.  For
intermediate values, we see two stable steady states, and an unstable
steady state separating them.  If we start with low cytosolic calcium and
low AC levels, the solution converges to the same steady state as before.  
However, by raising the level of L, only the high steady state remains.  
Changing the concentration of the ligand even very slightly changes the
steady state level of AC (and calcium) drastically.  This provides the
mechanism of the cAMP switch.

\begin{figure}[ht!]
\centering
\centerline{\includegraphics[width=0.8\textwidth]{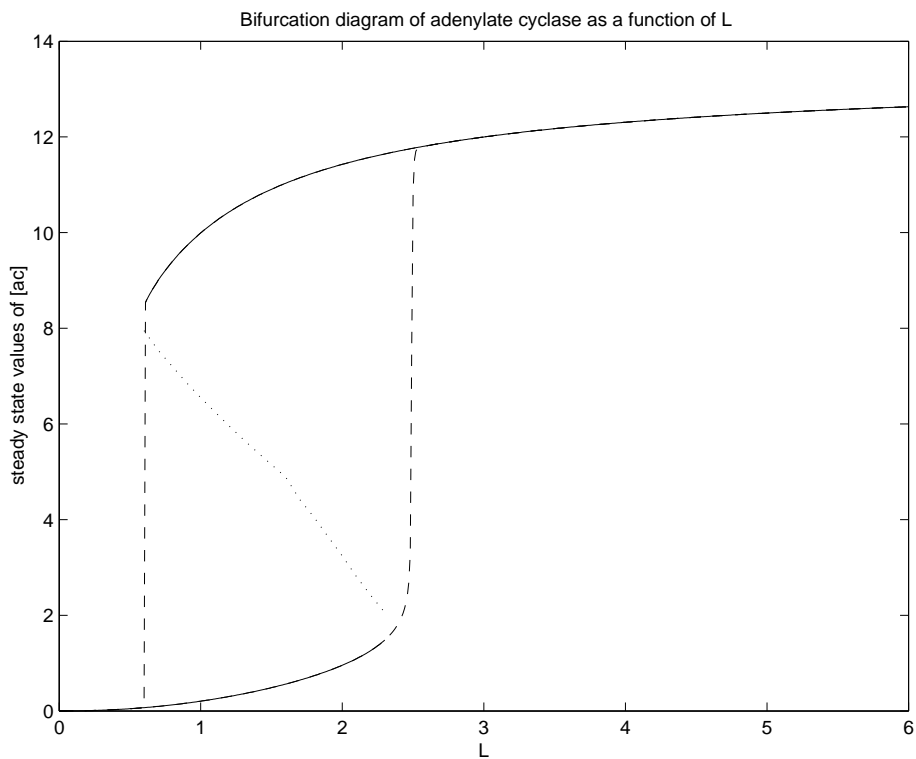}}
\caption[Bifurcation of AC] {Bifurcation diagram of adenylate cyclase as a 
function of L} 
\label{bifurcation6.eps}
\end{figure}

\begin{figure}[h!]
\centering
\centerline{\includegraphics[width=0.8\textwidth]{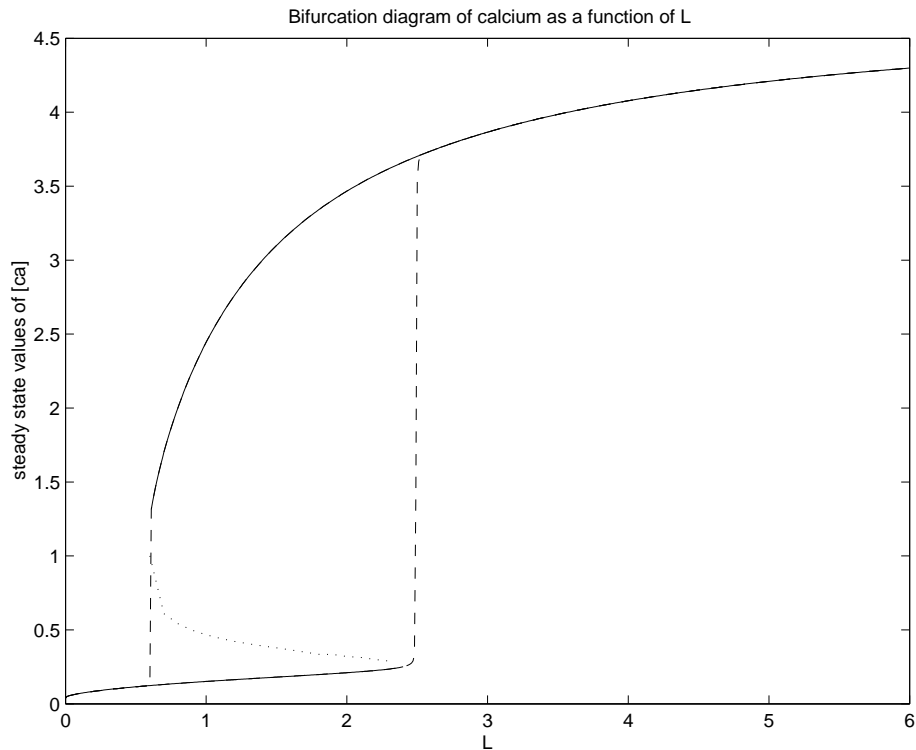}}
\caption[Bifurcation of Ca.] {Bifurcation diagram of calcium as a function 
of L}
\label{bifurcation5.eps}
\end{figure}

The bifurcations are further illustrated by Figure \ref{bifurcation6.eps},
showing the steady state of AC as a function of netrin-1 and Figure
\ref{bifurcation5.eps} with the calcium steady state as a function of
netrin-1.  In Figure \ref{bifurcation6.eps} we see that increasing ligand
concentrations slowly increases the low steady state value of AC.  This
continues until netrin-1 reaches the threshold value of L=2.3, at which
point the lower stable steady state disappears and only the higher steady
state value remains.  This results in a sudden drastic change in the
steady state AC concentration which jumps from the lower steady state
value of approximately 1.7 $\mu M$ to the higher steady state value of
approximately 12 $\mu M$.  This process corresponds to following the solid
line representing the small steady state values for AC, then making the
jump through the dashed line to the upper solid line which represents the
high AC steady state values.  In fact, there is a hysteresis here, because
now decreasing the ligand value below 2.3 (moving to the left along the
upper solid line) will not change the steady state value of AC
dramatically.  The high steady state value decreases gradually until the
ligand concentration reaches about 0.6.  At this point we drop from the
high steady state value to the low one instantaneously.  Similar behavior
is shown for calcium.

\begin{figure}[h!]
\centering
\centerline{\includegraphics[width=0.8\textwidth]{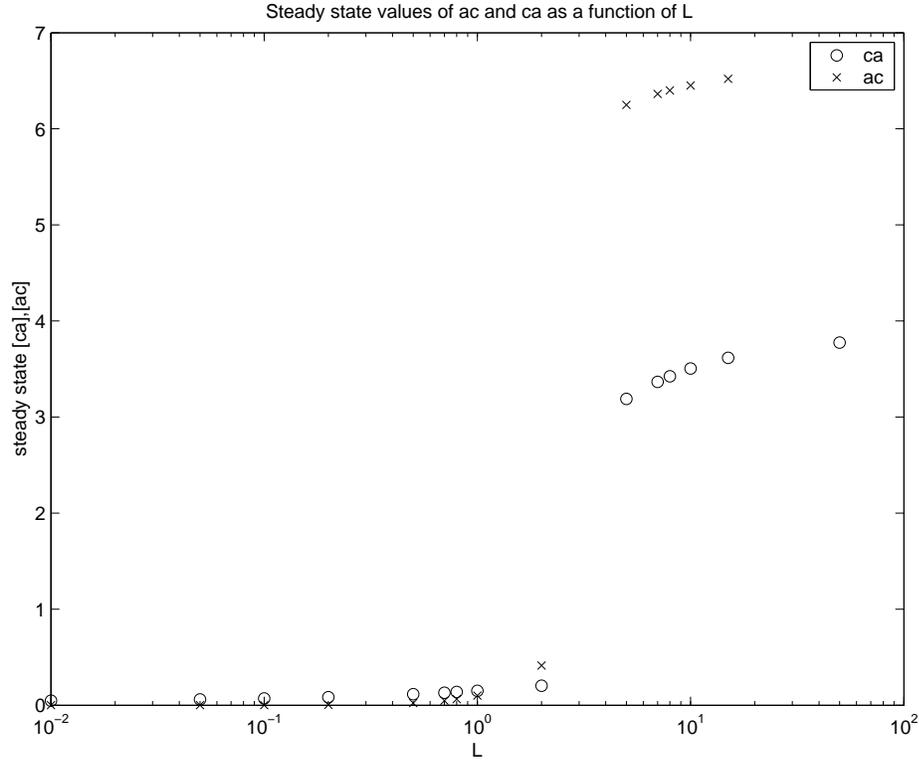}}
\caption[Steady states of Ca and AC.] {Steady state values of calcium and 
adenylate cyclase as a 
function of the netrin-1 concentration.  Initial conditions for all 
simulations are given by Ca(0)=$Ca_b$ ; AC(0)=0.}
\label{ca_acss7.ps}
\end{figure}

In terms of the chemotactic sensing, we have modeled a cAMP switch that
signals unambiguously in certain ligand gradients.  We can consider a cell
divided into two internal compartment, both of which contain the same
signal transduction pathways.  We also note that the concentration of cAMP
is proportional to the concentration of adenylate cyclase.  In some ligand
gradients that contain the threshold value of 2.3 $\mu M$, a small steady
state value of AC is obtained at one side of the cell where the ligand
concentration remains below 2.3 $\mu M$, and there is a sharp jump in the
steady state of the cytosolic AC concentration at the part of the cell
membrane where the ligand concentration exceeds the threshold. If the
ligand concentration obtains the threshold value, the model does satisfy
the goals as follows.  The size of the cAMP response is nearly independent
of the stimulus, as the steady state value does not change much with the
netrin-1 concentration.  This is shown in Figure \ref{ca_acss7.ps} where
ligand concentrations are changed over five orders of magnitude.  The only
significant increase in the steady state values is where the system goes
through a bifurcation.  Because of the bifurcation of the AC
concentration, even a very small change in netrin-1 can produce a very
large change in the concentration of AC.  The model does not assume
adaptation, and it is built on realistic calcium-adenylate cyclase
dynamics.

It is also clear from the above discussion that the model has serious
limitations.  Even if we continue to assume that growth cones do not adapt
(although recent experiments by Poo suggest that they do), the model fails
in ligand gradients that does not contain the threshold value.  Clearly,
this is the most important shortcoming of the model, because this implies
that chemotactic sensing in growth cones is dependent on the absolute
concentration of the ligand, and sensing is only able to occur under
special circumstances.  This is clearly not true.  

There are a few other problems as well.  First of all, as it became clear
that the model did not properly explain chemotactic sensing, efforts to
calibrate the parameter values and to non-dimensionalize the model were
abandoned.  Therefore, the numerical value for the bifurcation parameter,
and the steady state values of $Ca^{2+}$ and AC are meaningless.  We
assume the appropriate ligand gradient necessary to orient the cell is
presented.  In this case, the cell is locked onto this direction, and it
will be unable to re-orient itself, unless the ligand concentrations fall
below the value of 0.6 $\ mu M$ again.  During all the previous discussion
only the temporal dynamics of calcium and adenylate cyclase are
considered, and the spatial behavior was greatly simplified by the
two-compartment model in which the signaling components of the two
compartments do not interact.  Diffusion of the second messengers always
occurs.  Even considering that the fast half-life of cAMP would limit the
range of diffusion, there must be diffusion between the two compartments,
leading to a diminished difference between the cAMP concentration of the
two compartments.  Finally, this model also fails to explain the
experimental observation that lowering the cytosolic calcium level changes
attractive turning to repulsive turning.

\subsection{Adaptation and diffusion model}

The model for the $Ca^{2+}$-AC switch does not account adaptation, i.e.
for transient signaling in uniform changes of the attractant
concentration, and for a cell's ability to choose new orientation in a
changed attractant gradient.  It also predicts that gradient sensing can
only occur in certain ligand concentrations. In this section we provide a
model which corrects these problems.  The model is fundamentally similar
to the Levchenko \& Iglesias model which allowed for both adaptation to
uniform increases and persistent signaling in attractant gradients.  The
basic concept, based on the Levchenko \& Iglesias model is summarized
below.  We include this heuristic explanation in order to give a general
idea of the approach, and the details are made more concrete in this
section.

Adaptation and persistent signaling are achieved by allowing the temporal
adaptation of some response element to any given ligand concentration.  
This is similar to the Barkai-Leibler bacterial chemotaxis model \cite{BL}
where the receptors are known to adapt.  Then, by letting other components
(in Levchenko \& Iglesias, the inhibitor) diffuse, the response element is
forced to adapt to slightly different levels throughout the cell whenever
a gradient of the "inhibitor" are produced.  In these cases the process
leads to an internal spatial gradient of the response element.  Thus,
uniform changes in the attractant concentration lead to spatially uniform,
transient changes inside the cell, whereas an attractant gradient leads to
a spatial gradient of the response element. The conditions under which
such a mechanism can function are discussed in this section.  Our model
does not consider the amplification of the graded response, because we
assume that once an internal spatial gradient exists, its amplification
can be achieved in several ways downstream from sensing, and this is
briefly discussed in \ref{Discussion}.  Although our approach is very
similar to that of Levchenko \& Iglesias, we create a minimal model in
which it is sufficient to consider two components of the signaling
pathway, and we also offer a more detailed mathematical treatment of the
model than the original paper.

In summary, the goals of the model are adaptation to uniform attractant
increments; persistent signaling in gradients; and the ability to reorient
the cell when new stimulus is presented.  A model with such features must
be able to sense gradients in a wide range of attractant concentrations.  
We assume motility and sensing are independent, and focus on the
description of sensing only.

\subsubsection{\large Perfect adaptation} \label{perfect}

\begin{figure}[h!] \centering
\centerline{\includegraphics[width=0.2\textwidth]{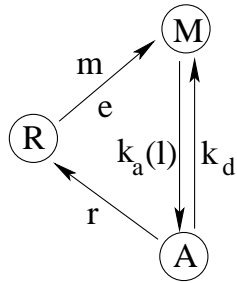}}   
\caption[Signal transduction pathway in one compartment.]
{Hypothetical signal transduction pathway.  R: stored, or "recycled"
substance; A:  "activated substance" which adapts; M: "modified 
substance",part of the signaling pathway that mediates adaptation.}
\label{chem_cycle.eps} \end{figure}

Let us consider the chemical pathway represented in Figure
\ref{chem_cycle.eps}.  We can describe the pathway by the following system
of equations with the appropriate initial conditions that we impose later.  
\begin{eqnarray}
\frac {dM} {dt} = m + \lambda (-k_a (l) M + k_d A) \nonumber \\
\frac {dA} {dt} = -r A + \lambda ( k_a (l)  M - k_d A)
\label{two_odes}
\end{eqnarray}
 
The variables are $A(t)$, the activated substance and $M(t)$, the modified
substance.  (In the original bacterial chemotaxis model, $R$, $A$ and $M$
are the number of receptors in the "refractory", "activated" and "modified"
state, respectively.)  The rate $k_a$ depends on the ligand concentration,
$l$ and $\lambda$ is a large nondimensional ratio of the time scales of the
slow and fast reactions.  In bacterial chemotaxis the fast and slow
reactions correspond to phosphorylation and methylation, respectively.  In
animal cell chemotaxis there is no reason to assume a priori that
similarly, a quick biochemical response is followed by slow adaptation,
therefore our analysis also considers the case when $\lambda \approx O(1)$.  
The equations assume that the production of $M$ does not depend on the
concentration of $R$, the stored substance.  This is true when the
concentration of $R$ is very large compared to concentration of the enzymes
mediating the transition from $R$ to $M$.  Another assumption is that only
the activated substance, $A$ is able to become "recycled", or $R$, and the
modified substance is not.  Furthermore, we assume that downstream effects,
such as turning, depend on the concentration of the active substance.

The steady state of the system is given by $$A_s= \frac {m} {r}$$ and
$$M_s = \frac {m(r+ \lambda k_d)} {k_a(l) \lambda r}\simeq \frac {k_d m}
{k_a(l) r}$$ so the steady state of A, the activated substance does not
depend on the ligand concentration.  $M_s$ assumes that $\lambda \gg 1$.  
The steady state of $A$ implies that regardless of the ligand
concentration, l, the concentration of the activated substance will always
adapt, i.e. return to the same value which is intrinsic to this system, as
it depends on two fixed rates, $m$ and $r$.  In the simplest case $k_a(l)
= k l$.  Let us assume that the ligand concentration $l_0$ jumps to $l_1$
at time $t=0$, and we start from the steady state of the system at the
ligand concentration $l_0$, so $A(0)=A_s = \frac {m} {r}$ and $M(0)=M_s =
\frac {m} {r} \frac {k_d} {k l_0}$.  The analytical solution of the system
of equations is given by (see Appendix \ref{approx_sol} for details):
\begin{eqnarray} 
M(t) \simeq M_2 + (M_0 - M_1) e^{-r_f t} + (M_1 - M_2)e^{-r_s t} \nonumber 
\\
A(t) \simeq A_s + (A_1 - A_s)(e^{-r_s t}- e^{-r_f t}) 
\label{anal_sol_2}
\end{eqnarray}

where we have defined
\begin{eqnarray}
r_f \simeq \lambda (k_d+kl_1)\ \ \ \ r_s \simeq r \frac {kl_1} {k_d+kl_1},
\nonumber \\
A_s = \frac {m}{ r}, \ \ \ \
A_1 = \frac {m} {r} \frac {1+(k_d/kl_0)} {1+(k_d/kl_1)}, \nonumber \\
M_0 = \frac {m} {r} \frac {k_d} {kl_0}, \ \ \ \
M_1 = \frac {m} {r} \frac {k_d} {kl_1} \frac {1+(k_d/kl_0)}{1+(k_d/kl_1)}, 
\ \
\ \
M_2 = \frac {m} {r} \frac {k_d} {kl_1}.
\label{anal_sol_2_param}
\end{eqnarray}

When the ligand concentration increases, then the concentration of $A$
grows from the base line $A_s$ to $A_1 > A_s$ on the fast time scale $T_f
\sim 1/r_f$.  Meanwhile, $M$, concentration of the modified substance
drops from $M_0$ to $M_1$. On the fast time scale, the levels of modified
and activated substances change, but their sum is not altered.  On the
slow time scale $T_s \sim 1/r_s$, the concentration of activated substance
returns to the base line $A_s$ from $A_1$, while the concentration of the
modified substance changes from $M_0$ to $M_2 \neq M_0$.  Now we have a
set of equations (eqn. \ref{two_odes}) that describes the perfect
adaptation of a system to a given ligand concentration.

\begin{figure}[h!] \centering
\centerline{\includegraphics[width=0.8\textwidth]{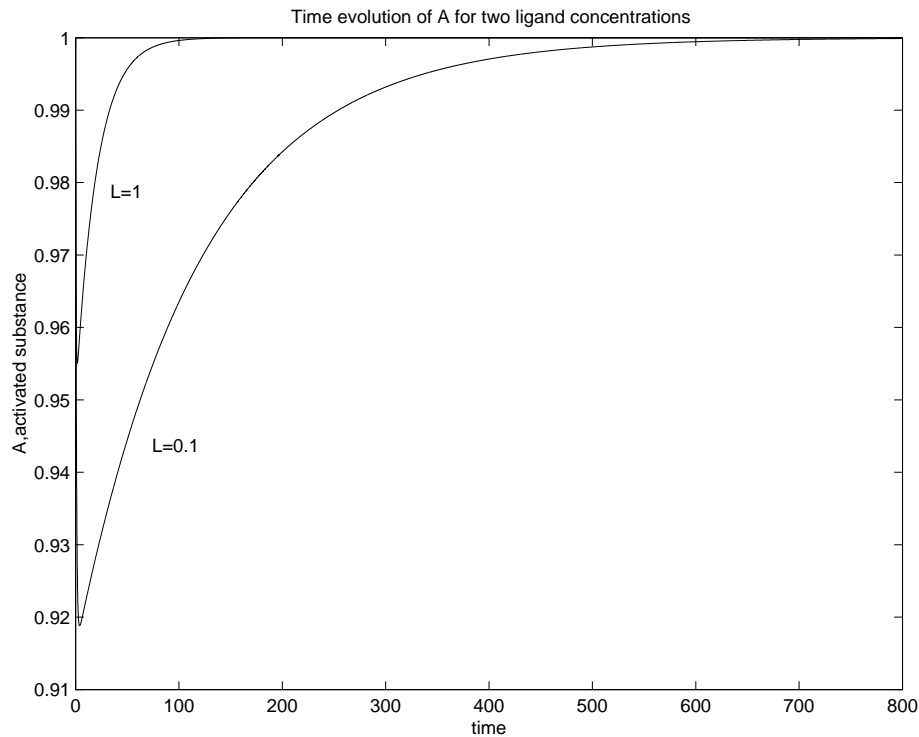}}
\caption[Time evolution of the active substance.]{Time evolution of the 
active substance, A for two different 
ligand concentrations, l=0.1 and l=1.  The concentration of the active 
substance returns to the baseline in each case.} \label{adapt.eps} 
\end{figure}

Figure \ref{adapt.eps} illustrates how the model works.  This system of
differential equations equations and all the following ones, are solved
with the Euler method.  The simulation is run for 800 seconds, and it
shows that after a transient drop in the concentration of the active
substance it returns to the baseline level.  The time step is chosen to be
$h=0.1$.  As it is clear from the analytical solution, the adaptation is
slower when the ligand concentration is small.  The parameter values in
the simulations are $m=0.1$, $\lambda = 5$, $k=0.2$, $k_d = 0.2$ and
$r=1$.  $k_a(l) = k \cdot l$ in these and all the following simulations.  
The same parameter values are used for $m$, $\lambda$, $k$, $k_d$ and $r$
in all of the following simulations, unless otherwise stated.

\subsubsection{\large Spatial models}

We are interested in spatial gradient sensing, so we must look at what
happens when both $A$ and $M$ are functions of time and space.  In order
to gain some insight on the behavior of the system, let us first consider
an axon consisting of two compartments.  We want to investigate how a
spatial gradient of $A$ can develop in this model.  In each compartment
the same set of reactions happens, and in addition, there is a flux
between the compartments.  We consider the case when the ligand
concentration jumps from $l_0$ to $l_1$ in the first compartment and to $
l_2$ in the second compartment. The equations in the two compartment model
are:
\begin{eqnarray*}
\frac {dM_1} {dt} = m + \lambda (-k_a(l_1) M_1 + k_d A_1) + k_1(M_2-M_1)  
\\
\frac {dA_1} {dt} = -r A_1 + \lambda ( k_a(l_1) M_1 - k_d A_1) + 
k_2(A_2-A_1)  \\
\frac {dM_2} {dt} = m + \lambda (-k_a(l_2) M_2 + k_d A_2) - k_1(M_2-M_1) 
\\ \frac {dA_2} {dt} = -r A_2 + \lambda ( k_a(l_2) M_2 - k_d A_2) - k_2 
(A_2-A_1)  \\
M_1(0) = M_2(0) = \frac {m} {r}  \frac {k_d} {k(l_0)} \ \
A_1(0) = A_2(0) = \frac {m} {r} 
\end{eqnarray*}

The steady state solutions are given by  
\begin{eqnarray*} 
A_{1s} = \frac {m} {r} \cdot \Big[ 1 + \frac {r_1 k_1 k} {\lambda 
r_2 k_p + k_1 k_s(r_2 + \lambda k_d)} \Big] \\
A_{2s} = \frac {m} {r} \cdot \Big[ 1 + \frac {- r_1 k_1 k } {\lambda r_2
k_p + k_1 k_s(r_2 + \lambda k_d)} \Big] \\
M_{1s} = \frac {m r_1} {\lambda r} \cdot \Big[\frac {r_2 ( \lambda
k_a(l_2) + 2 k_1) + 2 \lambda k_d k_1} {r_2[\lambda k_p + k_1 k_s] +
\lambda k_d k_1 k_s} \Big] \\
M_{2s} = \frac {m r_1} {\lambda r} \cdot \Big[\frac {r_2( \lambda
k_a(l_1) + 2 k_1) + 2 \lambda k_d k_1} {r_2 [\lambda k_p + k_1 k_s] +  
\lambda k_d k_1 k_s} \Big]
\end{eqnarray*}

where we have defined 
$$r_1= r+ \lambda k_d, \ \ r_2 = r + 2 k_2 $$
$$ k = k_a(l_1)-k_a(l_2), \ \ k_s = k_a(l_1) + k_a(l_2)$$
$$k_p = k_a(l_1)k_a(l_2)$$

The calculations are shown in Appendix \ref{ss_calc}.  If $l_1 = l_2 =
\bar{l}$, or the ligand concentration on the two sides of the growth cone
are the same, then $$A_{1s}=A_{2s} = \frac {m} {r}.$$ $$M_{1s} = M_{2s} =
\frac {m} {r} \frac {r + \lambda k_d} {\lambda k_a(\bar{l})}$$ and if
we let $$\lim_{\lambda\to\infty} M_{1s} =\lim_{\lambda\to\infty} M_{2s} =
\frac {m} {r} \frac {k_d} {k_a(\bar{l})}$$ so we have recovered the steady
state solution to equations \ref{two_odes}. This means that a spatially
uniform increment in the ligand concentration still results in adaptation
in the system, and no spatial gradient of the activated substance
develops.  These calculations are also included in Appendix
\ref{ss_calc}.

Now let us investigate how the flux of the modified substance, $k_1$ and
the flux of the activated substance , $k_2$ will influence the steady
state of the system.  It is important to note that the results of the
following calculations and the qualitative behavior of the system remains
the same under the assumption that $\lambda \gg 1$.  See Appendix
\ref{ss_calc} for details.

First, assume that $M$ is non-diffusible, so $k_1=0$.  Then
\begin{eqnarray*}
A_{1s} = A_{2s} = \frac {m} {r}  \\  
M_{1s} = \frac {m r_1} {\lambda r} \cdot \frac {\lambda
k_a(l_2)r_2} {\lambda k_p r_2 } 
=\frac {m(r + \lambda k_d)} {r \lambda k_a (l_1)} \simeq \frac {m k_d}
{r k_a(l_1)} \\
M_{2s} \simeq \frac {m k_d } {r k_a(l_2)} 
\end{eqnarray*}

The steady state of the system reveals that if the modified substance
does not diffuse between the two compartments, then it is as if the two
compartment were entirely separated.  $A$ and $M$ settle into the same
steady state values that they would have if the two compartments were not
connected at all.

We do not gain additional information from assuming that $k_2=0$, because
the qualitative behavior of the system remains the same in this case.  
(See Appendix \ref{ss_calc}.)  However, assuming that $k_2 \gg 1$ changes
the qualitative behavior.
\begin{eqnarray*}  
\lim_{k_2\to\infty} A_{1s} = \lim_{k_2\to\infty} \frac {m} {r} \Big[ 1+
\frac {r_1 k_1 k} {\lambda (r + 2k_2)k_p +k_1 k_s (r_1+2 k_2) } \Big] =
\frac {m}{r} \\
\lim_{k_2\to\infty} A_{2s} =  \frac {m}{r} \\
\lim_{k_2\to\infty} M_{1s} = \frac {m r_1} {\lambda r} \cdot \Big[ \frac
{\lambda k_a(l_2 )  + 2 k_1} {k_a(l_1)(\lambda k_a(l_2) + 2 k_1) - k_1
k} \Big]\\
\lim_{k_2\to\infty} M_{2s} = \frac {m r_1} {\lambda r}
\cdot \Big[\frac {\lambda k_a(l_1 ) + 2 k_1} {k_a(l_2)(\lambda k_a(l_1) +
2 k_1) + k_1 k} \Big]
\end{eqnarray*} 

These calculations show that the steady state of the active substance will
be independent of the external ligand concentrations in the limit as $k_2
\rightarrow \infty$.  Therefore $k_2 \gg 1$ leads to a diminished ability
of the system to respond to stimulus.  The steady state of the modified
substance still depends on the external stimulus for $\lambda \approx
O(1)$.  However, if $\lambda \gg 1$, then $M_{1s} = \frac {m k_d} {r
k_a(l_1)}$, and similarly, $M_{2s} = \frac {m k_d} {r k_a(l_2)}$, so in
the limit as $\lambda \rightarrow \infty$, our system again responds to
stimulus as two unconnected compartments would.

These simple calculations above suggest that the flux of $M$, the modified
substance must be larger than the flux of the activated substance, $A$, or
$k_1 >> k_2.$ A heuristic explanation of this is as follows.  If $k_2 >
k_1$, i.e if $A$ were allowed to diffuse faster than $M$, then regardless
of the external ligand concentration the amount of activated substance in
the two compartments becomes the same.  In this case the steady state of
$M$ still depends on the ligand concentration, unless the reactions
between $M$ and $A$ are much faster than all other rates, i.e. if $\lambda
\gg 1$.  If $\lambda$ is large, then $M$ and $A$ are allowed to exchange
quickly, and the concentration of the modified substance will depend
mainly on $k_a(l)$ and $k_d$, therefore $M_{1s}$ and $M_{2s}$ return to
the values they would have in case of two unconnected compartments.
  
\begin{figure}[h!]\centering
\centerline{\includegraphics[width=0.45\textwidth]{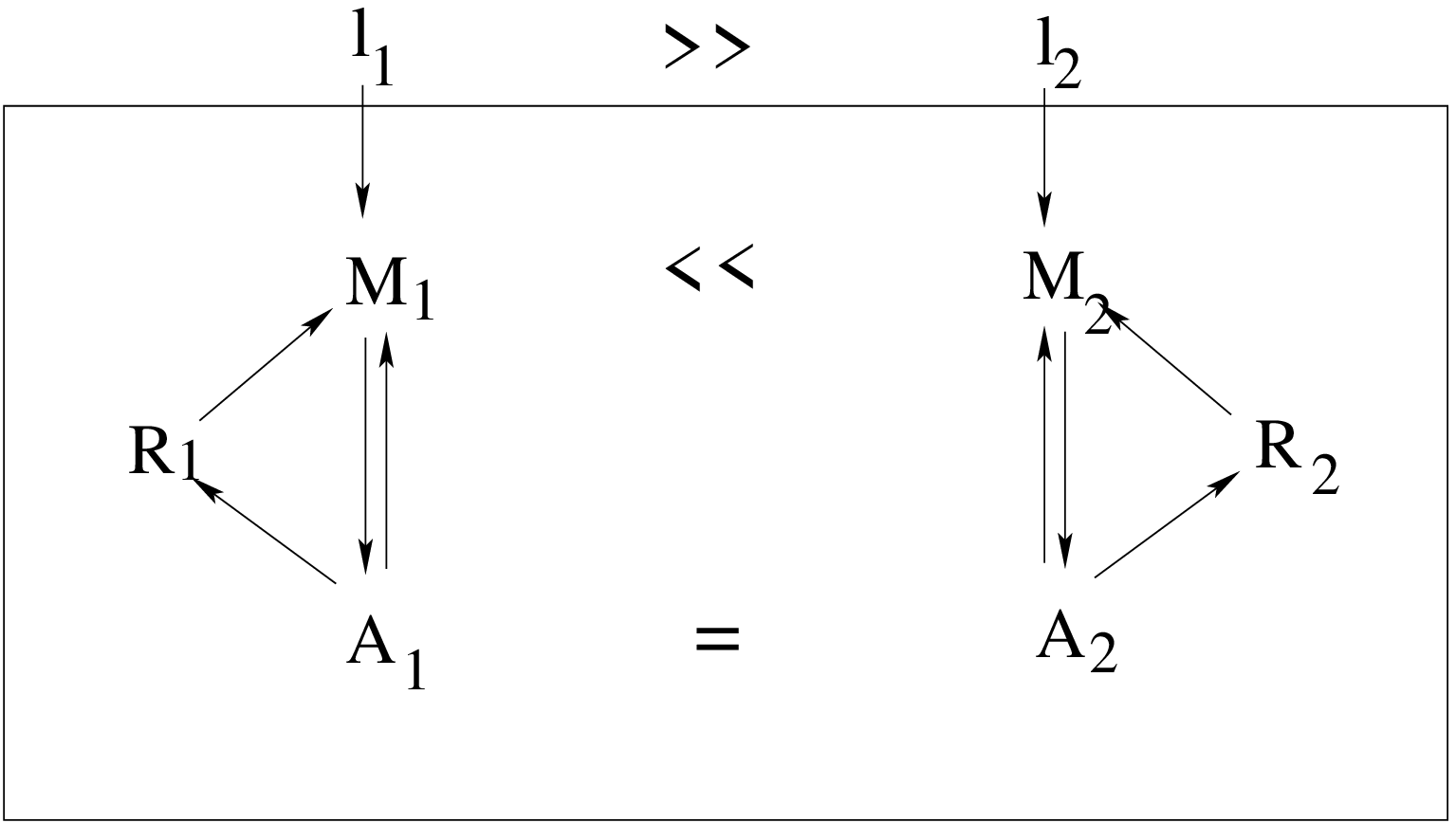} \quad
\includegraphics[width=0.45\textwidth]{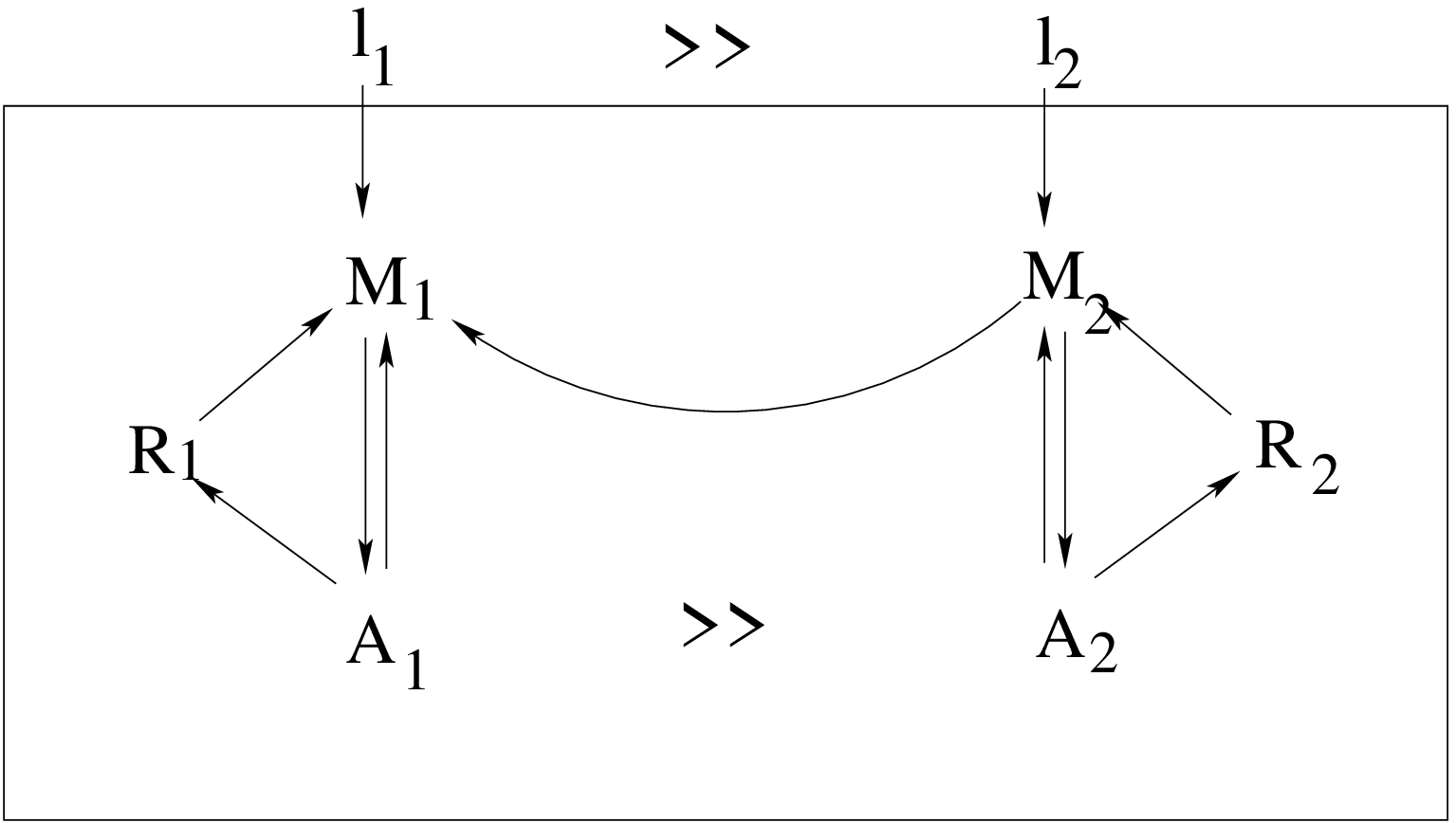}}
\caption[Explanation of the two compartment model.]{Illustration of the 
two compartment model when $k_1 \gg k_2$.  
The first figure shows the two compartments without any connection.  In
the second figure we assume that $M$ is allowed to diffuse between the
compartments.  This results in the creation of a spatial gradient of $A$.
} \label{heuristic} \end{figure}

Now let us consider the case that $A$ is non-diffusive, but $M$ is.  Let
us assume that the ligand concentration at the first compartment is larger
than in the second compartment, i.e. $l_1 > l_2$.  As we had seen from
eqn.  \ref{two_odes}, the steady state of $M$ is inversely proportional to
$k_a(l)$ which we are taking to be proportional to the ligand
concentration, i.e, the steady state of $M$ is inversely proportional to
the ligand concentration.  This implies that $M_1$, the concentration of
the modified substance in the first compartment is smaller than $M_2$.  
The flux between the two compartments will increase the value of $M_1$
(and decrease $M_2$), therefore more $A_1$ is produced in compartment one
than $A_2$ in compartment two.  Because the flux, $k_2$ between $A_1$ and
$A_2$ is negligible compared to $k_1$, the difference $A_1 - A_2$ is
maintained.  Figure \ref{heuristic} illustrates this mechanism.

We are also interested in how the difference between the levels the
activated substance in the two compartments, $A_1 -A_2$ changes with the
ligand concentration.  Based on our previous calculations, we assume that
$k_2 \approx 0$.  Furthermore, we assume that $k_a(l)$ is a monotonically
increasing function of $l$, and that $k_a(l) > 0$ for all $l$.  When two
different ligand concentrations are presented for compartments one and
two, we see a difference in the rates $k_a(l_1)$ and $k_a (l_2)$, thus a
spatial gradient in the ligand produces the spatial gradient of $A$.  We
notice that $|A_{1s}-A_{2s}|=0$ if $k_a(l_2)=k_a(l_1)$, and this is the
minimum value $|A_{1s}-A_{2s}|$ can obtain. $$0< |A_{1s}-A_{2s}| = \Big|
\frac {2m} {r} \cdot \frac {k_1 r_1 k} {\lambda r k_p + k_1 k_s r_1}
\Big|< \frac {2m} {r} \Big| \frac {k_1 r_1 k} {k_1 r_1 k_s} \Big| < \frac
{2m} {r} $$

The absolute difference, $|A_{1s}-A_{2s}|$ is bounded below by 0, and
above by $\frac {2m} {r}$.  It is also easy to show that for arbitrarily
large differences in the ligand concentrations the difference between
$A_1$ and $A_2$ approaches $\frac {2m} {r}$.  (We show this by taking
$\lim_{k_a(l_1)\to\infty}$ and $\lim_{k_a(l_2)\to 0}$. Because of the
symmetry of the expression, $\lim_{k_a(l_2)\to\infty}$ and
$\lim_{k_a(l_1)\to 0}$ leads to the same result.)

We want to know how $|A_{1s}-A_{2s}|$ depends on the difference in
ligand concentrations, $k$ and on absolute size of the ligand
concentrations, $k_s$.  We introduce new constants: $a_1 = \frac {2m}
{r}$, $a_2 = k_1 (r + \lambda k_d)$ and $a_3 = \lambda r.$
$$A_{1s}-A_{2s}= f(k,k_s) = a_1 \frac {a_2 k} {a_2 k_s + \frac {a_3} {4}
(k_s^2 -k^2)}.$$

We are interested in $\frac {\partial f} {\partial k}$, how the absolute
value of the difference in steady state of $A$ in our two compartment
depends on the ligand difference, and $\frac {\partial f } {\partial
k_s}$, how the difference depends on the size of the ligand
concentrations.
\begin{eqnarray*} \frac {\partial f} {\partial k} = a_1 a_2 \frac {a_2 k_s
+ \frac {a_3} {4} (k_s^2 +k^2)} {[a_2 k_s + \frac {a_3} {4}(k_s^2
-k^2)]^2} > 0 \ \ \forall k, k_s \end{eqnarray*}

The difference, $A_{1s}-A_{2s}$ always increases with the increasing 
difference in the ligand concentration.  (This is clear from the formula 
$A_{1s}-A_{2s}$ as well.)  Now we look at how the absolute concentration 
level of the ligand changes $A_{1s}-A_{2s}$:
\begin{eqnarray*}
\frac {\partial f} {\partial k_s} = a_1 \frac {-a_2 + \frac{a_3} {2} k_s} 
{[a_2 k_s + \frac {a_3} {4} (k_s^2 -k^2)]^2} \\
\frac {\partial f} {\partial k_s} < 0 \qquad\mbox{if}\qquad a_2 > \frac 
{a_3} {2k_s} 
\end{eqnarray*}

In our original notation this expression means that $|A_{1s}-A_{2s}|$ is a
decreasing function of the sum of the ligand concentrations if $$k_1 (r
+\lambda k_d ) > \frac {\lambda r} {k_a (l_1) + k_a(l_2)}$$ or $k_1 (r
+\lambda k_d) (k_a (l_1) + k_a(l_2)) > \lambda r$.  This relationship
shows that (depending on the explicit form of $k_a(l)$), there is an
absolute ligand concentration at which the difference between the
activated substance in the two compartments will be the largest.  Let us
assume that we fix the difference in the ligand concentrations
$k_a(l_1)-k_a(l_2)$, and only change their sum.  Increasing the absolute
concentration of the ligand beyond the point where $$k_a (l_1)  + k_a(l_2)
= \frac {\lambda r} {k_1 (r +\lambda k_d)}$$ will result in decreased
sensitivity in sensing, and similarly, smaller ligand concentrations also
result is a loss of sensitivity.  The existence of a range of 
ligand concentrations in which sensing is optimal corresponds to 
experimental observations.  If the ligand concentrations are too low, then 
the noise in the receptor occupancy leads to errors in gradient sensing.  
If the ligand concentrations are too high, then receptors are saturated, 
and the cell's ability to to detect gradients is compromised again.  

\begin{figure}[h!] \centering
\centerline{\includegraphics[width=0.8\textwidth]{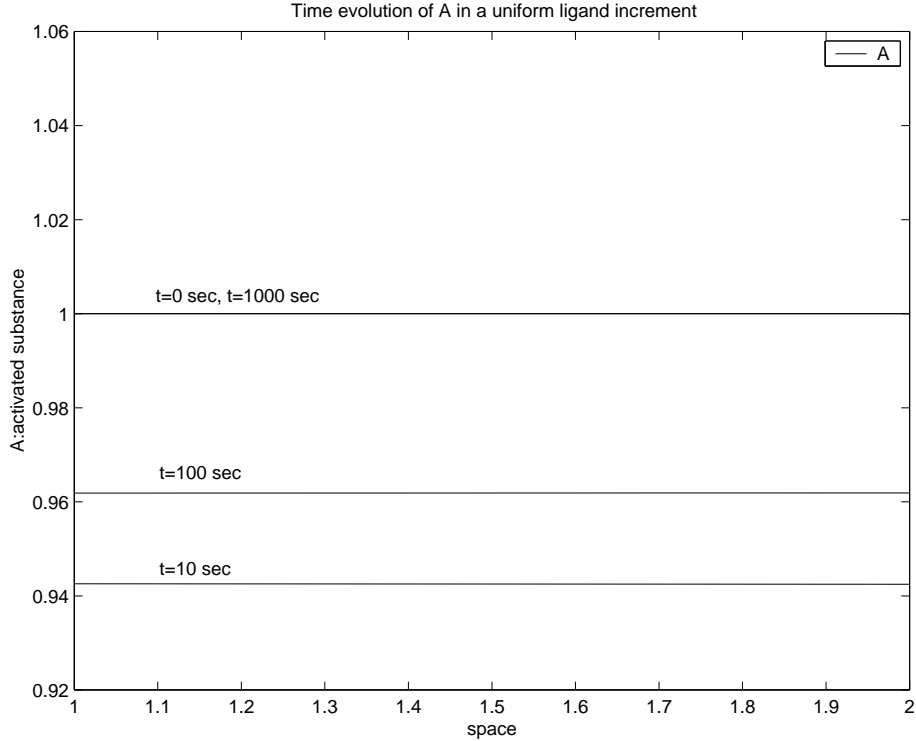}}
\caption[Time evolution of A in a ligand increment.]{Time evolution of the 
active substance, A in a ligand increment.}
\label{two_comp4.eps}
\end{figure}

\begin{figure}[h!] \centering
\centerline{\includegraphics[width=0.8\textwidth]{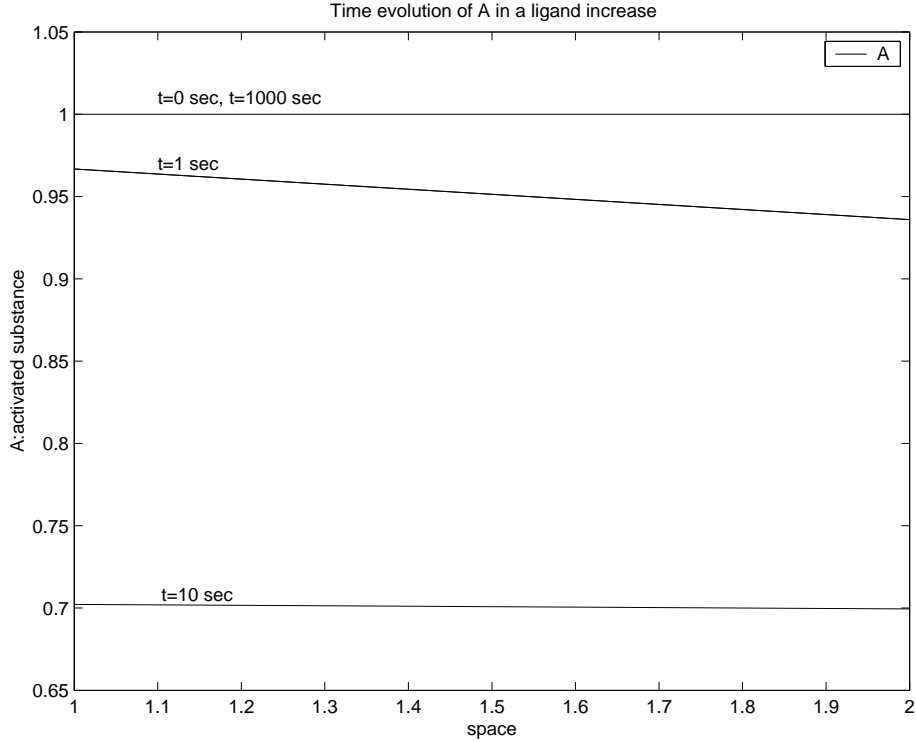}}
\caption[Time evolution of A from gradient to spatially uniform 
ligand.  ]{Initial condition: $A$ is adapted to a ligand gradient.
Temporal dynamics of $A$ when the two compartments are exposed to the same
ligand concentration.} \label{two_comp5.eps} \end{figure}

\begin{figure}[h!] \centering
\centerline{\includegraphics[width=0.8\textwidth]{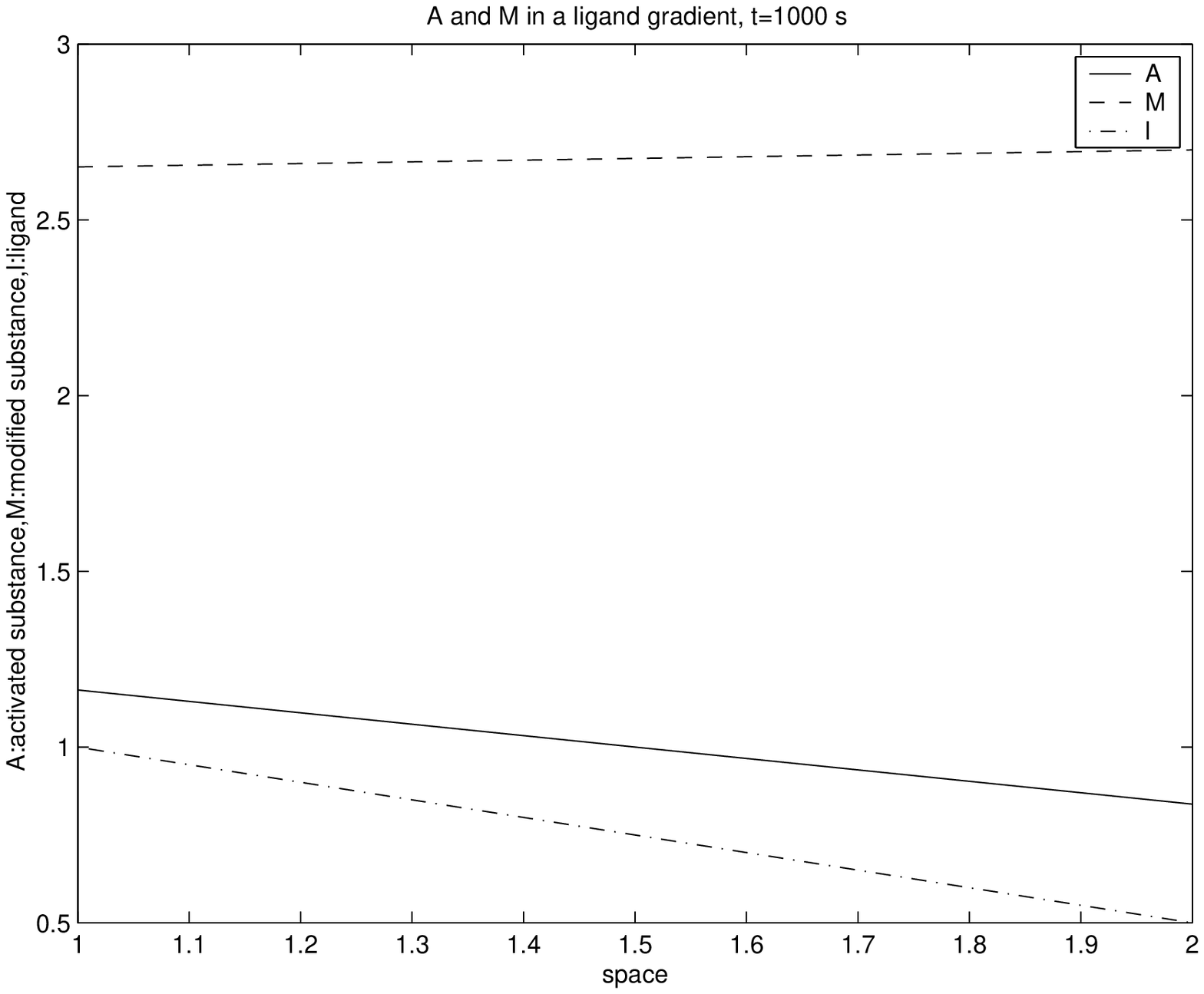}}
\caption[Steady state of the system in a gradient.]{Steady state of the 
two compartment system in a ligand
gradient.} \label{two_comp2.eps}
\end{figure}

The following numerical simulations confirm the behavior of the system.  
The first figure, \ref{two_comp4.eps} shows the response of the activated
substance, $A$ in a spatially uniform change of the ligand concentration.  
The initial condition of $A$ is $A(0)= \frac {m} {r}$ which is independent
of the ligand concentration.  We start with our two compartment model
adapted to a spatially uniform ligand concentration, $l_0 = 0.1$.  Figure
\ref{two_comp4.eps} shows the temporal dynamics of the system in a
spatially uniform ligand step to $l_1 = 1$. There is a quick drop in the
concentration of $A$, then a slow adaptation to the steady state level,
$A_{ss}=\frac {m} {r}$.

We see similar behavior if we start with our system adapted to a ligand
gradient (which does not reflect in the initial conditions for $A$), and
then let the ligand concentration jump to the same uniform level inside
both compartments.  In compartment one the initial ligand concentration is
0.1, and in compartment two the ligand concentration is 0.5, and the
ligand concentration jumps to $l=1$ in both compartments.  Although the
transient levels of $A$ are not the same in the two compartments, the
steady state levels are.  Figure \ref{two_comp5.eps} shows these
simulations.  The values for $k_1$ and $k_2$ are 1 and 0.1, respectively.  
The qualitative behavior of the system remains the same if $k_2=0$ is
used.  In these figures the active substance in compartment one, $A_1$ is
always given by $A(1)$, and $A_2$ by $A(2)$.  Matlab, instead of plotting
the value $A_1$ from 1 to 1.5 and $A_2$ from 1.5 to 2 connects $A_1$ and
$A_2$.  These numerical simulations illustrate that our two compartment
system responds transiently to changes in the ligand concentration, and it
settles into a ligand-independent steady state in uniform ligand
concentrations.

Now we must examine the other main claim of our model, namely, that it
sets up an internal gradient of the active substance when presented with a
spatial attractant gradient.  The numerical simulation for this is shown
in Figure \ref{two_comp2.eps}.  This figure depicts the steady state of
the modified and the active substance when the ligand concentration in the
first compartment is $l_1=1$, and in second compartment $l_2=0.5$.  The
simulation is run for t=1000 seconds.  $k_1=1$ and $k_2=0.1$, and again,
$k_2=0$ does not change the qualitative behavior.  The spatial gradient of
the activated substance is maintained while the ligand gradient remains 
unchanged.  This results in the persistent signaling of the system in 
ligand gradients.  

Appendix \ref{anal_sol} contains notes and comments on the analytical
solution of the two compartment model, as well as approximate solutions to 
the problem if $\lambda \gg 1.$ 

In a general case we can assume that both $A$ and $M$ are functions of
time and space, and they are both allowed to diffuse.  We obtain the
following equations:
\begin{eqnarray*}
\frac {\partial M} {\partial t} = m + \lambda (-k_a(l) M + k_d A) + D_1
\frac {\partial^2 M} {\partial x^2}  \\
\frac {\partial A} {\partial t} = -r A + \lambda ( k_a (l) M - k_d A) +
D_2 \frac {\partial^2 A} {\partial x^2} \\
A(x,0) = \frac{m} {r} \ \ M(x,0) = \frac {m} {r} \frac {k_d} {k_a(l)}  \\
\frac {\partial A} {\partial x}|_{x=0} = \frac {\partial A} {\partial 
x}|_{x=L} = 0 \\
\frac {\partial M} {\partial x}|_{x=0} = \frac {\partial M} {\partial 
x}|_{x=L} = 0 \\
\end{eqnarray*}

\begin{figure}[h!]
\includegraphics[width=0.45\textwidth]{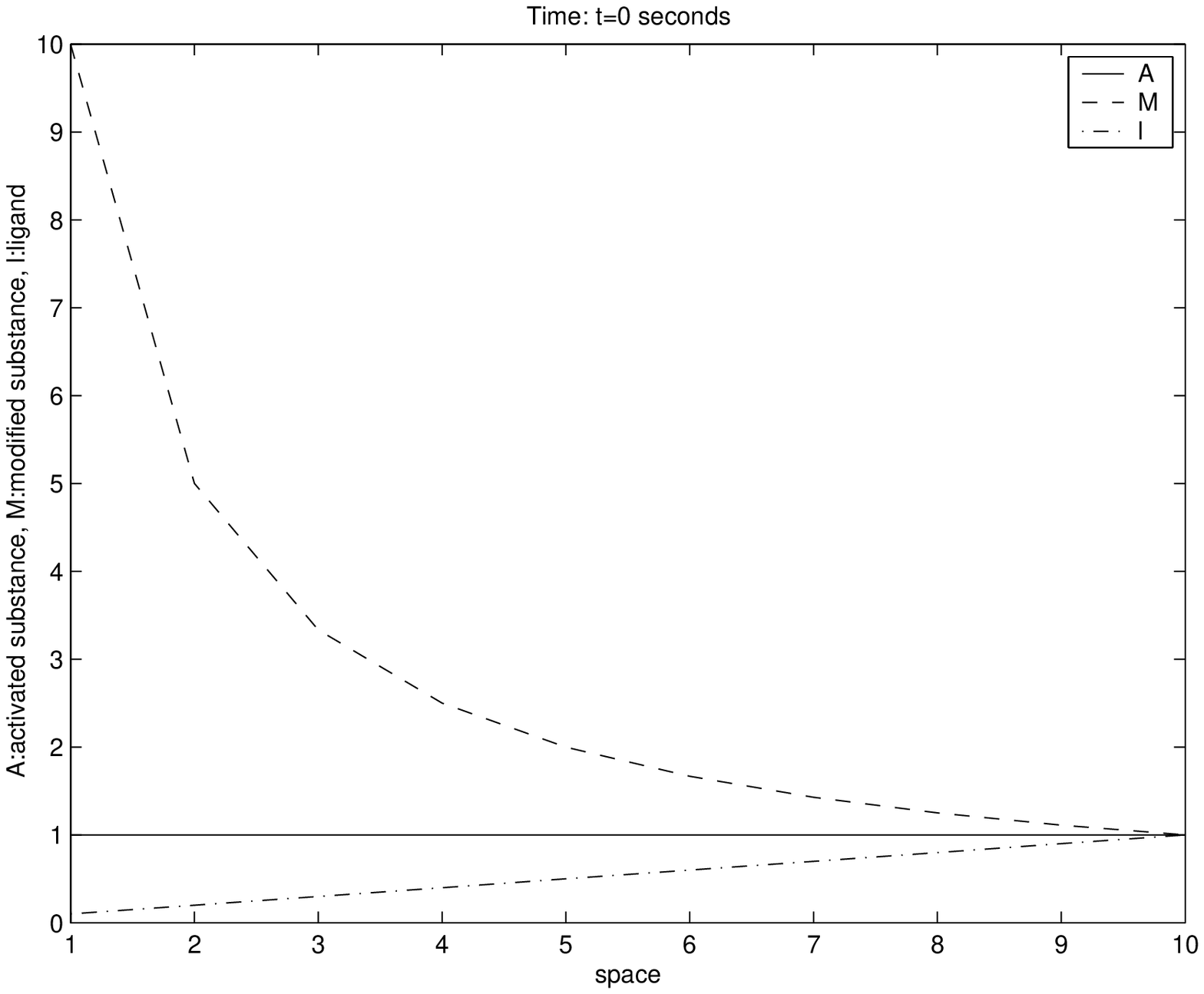}\qquad
\includegraphics[width=0.45\textwidth]{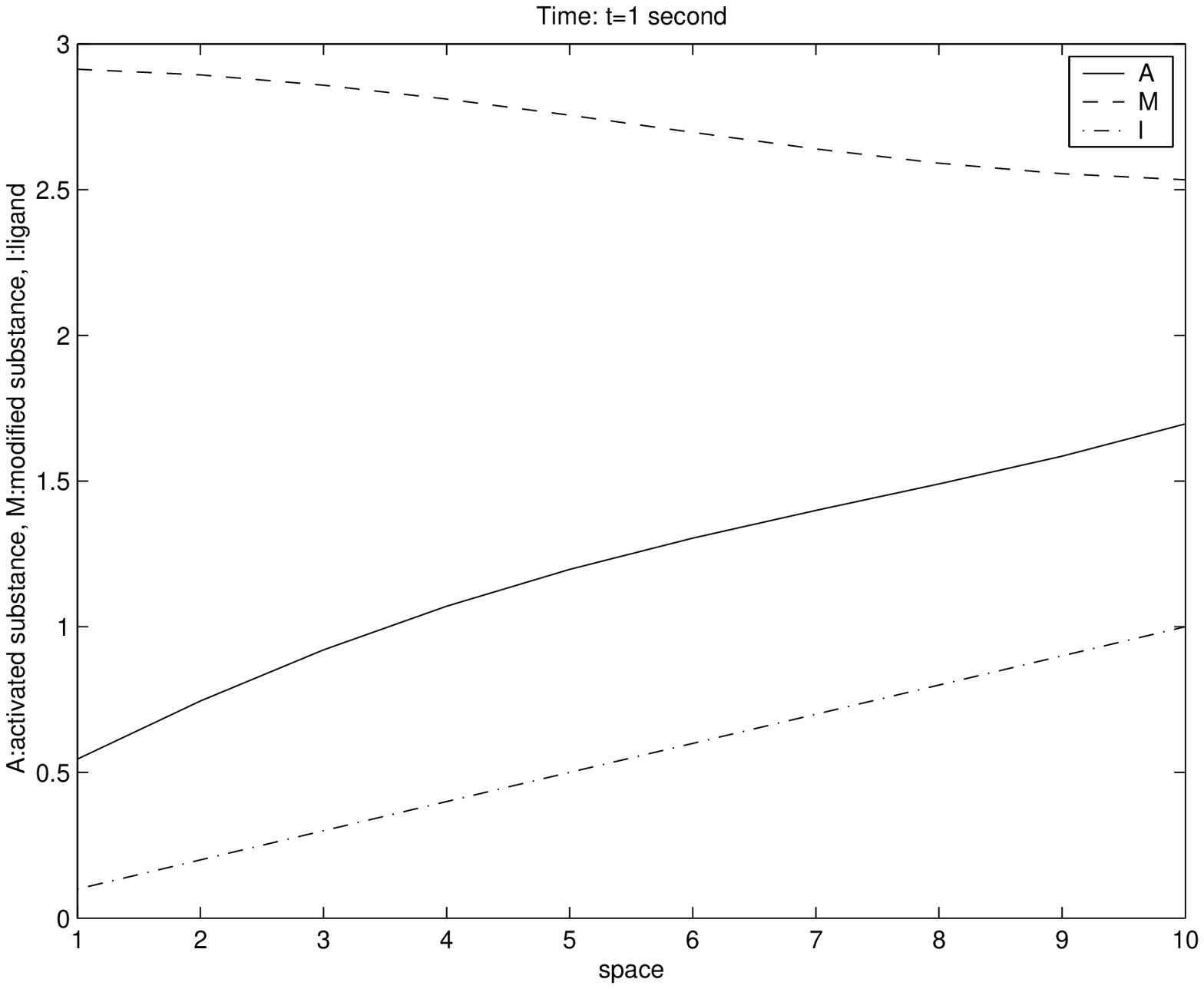}\\
\includegraphics[width=0.45\textwidth]{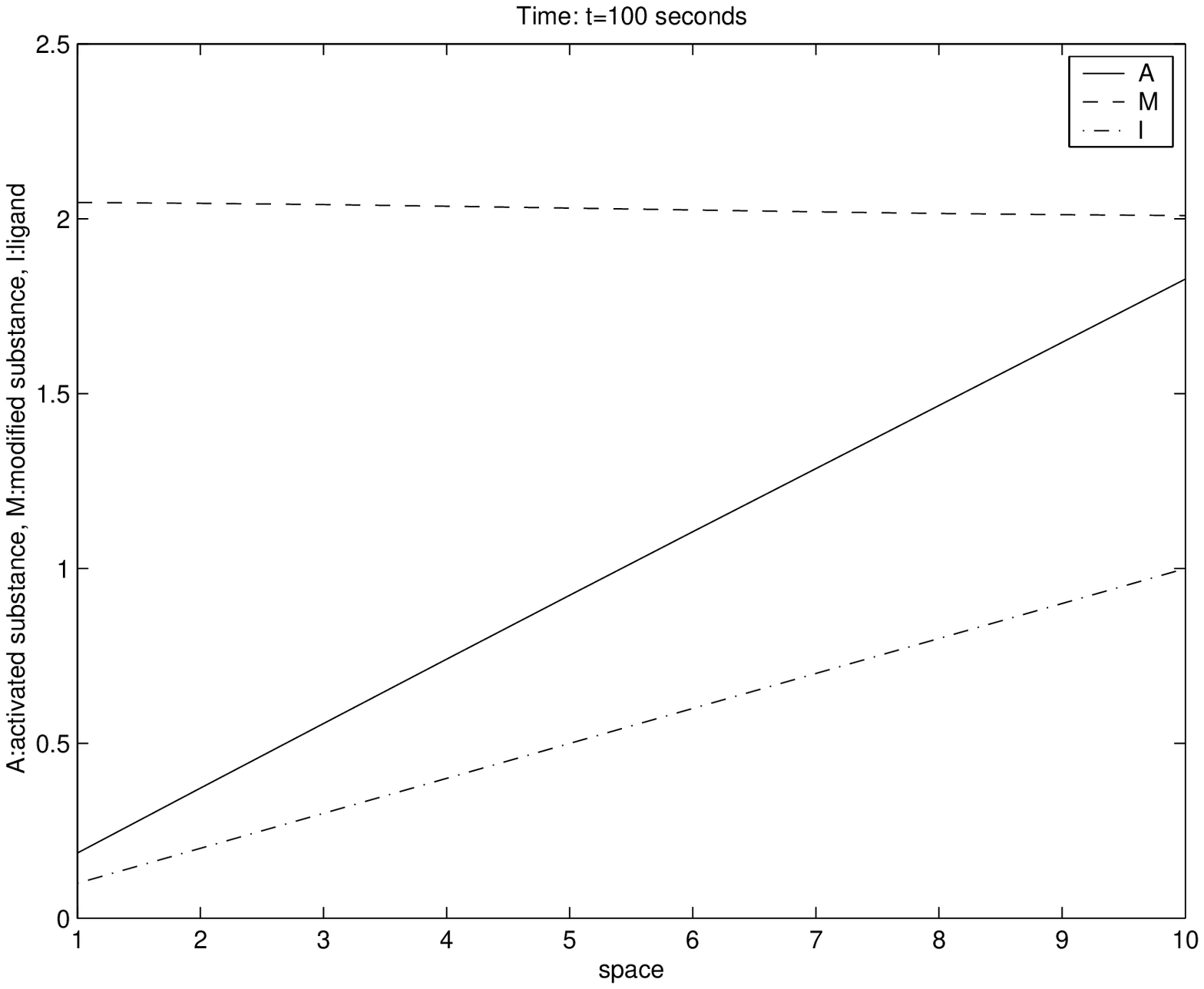}\qquad
\includegraphics[width=0.45\textwidth]{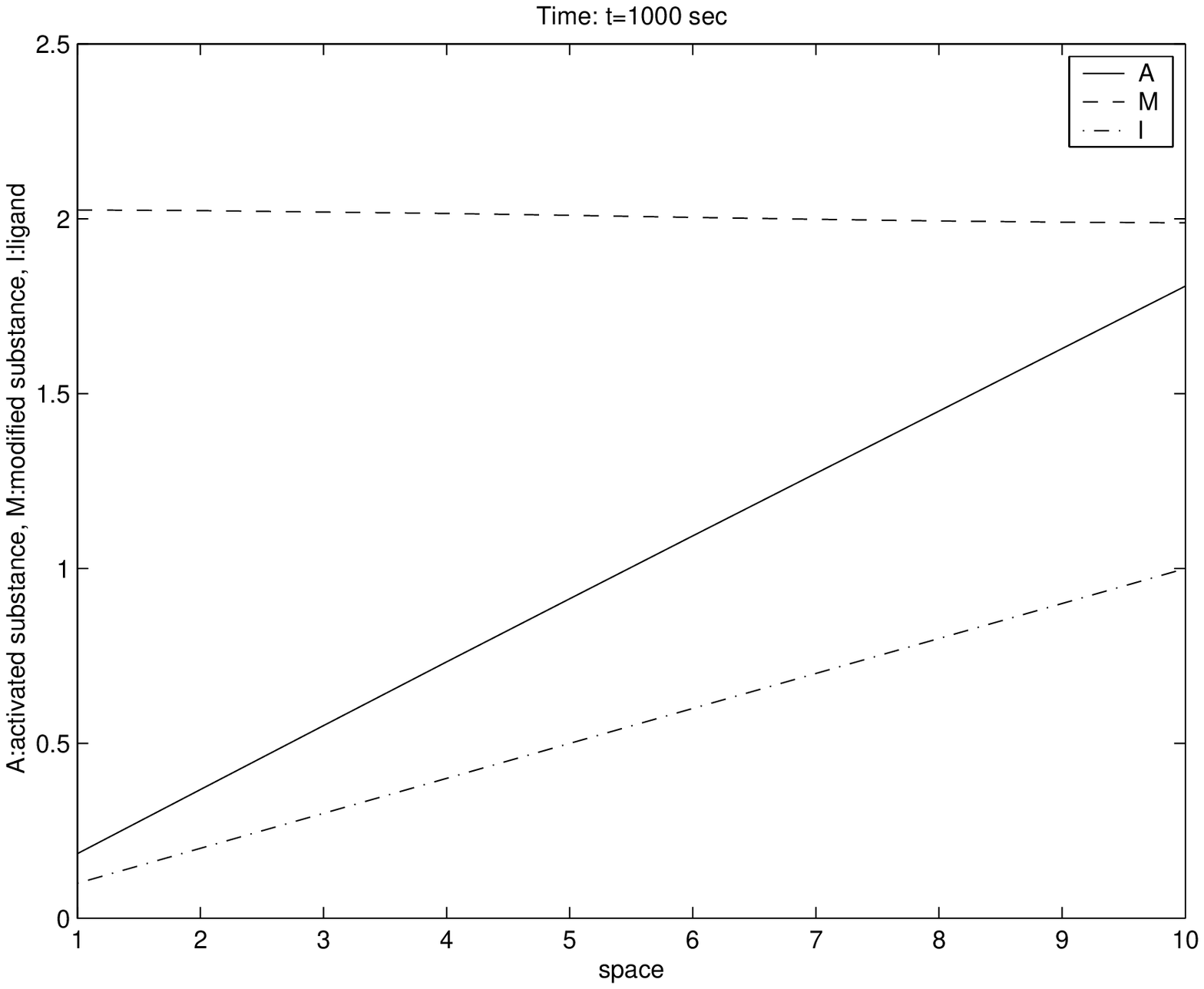}
\caption[Temporal dynamics of A in a linear gradient, reaction-diffusion 
equation. ]{The figures show the temporal dynamics of 
the 
reaction-diffusion system in a linear attractant gradient.  The first 
figure shows the initial condition, the second one the system at t=10 
seconds, and the third one shows the steady state of the system at t=100
seconds, and the fourth one at t=1000 seconds.} \label{diff_l_lin} 
\end{figure}

\begin{figure}[h!]
\includegraphics[width=0.45\textwidth]{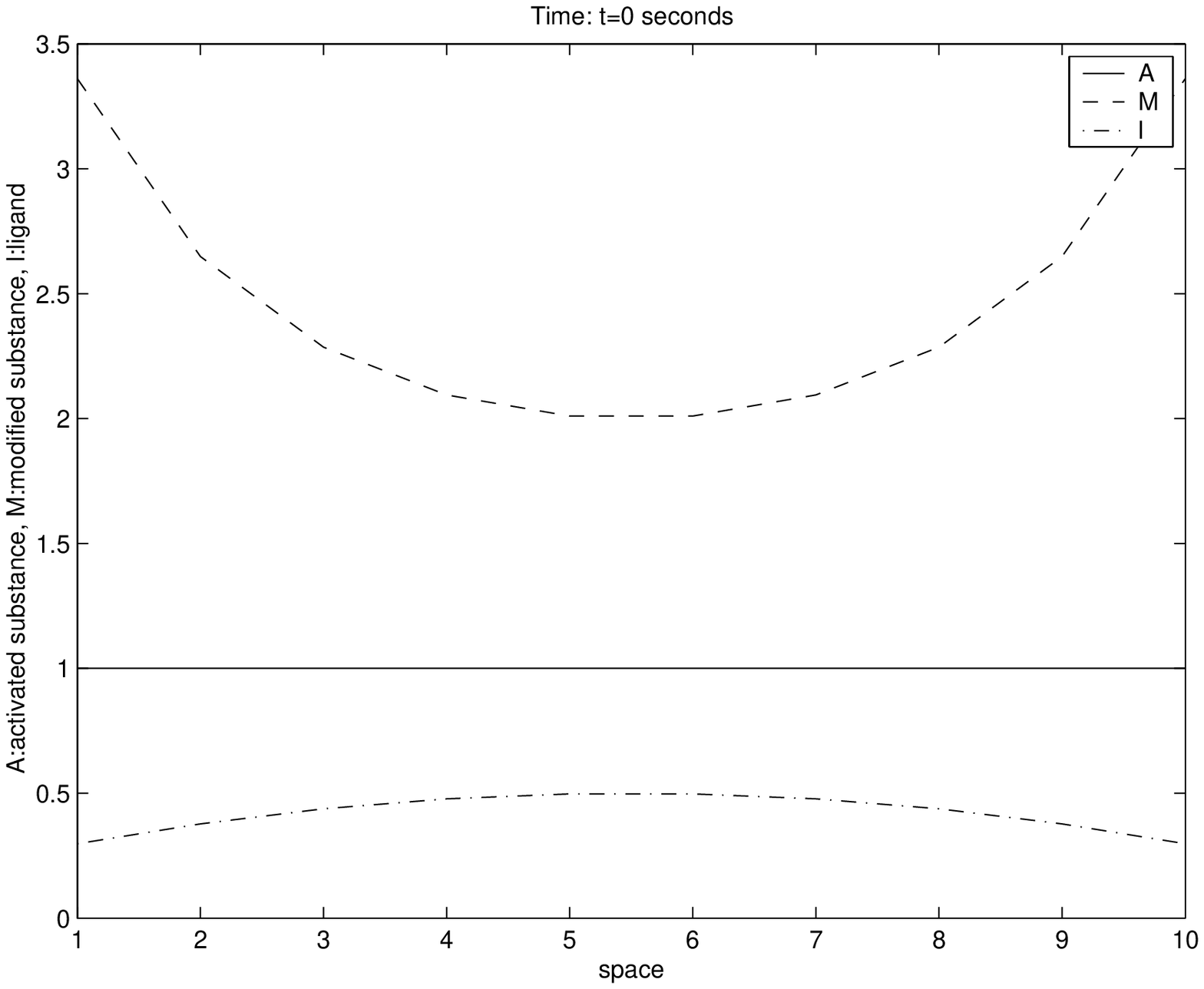}\qquad
\includegraphics[width=0.45\textwidth]{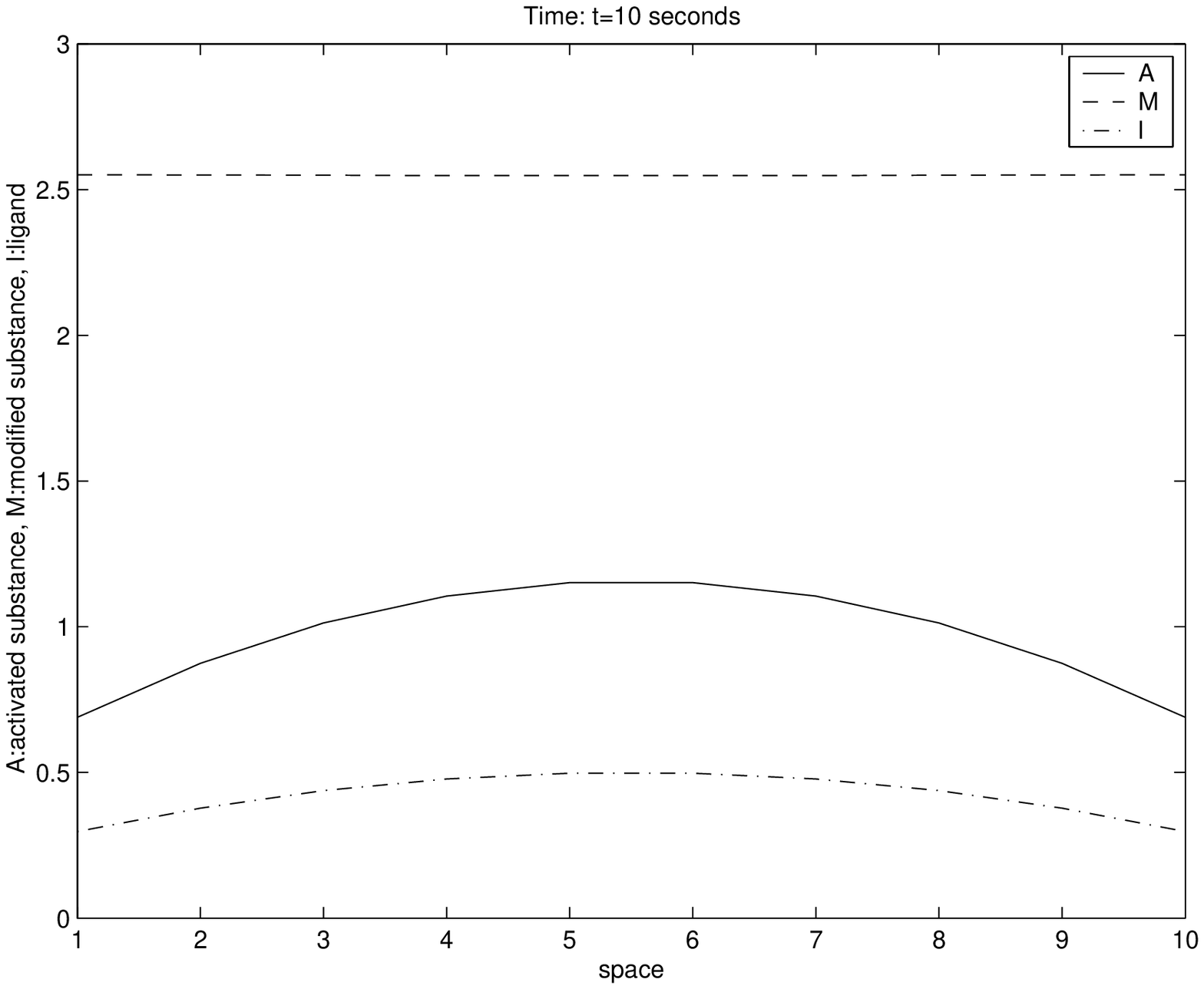}\\
\centerline{\includegraphics[width=0.45\textwidth]{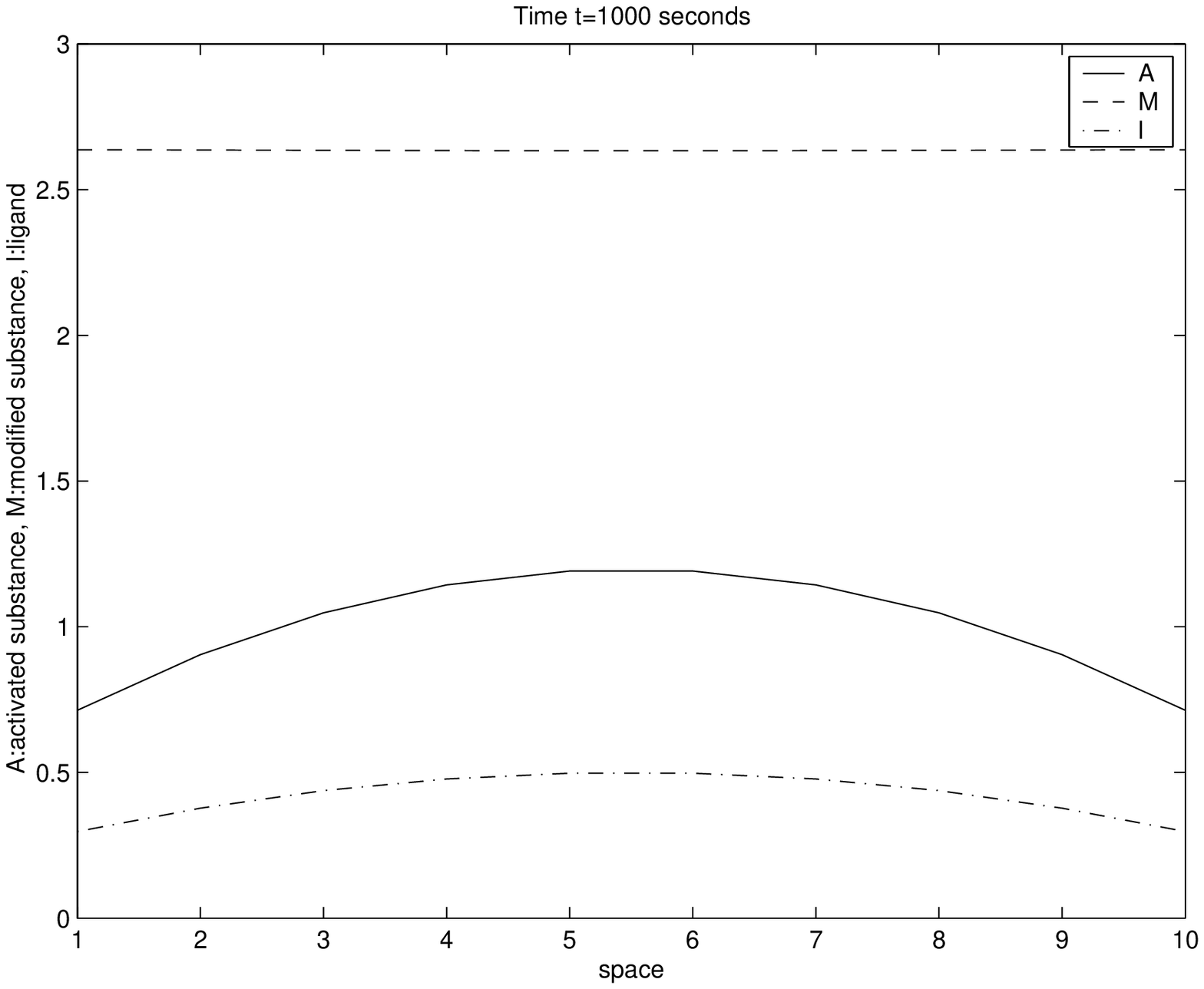}}
\caption[Temporal dynamics in a nonlinear gradient, 
reaction-diffusion system.  ]{The figures show the temporal dynamics of 
the
reaction-diffusion system in a quadratic attractant gradient.  The first
figure shows the initial condition, the second one the system at t=10
seconds, and the third one shows the steady state of the system at t=1000
seconds.}
\label{diff_l_quad}
\end{figure}

We provide numerical solutions to this system on the interval
$(0,L)=(0,10)$ in Figures \ref{diff_l_lin} and \ref{diff_l_quad}.  Based
on our analysis of the two compartment model, we assume $D_1 \gg D_2$ and
$\lambda \gg 1.$ The diffusion coefficient of $A$, $D_2$ is chosen to be
zero, but as before, $D_2 \ll 1$ would also provide qualitative similar
results.  The initial conditions are chosen to be the same as in the
differential equation with $\lambda \gg 1.$ The equations are solved with
a FTCS (forward time center space) method.  The time step is chosen to be
$\Delta t = 0.01$, and the grid size is $\Delta x = \frac {1} {9}= 0.11 $.  
We show three sets of simulations verifying that in uniform ligand 
concentrations there is a transient response, and that ligand gradients 
elicit a persistent graded response.  

\begin{figure}[h!] 
\includegraphics[width=0.45\textwidth]{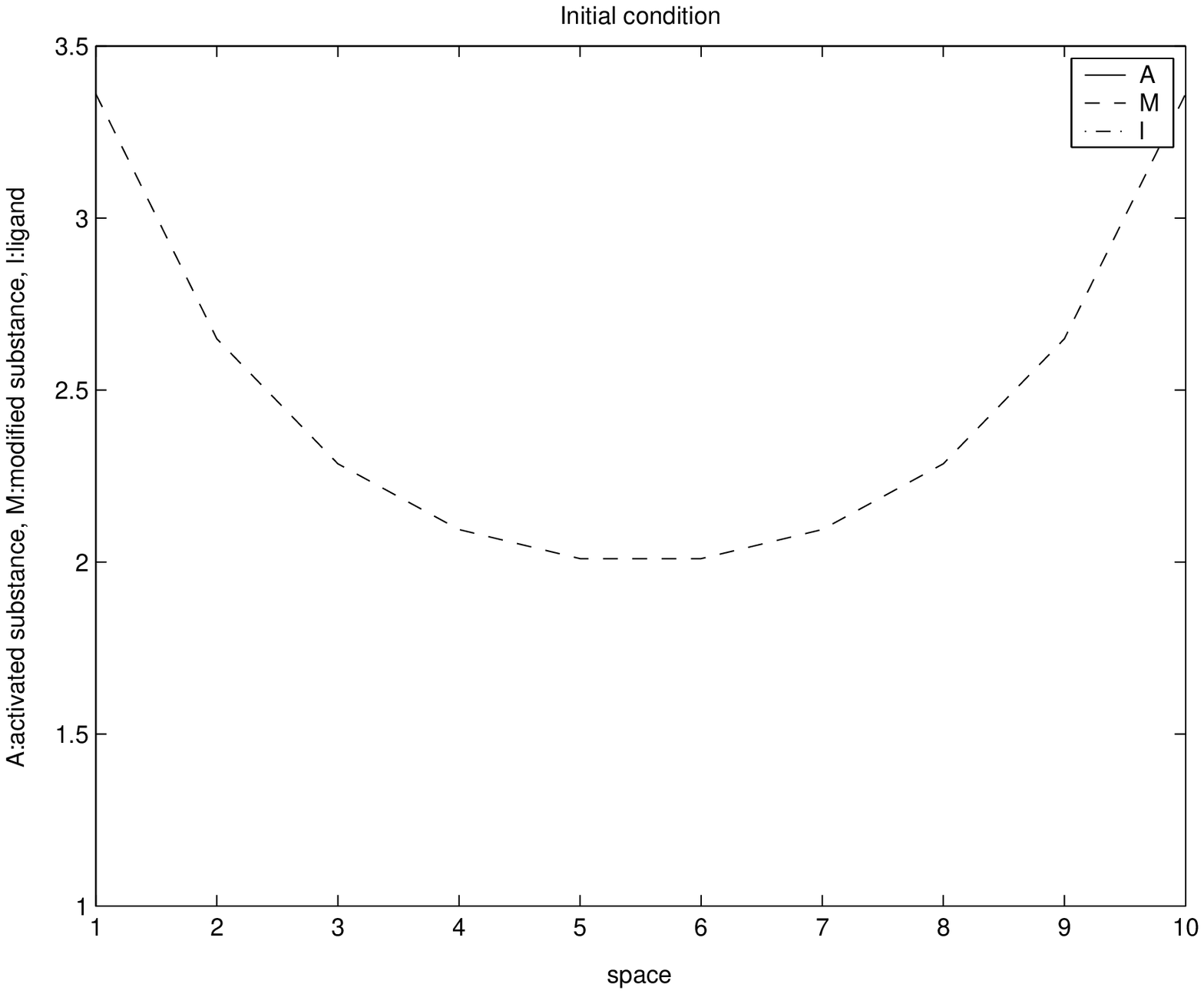}\qquad
\includegraphics[width=0.45\textwidth]{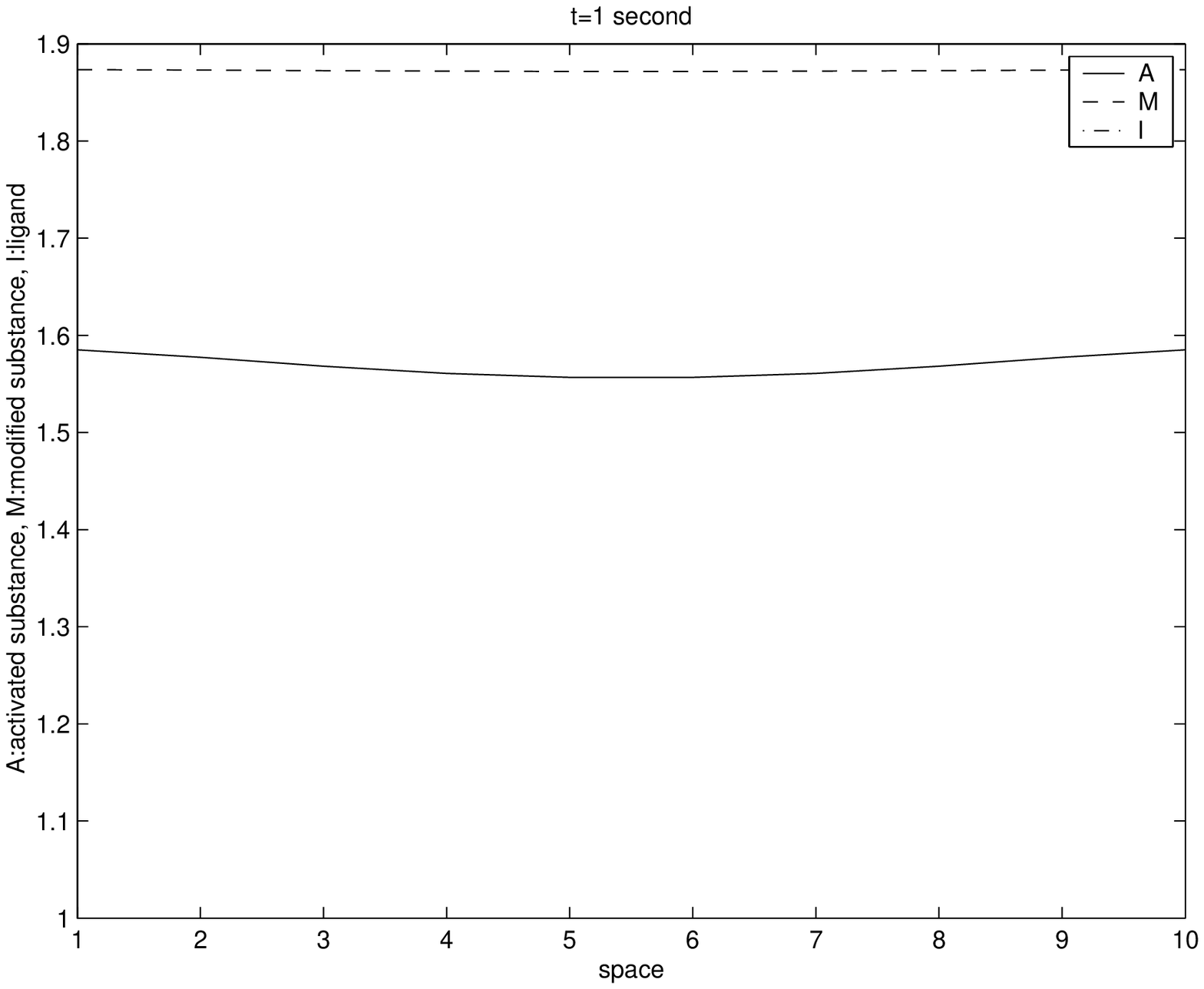}\\
\centerline{\includegraphics[width=0.45\textwidth]{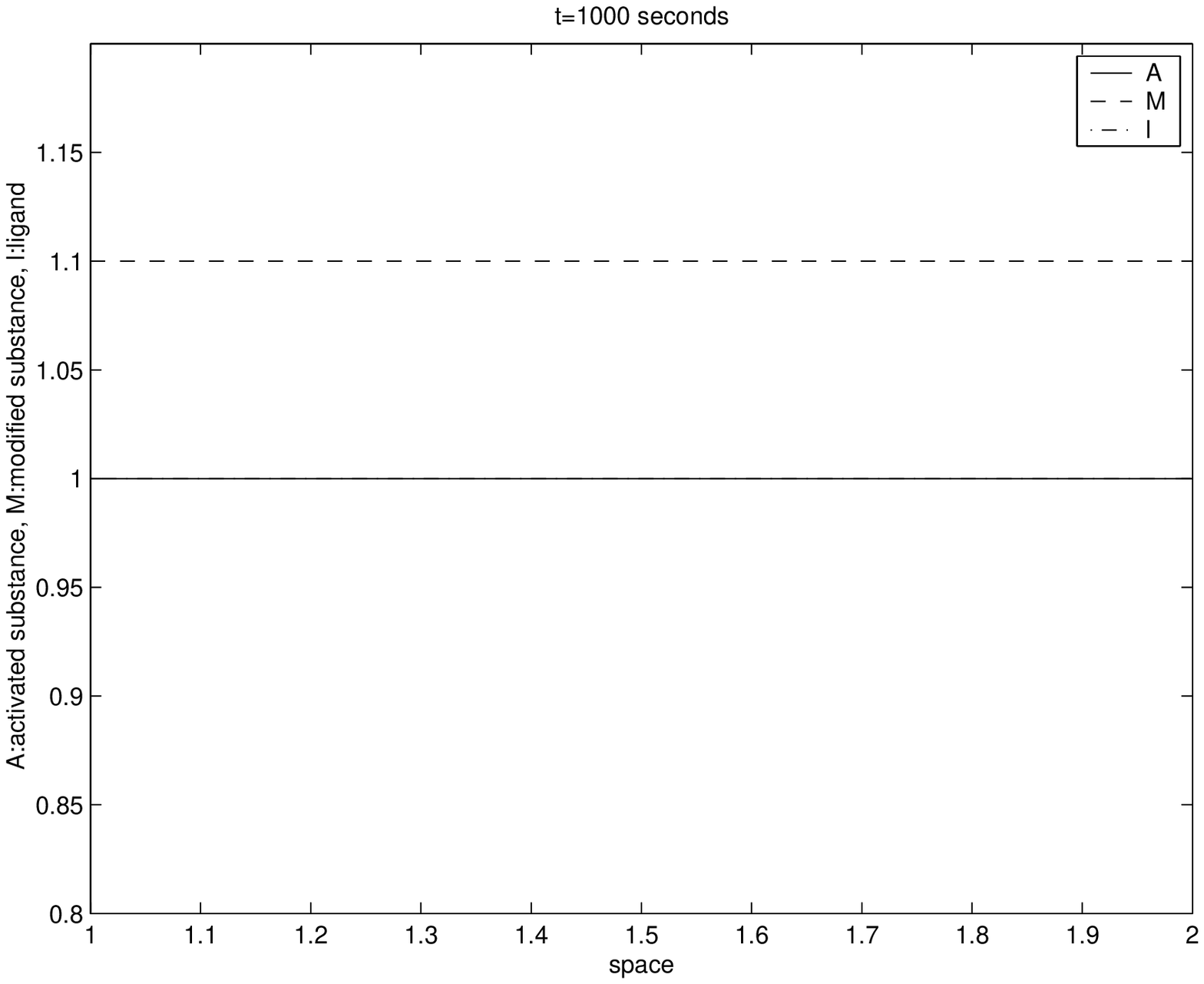}}
\caption[Temporal dynamics of A in a ligand increment, reaction-diffusion 
equation. ]{The figures show the temporal dynamics of the reaction-diffusion
system in a uniform attractant concentration.  The first figure shows the
initial condition, the second one the system at t=1 second, and the third
one shows the steady state of the system at t=1000 seconds. }
\label{udiff} \end{figure}

In figure \ref{udiff} we examine the behavior of the system which has
adapted to a ligand gradient, but at time t=0 we present to the cell a
spatially uniform ligand concentration, $l=1$.  The middle panel shows
that the cell responds very quickly, and in one second the both $A$ and
$M$ are almost uniform.  (The ligand concentration is equal to one
everywhere which is difficult to see in the figures.)  Finally, the panel
on the right shows the steady state of the system at t=1000 seconds.  $A$
returns to the baseline value of one, and $M$ also becomes a constant in
space.

The following two figures show the response of a cell to an attractant
gradient. In the first one, in Figure \ref{diff_l_lin} the attractant
concentration, $l$ is linear, as before, and in the second one, in Figure
\ref{diff_l_quad}, $l$ is quadratic.

As before, the active substance, $A$ is spatially uniform initially, while
the modified substance, $M$ is inversely proportional to the ligand
concentration.  The system is shown at t=1 second on the second figure, at
t=100 seconds in the third figure and at the steady state, in the last
figure.  As in the two compartment model with a ligand gradient, in the
steady state the concentration of $A$ is proportional to the ligand
concentration, and the concentration of $M$ is inversely proportional to
it.  This spatial profile persists, representing the persistent signaling
of the system in a spatial gradient.

Figure \ref{diff_l_quad} shows dynamics in a nonlinear ligand gradient.  
Similarly to Figure \ref{diff_l_lin}, initial conditions of $A$ are
spatially uniform, and those of $M$ are inversely proportional to the
ligand concentration.  The second subfigure shows the system after t=10
seconds, and finally, the third subfigures shows the steady state after
t=1000 seconds.  At the steady state the modified and the activated
substances are quadratic.  As expected, the highest value of $M$
corresponds to the lowest value of the ligand, and the highest value of
$A$ corresponds to the highest value of the ligand.  Persistent signaling
is predicted again.

\subsubsection{\large Discussion}
\label{Discussion} 

Our model aimed at recreating a few key features of chemotactic sensing.  
First, we wanted the cell model to respond with transient signaling to a
spatially uniform ligand concentration and with a persistent signal in a
spatial ligand gradient.  The numerical simulations for both the two
compartment model (which can be considered a discretized version of the
reaction diffusion equations) and the reaction diffusion equations show,
in Figures \ref{two_comp4.eps},\ref{two_comp5.eps} and \ref{udiff} that a
spatially uniform ligand concentration elicits a transient response from
the system, but the steady state of the active substance is spatially
uniform.  Persistent signaling is also achieved in ligand gradients,
illustrated by Figures \ref{two_comp2.eps}, \ref{diff_l_lin} and
\ref{diff_l_quad}.  Both transient and persistent signaling can also be
deduced from the steady state analysis of the two compartment model where
the steady state of the active substance depends on the difference of
ligand concentrations.

Secondly, the model must allow the cell to choose new orientation in a 
changed attractant gradient.  This is clearly the case, as both the active 
and the modified substance depend on the ligand concentration through 
$k_a(l)$, the rate of production of $A$ from $M$.  

Finally, this model does not have the limitation of the calcium-cAMP
model, as it allows sensing in all ligand concentrations.  The
adaptation-diffusion model predicts that there is an optimal range of
ligand concentrations when $k_a (l_1)  + k_a(l_2) = \frac {\lambda r}
{k_1 (r +\lambda k_d)}$.  A fixed ligand difference will result in the 
largest signal when this condition is satisfied.  The model also predicts 
that increasing the ligand difference results in an increase in the 
signal.  

We can also improve the model by showing the simple modifications can 
result in the amplification of the signal.  In all previous simulations 
parameters were chosen so that the ligand concentration is amplified 
moderately by the concentration of the active substance, however, this 
need not be the case for our current system with other parameter values.  
It is also important to note that huge amplification of the external 
signal is characteristic of chemotaxis, and such amplification is 
impossible to produce with our current model.  As in the Levchenko \& 
Iglesias model, we assume that amplification takes place downstream from 
$A$, because this assumption guarantees that when the external signal, $l$ 
is terminated, the internal signal stops as well.  The nature of this 
amplification could be similar to what Levchenko \& Iglesias has assumed, 
shown in Figure \ref{lev_igl2.eps} where $A$ would play the role of the 
activated response element, $R*$.   

\begin{figure}[h!] \centering
\centerline{\includegraphics[width=0.2\textwidth]{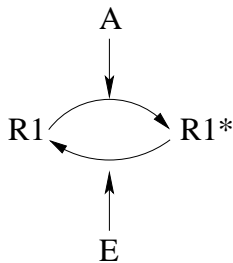}}
\caption[Signal amplification based on Goldbeter and Koshland.]{$A$ as an 
enzyme in the production of an activated response 
element, $R1*$.} \label{gk.eps}
\end{figure}

Another possibility for amplification is based on Goldbeter \&
Koshland \cite{GK} who show that amplification can result not only
from nonlinear interactions, but also form covalent modifications
under certain conditions.  Figure \ref{gk.eps} shows a reaction in
which $A$ is the enzyme in the production of an active response
element, $R1*$.  Based on Goldbeter \& Koshland, there can be a
drastic amplification in the amount of $R1*$ produced when the enzymes
$A$ and $E$ are saturated, so the total amount of free $A$ and $E$ are
negligible when compared with their concentrations in the complexes
produced with $R1$ and $R1*$.  This is one possible mechanism for
amplification of $A$.

It is clear both from the analytical solution to the perfect adaptation
model \ref{anal_sol_2} and the analysis of the steady state solution of
the two compartment model that the size of $\lambda$ is very important. In
the original equations \ref{two_odes} $\lambda$ is the non-dimensional
ratio of the fast and slow time scales which means that $\lambda$
determines how fast the initial transient response is in comparison to the
adaptation.  In animal cell chemotaxis the time scales for a transient
response and for adaptation to a signal are not well known.  We have run
the numerical simulations for the two-compartment model and the
reaction-diffusion model for several values of $\lambda$.  The only
observation that can be made based on these simulations is that for small
values of $\lambda$ the initial transient of $A$ is longer in comparison
with the time it takes for $A$ to return to its baseline value.  Changing
$\lambda$ by three orders of magnitude did not appear to alter the
qualitative behavior of the system.  Further analytical treatment of the
problem is necessary to understand the role of and the constraints on the
size of $\lambda$.

Finally, we must discuss the limitations of the adaptation-diffusion
model.  So far we have not linked the activated or the modified substance
of the model to any particular components of the biochemical pathway,
because there are no clear candidates for such connections.  It is
possible that one might find chemical species that adapt to a certain
ligand concentration, but based on the current experimental data one
cannot say with certainty what it might be.  Related to this, our original
aim was to explain the results of the calcium and cAMP experiments on {\it
Xenopus} neurons in the framework of chemotactic sensing.  We have not
addressed this, although in Appendix \ref{switch} we propose an extension
of the current model that accounts for how changing the absolute calcium
concentration changes turning behavior in growth cones.

\clearpage
\section{Conclusions and further direction}

We presented two mathematical models that attempted to describe
chemotactic movement of growth cones in response to a netrin-1 gradient.  
In our first model we assumed that the signal transduction pathways would
not be able to adapt to a constant signal.  In addition, we aimed at
basing our model on experimentally observed phenomena, in particular, the
way turning response is determined by cytosolic cAMP and $Ca^{2+}$
concentrations.  This model described a cAMP switch which is very
sensitive at a certain threshold concentration of netrin-1.  
Mathematically, the threshold value is a bifurcation parameter.  The
system goes through a bifurcation as the two stable steady states
separated by an unstable steady state change to one stable steady state.
This model did not successfully explain how gradient sensing is possible 
in concentration ranges that do not include the threshold concentration.  
This is a severe limitation of the model, because, a growth cones 
moving past the threshold level concentration would permanently lock onto 
the same direction, even if the ligand gradient changed.  In addition, it 
is possible that the assumption of no adaptation is also flawed.   

The second model attempts to satisfy most of the criteria set for
chemotactic sensing, except internal signal amplification.  This model
demonstrably explains how a graded internal response develops in ligand
gradients, and how a uniform increase in the ligand concentration leads to
a transient internal response.  Although the model is theoretically more
sound than the first one, it is equally limited to that, because its lack
of connection to experimentally observable signaling pathways.  The
variables of the model may represent particular chemicals, or they might
represent many components of the signal transduction pathway which acts as
a unit on a time scale faster than what is considered in the model.

Much of the model builds on the understanding how perfect adaptation
occurs in bacterial chemotaxis.  However, there are important differences
between bacterial and growth cone chemotaxis which might influence the
mathematical analysis of the problem.  One such example is the known
separation of time scales between phosphorylation and methylation in
bacterial chemotaxis which allows simplifications in the mathematical
treatment of the problem.  The same simplifications may be incorrect in
the description of growth cone chemotaxis.

Clearly, this topic is open for further theoretical and experimental
research.  Experimental observations of Song \& Poo \cite{SP1} on
adaptation might be a promising starting point for further work.  Song \&
Poo believe that growth cones periodically lose and regain their
sensitivity to gradients.  Such adaptation (which is different from the
adaptation assumed in our second model) appears similar to an idea
proposed by Meinhardt \cite{Me}.  He states that in a reaction-diffusion
system, if the half-life of the inhibitor is shorter than that of the
activator, then oscillations occur.  The accumulation of the activator is
overtaken by the inhibitor, and for the period of time that only the
inhibitor is present, the cell is unable to respond to new gradients.  
However, in order to pursue a mathematical model based on this idea, more
biological data is necessary to formulate a hypothesis.

        
\newpage
\pagestyle{myheadings}
\markright{  \rm \normalsize CHAPTER 4. \hspace{0.5cm}
  MATHEMATICAL MODELS IN BIOLOGY}
\chapter{Endothelial cell deformation} \label{endothelial}
\thispagestyle{myheadings}

\section{Introduction}

The present chapter concerns a topic very distinct from gradient sensing,
the subject of the previous sections.  In the model developed in this
section, we investigate the mechanical effects of blood flow on endothelial
cells.  As opposed to previous models in which we aimed to understand
biochemical signal transduction pathways, here we investigate a
mechano-transduction pathway, i.e. a pathway which transmits physical
signals such as forces and deformations.  

Endothelial cells form a monolayer inside blood vessels, acting as a
boundary between the blood flow and the vessel walls.  The endothelial
layer is exposed to various mechanical stresses, such as pressure,
circumferential stretch, and tangential shearing forces due to blood flow.  
Endothelial cells respond in a wide variety of ways to these forces, and
their response is thought to protect the arterial system from potential
damages.  Evidence supporting this hypothesis is offered by experiments
which have demonstrated that the development of certain vascular diseases,
such as atherosclerosis, coincides with the failure of the proper responses
of the endothelium to flow.

There is an array of events that takes place in endothelial cells exposed
to flow, some of which are immediate upon the application of shear
stresses, and some, which occur on the time scale of many hours.  Among the
fast responses are activations of flow-sensitive ion channels and
activation of G-proteins.  Changes in gene expression are also observed,
although on these events happen on a slower time scale, and finally, cells
go through morphological changes after they have been exposed to shear
stress approximately a day.  In steady and in pulsatile flow (a
superposition of oscillatory and steady flow) endothelial cells tend to
elongate and align with the direction of the flow.  Such enormous changes
require cells to extensively rearrange their cytoskeletal structure.  
Although the cytoskeletal reorganization in response to flow induced shear
stress has been studied, it is still an open area of research.

The morphological changes in the cytoskeleton are well documented
(\cite{Ba} original source: \cite{H2, N}), but it is unknown how the cells
are able to sense shear stress, and once it is detected, how shear stress
is transmitted from the cell membrane to the cytoskeleton. There have been
quite a lot of previous theoretical and numerical investigations of the
cytoskeletal changes produced by shearing forces.  Our aim is to
incorporate the effects of flow on the cytoskeleton, but in addition, we
want to examine the viscoelastic behavior of other structures, such as the
nucleus, cell-cell adhesions, and focal adhesion sites.  However, as the
cytoskeleton is thought to be the main force-bearing structure of
endothelial cells, much of the deformation and mechanical response must
come from it.  For this reason, we summarize some previous results
regarding the reorganization of cytoskeleton without attempting to provide
an exhaustive review of theoretical work in this area.

Theoretical models focusing on how shear stress is transmitted often make
one of two assumptions.  They either propose that stress is transmitted by
producing deformations in filaments, or, that there exists an internal
mechanical tension inside cells independently of the external shearing
forces, and that instead of significant deformations, there is simply a
rotation or change in spacing in order to respond to stress.  The first
assumption is consistent with a model by Satcher \& Dewey, called an
"open-cell foam" model, the second assumption is used by Wang \& Ingber
(\cite{St}, original source: \cite{W}), and by Stamenovic et al.  
\cite{St} in so called "tensegrity" models as well as "cable net" models.

Satcher \& Dewey \cite{SD} investigate how shearing forces distort the
polymers of the cytoskeleton:  F-actin, intermediate filaments and
microtubules.  The thesis of their work is that F-actin stabilizes cells by
decreasing deformability.  The role of stress fibers, which are
microfilaments connected into bundles, is also considered.  Satcher \&
Dewey model the F-actin filaments as open lattices.  (The structure of the
actin filaments is similar to other material, such as glass foams, etc,
hence the name "open-cell foam model".)  The advantage of such a model
formulation is that the Young's modulus, the measure of the material's
ability to resist distortion, can be expressed in terms of filament
properties, instead of having to find the density and moment of inertia for
the entire F-actin cytoskeleton.  With the open-lattice model properties of
the cytoskeleton (shear modulus, which is the coefficient of the rigidity
of a material, and the modulus of elastic deformation) are computed.  The
obtained values are the same order of magnitude as experimental data.  The
article also shows that although stress fibers increase the rigidity of the
actin network, their elastic deformation is too small to effect the
network, therefore their role is unclear based on the model.

Stamenovic \& Coughlin \cite{SC} compare predictions of three different
types of models of the cytoskeleton.  The first type is the "open-cell
foam" model, the second is the "cable net" model, and the third is the
"tensegrity" model.  The latter two share the assumption that pre-existing
tensions are present in the cytoskeleton even in absence of external
stress.  In the cable net model the actin filaments are represented by
elastic cables that are pulled tight by various forces generated inside the
cell.  The tensegrity model is similar: it assumes cables connected to
rigid beams or "struts", and the network of struts and cables represents
the cytoskeleton.  The basis of comparison of the three models is the
model's estimate of Young's modulus for the cytoskeleton.  This question
if further complicated by the large range ($10^0-10^5 Pa$) given for the
Young's modulus attained by different experimental techniques, such as
magnetic bead microrheometry, magnetic twisting cytometry,micropipette
aspiration and atomic force microscopy.  The paper discusses possible ways
experimental procedures might bias the obtained values.

The Young's modulus predicted by the open-cell foam model were much higher
($10^3-10^4$ Pa) than the Young's modulus predicted by the cable net and
tensegrity models ($10-10^2$ Pa).  The article concludes that different
models may be appropriate for modeling cytoskeletal changes under different
conditions.  For example, the cable net and tensegrity models may be
applicable in low stress whereas under large stress the cytoskeletal
filaments bend, and the open-cell foam model provides a better description.  
The models compared by Stamenovic \& Coughlin describe only the elastic
properties of the cytoskeleton.  However, Satcher \& Dewey in \cite{SD}
note that the open-cell foam model lends itself easily to the description
of viscoelastic properties, if the open space between the lattices is
assumed to be filled with a liquid.

Now we return to the larger goal of understanding how shear stress sensing
and transduction occur in endothelial cells.  One hypothesis is that shear
stress deforms flow-sensitive parts of the membrane, such as certain
ion-channels or receptors.  The deformation of these structures could then
be immediately transmitted to cytoskeletal elements connected to them.  
Such a mechano-transduction pathway may complement other, biochemical
signaling pathways.  A realistic model of the mechanical signaling pathway
would allow quantitative tests of whether deformations and stresses
generated in the cell would provide a sufficient signal.  A previous model
by Barakat \cite{Ba} shows that flow sensors, modeled by a viscoelastic
body, respond differently in oscillatory, pulsatile and steady flows, and
the differences in the response could provide the necessary signal for
downstream components of the pathway.

Our current work extends this model, and we represent other parts of the
endothelial cell, namely actin filaments, the nucleus, and transmembrane
proteins as viscoelastic Kelvin bodies as well.  The final goal of this
line of investigation is the development of a complex network of
viscoelastic bodies where each body is described by experimentally obtained
parameters.  This dissertation is only concerned with two small networks,
one consisting four, and the other one of seven viscoelastic bodies.  Our
work derives the equations for single bodies in series and in parallel,
however, further work is necessary to describe more complicated connections
(for example, n bodies connected in series which are connected to in
parallel to n bodies in series again).  Numerical simulations of the model
are generated to test the model's dependence on parameter values, and
numerical solutions of the deformation due to shear stress of the two
simple model networks are also given.  We discuss the implications of our
results regarding endothelial cell behavior.

\clearpage
\section{Mathematical Model}

In this section, first we derive the equations describing the deformations
of coupled Kelvin bodies and discuss solutions to the equations.  Next, we
describe the networks we use to model endothelial cells, and finally we
discuss how the parameter values were obtained for parts of the network.  

\subsection{Kelvin bodies in series}

Our goal is to develop a mathematical framework to describe the
deformations of the cell surface and intracellular structure within
endothelial cells with respect to steady and oscillatory flow.  Kelvin
bodies are the most general models for viscoelastic materials, and they
have frequently been used to model how the deformation of cell tissues
depends on the forcing \cite{Ba}, \cite{Fu}.  In addition, experimental
data is also available which describes the viscoelastic properties of
the cell nucleus \cite{GT}, cytoskeletal structures \cite{Sa} and
transmembrane proteins \cite{BZ} in terms of the parameters of a Kelvin
body.  This makes the Kelvin body a very effective tool for theoretical
modeling of endothelial cell deformations.

First, we give the equation relating the deformation and the force exerted
by the flow for one Kelvin body as derived by Fung \cite{Fu}. (The
derivation is very similar to the case where two Kelvin bodies are
connected in parallel which is shown in detail below.) 

\begin{figure}[h]
\centering
\includegraphics[width=0.3\textwidth]{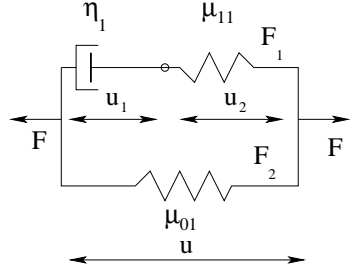}
\caption[Diagram to illustrate one viscoelastic  Kelvin body.]{Diagram to 
illustrate one 
Kelvin body.  The
parameters to characterize the body are as follows.  Dashpot viscosity: 
$\eta_1$, spring constant in upper branch: $\mu_{11}$, spring constant in 
lower branch: $\mu_{01}$.}  
\label{one_body}
\end{figure}

The deformation $u(t)$ of one Kelvin body as a function of a given
forcing, $F(t)$, is obtained by solving a first order linear equation
(Fung, \cite{Fu}):
\begin{eqnarray}
F + \frac {\eta_1} {\mu_{11}} \dot{F} = \mu_{01} u + \eta_1 (1 + \frac 
{\mu_{01}} {\mu_{11}} ) \dot{u} \label{1_ode} \\
u(0) = \frac {F(0)} {\mu_{01}+ \mu_{11}}  
\label{1_ode_ic}
\end{eqnarray}

The solution in steady flow, $F=F_0$ is (Barakat, \cite{Ba}): 
\begin{eqnarray} u(t) = \frac {F_0} {\mu_0} \Big[1 - \Big( 1-\frac
{\tau_{\epsilon}} {\tau_{\sigma}}\Big) \exp \Big(\frac {-t} {\tau_{\sigma}}
\Big)\Big] \label{u_steady} \end{eqnarray} where $$\tau_{\sigma} = \frac
{\eta_1} {\mu_{01}} \Big( 1 + \frac {\mu_{01}} {\mu_{11}} \Big)$$ and
$$\tau_{\epsilon}= \frac {\eta_1} {\mu_{11}}.$$ $\tau_{\epsilon}$ is the
relaxation time for a Kelvin body under constant strain (i.e. the stretch
per unit length), and $\tau_{\sigma}$ represents the relaxation time for
constant stress (i.e.  the force per unit area).  The deformation of the
two springs in the Kelvin body is instantaneous, as it is seen from the
initial condition for the deformation, \ref{1_ode_ic} while the dashpot
slowly creeps to the steady state of its deformation.  Examining the
expression for the deformation in steady flow (Eqn. \ref{u_steady}), it is
clear that large coefficients of viscosity lead to longer relaxation times
under constant stress while large spring coefficients result in decreasing
relaxation time.  The deformation in oscillatory and pulsatile flow is also
given by Barakat, \cite{Ba}.

\begin{figure}[h]
\centering
\includegraphics[width=1\textwidth]{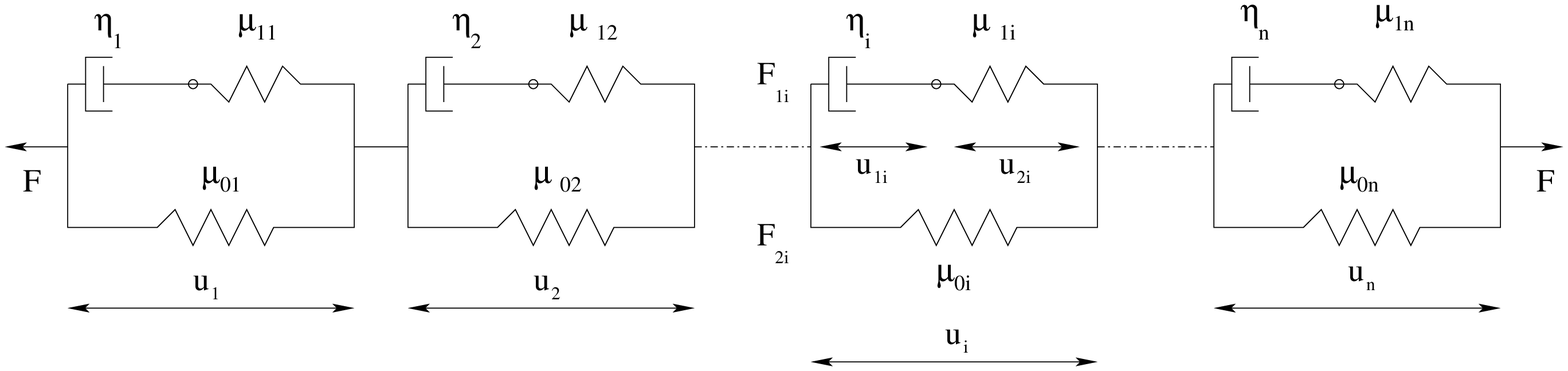}
\caption[n Kelvin bodies in series.]{Diagram to 
illustrate n Kelvin bodies in series.  The
parameters to characterize the ith body are as follows.  Dashpot 
viscosity: $\eta_{i}$, spring constant in upper branch:
$\mu_{1i}$, spring constant in lower branch: $\mu_{0i}$.}
\label{nbodies_series}
\end{figure}

When Kelvin bodies are coupled in series, the deformation for each body
can be found individually, so in essence it is exactly like the one-body
case.  The force acting on each Kelvin body is the same, so the overall 
deformation can be calculated as the sum of deformations, $u_i$ for 
i=1,...,n where each deformation is the function of the same forcing 
function, $F$:
\begin{eqnarray*}    
F + \frac {\eta_i} {\mu_{1i}} \dot{F} = \mu_{0i} u + \eta_i (1 + \frac   
{\mu_{0i}} {\mu_{1i}} ) \dot{u}_i \\
u_i(0) = \frac {F(0)} {\mu_{0i}+ \mu_{1i}} \\
u(t) = \sum_{i=1}^n u_i(t)
\end{eqnarray*}

\subsection{Kelvin  bodies in parallel}

\begin{figure}[h]
\centering
\includegraphics[width=0.5\textwidth]{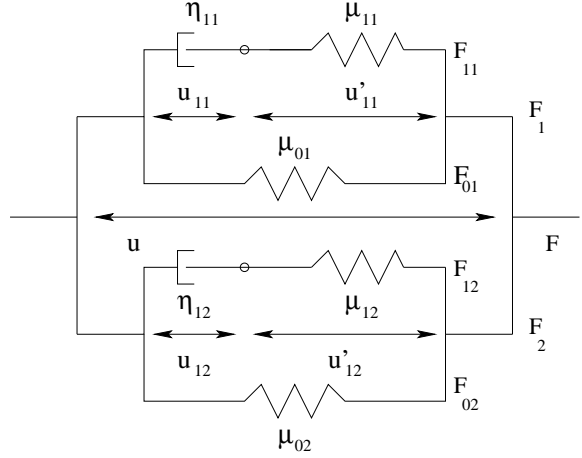}
\caption[Two Kelvin bodies in parallel.]{Diagram to illustrate two Kelvin 
bodies in parallel.  The
parameters to characterize the bodies are as follows.  Upper Kelvin body:
dashpot viscosity: $\eta_{11}$, spring constant in upper branch:
$\mu_{11}$, spring constant in lower branch: $\mu_{01}$.  Lower Kelvin
body: dashpot viscosity: $\eta_{12}$, spring constant in upper branch:
$\mu_{12}$, spring constant in lower branch: $\mu_{02}$.}
\label{twobodies.}
\end{figure}

We want to find the deformation for Kelvin bodies in more complicated
networks.  We begin by taking two Kelvin bodies coupled in parallel. Each
of the Kelvin bodies is described by three parameters: the viscosity of
the dashpot, and the two spring constants, as Figure \ref{twobodies.}
shows below.

In both the upper and the lower Kelvin body there are some relationships
that must hold, namely, the total deformation, $u$ must be a sum of the
deformations of the dashpot and the spring in the upper branch, and this
deformation is the same as the deformation of the spring in the lower
branch.  We must also note that the deformation of the upper body and the
deformation of the lower body must be identical.  These relationships give
us the following equations.
$$ u_{11} + u'_{11} = u $$
$$ u_{12} + u'_{12} = u $$

Another observation is that the total force, $F$ of the two bodies splits
into the force in the upper body, $F_1$ and the force in the lower body,
$F_2$.  The total force in the upper body is also given as a sum of the
force in the upper branch and the lower branch, and similarly for the
lower body.  This is described by the equations:
$$ F_1 + F_2 = F $$
$$ F_{01} + F_{11} = F_1 $$
$$ F_{02} + F_{12} = F_2 $$

Now let us consider the upper body only.  The same force in the upper
branch, $F_{11}$ is transmitted from the dashpot to the spring.  The force
acting on the spring is proportional to the deformation it produces, and
the force acting on the dashpot is proportional to the velocity of the
dashpot.  Using the variables of the diagram we can write this as:
$$ F_{11} = \eta_{11} \dot{u}_{11} = \mu_{11} u'_{11}.$$

Here $\dot{u}_{11}$ denotes the derivative of $u_{11}$, the deformation of
the dashpot, and $u'_{11}$ is the deformation of the spring.  The force in
the lower branch of this body acts entirely on the dashpot, and here the
deformation is going to be the sum of the deformation due to the dashpot
plus the deformation due to the spring in the upper branch, therefore
$$ F_{01} = \mu_{01} (u_{11}+ u'_{11}).$$

Similarly, we can write down the corresponding equations for the lower
body as well.
$$ F_{12} = \eta_{12} \dot{u}_{12} = \mu_{12} u'_{12}$$
$$ F_{02} = \mu_{02} (u_{12}+ u'_{12})$$

Using the nine equations above, one can derive two differential
equations that describe the the deformation of the two coupled Kelvin
bodies as a function of time and the force acting on the bodies, with
the appropriate boundary conditions.  The derivation is similar to
Fung's \cite{Fu}.  In the equations below we assume that $F(t)$ is
given, therefore we can also find $\dot{F(t)}.$ For the particular
forcing functions we are interested in, the derivative always exists,
and it is continuous.
\begin{eqnarray}
F_1 + \frac {\eta_{11}} {\mu_{11}} \dot{F}_1 = \mu_{01} u + \eta_{11}(1 +
\frac {\mu_{01}} {\mu_{11}}) \dot{u}  \nonumber \\
u(0) = \frac {F_1(0)} {\mu_{01}+\mu_{11}}
\end{eqnarray}

\begin{eqnarray}
F_2 + \frac {\eta_{12}} {\mu_{12}} \dot{F}_2 = \mu_{02} u + \eta_{12}(1 + 
\frac {\mu_{02}} {\mu_{12}}) \dot{u} \nonumber \\
u(0) = \frac {F_2(0)} {\mu_{02}+\mu_{12}}
\end{eqnarray}

Now we can use the fact that $F_1$ and $F_2$ sum to $F$ to substitute
$a(t) F(t) = F_1(t) $ and $ 1 - a(t) F(t) = F_2(t) $.  $a(t)$ is the
coefficient of force splitting, and it is an unknown function of time.
\begin{eqnarray}
a F + \frac {\eta_{11}} {\mu_{11}} \dot{ (a F)} = \mu_{01} u +
\eta_{11}(1 + \frac {\mu_{01}} {\mu_{11}}) \dot{u} \nonumber \\
u(0) = \frac {a(0)F(0)} {\mu_{01}+\mu_{11}}
\label{2bodeq1}
\end{eqnarray}

\begin{eqnarray}
(1-a)F + \frac {\eta_{12}} {\mu_{12}} \dot{((1-a)F)} = \mu_{02} u +
\eta_{12}(1 + \frac {\mu_{02}} {\mu_{12}}) \dot{u} \nonumber \\
u(0) = \frac {(1-a(0))F(0)} {\mu_{02}+\mu_{12}}
\label{2bodeq2}
\end{eqnarray}

Now we have a system of two differential equations with initial conditions
and two unknown functions, $u(t)$ and $F(t)a(t)$, and we would like to
solve for them.  Solving the equations for $u(t)$ and $a(t)$ is not
practical, for two reasons.  For particular choices of the flow, for
example, for oscillatory flow $F(t)=F_0 \cos(\omega t)$, the matrix A,
defined below, would depend on the forcing function, and $A^{-1}$ would
become singular periodically when $\cos(\omega t) = 0 $.  Also, matrices A
and D (also defined below)  would both depend on time, and this could
considerably slow down the computations.  In order to avoid these problems,
we compute $u(t)$ and $F(t)a(t)$.  $F(t)$ is a know function of time, so it
is always possible to find $a(t)$, if necessary. Now rearranging equations
\ref{2bodeq1} and \ref{2bodeq2} gives us the following.
\begin{eqnarray*}
\left[ \begin{array}{cc}
\eta_{11}(1 + \frac {\mu_{01}} {\mu_{11}}) & - \frac {\eta_{11}}
{\mu_{11}}\\
\eta_{12}(1 + \frac {\mu_{02}} {\mu_{12}}) & \frac {\eta_{12}} {\mu_{12}}
\end{array} \right]
\left[ \begin{array}{c}
\dot{u} \\
\dot{(aF)}
\end{array} \right] =
\left[ \begin{array}{cc}
-\mu_{01} & 1 \\  -\mu_{02} & -1
\end{array} \right]
\left[ \begin{array}{c}
u  \\ (aF)
\end{array} \right]
\end{eqnarray*}
\begin{eqnarray*}
+\left[ \begin{array}{c}
0  \\  F + \frac {\eta_{12}} {\mu_{12}} \dot{F}
\end{array} \right]
\end{eqnarray*}

\begin{eqnarray}
\left[ \begin{array}{c}
u(0) \\ a(0)F(0)
\end{array} \right] =
\left[ \begin{array}{c}
\frac {F(0)} {\mu_{01} + \mu_{11} + \mu_{02} + \mu_{12}}\\
\frac {F(0)(\mu_{01}+\mu_{11})} {\mu_{01} + \mu_{11} + \mu_{02} +
\mu_{12}} \end{array} \right]
\end{eqnarray}

We can let
\begin{eqnarray*}
A = \left[ \begin{array}{c c}
\eta_{11}(1 + \frac {\mu_{01}} {\mu_{11}}) & - \frac {\eta_{11}}
{\mu_{11}} \\
\eta_{12}(1 + \frac {\mu_{02}} {\mu_{12}}) & \frac {\eta_{12}} {\mu_{12}}
\end{array} \right]
\end{eqnarray*}

\begin{eqnarray*}  
D = \left[ \begin{array}{c c}
-\mu_{01} & 1  \\ -\mu_{02} & -1
\end{array} \right]
\end{eqnarray*}

\begin{eqnarray*}
\vec{c} = \left[ \begin{array}{c}
0  \\F + \frac {\eta_{12}} {\mu_{12}} \dot{F}
\end{array} \right]
\end{eqnarray*}

\begin{eqnarray*}
\vec{u} = \left[ \begin{array}{c}
u(t) \\ a(t)F(t)
\end{array} \right]
\end{eqnarray*}

Now we can write the system of equations as $$ A \frac{d \vec{u}} {dt} = D
\vec{u} + \vec{c}. $$

In order to compute the solution for $\vec{u}$, we can express the system
of differential equations in the form
\begin{eqnarray}
\frac {d \vec{u}} {dt} = A^{-1} D \vec{u} + A^{-1} \vec{c} 
\label{maineqs1} \\
\vec{u}(0) = \vec{u}_0
\label{maineqs2}
\end{eqnarray} 

This is the equation with the appropriate initial conditions that
describes the dynamics of two Kelvin bodies coupled in parallel.

\begin{figure}[h]
\centering
\includegraphics[width=0.4\textwidth]{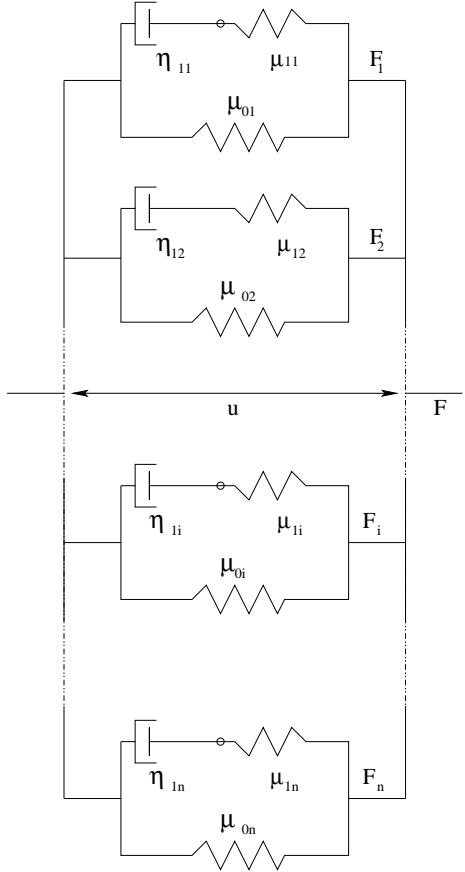}
\caption[n Kelvin bodies in parallel.]{Diagram to illustrate n Kelvin 
bodies in parallel. The parameters
to characterize the bodies are as follows.  ith body (for i=1,...,n):
dashpot viscosity: $\eta_{1i}$, spring constant in upper branch: 
$\mu_{1i}$, spring constant in lower branch: $\mu_{0i}$.}
\label{nbodies.}
\end{figure}

Now we can look at a generalization of the two body problem to deriving the
differential equations governing n Kelvin bodies coupled in parallel. Let
us start again with the diagram of the bodies with their appropriate
parameters, shown in Figure \ref{nbodies.}.  As before, the sum of the
forces in the branches has to be the total force, F.

$$ \sum_{i=1}^n F_i = F$$

In general, we do not know how the forces split into these branches
because this depends on the particular parameter values of the Kelvin
bodies.  Therefore we can call the force splitting coefficient for the ith
branch to be $a_i(t)$, and the force in this branch $a_i(t)F(t)$.  Using
the above relationship we get the following:
$$ (1 - \sum_{i=1}^{n-1} a_i) F = F_n $$

The equation for the the ith body for $ i = 1, ..., n-1 $ is:
\begin{eqnarray}
a_i F + \frac {\eta_{1i}} {\mu_{1i}} \dot{ (a_i F)} = \mu_{0i} u +
\eta_{1i}(1 + \frac {\mu_{0i}} {\mu_{1i}}) \dot{u} \nonumber \\
u(0) = \frac {a_i(0)F(0)} {\mu_{0i}+\mu_{1i}}
\end{eqnarray}  

And for the nth body we have a similar expression:
\begin{eqnarray}
(1 - \sum_{i=1}^{n-1} a_i )F + \frac {\eta_{1n}} {\mu_{1n}} \cdot \frac
{d} {dt}((1 - \sum_{i=1}^{n-1} a_i ) F)  = \mu_{0n} u +
\eta_{1n}(1 + \frac {\mu_{0n}} {\mu_{1n}}) \dot{u} \nonumber \\
u(0) = \frac {(1 - \sum_{i=1}^{n-1} a_i(0))F(0)} {\mu_{0n}+\mu_{1n}}
\end{eqnarray}

Now we must rearrange the equations so we are solving for $u(t)$ and
$a_i(t)F(t)$ again.  When we rearrange differential equations we
get for $i=1,...,(n-1)$:
$$ \eta_{1i}(1 + \frac {\mu_{0i}} {\mu_{1i}}) \dot{u} - \frac {\eta_{1i}}
{\mu_{1i}} \dot{(a_i F)} = - \mu_{0i} u + a_i F $$

for n:
$$ \eta_{1n}(1 + \frac {\mu_{0n}} {\mu_{1n}}) \dot{u} + \frac {\eta_{1n}}
{\mu_{1n}} \sum_{i=1}^{n-1} \dot{(a_i F)} = -\mu_{0n} u - \sum_{i=1}^{n-1}
a_i F  + F + \frac {\eta_{1n}} {\mu{1n}} \dot{F}$$

We also need to find the appropriate expression for the initial
conditions. We have, first for $i=1,...,(n-1)$:
\begin{eqnarray}
u(0) (\mu_{0i}+\mu_{1i}) = a_i(0) F(0) 
\label{ic1}
\end{eqnarray}

and for i=n:
\begin{eqnarray}
u(0) = \frac {F(0)- F(0) \sum_{i=1}^{n-1}a_i(0)} {\mu_{0n}+\mu_{1n}}
\label{ic2}
\end{eqnarray}

$u(0)$, the initial deformation of all of the bodies is the same, so we
get n equations and n unknowns: $u(0)$ and $a_i(0)$ for i=1,...,n-1.  The
force splitting coefficient of the nth body is already determined from
this to be $a_n(0)=1 - \sum^{n-1}_{i=1} a_i(0)$.  We must rearrange the
equations \ref{ic1} and \ref{ic2} to solve for $u(0)$ and $a_i(0)$.
$$ u(0) (\mu_{0n}+\mu_{1n}) = F(0) - \sum_{i=1}^{n-1}a_i(0) F(0) $$
$$ u(0) (\mu_{0n}+\mu_{1n}) = F(0) - \sum_{i=1}^{n-1} u(0) (\mu_{0i} +
\mu{1i}) $$
$$ u(0) \sum_{i=1}^{n}(\mu_{0i}+\mu{1i}) = F(0) $$

Therefore the initial conditions are:
$$ u(0) = \frac {F(0)} {\sum_{i=1}^{n}(\mu_{0i} + \mu_{1i})} $$
$$ a_i(0) = \frac {F(0)(\mu_{0i} + \mu_{1i})} {\sum_{i=1}^{n}(\mu_{0i} +
\mu_{1i})} $$

Now we can look at the equations in matrix form:
\begin{eqnarray*}
A = \left[ \begin{array}{c c c c c c c}
\eta_{11}(1 + \frac {\mu_{01}} {\mu_{11}}) & - \frac {\eta_{11}}
{\mu_{11}} & 0 & \cdots & \cdots & \cdots &  0 \\
\eta_{12}(1 + \frac {\mu_{02}} {\mu_{12}}) & 0 & -\frac {\eta_{12}}
{\mu_{12}} & 0 &\cdots & \cdots &  0 \\
\vdots & \vdots &   & \ddots  &  & \ & \vdots \\  
\eta_{1i}(1 + \frac {\mu_{0i}} {\mu_{1i}}) & 0 &  & &  -\frac  
{\eta_{1i}} {\mu_{1i}} &  & 0 \\
\vdots & \vdots &   &  &  &\ddots  & \vdots \\
\eta_{1(n-1)}(1 + \frac {\mu_{0(n-1)}} {\mu_{1(n-01)}}) &  0  & \cdots  &
\cdots  &\cdots & 0 &  - \frac {\eta_{1(n-1)}} {\mu_{1(n-1)}} \\
\eta_{1n}(1 + \frac {\mu_{0n}} {\mu_{1n}}) & \frac {\eta_{1n}} {\mu_{1n}}
& \cdots & \cdots & \cdots & \cdots  &  \frac {\eta_{1n}} {\mu_{1n}}
\end{array} \right]
\end{eqnarray*}

\begin{eqnarray*}
D = \left[ \begin{array}{c c c c c c c}
-\mu_{01} & 1 & 0 & \cdots & \cdots & \cdots & 0 \\
-\mu_{02} & 0 & 1 & 0 & \cdots & \cdots &  0 \\
\vdots & \vdots & & \ddots  & &  & \vdots \\
-\mu_{0i} & 0 &  &  & 1 &  & 0 \\
\vdots & \vdots &  &  &  & \ddots & \vdots \\
-\mu_{0(n-1)} & 0 & \cdots & \cdots &\cdots  & \cdots & 1 \\
-\mu_{0n} & - 1 & \cdots & \cdots & \cdots & \cdots & -1
\end{array} \right]
\end{eqnarray*}

\begin{eqnarray*}
c = \left[ \begin{array}{c}
0  \\  \vdots \\ 0 \\
F + \frac {\eta_{12}} {\mu_{12}}
\dot{F}
\end{array} \right]
\end{eqnarray*} 

Just like before, we have the differential equation for $\vec{u}(t)$, an 
$n \times 1 $ vector whose entries are $u(t)$ and $a_i(t)F(t)$ for $ i = 
1, ..., (n-1)$ with the appropriate initial conditions:
\begin{eqnarray*}
\frac {d \vec{u}(t)} {dt} = A^{-1}D \vec{u} + A^{-1} c \\
\vec{u}(0) = \vec{u}_0
\end{eqnarray*}

This is the same linear equation as \ref{maineqs1} and \ref{maineqs2} with
the matrix $A^{-1}D$ and vector $A^{-1}c$ defined appropriately, and its
solution is given by
\begin{eqnarray} 
\vec{u} = (\vec{u}_0 + D^{-1} c ) e^{\Lambda t} - A^{-1}c 
\end{eqnarray} 

where $\Lambda$ is a diagonal matrix whose eigenvalues are the same as the
eigenvalues of the matrix $M=A^{-1} D$.  (Obtaining the solution to
equations \ref{maineqs1} and \ref{maineqs2} is similar to the derivation
shown in Appendix \ref{anal_sol}.)Let us first discuss how to find $D^{-1}
c$ and $A^{-1} c$, then turn to finding $\Lambda$.  If $D^{-1} c = y$,
then $c = D y$, in other words,
\begin{eqnarray*}
\left[ \begin{array}{c}
0 \\  \vdots \\  \vdots  \\  0 \\
F + \frac {\eta_{1n}} {\mu_{1n}}
\dot{F}
\end{array} \right]
=
\left[ \begin{array}{c c c c c}
-\mu_{01} & 1 & 0 & \cdots & 0 \\
-\mu_{02} & 0 & 1 &  & 0 \\
\vdots & \vdots  &  & \ddots  & \vdots \\
-\mu_{0(n-1)} & 0 &  \cdots &\cdots & 1 \\
-\mu_{0n} & - 1 & \cdots  & \cdots  & -1
\end{array} \right]
\left[ \begin{array}{c}
y_1 \\ y_2 \\ \vdots \\  y_{n-1} \\ y_n   
\end{array} \right]   
\end{eqnarray*}  

In the ith row we have $-\mu_{0i} y_1 + y_{i+1} = 0$, and this implies
that for i=1,...,n-1 $$y_{i+1} = y_1 \mu_{0i}.$$ In the nth row
we get $-\mu_{0n} y_1 - \sum^{n}_{i=2} y_i = F + \frac {\eta_{1n}}
{\mu_{1n}} \dot{F}$.  By substituting the expression for $y_{i+1}$ in
terms of $y_1$ into this

$$ y_1 ( -\mu_{0n} - \sum_{i=2}^{n} \mu_{0i-1} ) = F + \frac
{\eta_{1n}} {\mu_{1n}} \dot{F}. $$
$$ y_1 = \frac {F + \frac{\eta_{1n}}{\mu_{1n}} \dot{F}} { -\mu_{0n} -
\sum_{i=2}^{n} \mu_{0i-1} } $$
This gives an explicit formula for $y_1$ and based on our expression for
$y_i$ in terms of $y_1$ we can find all the other components of $y$.

We can use a similar argument to find $x=A^{-1}c$.  If we let $h_i =
\eta_{1i} (1 + \frac {\mu_{0i}} {\mu_{1i}} )$ and $d_i = - \frac 
{\eta_{1i}}
{\mu_{1i}}$ then we get that
$$x_1 = \frac {F + \frac {\eta_{1n}} {\mu_{1n}} \dot{F} } {h_n - d_n
\sum^{n-1}_{i=1} \frac {h_i} {d_i}} $$
$$ x_{i+1} = \frac {h_i x_1} {d_i}.$$

Now let us return to the matrix of eigenvalues, $\Lambda$.  We are looking
for the diagonal matrix whose entries $\lambda_1$, ..., $\lambda_n$
satisfy the equation $A^{-1}D x = \lambda_i x $ for i=1,...,n and for $x
\neq 0$. This is equivalent to $(D-\lambda_i A) x = 0$ which is called the
generalized eigenvalue problem.  It is easy to see that the matrix
$(D-\lambda_i A)$ has the same very nice and sparse structure that both A
nd D have, but solving the generalized eigenvalue problem leads to having
to find the roots of an nth degree polynomial.  In the most general case
this can only be solved numerically, but some special cases of the problem
can be solved analytically (for example, if all n Kelvin bodies have the
same parameter values).  Other methods can also be applied to solve the
original system of linear differential equations, for example Laplace
transform methods, but they lead to the same problem of having to find
roots of an n-degree polynomial.

In this dissertation we only use networks with two Kelvin bodies in
parallel, and in these cases the deformation can be found analytically as
well as numerically, because it only requires finding solutions to a
quadratic equation. Extensions of the model to a large number of bodies in
parallel would have to be done numerically, and even analytical solutions
can only be found by numerically computing the roots of an nth-degree
polynomial.

\subsection{Model networks}

The aim of our model is to gain further understanding of how flow over the
surface of endothelial cells leads to regulation of the cell shape.  It is
known that the cytoskeletal structure is re-organized in a period of
approximately a day, but there is no clear evidence to what extent  
biochemical events, and to what extent purely mechanical processes
contribute to this.  We want to examine how the flow-induced shear stress
which deforms flow sensors is transmitted through the cytoskeleton to the
nucleus and to other transmembrane proteins such as ion channels and
attachments to the substrate. We want to use networks of coupled Kelvin
bodies to model endothelial cells, and we want to investigate how
deformations of the individual parts will contribute to the overall
deformation of the cell.

The cytoskeleton consists of three types of polymers: actin filaments,
microtubules and intermediate filaments each of which deform if shear
stress is applied to the cell \cite{SD}.  Qualitatively, actin filaments
rupture at relatively low strain, but actin can be rapidly recycled and
filaments re-formed as it is required for cell motility, among its other
functions. Below a critical strain actin networks show the greatest
rigidity \cite{J}. Microtubule networks can withstand very high strain,
and the greatest deformability.  This is consistent with their role as
structural support of the actin filaments \cite{J}.  Vimentin networks,
which mostly make up the intermediate filaments, tend to be less rigid at
low shear strain, but they harden at high strains.  These responses make
them ideal for the support of nucleus \cite{J}.  Because of this, we
expect very significant differences in the responses of the actin
cytoskeleton and the nucleus of the cell.

\begin{figure}[h!]
\centering
\includegraphics[width=0.25\textwidth]{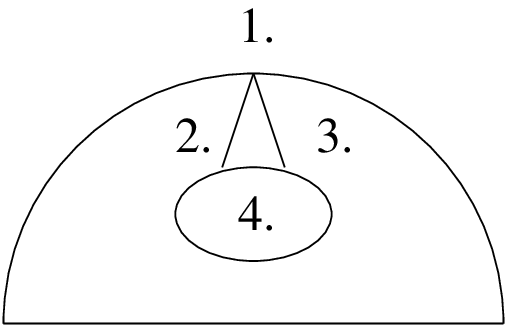}\quad
\includegraphics[width=0.45\textwidth]{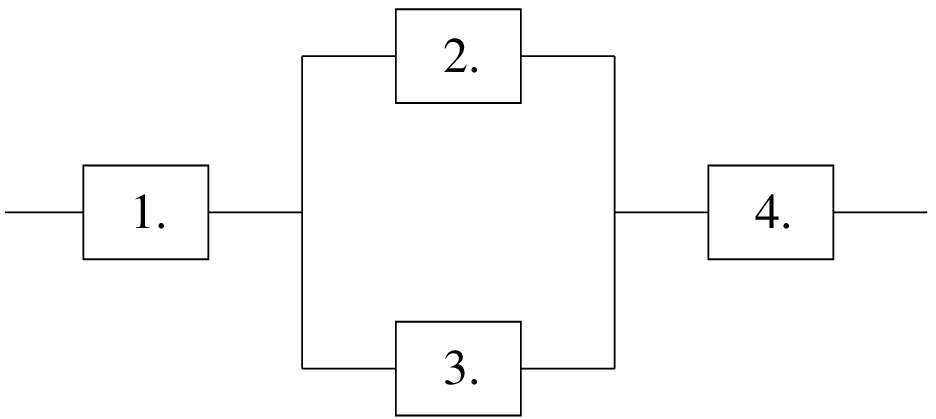}
\caption[Network I: model of an endothelial cell.]{Model I. of an 
endothelial 
cell.}
\label{Network I.}
\end{figure}

Based on the above, we choose a simple network to model an endothelial 
cell.  This model is shown in Figure \ref{Network I.}. The flow sensor,
body 1, is attached to the actin cytoskeleton which is represented by
bodies 2 and 3 in parallel.  The actin cytoskeleton then attaches to the
nucleus, body 4 in the diagram.

Next, we can elaborate on our initial diagram and add the connections  
between the nucleus and the attachments to the substrate.  Part of this
second model is the same as the first network, but now the nucleus (body
4) is further connected to actin bundles (bodies 5 and 6) that end at
transmembrane proteins (body 7).

\begin{figure}[h!]
\centering
\includegraphics[width=0.25\textwidth]{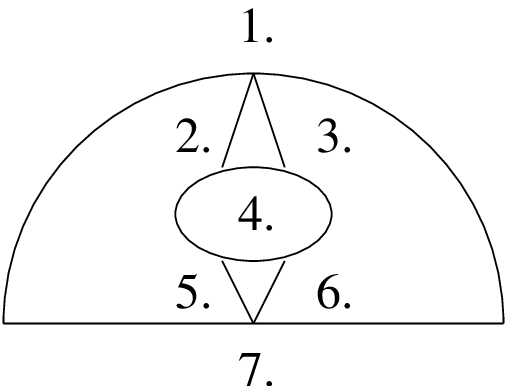}\quad
\includegraphics[width=0.7\textwidth]{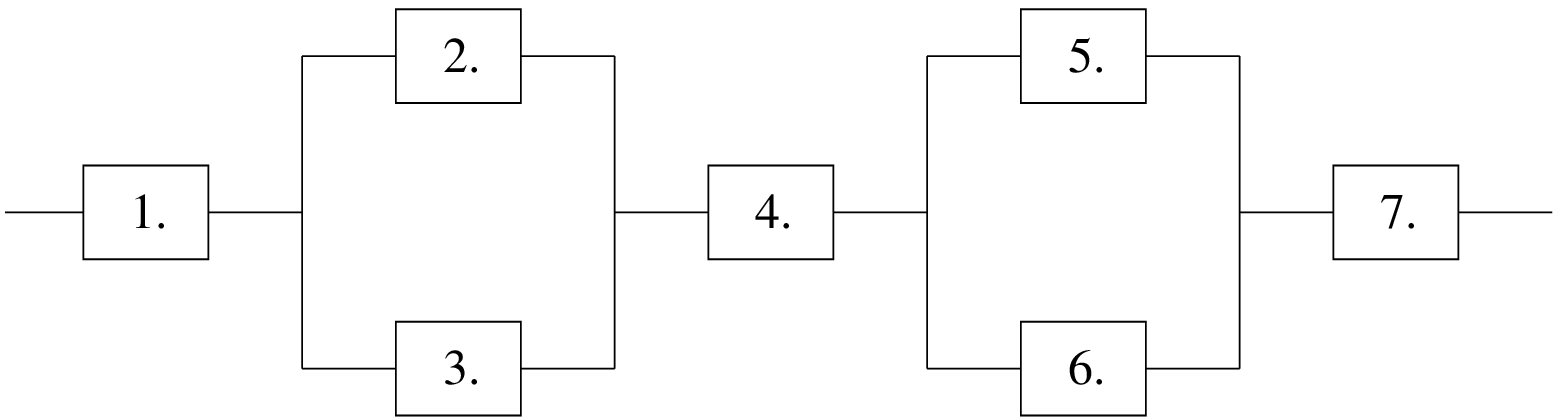}
\caption[Network II: model of an endothelial cell.]{Model II. of an 
endothelial cell.}
\label{Network II.}
\end{figure}

\subsection{Parameter values}
\label{parameter_calc}

Appropriate parameter values must be chosen for all of the bodies.  The
parameter values for the actin filaments are taken from Sato et al.
\cite{Sa} who measure viscoelastic properties of endothelial cells with 
a
micropipette technique.  Guilak et al. give the parameter values for the
nucleus based on a study also with a micropipette aspiration, and they
conclude that the nucleus is about 3-4 stiffer and approximately twice as
viscous as the cytoplasm \cite{GT}.  Finally, the parameter values for
transmembrane proteins is found by Bausch et al. \cite{BZ} who use the
novel technique of magnetic bead microrheometry.  All the parameter values
are summarized in Table \ref{parameter_values}.

\begin{table}[h!]
\centering
\begin{tabular}{|c|c|c|c|c|}
\hline
 & $\eta_1$ (Pa s) & $\mu_{01}$ (Pa) & $\mu_{11}$ (Pa) & Ref. \\
\hline
\hline
Actin filaments & 5000 & 50 & 100 & \cite{Sa} \\
\hline
Nucleus & 10 000 & 200 & 400 & \cite{GT} \\
\hline
Transmembrane proteins & 7.5 & 100 & 200 & \cite{BZ} \\
\hline
\end{tabular}
\caption{Parameter values for the endothelial cell models}
\label{parameter_values}
\end{table}

The parameter values for transmembrane proteins were given in the original
paper in units of Pa m for the spring constants and Pa s m for the dashpot
viscosity, and had to be converted to Pa and Pa s, respectively.  These   
calculations are given below.

Sato et al. give the formula for the deformation as $$ L(t) = \frac {2 a
\Delta p} {\pi \mu_{01}} (1 - \frac {\mu{11}} {\mu_{01} + \mu_{11}} e^{-
\frac {t} {\tau}}).$$ The dimensions are as follows: $[a] =$ m, $ [\Delta
p]=$ Pa and $[\mu_{01}, \mu_{11}] = $ Pa.  The formula for the deformation
in Bausch et al. is $$ L(t)  = \frac {F} {\mu_{01}} (1 - \frac {\mu_{11}}
{\mu_{01} + \mu_{11}} e^{- \frac {t} {\tau}})$$ with dimensions $[F]=$ N,
$[\mu_{01}, \mu_{11}]=$ Pa m = N/m.  In order to compare the spring
constants given in the two papers, we must have them in the same
dimensions.  First we will find what the applied force is in Sato et al.,
then we use the expression for the initial deformation and final
deformation to determine the spring constants.  The force is given by $ F
= \Delta p \pi a^2 \simeq 2500$ pN. The initial deformation, $$L_0 (t)=
\frac {F} {\mu_{01}} (1 - \frac {\mu_{11}} {\mu_{01} + \mu_{11}}) = \frac
{F}{\mu_{01}+ \mu_{11}}$$ and the deformation at steady state is $$L_s(t)
= \frac {F} {\mu_{01}}.$$ Using the data in Sato et al. this gives us
$\mu_{01} = 6.35 \times 10 ^{-4} $ Pa m and $\mu_{11} = 9.38 \times
10^{-4}$ Pa m. The same calculations with the data in Bausch et al. leads
to $\mu_{01} = 1.25 \times 10^{-3}$ Pa m and $\mu_{11} = 1.61 \times
10^{-3}$ Pa m. Finally, we must obtain the viscosity of the dashpot,
$\eta_1.$ $$ \eta_1 = \frac {\tau \mu_{01} \mu_{11}} {\mu_{01} +
\mu_{11}}.$$ In Sato, $\eta_1 = 4.125 \times 10^{-2}$ Pa m s and in Bausch
$\eta_1 = 6.33 \times 10^{-5}$ Pa m s. Now we can compare the parameter
values.  $\mu_{01}$ and $\mu_{02}$ for the nucleus is approximately twice
the value of $\mu_{01}$ and $\mu_{11}$, respectively, for actin filaments,
and the viscosity in the nucleus is approximately $1.5 \times 10^{-3}$ the
dashpot viscosity of actin filaments.  Based on this we arrive at:
$$\mu_{01} = 2 (50) = 100 \: Pa $$ $$ \mu_{11} = 2(100) = 200 \: Pa $$ $$
\eta_1 = (1.5 \times 10^{-3}) (5000)  = 7.5 \: Pa \: s.$$

\clearpage
\section{Results}

This section contains two sets of numerical simulations.  The first set of
simulations examines the relationships between the parameters of the
Kelvin bodies and the deformation, $u(t)$ and the force splitting,
$a(t)F(t)$ in a the two-body problem.  Both steady and oscillatory flow
are considered.

The second set of simulations takes the four-body and the seven-body model
networks with realistic parameter values and finds the temporal dynamics
of the deformation, and the force splitting coefficients.  

The differential equations are solved with a four-stage fourth-order 
Runge-Kutta method.  The Matlab code used to solve the equations is 
presented in Appendix \ref{num_sol}.  The time step chosen for the 
simulations is $h=0.1$.  This method is of order 4, so the error is  
$O(h^5)= 0.00001$. 

\subsection{Parameter sensitivity analysis}

We investigate the behavior of two Kelvin bodies coupled in parallel.  We
focus on three questions: \begin{itemize} \item temporal dynamics of the
system with one parameter value changed over three orders of magnitude;
\item the steady state behavior of the system as a parameter value is
changed over several orders of magnitude;  \item the behavior of the
system when all three parameters are changed over three orders of
magnitude \end{itemize} This analysis allows us to understand how the
parameter values determine the viscoelastic properties of a material which
gives us intuition into the behavior of the more complicated model
networks involving materials of different properties.

\begin{figure}[h!]
\centerline{\includegraphics[width=0.7\textwidth]{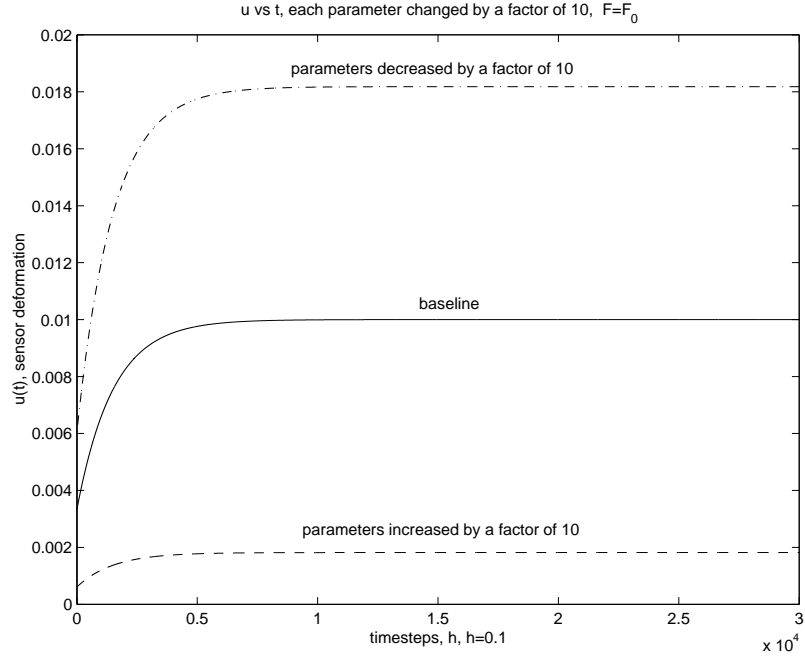}}
\caption[Deformation with all parameters changed.  Steady flow.]
{Deformation when every parameter of body 2 is changed by a  
factor of 10.  Steady flow.}
\label{factor10_u_s.ps}
\end{figure}

\begin{figure}[h!]
\centerline{\includegraphics[width=0.7\textwidth]{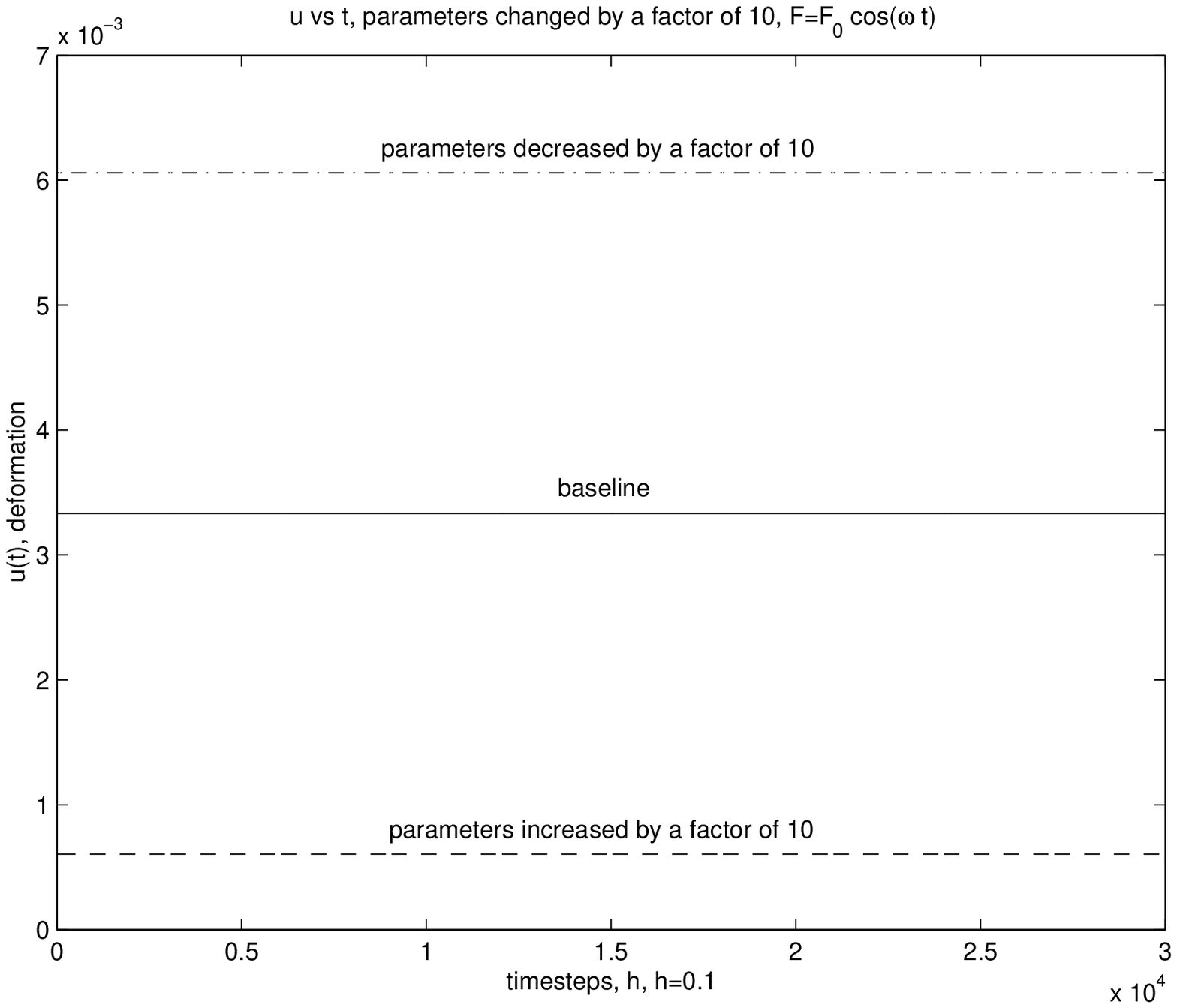}}
\caption[Deformation with all parameters changed.  Oscillatory flow.]
{Deformation when every parameter of body 2 is changed by a
factor of 10. Oscillatory flow.}
\label{factor10_u_o.ps}
\end{figure}

\begin{figure}[h!]
\centerline{\includegraphics[width=0.7\textwidth]{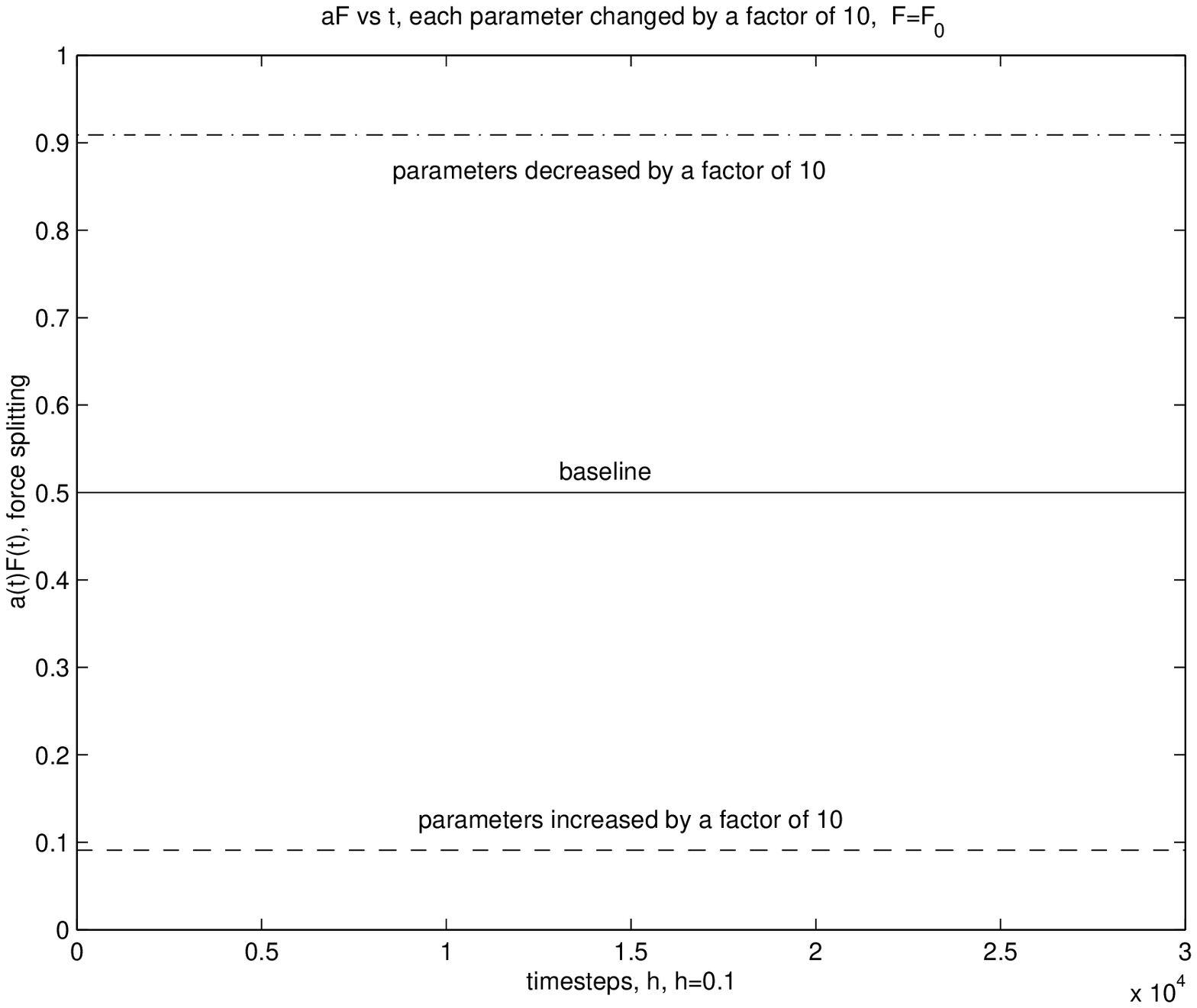}}
\caption [Force splitting with all parameters changed.  Steady
flow.]{Force splitting when every parameter of body 2 is changed by a  
factor of 10.  Steady flow.} \label{factor10_aF_s.ps} \end{figure}

\begin{figure}[h!]
\centerline{\includegraphics[width=0.7\textwidth]{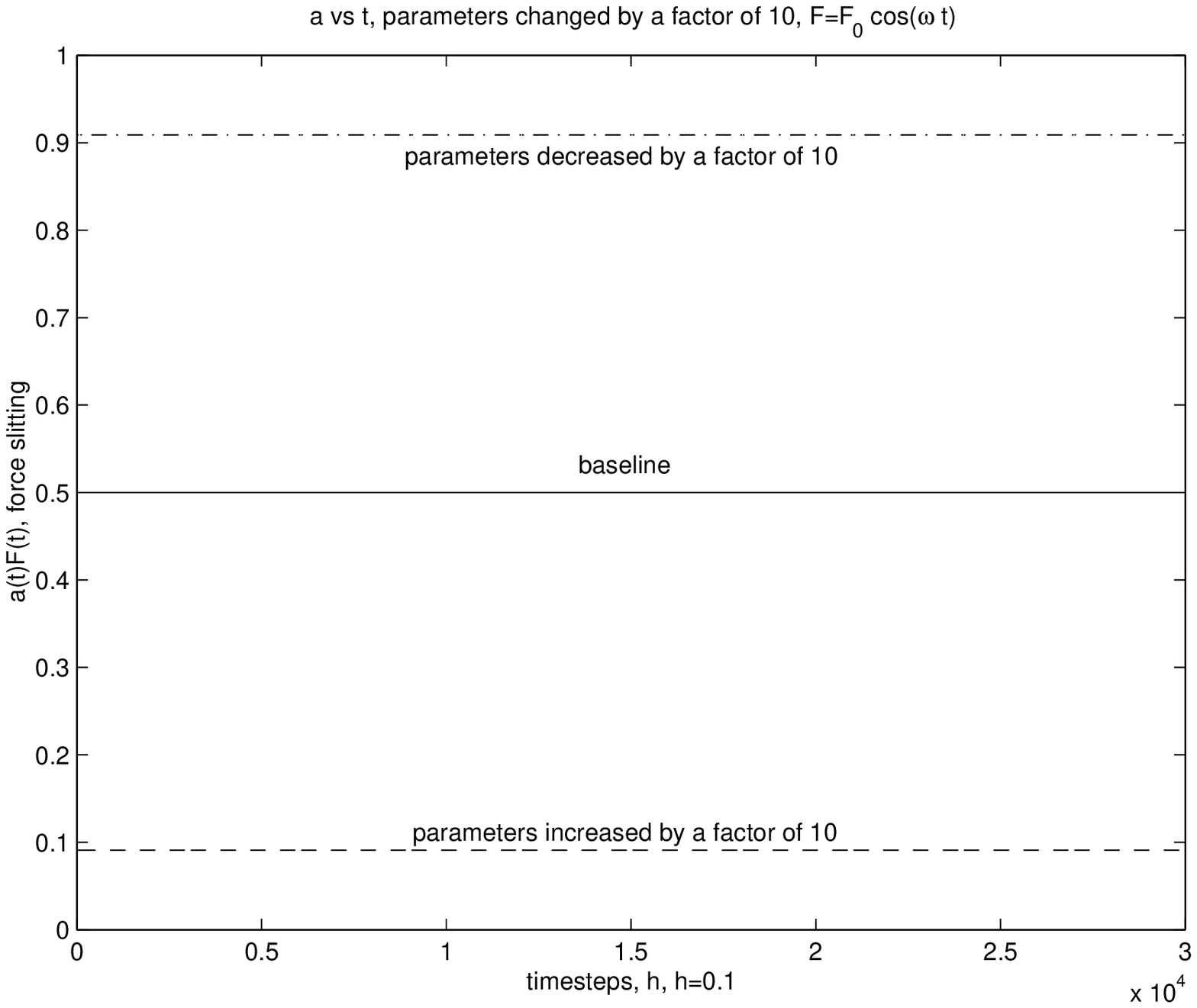}}    
\caption[Force splitting with all parameters changed.  Oscillatory flow.]
{Force splitting when every parameter of body 2 is changed by a factor of
10.  Oscillatory flow.} \label{factor10_aF_o.ps} \end{figure}

In all of the figures depicting the deformation $u(t)$ and the force
splitting, $a(t)F(t)$, the parameters for body one are set to baseline
levels, $\mu_{01}=50$, $\mu_{11}=100$ and $\eta_{11}=5000$.  With the
exception of Figures \ref{factor10_u_s.ps}-\ref{factor10_aF_o.ps} in which
all three parameters of body 2 are changed, in all other figures of this
section two parameters of body two are kept constant, and only one
parameter is changed.  When the steady state behavior of the two Kelvin
bodies is investigated, the following values are used:  $\mu_{02} =
\mu_{12}$ = 1, 2.5, 5, 7.5, 10, 25, 50, 75, 100, 250, 500, 750 and 1000
and $\eta_{12}$ = 1, 2, 5, 10, 20, 50, 100, 250, 500, 1000, 2500, 5000,
10000, 25000, 50000 and 100 000.  The baseline value for $\eta_{12}$ is
larger than the other two constants, and this is why we consider much
larger values when examining the steady state behavior.  The figures show
which parameter is altered, and the values used in the simulations.  The
time steps of the simulations are $h=0.1$, so one second of time is
depicted by ten time steps on the graphs.

\begin{figure}[hb!]
\centerline{\includegraphics[width=0.7\textwidth]{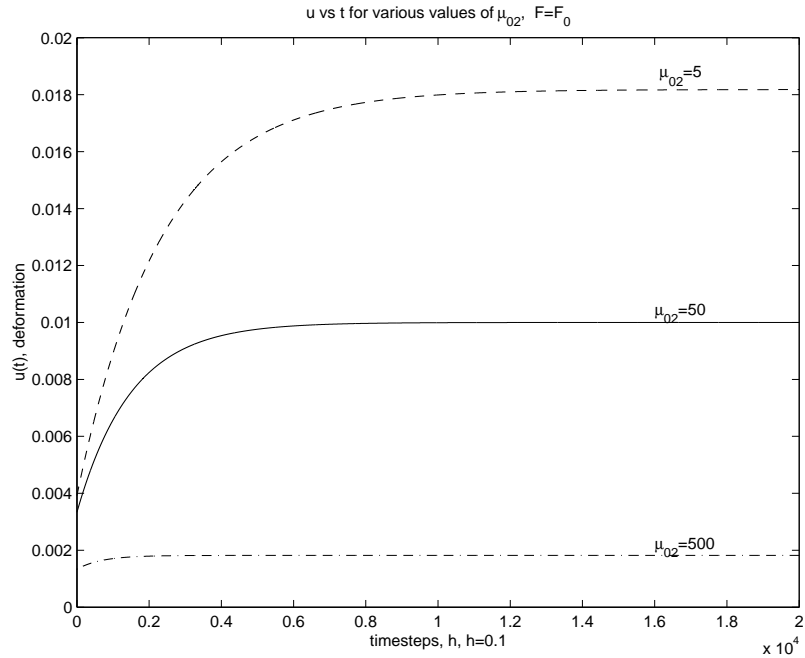}}
\caption[Dependence of deformation on $\mu_{02}$. Steady flow.]{Dependence 
of deformation on $\mu_{02}$. Steady flow.}
\label{2bodiesaF_u_m02_s2.ps}
\end{figure}

\begin{figure}[h!]
\centerline{\includegraphics[width=0.7\textwidth]{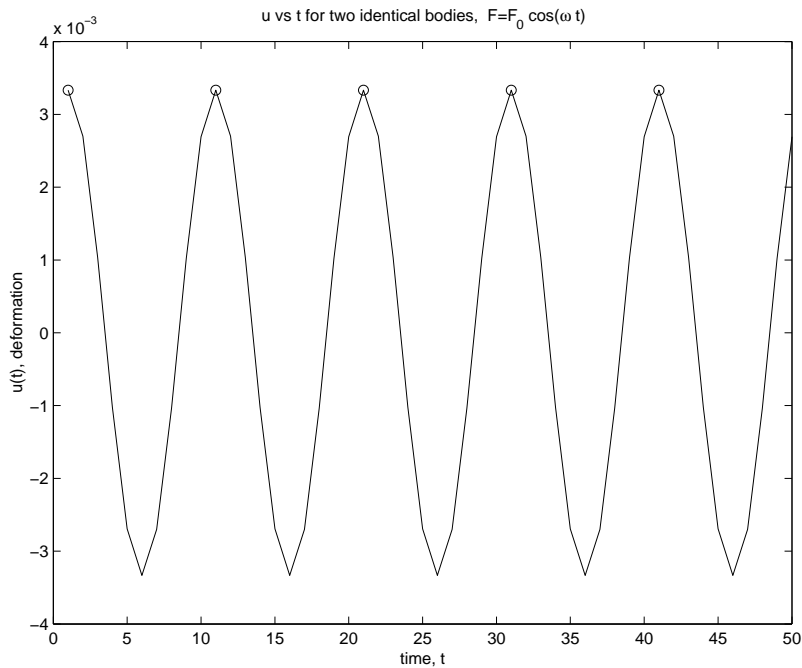}}
\caption[Dependence of deformation on $\mu_{02}$. Oscillatory 
flow.]{Oscillations in the deformation for a particular value of
$\mu_{02}$.  The circles show the only points which are
displayed in the subsequent figures displaying the dependence of
deformation on parameter values. Oscillatory flow.}
\label{2bodiesaF_u_peak_o.ps}
\end{figure}

\begin{figure}[h!]
\centerline{\includegraphics[width=0.7\textwidth]{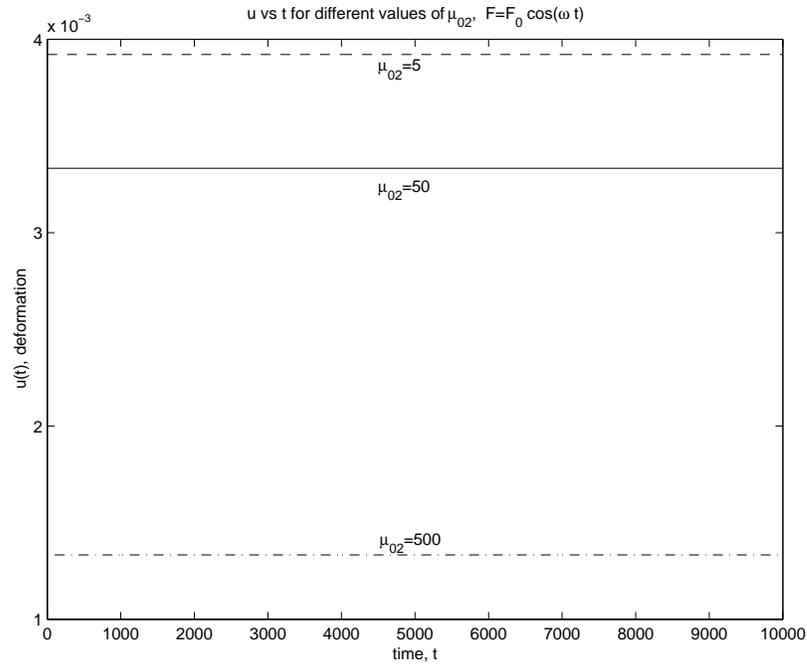}}
\caption[Dependence of peak deformation on $\mu_{02}$. Oscillatory 
flow.]{Dependence of deformation on $\mu_{02}$. Oscillatory. }
\label{2bodiesaF_u_m02_o.ps}
\end{figure}

\begin{figure}[h!]
\centerline{\includegraphics[width=0.7\textwidth]{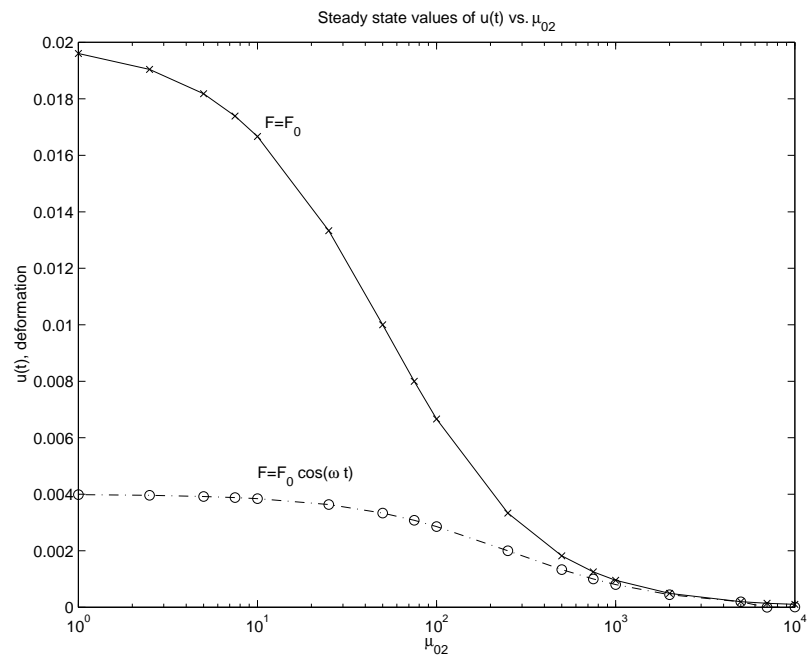}}
\caption[Dependence of steady state deformation on $\mu_{02}$.]{Dependence 
of steady state deformation on $\mu_{02}$.}
\label{peakvalues_u_m02.ps}
\end{figure}

Figures \ref{2bodiesaF_u_m02_s2.ps}-\ref{peakvalues_u_m02.ps} show the
deformation of the two body model as one of the spring constants,
$\mu_{02}$ changes. ($\mu_{02}$ is the constant that characterizes the
isolated spring.)  Figure \ref{2bodiesaF_u_m02_s2.ps} shows the
deformation in steady flow for three values of the spring constant
$\mu_{02}$ =5, 50, 500.  The spring becomes stiffer as $\mu_{02}$
increases, and this results in reducing the deformation, because 
the spring is more difficult to stretch.  The steady state is reached very
quickly for $\mu_{02}$ large.

Figure \ref{2bodiesaF_u_peak_o.ps} illustrates the treatment of
oscillatory flow.  Because the frequency of oscillations is high, we do
not want to display all the oscillations, only the peak values that the
deformation reaches.  These values are marked with an 'o' on this figure,
and in subsequent graphs of oscillatory flow only these are displayed.  
We also note in this figure that the steady state value in oscillatory
flow is obtained within the first few oscillations (2-3)  seconds.  This
observation is confirmed by figure \ref{2bodiesaF_u_m02_o.ps} which
depicts the time evolution of the deformation for oscillatory flow as the
spring constant $\mu_{02}$ is varied.  The steady state of the deformation
is reached almost immediately, and it is a very small value (< 0.004),
even if the spring constant is small $\mu_{02}=5$.  Stiffer springs result
in even smaller deformations, as expected.

We can compare how the steady state of deformation changes with the spring
constant $\mu_{02}$ in steady flow, and the steady state of the peak
deformation in oscillatory flow in Figure \ref{peakvalues_u_m02.ps}.  
Each peak value is obtained after $2 \times 10 ^4$ time steps, which
corresponds to 2000 time units (seconds).  Clearly, much larger
deformation is produced in steady flow than in oscillatory flow, then, as
the spring becomes increasingly stiff, the deformation tends to zero
regardless of the flow.

\begin{figure}[h!]
\centerline{\includegraphics[width=0.7\textwidth]{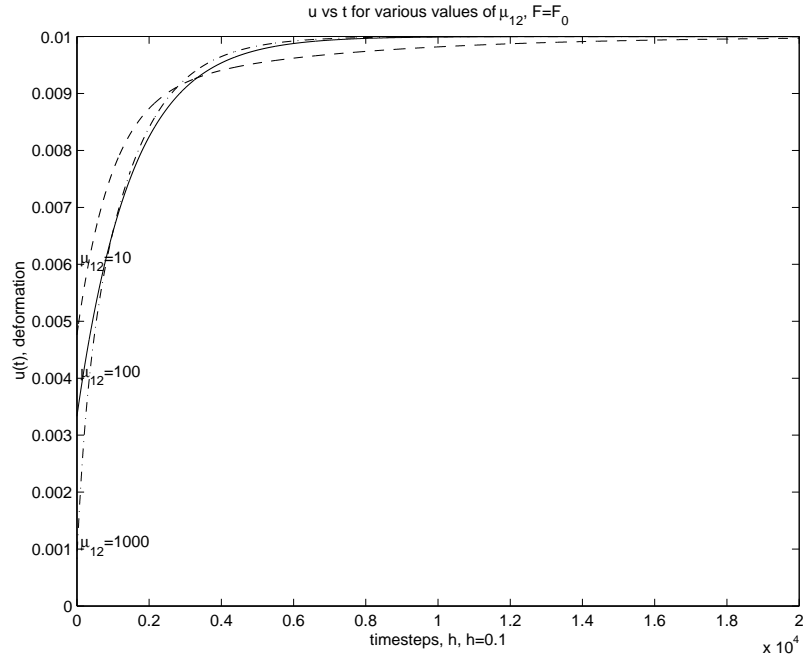}}
\caption[Dependence of deformation on $\mu_{12}$. Steady flow.]{Dependence 
of deformation on $\mu_{12}$. Steady flow.}
\label{2bodiesaF_u_m12_s.ps}
\end{figure}

\begin{figure}[h!]
\centerline{\includegraphics[width=0.7\textwidth]{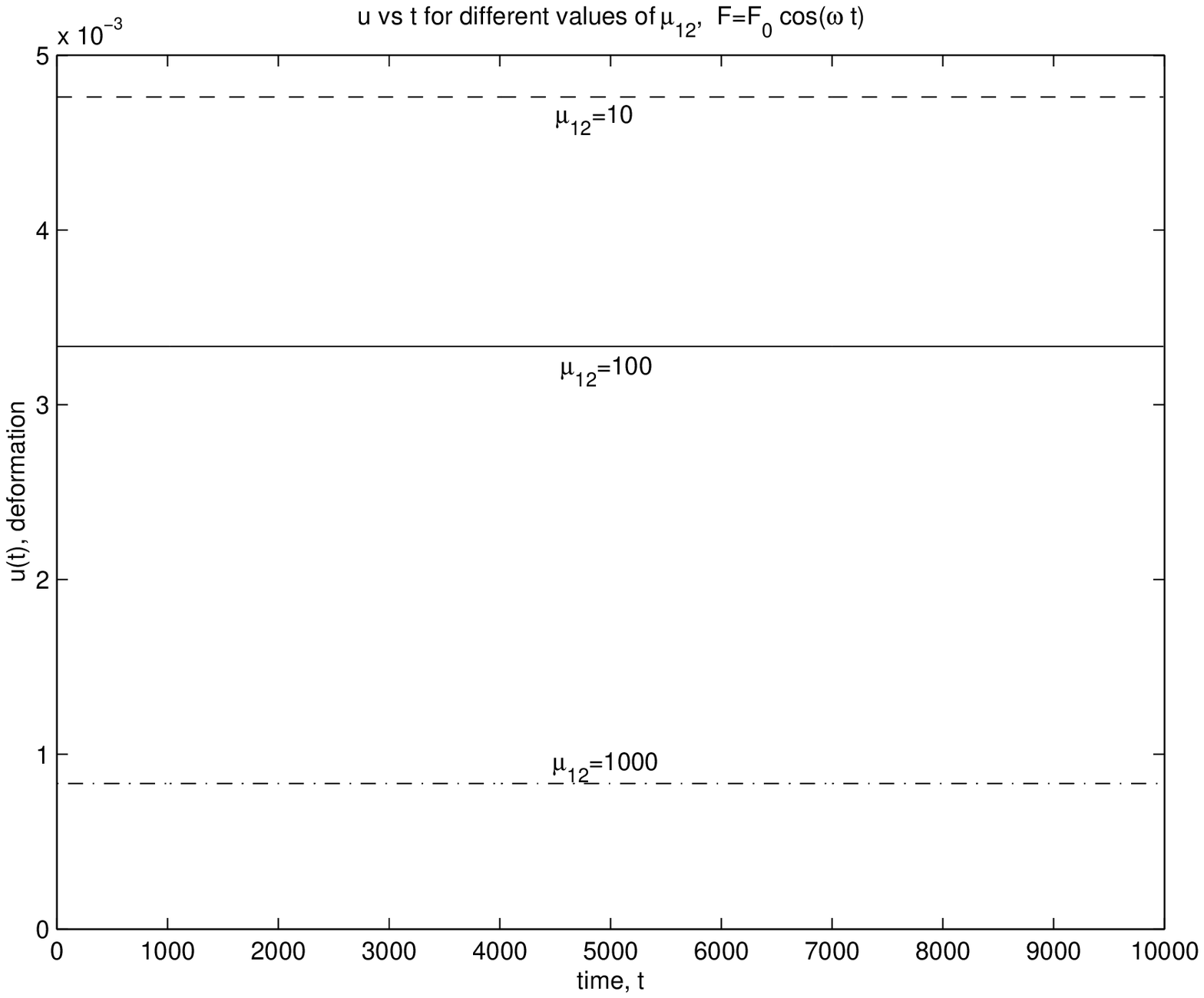}} 
\caption[Dependence of peak deformation on $\mu_{12}$. Oscillatory 
flow.]{Dependence of deformation on $\mu_{12}$. Oscillatory flow. }   
\label{2bodiesaF_u_m12_o.ps}
\end{figure}

\begin{figure}[h!]
\centerline{\includegraphics[width=0.7\textwidth]{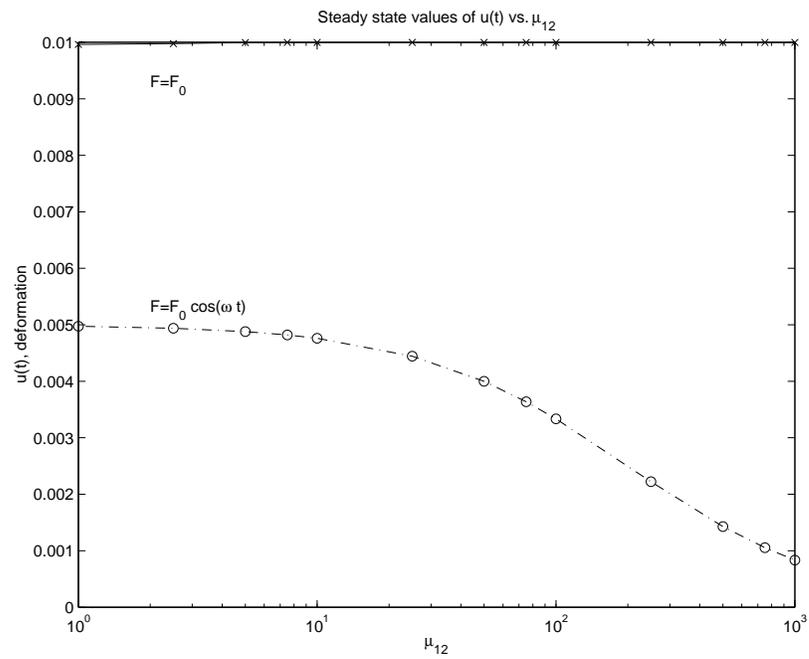}}
\caption[Dependence of steady state deformation on $\mu_{12}$.]{Dependence 
of steady state deformation on $\mu_{12}$.}
\label{peakvalues_u_m12.ps}
\end{figure}

The next three figures, Figures
\ref{2bodiesaF_u_m12_s.ps}-\ref{peakvalues_u_m12.ps} show how the
deformation of a cell changes as the other spring constant of body two,
$\mu_{12}$ changes.  (This spring constant characterizes the spring next
to the dashpot in the second Kelvin body.)  Figure
\ref{2bodiesaF_u_m12_s.ps} shows that although the initial deformation is
different for the three values of $\mu_{12}$ (as one expects, because the
initial condition depends on $\mu_{12}$), the steady state obtained is the
same for the three values.  Similarly, in Figure
\ref{2bodiesaF_u_m12_o.ps}, we see that in oscillatory flow the steady
state is obtained virtually immediately again (in 2-3 seconds), and that
the deformation for all three values of $\mu_{12}$ is very small, less
than 0.005.

In Figure \ref{peakvalues_u_m12.ps} most steady state values are obtained
after 2000 seconds again.  In oscillatory flow, as the spring becomes
stiffer, the deformation decreases, just like it did in Figure
\ref{peakvalues_u_m02.ps}.  However, in steady flow, the deformation stays
almost the same regardless of the spring constant.  This is due to more of
the force concentrating on the dashpot if the spring is very pliable. In
this case the steady state value is obtained slower too, in (3000 seconds)
which also suggests that the force is concentrated on the dashpot which
creeps to the steady deformation slower.  Just like before, the overall
deformation is always smaller in oscillatory flow, because the oscillatory
flow only produces a large force periodically and not continuously.  
Also, there is a force of equal magnitude but opposite direction acting on
the dashpot periodically, therefore the dashpot is unable to extend fully.  
(This claim is verified later in simulations where the frequency of
oscillations is changed, allowing the dashpot more time to deform.)  When
the spring next to the dashpot is very pliable, i.e for small values of
$\mu_{12}$, much of the force is allowed to act on the dashpot, therefore
the overall deformation is larger when the $\mu_{12}$ is changed ( Figure
\ref{peakvalues_u_m12.ps}) than if $\mu_{02}$ is changed
 (Figure \ref{peakvalues_u_m02.ps}).

\begin{figure}[h!]
\centerline{\includegraphics[width=0.7\textwidth]{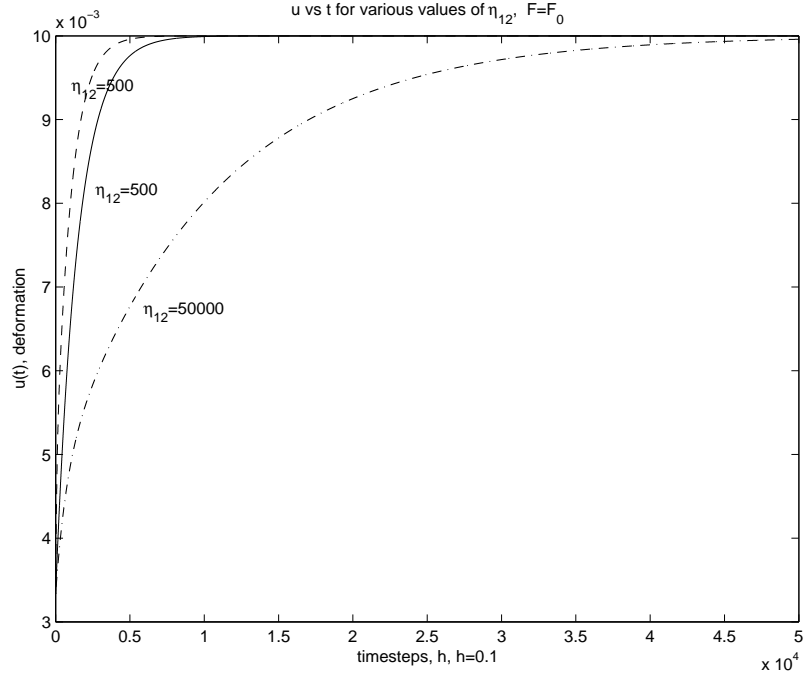}}
\caption[Dependence of deformation on $\eta_{12}$. Steady 
flow.]{Dependence of deformation on $\eta_{12}$. Steady flow.}
\label{2bodiesaF_u_eta_s.ps}
\end{figure}

\begin{figure}[h!]
\centerline{\includegraphics[width=0.7\textwidth]{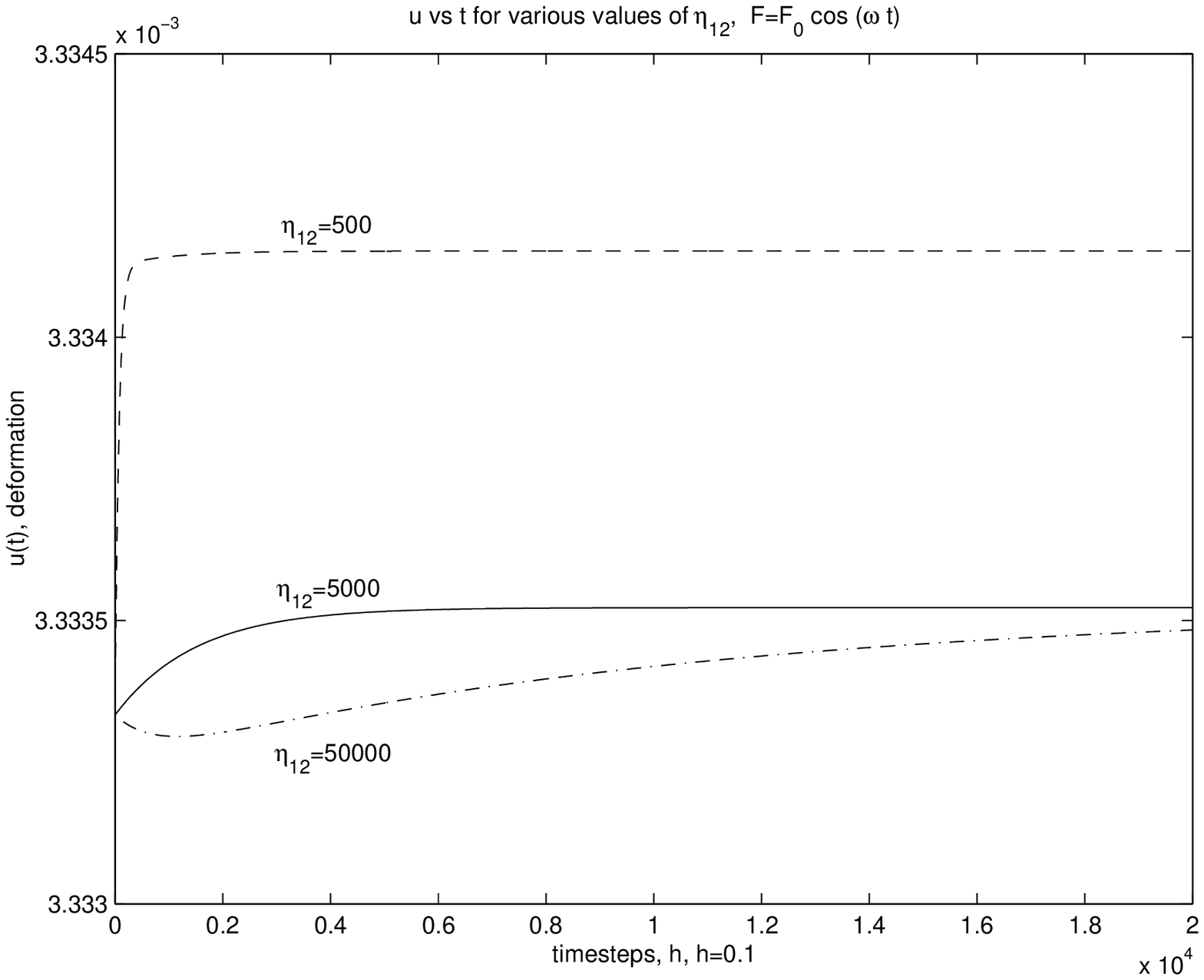}}
\caption[Dependence of peak deformation on $\eta_{12}$. Oscillatory 
flow.]{Dependence of deformation on $\eta_{12}$.  Oscillatory flow. }
\label{2bodiesaF_u_eta_peak_o.ps}
\end{figure}

\begin{figure}[h!]
\centerline{\includegraphics[width=0.7\textwidth]{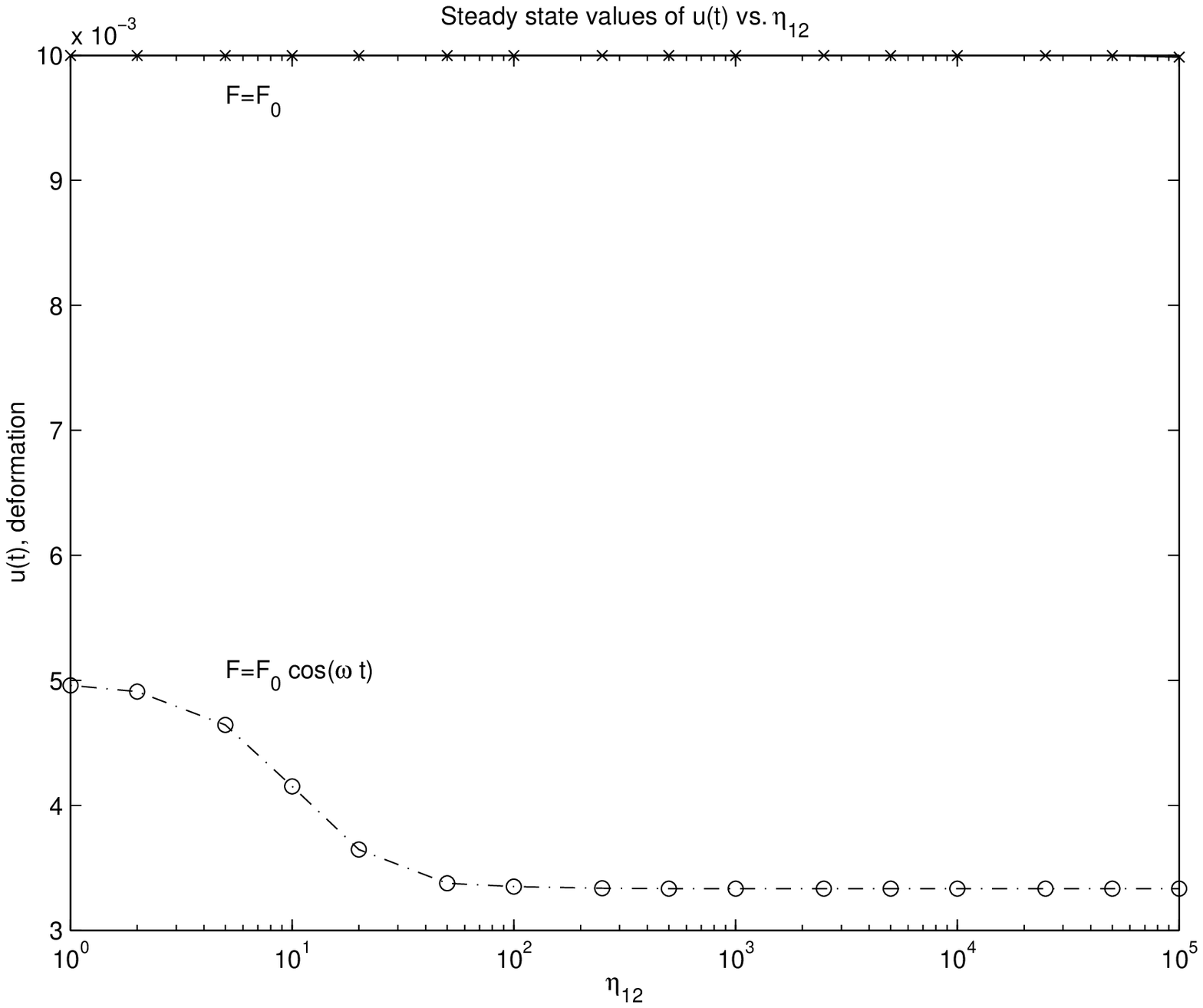}}  
\caption[Dependence of steady state deformation on $\eta_{12}$] 
{Dependence of steady state deformation on $\eta_{12}$.}
\label{peakvalues_u_eta.ps}
\end{figure}

Now we turn to examining how the dashpot viscosity influences the behavior
of the system.  Figures \ref{2bodiesaF_u_eta_s.ps} and
\ref{2bodiesaF_u_eta_peak_o.ps} depict the time evolution of the
deformation for $\eta_{12}$=500, 5000 and 50000.  In steady flow the
steady state of the deformation is the same regardless of the value of
$\eta_{12}$, but for large dashpot viscosities the steady state takes
longer to obtain, approximately 5000 seconds.  Figure
\ref{2bodiesaF_u_eta_peak_o.ps} shows the steady state of the deformation
in oscillatory flow.  For the three values of $\eta_{12}$ shown here, the
difference in deformations is negligible, on the order of $10^{-6}$.  The
largest deformation is obtained when $\eta_{12}$ is the smallest, and the 
bodies in this case deform very quickly, within the first 100 seconds.

Figure \ref{peakvalues_u_eta.ps} is created by running the simulations
longer than in the previous figures, to time $t=1.2 \times 10^{4}$ seconds
to ensure that the deformation reaches its steady state. In oscillatory
flow, the steady state of the deformation is slightly larger for small
dashpot viscosities, but after about $\eta_{12}=100$, the deformation is
independent of the viscosity. The interpretation of this is that when the
viscosity is sufficiently low, the force is able to deform the dashpot
quickly, then, as the forcing oscillates the deformation remains the same.  
For larger viscosities the initial deformation of the dashpot is
negligible, and the later oscillations are again unable to change the
deformation of the dashpot.  As noted before, the deformation in steady
flow is independent of the dashpot viscosity.

\begin{figure}[h!]
\centerline{\includegraphics[width=0.7\textwidth]{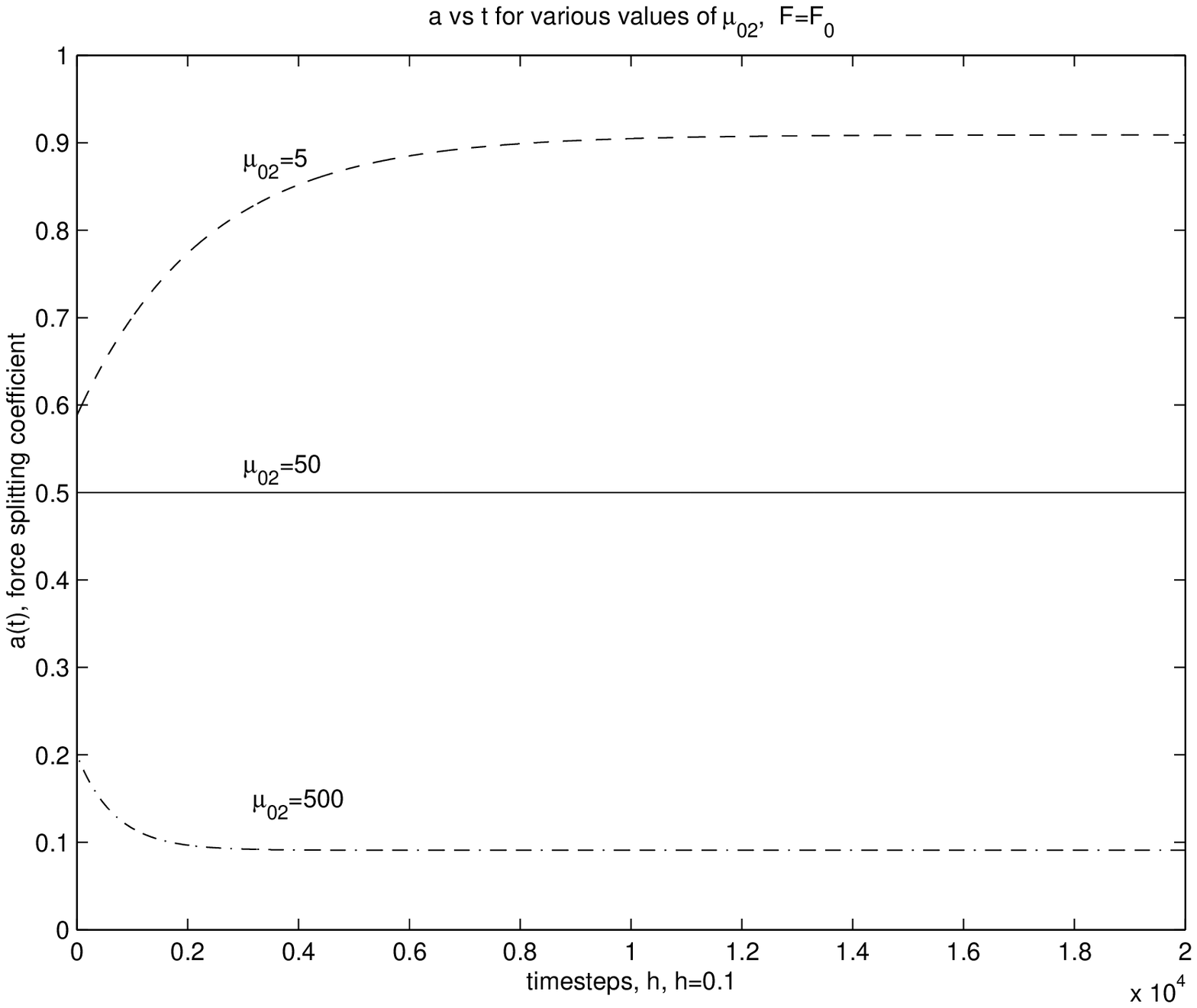}}
\caption[Dependence of force splitting on $\mu_{02}$. Steady
flow.]{Dependence of force splitting coefficient on $\mu_{02}$. Steady
flow.} \label{2bodiesaF_a_m02_s.ps} \end{figure}

\begin{figure}[h!]
\centerline{\includegraphics[width=0.7\textwidth]{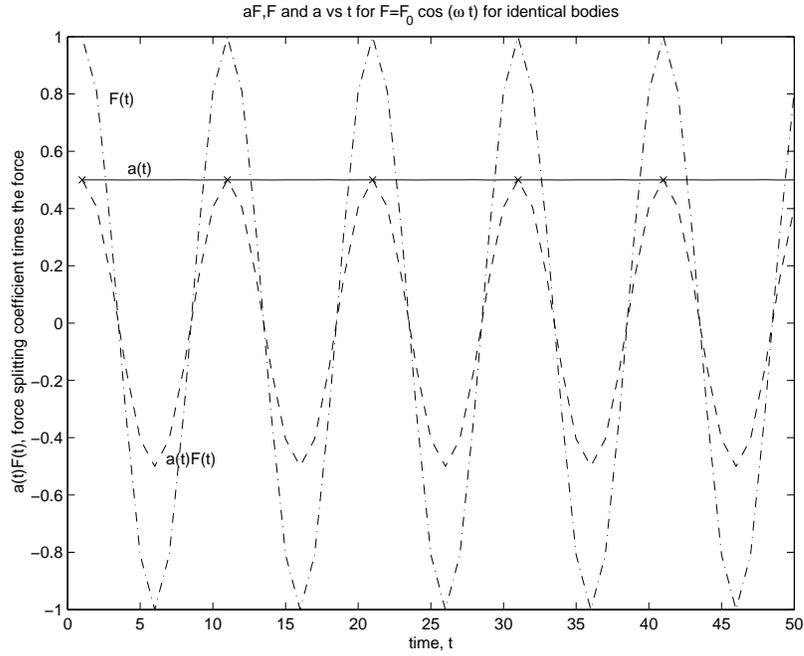}}
\caption[Dependence of force splitting on $\mu_{02}$. Oscillatory
flow.]{Oscillations of $a(t)F(t)$, $F(t)$ for identical bodies. $a(t)$
also displayed. 'x' marks the points which are used to display peak
oscillations of $a(t)$ in the following figures. Oscillatory flow. }
\label{2bodiesaF_aaFF_o.ps} \end{figure}

\begin{figure}[h!]
\centerline{\includegraphics[width=0.7\textwidth]{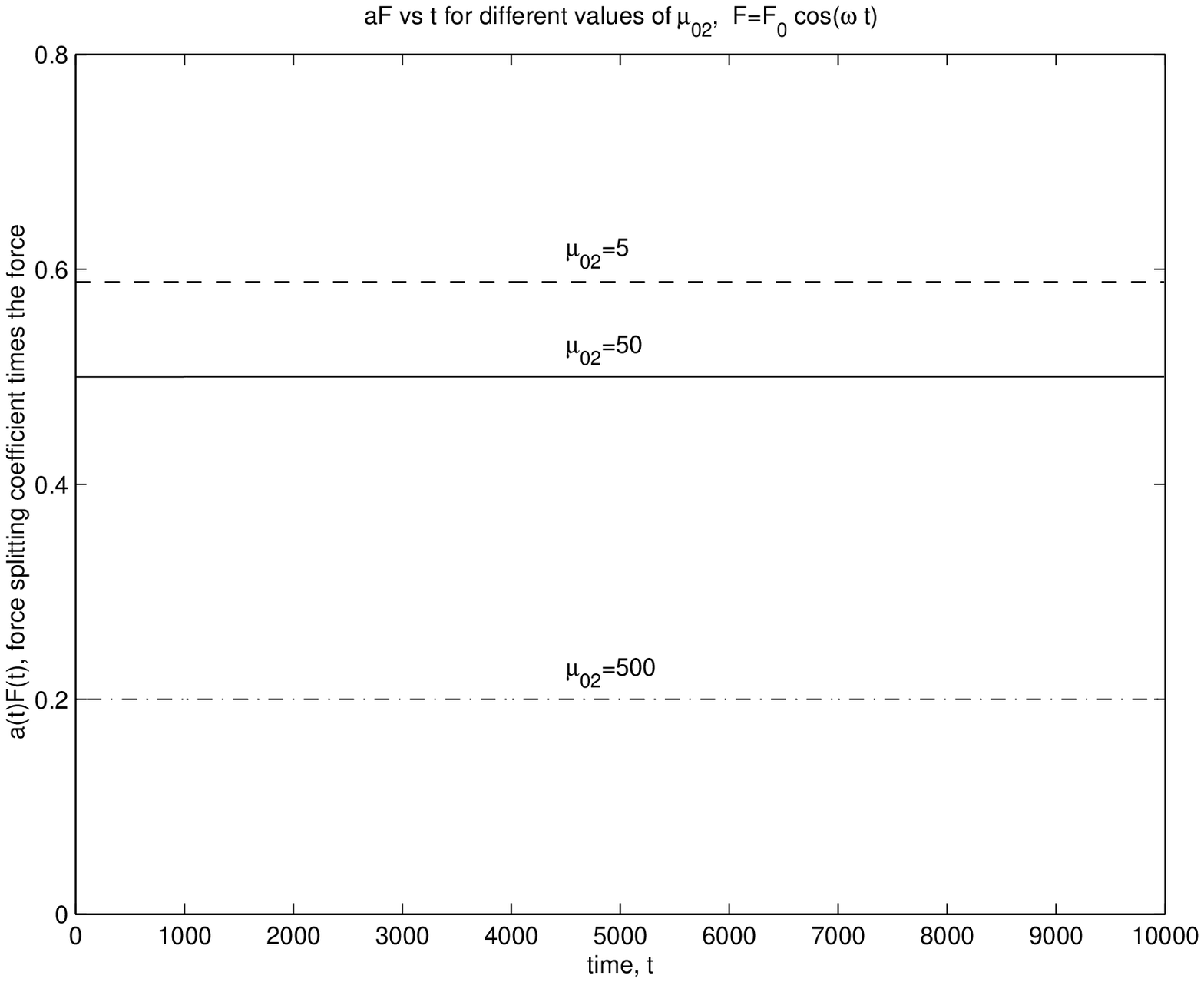}}
\caption[Dependence of peak force splitting on $\mu_{02}$. Oscillatory
flow.]{Dependence of $a(t)$ on $\mu_{02}$.  Graph shows peak values of
$a(t)F(t)$.  Oscillatory flow.} \label{2bodiesaF_aF_m02_o.ps} \end{figure}

\begin{figure}[h!]
\centerline{\includegraphics[width=0.7\textwidth]{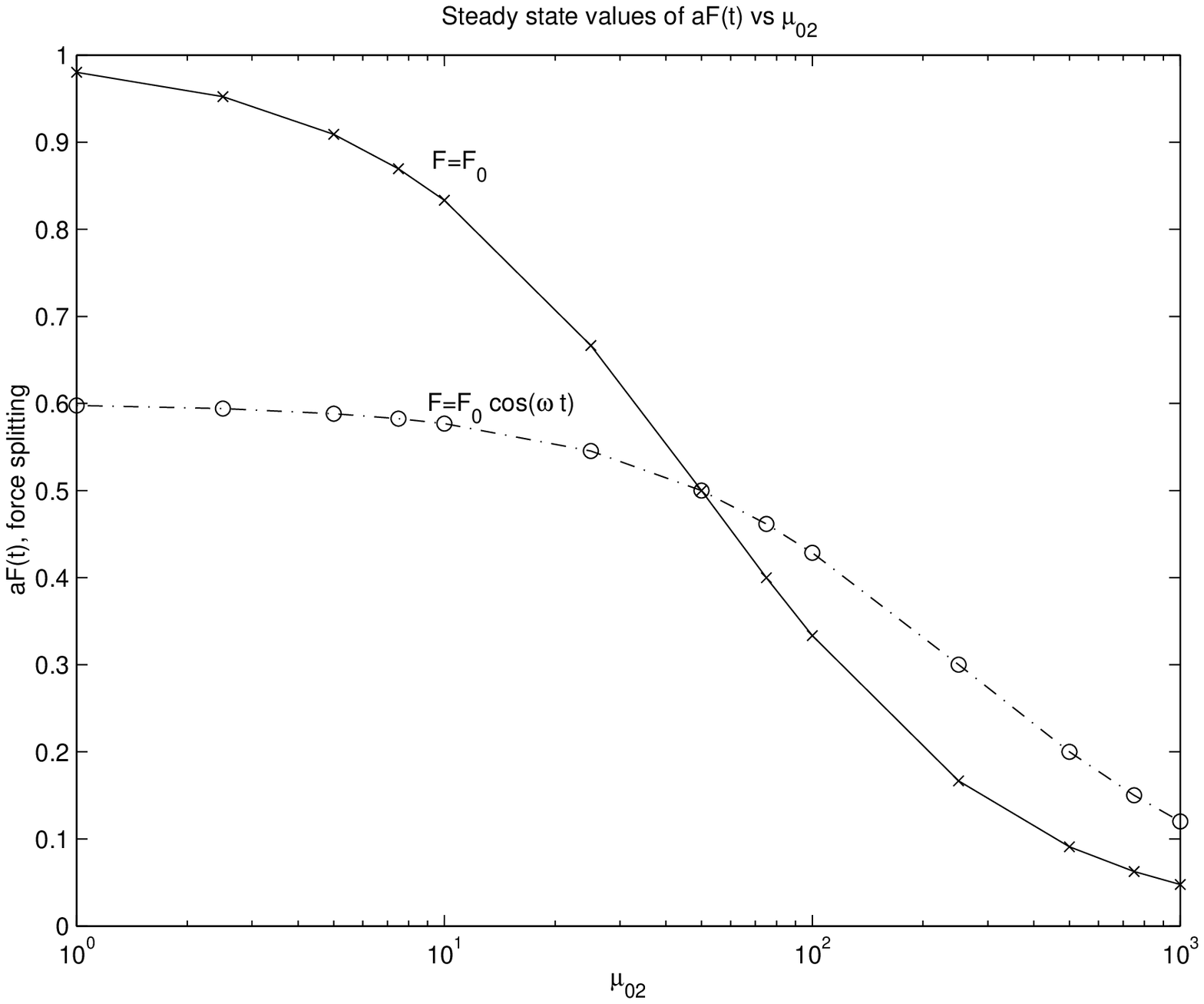}}
\caption[Dependence of steady state force splitting on
$\mu_{02}$.]{Dependence of steady state force splitting on $\mu_{02}$.}
\label{peakvalues_aF_m02.ps} \end{figure}

Now we can turn to investigating how the force splitting is effected by
the parameter values.  All the figures of the force splitting show the
force acting on body one, and the force on the second body can be found
from this by recalling that the sum of the two forces is one.  We use the
same values for $\mu_{02}$ as previously.  The simulations in steady flow
are shown in Figure \ref{2bodiesaF_a_m02_s.ps}.  As expected, if the
spring is very pliable in body two, then most of the force will have to
focus on body one.  For two identical bodies (i.e. if $\mu_{02}$ = 50),
the force splits equally between the two bodies.  A very stiff spring in
body two means that all the force has to concentrate here.  In oscillatory
flow we are only interested in the peak force acting on body one.  Figure
\ref{2bodiesaF_aaFF_o.ps} shows the oscillations of the overall force, the
force acting on body one, and the force splitting coefficient, $a$
displayed on the same graph.  Only the peak values of $aF$ (which are
identical to the value of $a$) are displayed subsequently.  As before, the
peak values of the force do not change with time, as shown in Figure
\ref{2bodiesaF_aF_m02_o.ps}.  In the oscillatory flow, just like in steady
flow, a pliable spring in body two leads to more of the force
concentrating on body one, in identical bodies the forces split evenly,
and stiff springs require more force to deform.

Figure \ref{peakvalues_aF_m02.ps} depicts the steady state values of force
splitting for a range of values of $\mu_{02}$.  The curves for steady flow
and oscillatory flow intersect when $\mu_{02}$=50, because this is the
baseline value at which the two bodies are identical, therefore the forces
split evenly regardless of the flow.  It is interesting to note that in
steady flow the force splitting is more extreme than in oscillatory flow.  
More specifically, a small $\mu_{02}$ (pliable spring in body two) in
steady flow allows much more of the force to act on body one than in
oscillatory flow, but a large value of $\mu_{02}$ in steady flow leads to
a smaller force on body one than in oscillatory flow.  The interpretation
of this is that in oscillatory flow the dashpot offers a constant
resistance, thus much of the overall force is always trying to stretch the 
dashpot.  

\begin{figure}[h!]
\centerline{\includegraphics[width=0.7\textwidth]{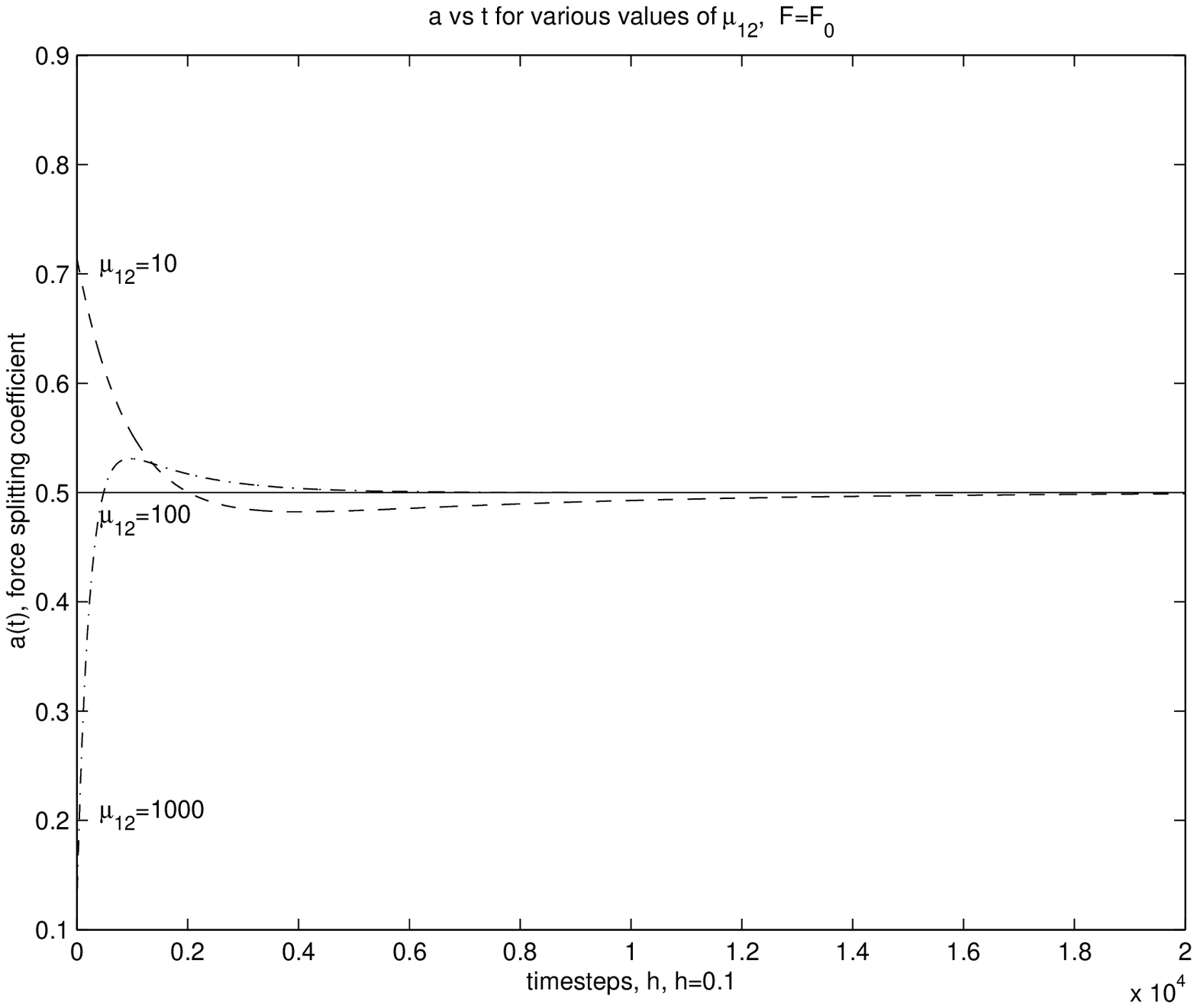}}
\caption[Dependence of force splitting on $\mu_{12}$. Steady
flow.]{Dependence of force splitting coefficient on $\mu_{12}$. Steady
flow.} \label{2bodiesaF_a_m12_s.ps} \end{figure}

\begin{figure}[h!]
\centerline{\includegraphics[width=0.7\textwidth]{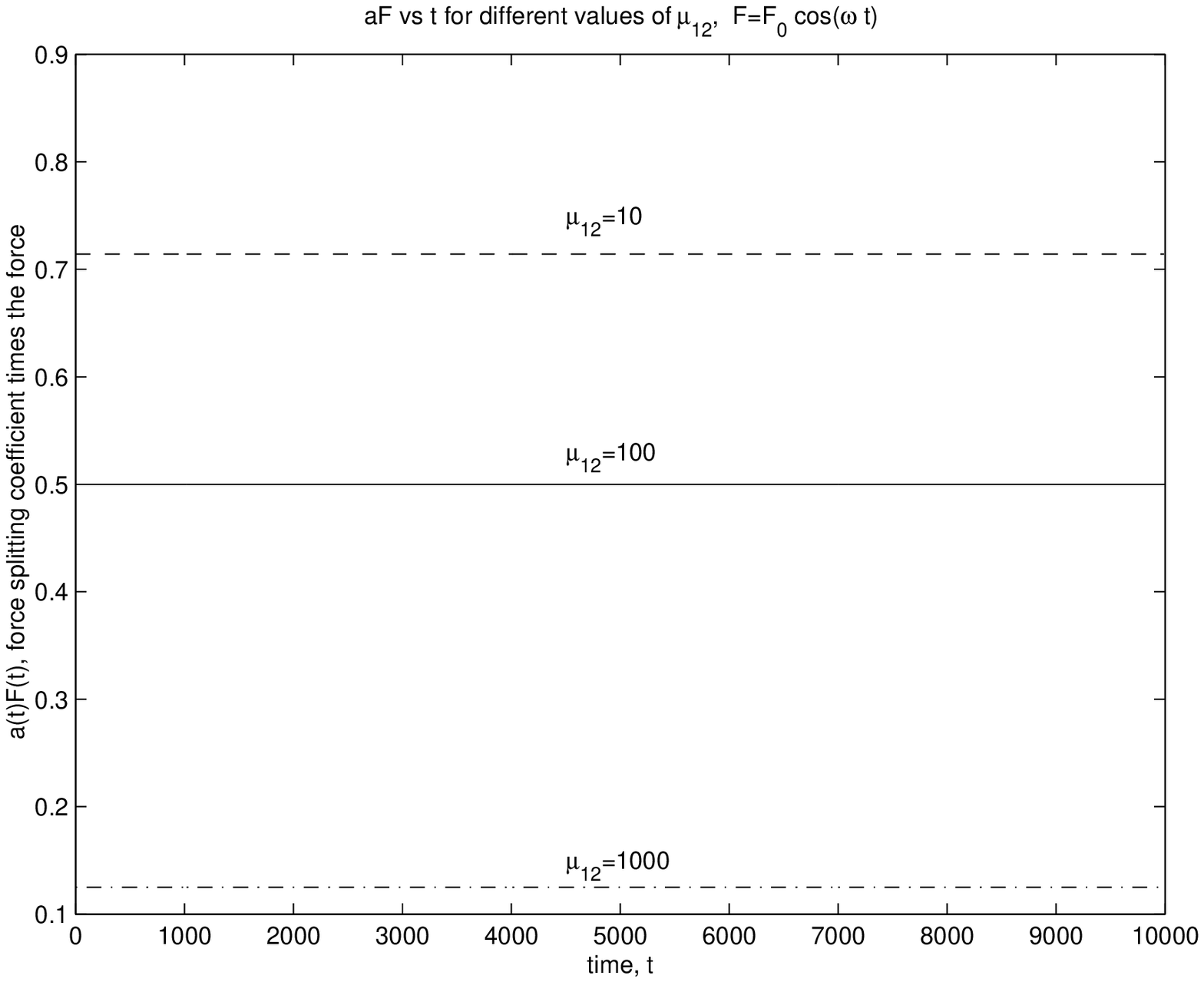}}
\caption[Dependence of peak force splitting on $\mu_{12}$. Oscillatory
flow.]{Dependence of $a(t)$ on $\mu_{12}$.  Graph on shows peak values of
$a(t)F(t)$.  Oscillatory flow.} \label{2bodiesaF_aF_m12_o.ps} \end{figure}

\begin{figure}[h!]
\centerline{\includegraphics[width=0.7\textwidth]{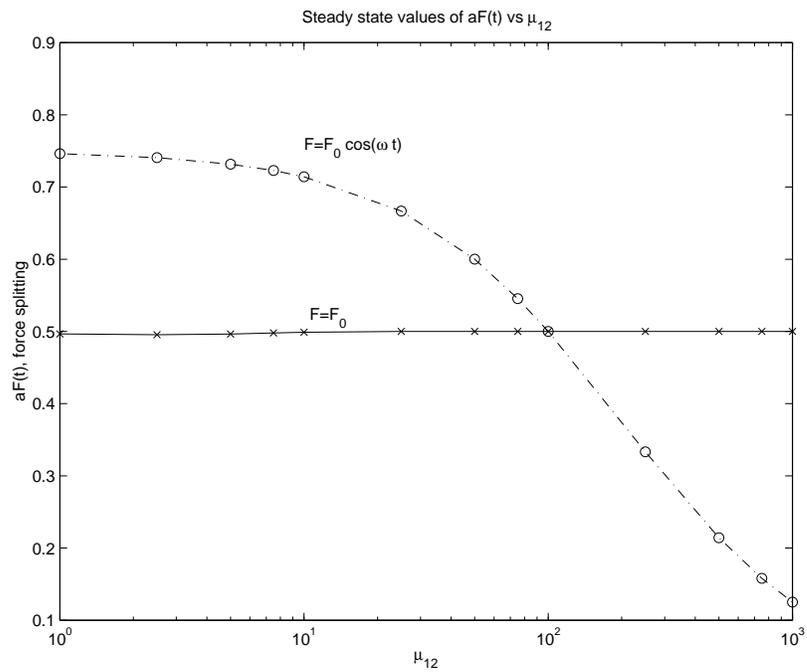}}  
\caption[Dependence of steady state force splitting on 
$\mu_{12}$.]{Dependence of steady state force splitting on $\mu_{12}$.}
\label{peakvalues_aF_m12.ps}
\end{figure}

Next, we examine the dependence of the force splitting on the other spring
constant, $\mu_{12}.$ Figure \ref{2bodiesaF_a_m12_s.ps} shows the force
splitting in steady flow.  The initial condition of the force in the bodies
depends on $\mu_{12}$, but clearly, the steady state of the force tends to
the same value, 0.5.  Regardless of the stiffness of the spring in body
two, eventually the forces acting on the two bodies become the same.  
When the force acts on a pliable spring, more force goes to the dashpot,
and when the spring is stiff, more of the force goes to it.  The mediating
effects of the dashpot lead the equal force splitting between the two
bodies.

In oscillatory flow, shown in Figure \ref{2bodiesaF_aF_m12_o.ps} the 
mediating effects of the dashpot are smaller, so there is a larger 
difference between force splitting between bodies one and two for 
different values of the spring constant, $\mu_{12}$.  As before, the peak 
force acting on body one reaches its steady state very quickly, so we see 
constants.  For small values of $\mu_{12}$ the more of the force acts on 
body one, and for stiff springs more force is concentrated on body 
two.  The comparison of steady states in oscillatory flow and steady flow, 
as depicted by Figure \ref{peakvalues_aF_m12.ps} brings no significant new 
information.        
   
\begin{figure}[h!]
\centerline{\includegraphics[width=0.7\textwidth]{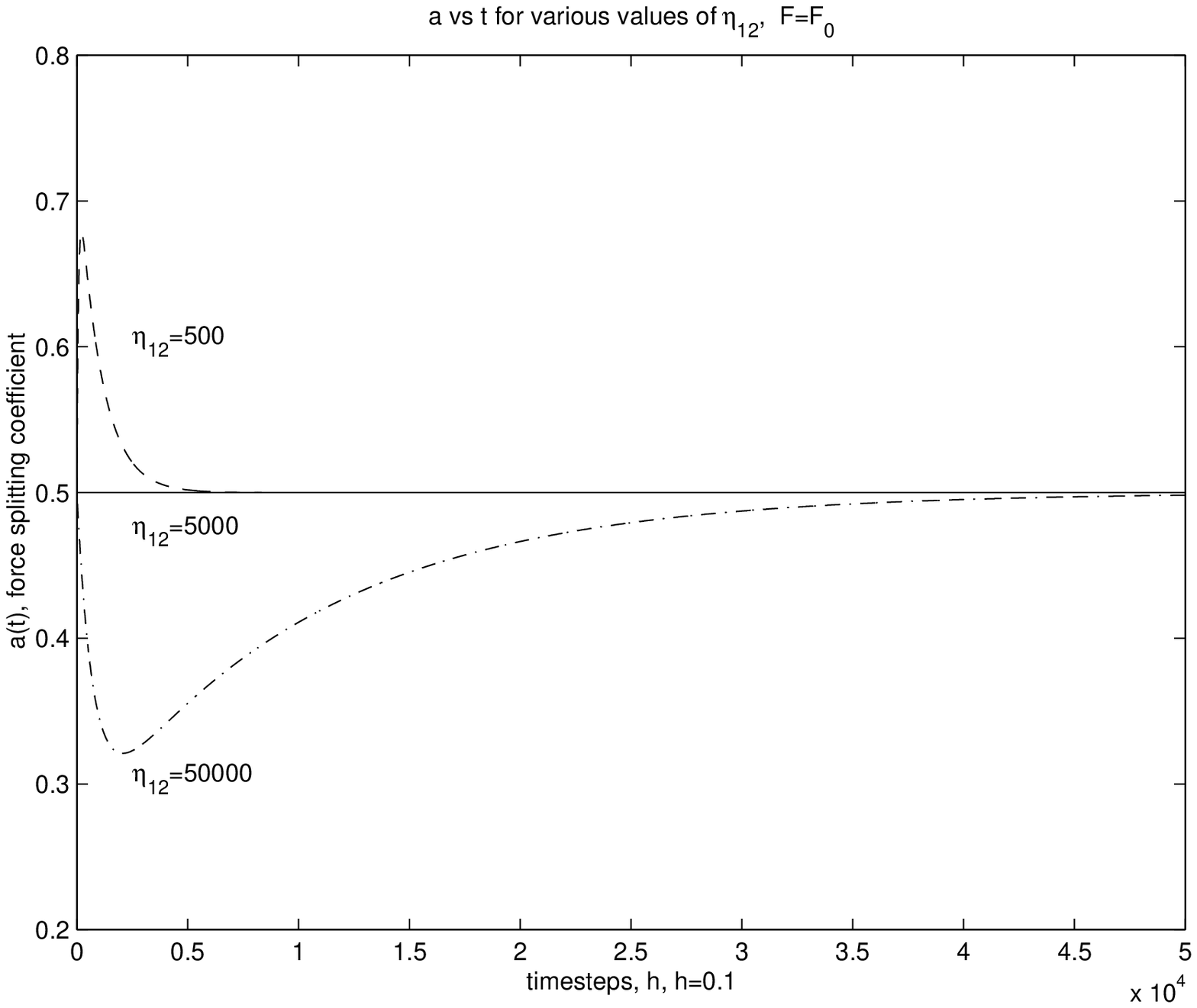}}
\caption[Dependence of force splitting coefficient on 
$\eta_{12}$. Steady flow.]{Dependence of force splitting coefficient on 
$\eta_{12}$.
Steady flow.}
\label{2bodiesaF_a_eta_s.ps}
\end{figure}

\begin{figure}[h!] 
\centerline{\includegraphics[width=0.7\textwidth]{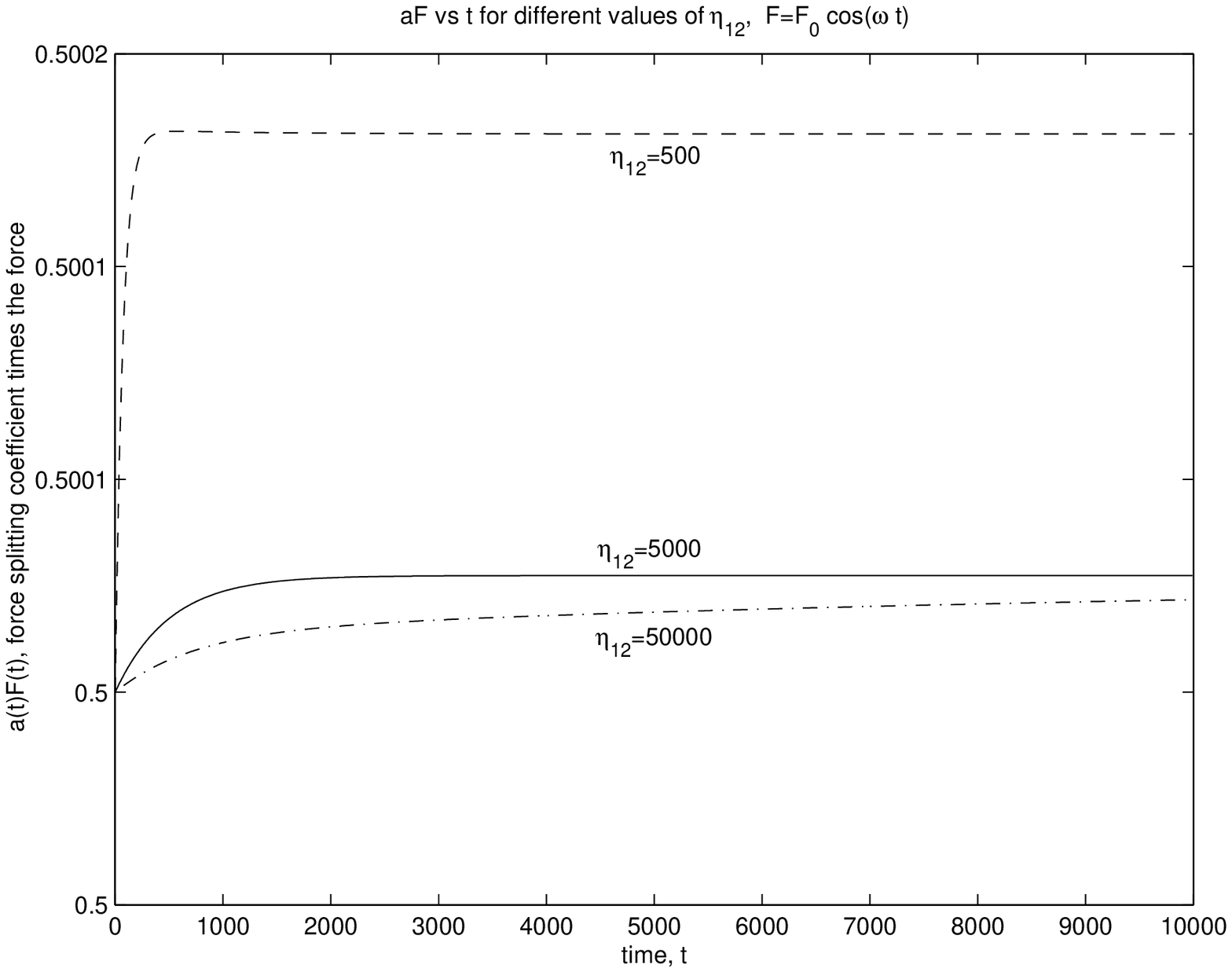}}
\caption[Dependence of peak force splitting on $\eta_{12}$. 
Oscillatory flow.]{Dependence of $a(t)$ on $\eta_{12}$.  Graph shows 
peak values of $a(t)F(t)$.  Oscillatory flow.} 
\label{2bodiesaF_aF_eta_o.ps}
\end{figure}

\begin{figure}[h!]
\centerline{\includegraphics[width=0.7\textwidth]{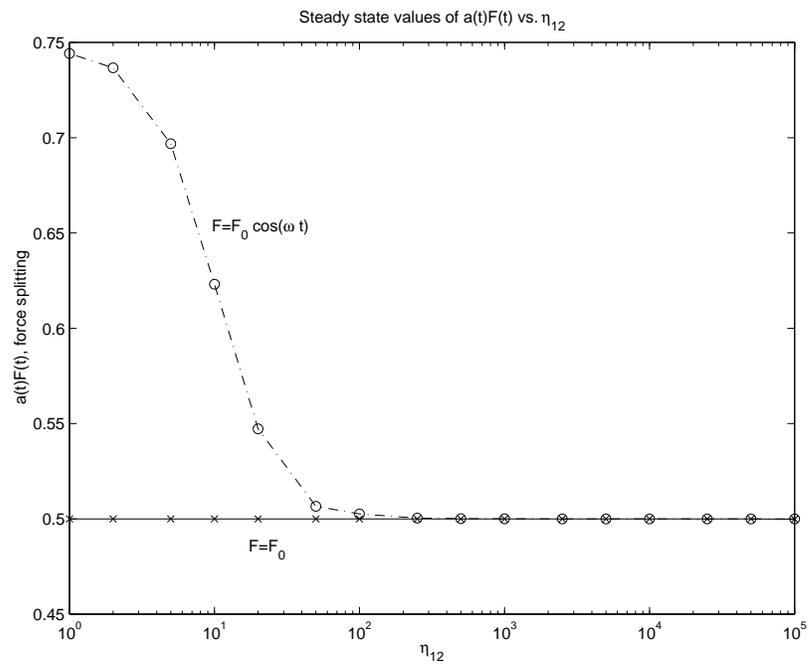}}
\caption[Dependence of steady state force splitting on 
$\eta_{12}$.]{Dependence of steady state force splitting on $\eta_{12}$.}
\label{peakvalues_aF_eta.ps}
\end{figure}

Now we turn to looking at the dependence of force splitting on the dashpot
viscosity.  Longer simulation times are necessary again to obtain the
steady state values.  In steady flow, shown in Figure
\ref{2bodiesaF_a_eta_s.ps} small dashpot viscosity in body two leads to a
transient increase of the force acting on body one, but at the steady
state the force splits equally between the bodies.  The transient increase 
is due to more force needing to deform the dashpot of body one.  When the 
dashpot viscosity of body two is large, there is a transient decrease in 
the force acting on body one.  The reason for this is similar to our 
argument above: now more force is necessary to deform the dashpot of body 
two.  It is also clear from the figure that for large dashpot viscosities 
the transient time is longer.  This is consistent with our previous 
observations on the role of the dashpot viscosity.  

\begin{figure}[h!]
\centerline{\includegraphics[width=0.7\textwidth]{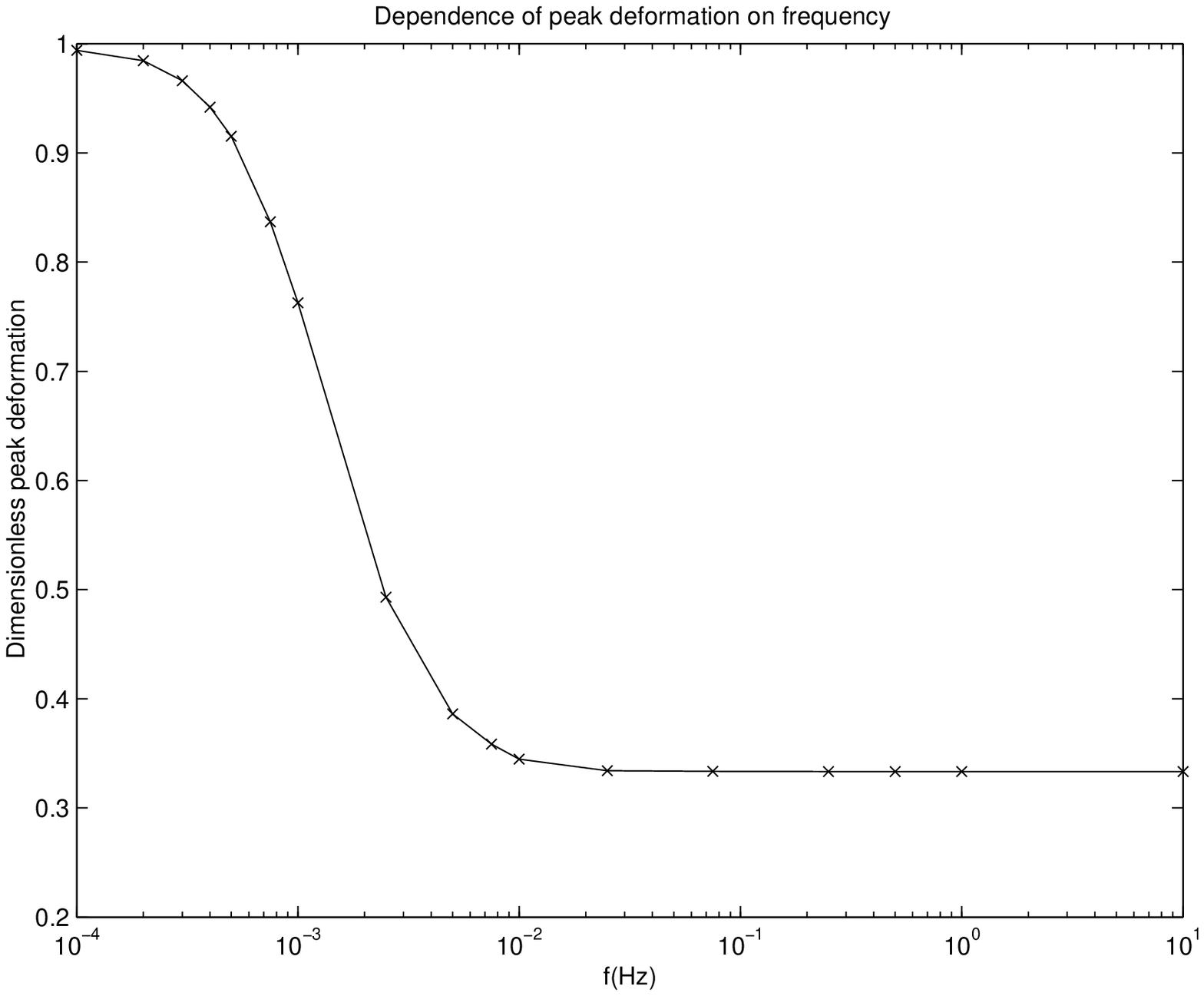}}
\caption[Peak steady state deformation as a function of frequency of  
oscillations.]{ Dependence of peak steady state deformation on the
frequency.
Peak steady state deformation of oscillatory flow is divided by the steady
state value of steady flow.}
\label{peakvalues_u_freq.ps}
\end{figure}

\begin{figure}[h!]
\centerline{\includegraphics[width=0.7\textwidth]{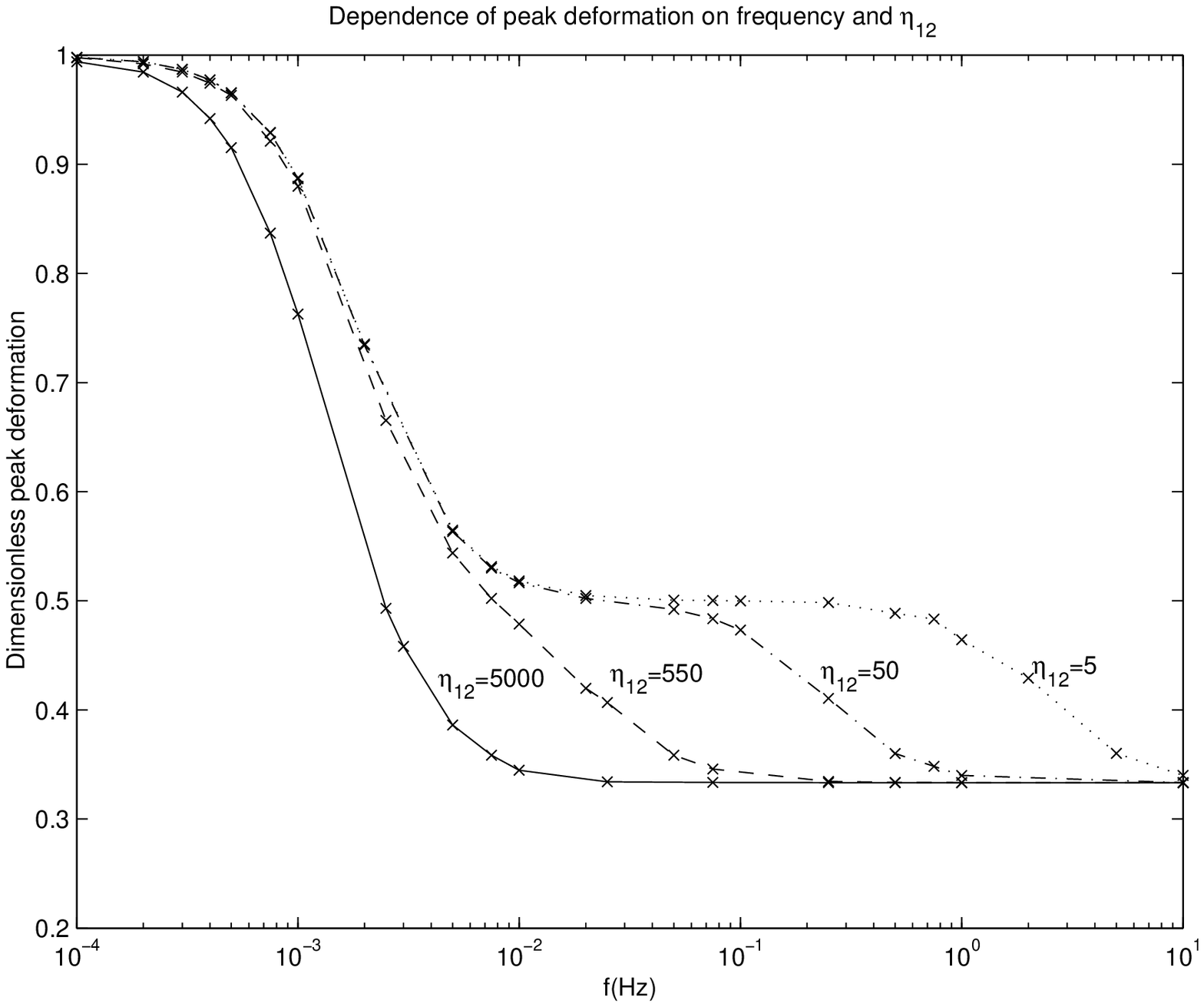}}
\caption[Peak steady state deformation as function of
frequency and $\eta_{12}$.]{Dependence
of peak steady state deformation on the frequency and
$\eta_{12}$.  Peak steady state deformation of oscillatory flow is divided
by the steady state value of steady flow.}
\label{peakvalues_u_freq_eta3.ps}
\end{figure}

\begin{figure}[h!]
\centerline{\includegraphics[width=0.7\textwidth]{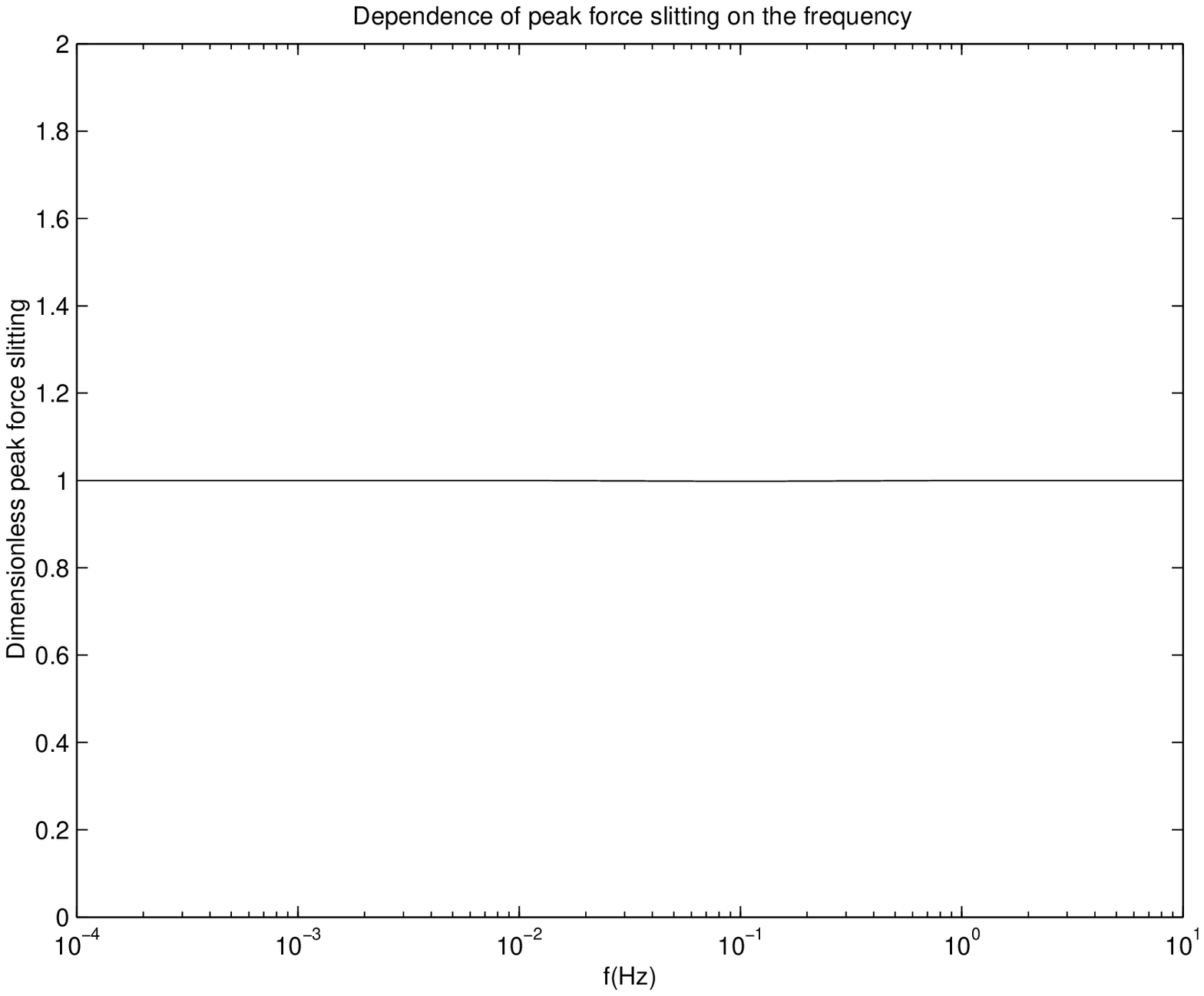}}
\caption[Dependence of peak steady state force splitting on the
frequency.]{Dependence of peak steady state force splitting on the
frequency.  Peak steady state force splitting of oscillatory flow is
divided by the steady state value of steady flow.}
\label{peakvalues_a_freq.ps} \end{figure}

In oscillatory flow, Figure \ref{2bodiesaF_aF_eta_o.ps}, the force split
is almost exactly equal independently of the dashpot viscosity.  Small
viscosities lead to a quicker and slightly larger deformation.  Figure
\ref{peakvalues_aF_eta.ps} shows how the steady state values of the force
in body one change with the dashpot viscosity in steady and in oscillatory
flow.  The most notable observation is that in oscillatory flow the forces
do not split equally.  This is the result of the initial quick dashpot
relaxation of body two (when $\eta_{12}$ is small) that leads to more
force being necessary for body one.  Similarly, if the dashpot viscosity
of the second body is large, then the initial dashpot relaxation occurs in
the first body, therefore more force will always be applied to the second
body.

So far all of our numerical simulations in oscillatory flow used the
frequency $f = 1 Hz = 2 \pi $ rad/sec.  This is a physiological value
which corresponds to the frequency at which the heart pumps the blood
through the blood vessels, however, this frequency may change during
exercise, or due to pathologies.  This raises the question of how the
frequency of oscillations may change our model, in particular how the
frequency effects the deformation and force splitting of the Kelvin
bodies.

Figures \ref{peakvalues_u_freq.ps} - \ref{peakvalues_a_freq.ps} show the
results of simulations where the frequency of oscillations is altered.  
In these plots we used the values $f=10^{-4}, 2 \times 10^{-3}, 3\times
10^{-3}, 4\times 10^{-3},5\times 10^{-3}, 10^{-2},2 \times 10^{-2},5
\times 10^{-2},7 \times 10^{-2}, 10^{-1}, 2.5 \times 10^{-1}, 7.5\times
10^{-1},1, 10$ Hz. All other parameter values for body 1 and body 2 are
baseline values.  The peak deformation (and the peak force splitting) for
each of the simulations is divided by the deformation (and the force
splitting) for steady flow.

Figure \ref{peakvalues_u_freq.ps} shows that if the frequency of
oscillations is very small, then the deformation in oscillatory flow is
essentially the same as in steady flow.  As the frequency increases, the
deformation decreases, because the dashpot is unable to fully deform once
the oscillations become fast enough.  This is the result of the sign of
the force changing quickly, therefore the force not acting for
sufficiently long periods of time to fully stretch the dashpot.  Once the
frequency reaches a critical value, $f_{crit} \approx 10^{-2}$, the
deformation is independent of the frequency of oscillations.  This signals
the frequency at which the deformation is entirely due to the springs.  
(Deformation of the springs is instantaneous when force is applied.)

Just like in the one Kelvin body case, $f_{crit}=10^{-2}$ \cite{Ba}, two
orders of magnitude below the physiological value.  This implies that
endothelial cells exposed to purely oscillatory flow will only deform to a
small fraction (approximately one third) of deformation possible in steady
flow.  If there is a threshold value of deformation that permits cells to
align and go through other significant physiological changes, it is very
likely that purely oscillatory flow will not be able to reproduce these
effects.  Also, it is interesting to note that the deformation for
frequencies $f>f_{crit}$ is the same for one Kelvin body as for two.  
This implies that any number of Kelvin bodies coupled in parallel will not
deform more than this value, therefore a model consisting of n Kelvin
bodies in parallel would retain the feature that $f_{crit}$ is well below
physiological values.

We have heuristically argued before that for all values of viscosity
$\eta_{12} > 100 $ the dashpot is not able to react in oscillatory flow,
because the direction of the force changes very quickly, and there is not
enough time for the dashpot to deform.  This suggests that the
relationship between frequency and deformation has to be examined.  
Figure \ref{peakvalues_u_freq_eta3.ps} shows these simulations.  We
observed in the previous plot, (Figure \ref{peakvalues_u_freq.ps}) that as
the frequency increases, the overall deformation decreases to a constant
value.  This can be explained, if very small frequency oscillations
allowed sufficient amount of time for the dashpot to react, but once the
frequency increased to about $f=0.01$ Hz, the dashpot did not have time to
respond at all.  Beyond this frequency all the deformation is due to the
springs.  Our figure confirms that by making the dashpot very inviscid
(decreasing $\eta_{12}$ below 100) allows the dashpot to react much more
quickly, and we see larger deformations for given frequencies. When the
frequency is very low, then the deformation is similar to the deformation
in steady flow for any value of the viscosity, therefore the dimensionless
value (peak deformation in oscillatory flow divided by the deformation in
steady flow) is near 1.  When the frequency is very large, then the
deformation is always a constant around 0.333 independently of the
frequency or the viscosity. Between these two regimes there is a range of
frequencies for which increasing the frequency means decreasing the
deformation.  Interestingly, for small values of the dashpot viscosity the 
normalized peak deformation is bi-sigmoidal.  

Figure \ref{peakvalues_a_freq.ps} describes the dependence of the peak
force splitting on the frequency of oscillations.  As before, the largest
amplitude at the steady state is taken for the appropriate value of the
frequency, and this amplitude is divided by the steady state value in
oscillatory flow.  As we can see, the frequency of oscillations does not
change the fact that the peak force split in oscillatory flow is always
the same as the steady state force split in steady flow.

We must also examine the case that all parameters in the second body are
changed, because contradictory predictions could be made based on changing
individual parameters only.  Figures \ref{factor10_u_s.ps}
-\ref{factor10_aF_o.ps} reveal that regardless of the type of flow, the
deformation decreases if the parameters are increased, and the deformation
is increased if all the parameters are decreased.  Similarly, in either
type of flow the force in body one decreases if the parameters are
increased in body two, and the force increases in body 1, if the
parameters in the other body are decreased.

\subsection{Network simulations} 

We represent endothelial cells as a network of viscoelastic bodies.  
Each part of the cell we model, namely, the transmembrane proteins, flow
sensors, cytoskeletal elements and the nucleus are thought of as Kelvin
bodies with different parameter values for the spring constant and the
dashpot viscosity.  The baseline values which we use for actin filaments
are: $\mu_{02}=50$ Pa, $\mu_{12}=100$ Pa and $\eta_{12}=5000$ Pa s.  
For the nucleus the spring constants are four times the baseline values
and the viscosity of twice the baseline: $\mu_{02}=200$ Pa,
$\mu_{12}=400$ Pa and $\eta_{12}=10000$ Pa s. These estimates are based
on experimental measurements by Guilak et al. \cite{GT}.  The parameter
values for the transmembrane proteins (flow sensors, ion channels,
attachments to the substrate) based on \cite{BZ} are: $\mu_{02}=100$ Pa,
$\mu_{12}=200$ Pa and $\eta_{12}=7.5$ Pa s.  The parameter values are
calculated in Section \ref{parameter_calc} and summarized in Table
\ref{parameter_values}.

\begin{figure}[h]
\includegraphics[width=0.47\textwidth]{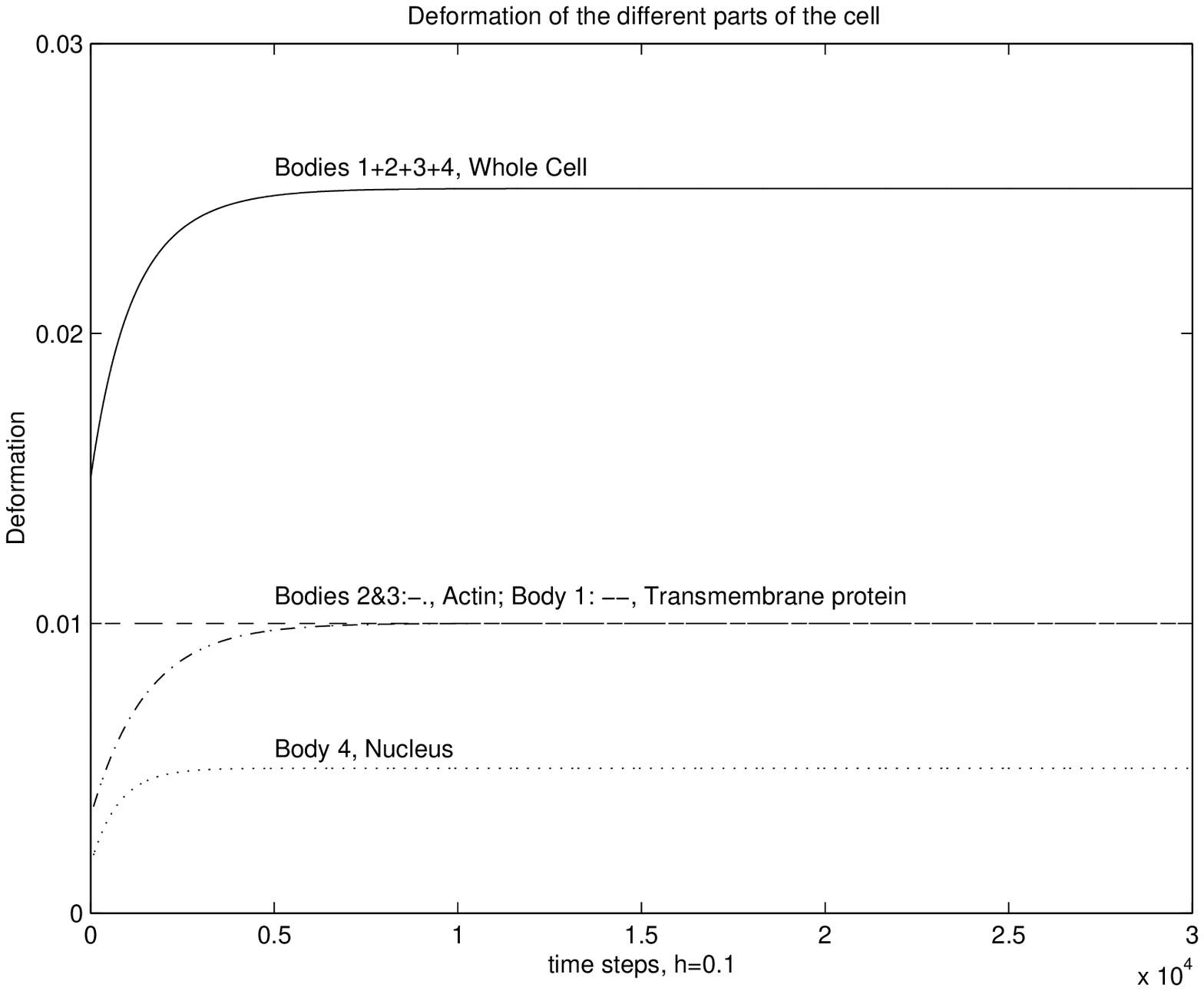} \quad
\includegraphics[width=0.47\textwidth]{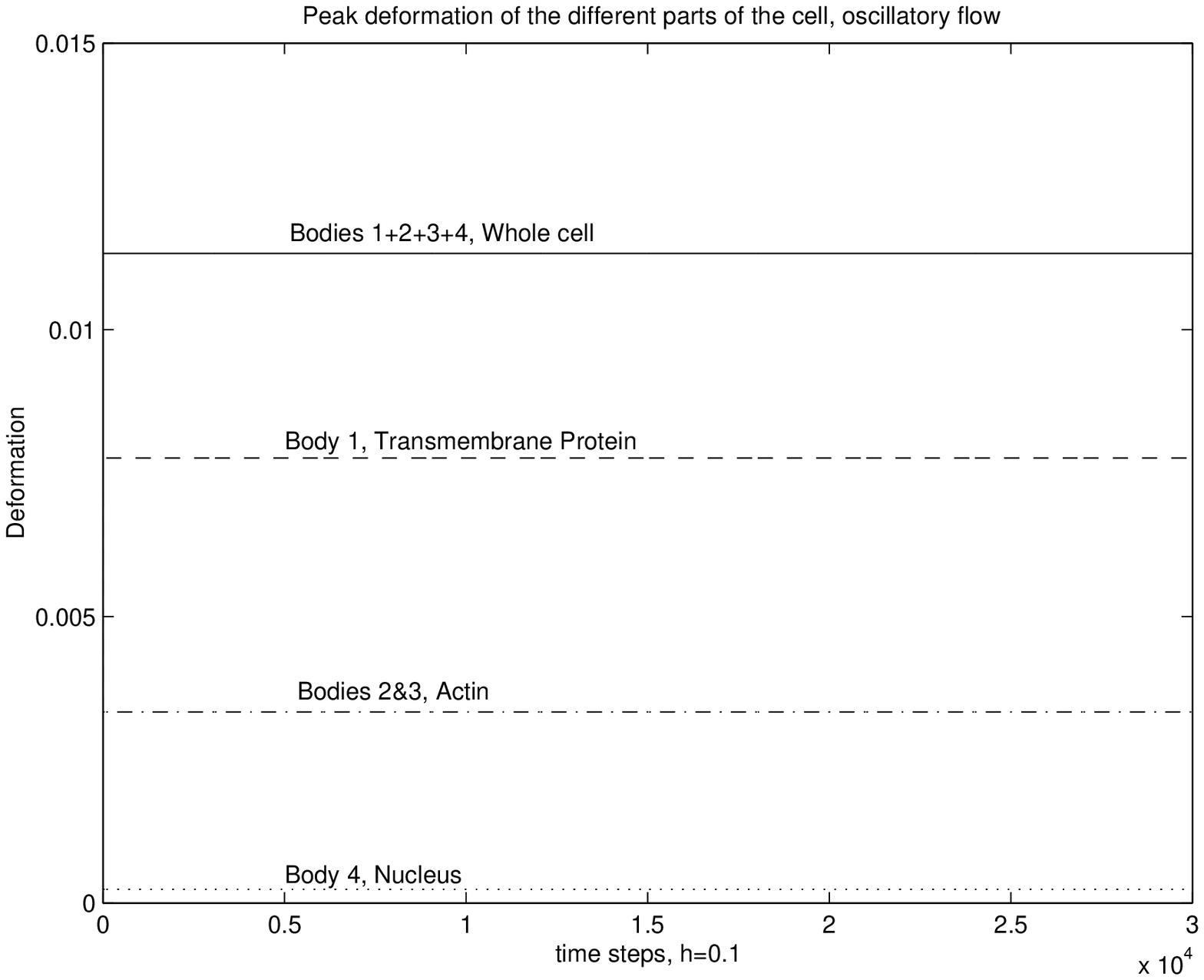} 
\caption[Deformation of network I.]{Deformation of
different parts of a simple endothelial cell in steady and oscillatory.}
\label{man_u_parts.ps} \end{figure}
 
Figure \ref{man_u_parts.ps} shows numerical simulations for the simple
four-body model of the endothelial cell (Figure \ref{Network I.}) in
steady and in oscillatory flow.  In steady flow transmembrane proteins
deform immediately, followed by the deformation of the nucleus, then the
deformation of the actin filaments.  This behavior is consistent with
transmembrane proteins having the smallest, and actin filaments having the
largest dashpot viscosities.  Transmembrane proteins and actin filaments
reach the overall deformation whereas the nucleus only deforms slightly.  
Because our model is linear, the overall cell deformation is the sum of
the deformation of the components.  The time constants are consistent with
experimental data in which transmembrane proteins such as flow sensors and
ion channels respond to flow very quickly, on the order of seconds, and
cytoskeletal reorganization is the slowest, in fact, it happens on the
time scale of many hours.

Next, we examine the same four-body model in oscillatory flow.  Only peak
values of the deformation are plotted, as before.  The overall deformation
is much smaller here than in steady flow, and the steady state of
deformation is attained immediately (within 2-3 seconds).  The nucleus
deforms the least amount again, and here the largest deformation is that
of the transmembrane protein. 

Because in our four-body model the nucleus and the flow sensor are both 
modeled as a single Kelvin body, the force acting on each of them is a 
constant, $F=1$.  The same parameter values characterize the two actin 
filaments which are modeled as two Kelvin bodies connected in parallel, 
therefore the forces acting on them are also equal, $aF = (1-aF) = 0.5$
This is true for steady as well as in oscillatory flow.  

\begin{figure}[h!]
\includegraphics[width=0.47\textwidth]{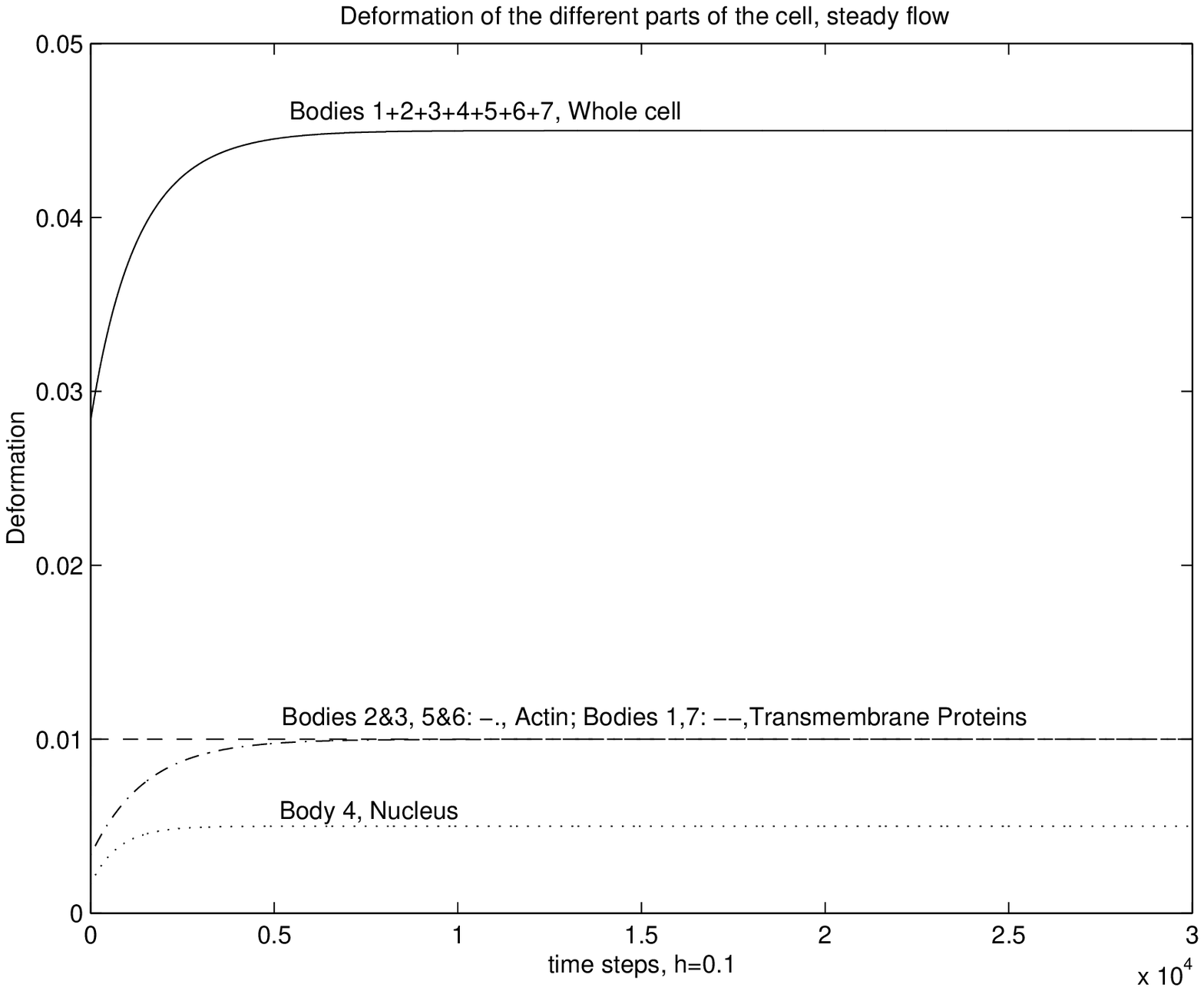} \quad
\includegraphics[width=0.47\textwidth]{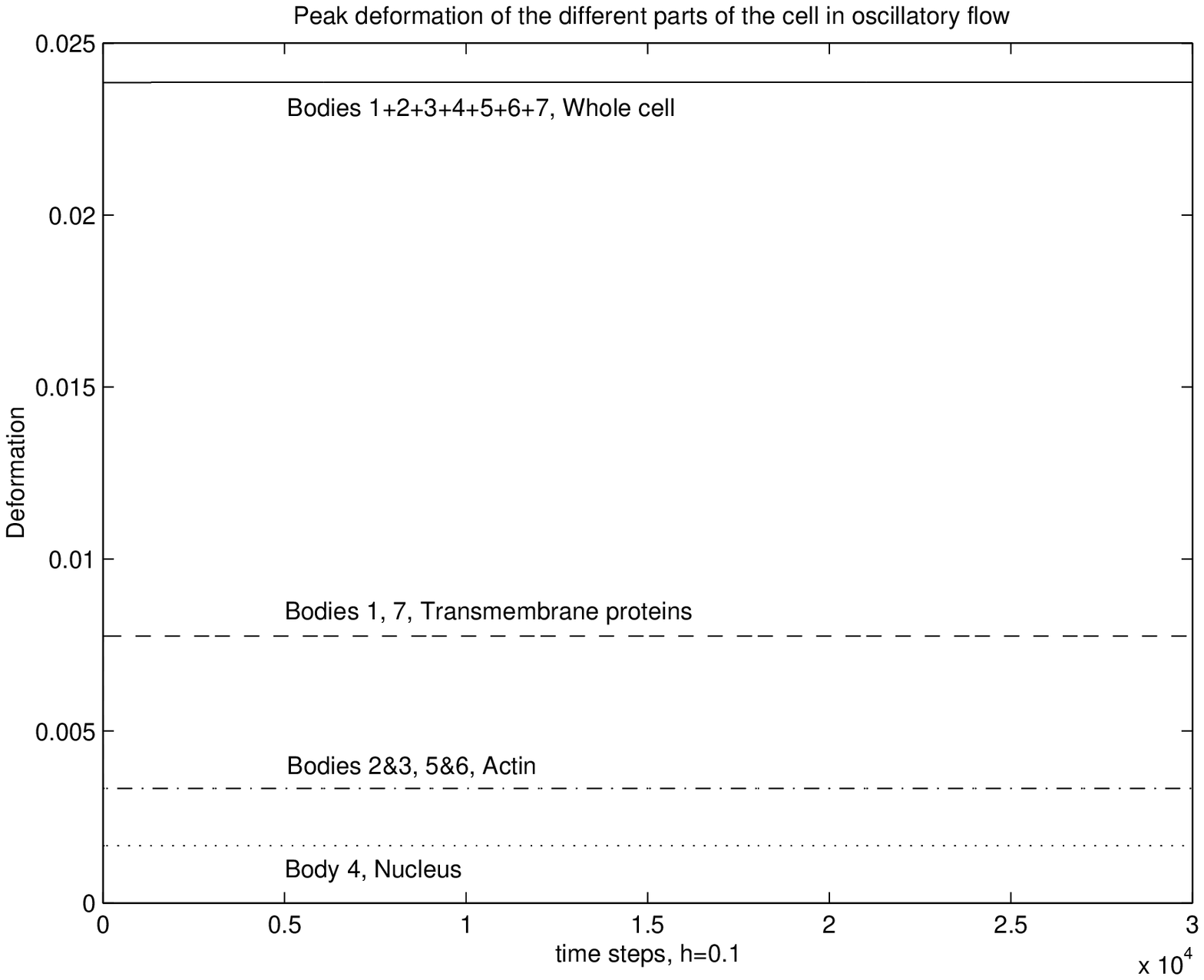}
\caption[Deformation of network II.]{Deformation of
different parts of a simple endothelial cell in steady and in 
oscillatory flow.}
\label{manas_u_parts.ps} \end{figure}

The final graph, Figure \ref{manas_u_parts.ps} and shows the response of
the seven Kelvin-body model (depicted in Figure \ref{Network II.})  in
steady and oscillatory flow.  The results are very similar
to the simulations of the four-body model. The forces for the bodies are
all constants: for the transmembrane proteins and nucleus (which are
modeled as single Kelvin bodies) the force is equal to one, and for acting
filaments, modeled as Kelvin bodies in parallel the force is split
equally.

\clearpage
\section{Discussion}

Our numerical results included a parameter sensitivity analysis and
simulations for two simple model networks of endothelial cells.  The
results of the simulations have been interpreted so far only in terms of
the elements of the Kelvin body, but it has not been discussed what their
implications are for endothelial cells.  This section describes the
conclusions we can draw regarding endothelial cell behavior from the
numerical results.  We also briefly mention further work that can be done
to improve the current model.

The most notable difference between the effect of steady and oscillatory
flow is that in oscillatory flow deformation tends to be much smaller than
in steady flow for the same viscoelastic materials, as it is depicted in
Figures \ref{peakvalues_u_m02.ps}, \ref{peakvalues_u_m12.ps} and
\ref{peakvalues_u_eta.ps}. The importance of the large deformation
difference in the two types of flow is clear when we recall experimental
results in which endothelial cells exhibit elongation and certain other
biological responses (such as activation of flow-sensitive ion channels
and mobilization of intracellular calcium) depending on the specific form
of the shear stress cells are exposed to.  Dewey et al. (\cite{Ba},
original source:\cite{D})  demonstrate that endothelial cells exhibit some
morphological responses when exposed to large shear stress, for example in
steady or pulsatile flow, but not when their exposure to shear stress is
below a certain level, for example, in oscillatory flow.

According to our model, the forces acting on each part of an endothelial
cell, regardless of its viscoelastic properties elicits a much smaller
response in oscillatory flow.  This can be particularly important for a
flow sensor, as discussed by Barakat \cite{Ba}, which in oscillatory flow
may not be stimulated sufficiently to initiate a signaling cascade for
downstream responses to the flow.  Our model shows, that even given an
signal from the flow sensor, the deformation of materials whose
viscoelastic properties are consistent with those of the nucleus and
cytoskeletal elements would be much smaller in oscillatory flow.

Another observation we can make (still based Figures
\ref{peakvalues_u_m02.ps}, \ref{peakvalues_u_m12.ps} and
\ref{peakvalues_u_eta.ps}) is that the viscoelastic materials which display
the largest difference in deformation can be characterized by small spring
constants but a viscosity coefficient which is at least $\eta \approx 10^2$.  
This characterization would allow us to estimate parameter values for
microtubules and intermediate filaments for future simulations, because we
have experimental observations of the relative behavior of actin filaments,
microtubules and vimentin \cite{J}.

We have investigated how the frequency of oscillations effects the
deformation.  Our simulations show that within physiological conditions
(even accounting for conditions where the frequency of the flow changes one
order of magnitude) the normalized deformation is small, unless the
viscosity coefficient becomes very small.  This implies that materials
whose viscosity is small, for example various transmembrane proteins, such
as attachments to the substratum or cell-cell adhesions, are less able 
to differentiate between oscillatory and steady flow than the cytoskeleton 
or the nucleus.  (As discussed above, the cytoskeleton and the nucleus 
respond very distinctively to steady and to oscillatory flow.)  

Now let us examine the network simulations of endothelial cells.  Both
networks hugely oversimplify the complex connections between endothelial
cells, therefore no detailed conclusions should be drawn from our results,
rather we ought to make simple qualitative observations.  More
sophisticated networks need to be created in the future, and the main
significance of our current simulations is to demonstrate the sort of
models we are able to create with our elements now (namely n Kelvin bodies
in parallel, coupled in series with single Kelvin bodies).  

Qualitatively, the time scale of deformations as predicted by our model
fits with experimental observations.  Instantaneous response can be seen
from transmembrane proteins (flow sensors, focal adhesions, cell-cell
adhesions), followed by the deformation of the nucleus, and finally, the
changes in the cytoskeleton.  Flow sensor and ion channel deformation can
occur within seconds whereas changes in the gene expression takes hours,
and the cytoskeletal re-organization takes place over the span of about a
day.  These time scales do not compare well with the prediction of our
model, which does predict flow sensor and ion channel deformations within
seconds, but it also predicts the response of the nucleus and the actin
filaments to be on the order of minutes.  Clearly, our model must be 
modified, and in order to obtain quantitatively accurate results, the 
effects of the biochemical signaling pathways must be considered as well.  

Our model is the first step in creating more realistic model networks of
endothelial cells.  The current project can be extended in a number of ways
to lead to better approximations of the cellular response to flow.  The
extension promising the most extensive changes in the results is the
assumption that the parameters describing the Kelvin bodies can alter in
time too.  This new assumption would turn our linear system of equations to
a nonlinear system with possibly much more complicated dynamics.  Other
changes to incorporate in the model are: formulating the equations for new
networks based on combinations of bodies which our system does not
describe, assuming that the forcing function is applied gradually (there is
a "ramp" for the force).  Allowing the force to act gradually on the bodies
could be compared to experimental results in which cell responses are
altered by applying shearing forces instantaneously or gradually.  
Although even the understanding of viscoelastic Kelvin bodies in parallel
lead to new insights on the hypothesized mechano-transduction of flow
induced shear stress in endothelial cells, further improvements and
extensions are necessary to clarify the role of this mechanism.


\appendix


\newpage
\thispagestyle{myheadings}
\markright{  \rm \normalsize APPENDIX. \hspace{0.5cm}
  MATHEMATICAL MODELS IN BIOLOGY}
\addcontentsline{toc}{chapter}{\protect\numberline{}{\bf Appendix}}
\chapter{Receptor model} \label{app_aerotax}
\thispagestyle{myheadings}

An important question is whether there is biological evidence of the
signal transduction pathway allowing the existence of the proposed
turning rates that depend on the oxygen concentration and the
spatial gradient of oxygen. In \cite{TJ}, Taylor and Johnson propose
a model that aims to explain experiments in which receptors are
rewired to produce inverse responses.  Using their model for a
receptor, it is easy to explain how the turning rates of our model
could be a result of a simple two-state receptor.

\begin{figure}[h]
\centerline{\includegraphics[width=0.6\textwidth]{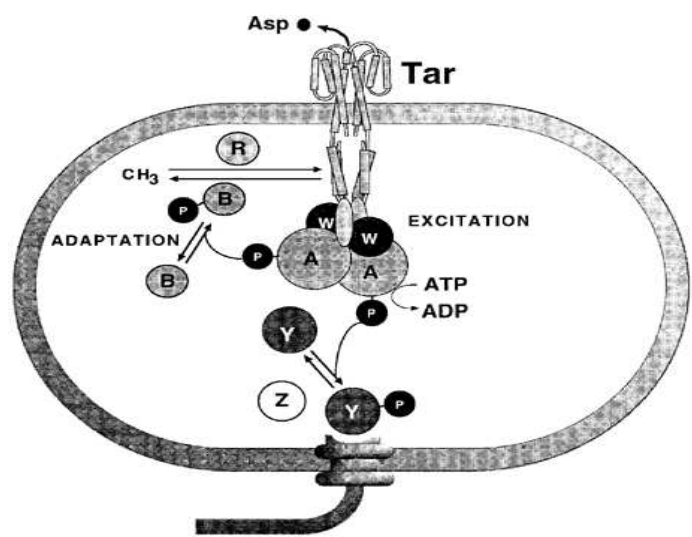}}
\caption[Tar receptor.]{The four-helix bundle of Tar receptor. (Figure
from Taylor and Zhulin, \cite{TZ}.)}
\end{figure}
  
It is known that part of the Tar (and Tsr) receptors is a four-helix
bundle outside the cell, and as ligands bind to the receptor, the
bundle goes through a conformational change.  By representing the
receptor as two independently moving parts, we have a simplified
model ('piston model') of the receptor.  Each part of the receptor  
is able to move in the $z$ direction as a function of the proton  
motive force, PMF.  One part of the receptor is able to change its
position fast in response to PMF, while the other part reacts
slower.  When the two parts of the receptor in the model are at the
same position, they lock (representing the conformation change in
the real receptor), and this initiates tumbling.  When the two parts
of the receptor are not in the same position, the cell swims
smoothly.

\begin{figure}[h]
\centering
\includegraphics[width=0.5\textwidth]{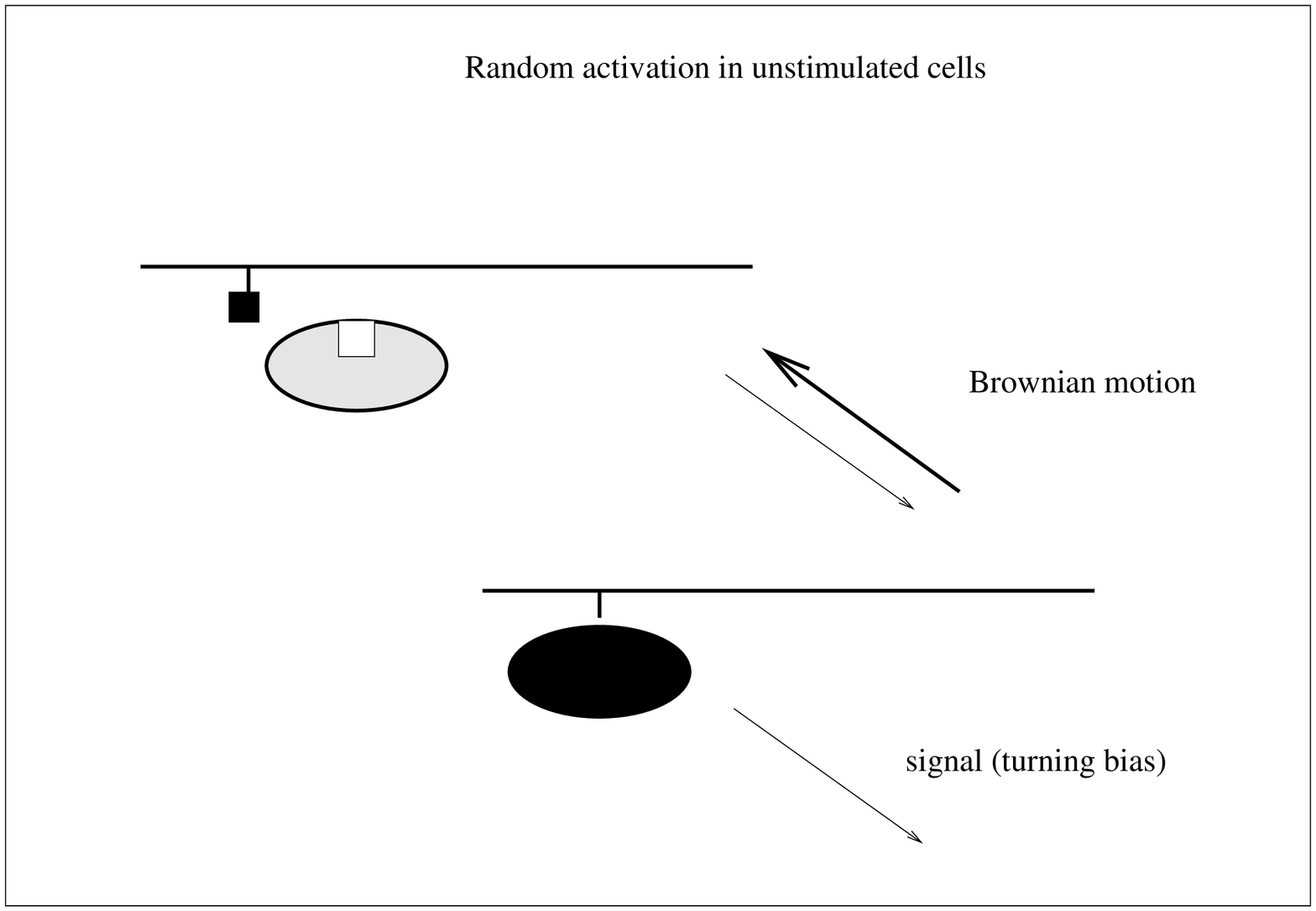} \\
\includegraphics[width=0.5\textwidth]{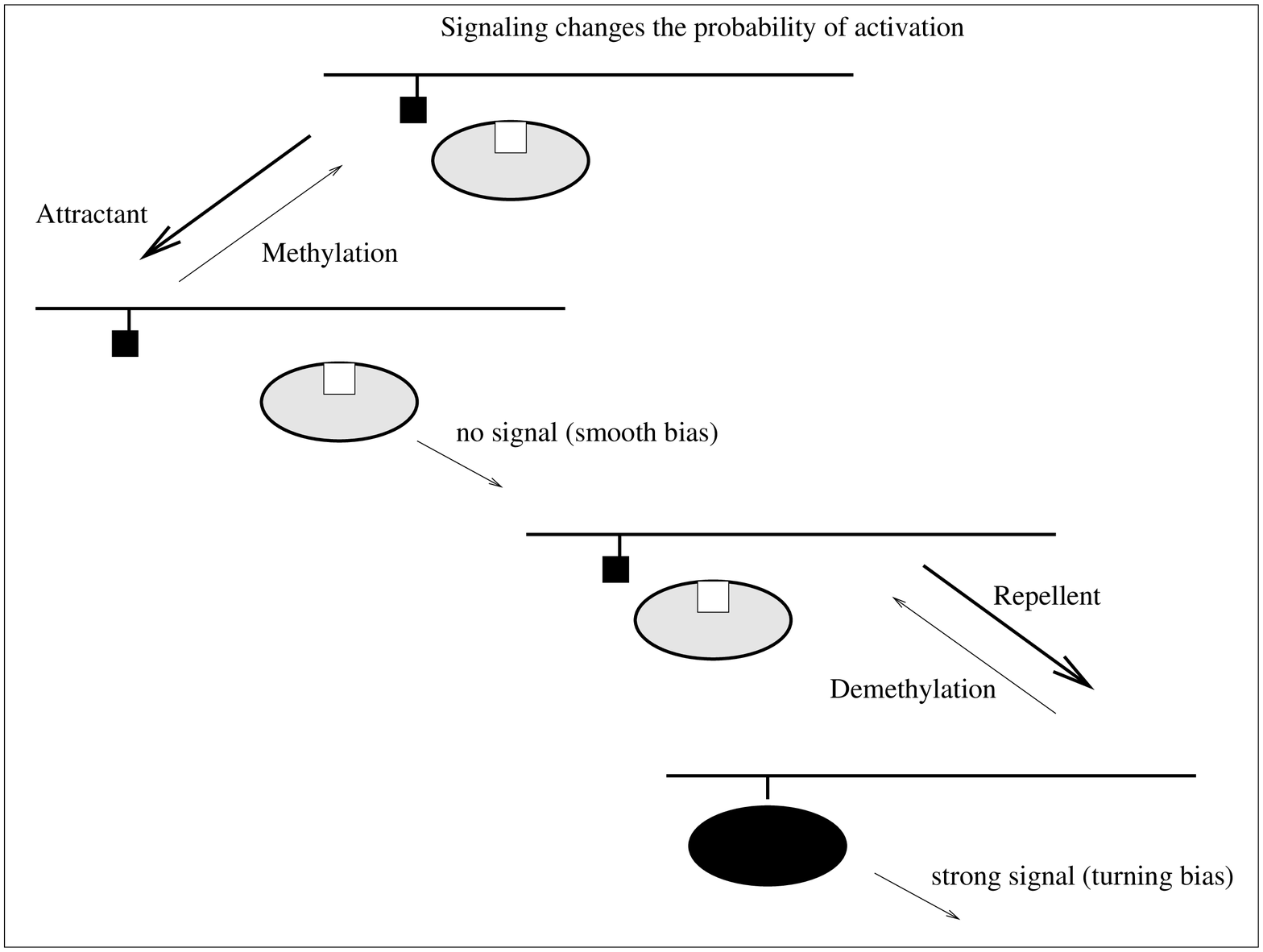}
\caption[Piston model of a receptor.]{Piston model of a receptor. 
(Figure
based on Taylor and Johnson, \cite{TJ}.)}
\label{tay_john}
\end{figure}

In ground state, the two parts are displaced from one another.
Thermal energy is able to randomly move the parts, which causes   
tumbling in unstimulated cells if the parts are moved together, and
it promotes smooth swimming if the parts are moved further apart.
Attractant binding to the receptor also moves the fast and slow
parts further, which makes smooth swimming more likely.  Adaptation
returns the receptor to the ground state.

In our mathematical model such a receptor could be responsible for
the same oxygen concentrations triggering different signals.  When a
bacterium is outside the optimal oxygen concentration, the cell is
in ground state; therefore, it has its baseline turning frequency.
As the cell enters the favorable oxygen concentration, the
fast-changing part moves; thus, the two parts get further from each
other, and the cell swims smoothly inside the band.  By the time the
cell gets across the band, the slow-changing part moves too, and the
two parts lock which causes tumbling.  The tumbling turns the cell
back into the favorable concentration, and the fast-moving part
changes position again, leading to smooth swimming.

\begin{figure}[h!]
\centerline{\includegraphics[width=0.5\textwidth]{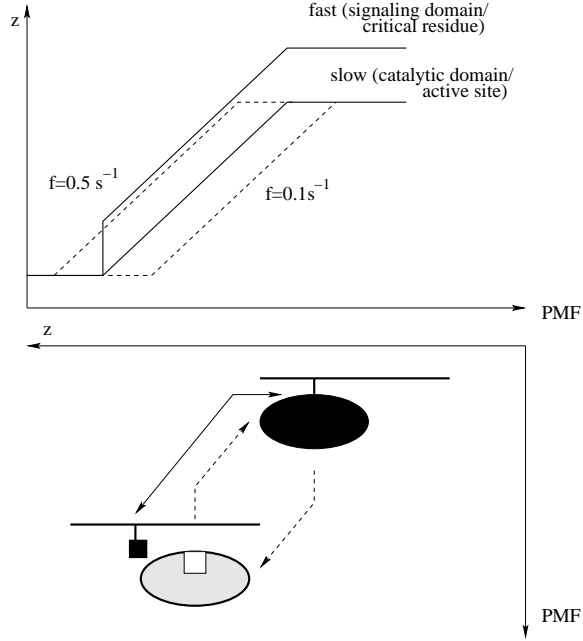}}
\caption[Fast and slow states of a receptor.]{Fast and slow states of
a receptor. The top figure shows how the receptor switches from high
turning frequency to low turning frequency and vice versa.  The bottom
figure shows a different representation of the same process. (Bottom
portion of the figure based on Taylor and Johnson, \cite{TJ}.)}
\label{f-s_rec.eps}
\end{figure}

Consider $z_f$: $z$ coordinate of fast-changing part, $z_s$:
$z$ coordinate of slow-changing part.
\begin{eqnarray}
\tilde{z_f} = z_f^{(0)} + c_1 p \label{2stateqn1} \\
\tilde{z_s} = z_s^{(0)} + c_1 p = z_f^{(0)} - \Delta z + c_1 p
\label{2stateqn2}
\end{eqnarray}

In these equations $\tilde{z_f}$ and $\tilde{z_s}$ are the
equilibrium positions of the fast-moving and slow-moving parts,
respectively.  In ground state, the equilibrium positions are not
the same, as discussed above. $z_f^{(0)}$ and $z_s^{(0)}$ are
constants.  $p$ is the proton motive force (which is proportional to
the oxygen concentration).  $$ p = \pm k t + c_0 $$ We can assume
that the PMF is a linear function of time, since the fully developed
oxygen gradient is linear in space, and it is observed by the  
bacteria swimming through it as oxygen changing linearly in time.

We can assume that the fast-moving part is just a function of the
proton motive force (i.e. its position changes immediately as it
detects changes in the proton motive force.)  We assume that the
slow-moving part changes its position after a delay.  This results
in the following two equations:
\begin{eqnarray} 
z_f \approx \tilde{z_f} (p(t)) \label{fast} \\
\frac {d z_s} {d t} = \frac {1} {\tau} (\tilde{z_s}(p(t)) - z_s)
\label{slow}
\end{eqnarray}

The first equation is an algebraic equation, and we can solve it by
substituting the expression for the proton motive force into
(\ref{fast}). This gives the equation: $$ z_f = [z_f^{(0)} + c_0 c_1]
\pm c_1 kt $$

This shows that the coordinate of the fast part will be increasing as it
moves up the concentration gradient, and it will be decreasing as the
cell swims down the gradient.  Now we can solve (\ref{slow}).
$$ \frac {d z_s} {d t} = \frac {-z_s} {\tau} + \frac {z_f^{(0)} -
\Delta z + c_1 p} {\tau} $$

We have to solve the homogeneous equation on the left-hand side,
then find the particular solution, $z_{sp} = A + B t$, using the
right-hand side of the equation:
$$ \frac {d z_s} {d t} + \frac {-z_s} {\tau} = \frac {z_f^{(0)}-
\Delta z + c_1( \pm kt + c_0)} {\tau} $$

Solving the homogeneous equation we get: $z_{sh} = z_0 e^{- t/ \tau}$.  
In the particular solution, we get $A= z_f^{(0)} + c_1 c_0 \mp c_1 k
\tau - \Delta z$; $ B = \pm c_1 k$.  This gives the full particular
solution, $z_{sp} = z_f^{(0)} + c_1 c_0 \mp c_1 k \tau - \Delta z \pm
c_1 k t$.  By making the appropriate substitution, this results in the
full solution, $$ z_s = z_f + [\mp c_1 k \tau - \Delta z] + z_0 e^{\frac
{-t} {\tau}} $$

If the adaptation time, $\tau$, is small compared to the run (which is
our assumption in the model), then the last term is negligible. Now we
can see how the receptor will work by examining the position of the
fast- and slow-moving parts.  The position of the slow-moving part,
$z_s$ is behind the position of the fast part.  To see the dynamics of
the two parts, the solution to $z_s$ can also be re-written as: $$ z_s =
z_f^{(0)} + c_1 c_0 - \Delta z \pm c_1 k (t-\tau) + z_0 e^{\frac {-t}
{\tau}} $$

As long as the time, $t$ is less than the adaptation time $\tau$, the
slow part will be behind the fast one, but after the time reaches
$\tau$, the sign of $c_1 k (t - \tau)$ changes, and the difference in
the position of the fast and slow part decreases.  This difference is
maintained until the cell turns down the gradient, at which point the
z-coordinate of the fast moving part will change quickly again, followed
by the z-coordinate of the slow moving part.

This simple model of the receptor based on proton motive force is
sufficient to explain the turning frequencies assumed in the full
mathematical model of aerotaxis.


\newpage
\thispagestyle{myheadings}
\markright{  \rm \normalsize APPENDIX. \hspace{0.5cm}
  MATHEMATICAL MODELS IN BIOLOGY}
\chapter{Analytical calculations} 
\thispagestyle{myheadings}

\section{Asymptotic approximation} \label{approx_sol}

We want to solve the system of equations
\begin{eqnarray*}
\frac {dM} {dt} = m + \lambda (-k_a (l_1) M + k_d A)  \\
\frac {dA} {dt} = -r A + \lambda ( k_a (l_1)  M - k_d A)
\end{eqnarray*}

with initial conditions:
\begin{eqnarray*}
M(0) = \frac{m k_d} {r k_a(l_0)}, \ \ \ A(0)=\frac {m} {r}
\end{eqnarray*}

On the fast time scale, i.e. when $\tau = \lambda t$, $$\frac 
{dM} {d\tau} + 
\frac{dA} {d\tau} = 0 $$ which implies 
$A+M=C$ for some constant, $C$.  We substitute this expression into the 
original system of equations to get:
\begin{eqnarray*}
\frac {dM} {d\tau} = -k_a (l_1) M + k_d(C- M) \\ 
\frac {dA} {d\tau} =  k_a (l_1)  (C-A) - k_d A 
\end{eqnarray*}

The initial conditions on the fast time scale are the same as the initial 
conditions of the original equations, 
\begin{eqnarray*}
M_0=M(\tau=0) = \frac{m k_d} {r k_a(l_0)} \\A_s= A(\tau=0)=\frac {m} {r}
\end{eqnarray*}

The solution to these equations is given by:
\begin{eqnarray*}
A(\tau) = \big(A_s - \frac {k_a(l_1) C} {r_f} \big) e^{-r_f \tau} + 
\frac {k_a (l_1) C} {r_f} \label{a_fast}\\
M(\tau) = \big(M_0 - \frac {k_d C} {r_f} \big) 
e^{-r_f \tau} + \frac {k_d C} {r_f} \label{m_fast}
\end{eqnarray*}
Where we defined the fast time scale to be $$r_f = k_a(l_1) + k_d.$$

Now we use the fact that the sum of $A$ and $M$ is always a constant, in 
particular, 
$A(0) + M(0) = C.$  
\begin{eqnarray*}
\big(A_s - \frac {k_a(l_1) C} {r_f} \big) + \frac {k_a (l_1) C} {r_f} + 
\big(M_0 - \frac {k_d C} {r_f} \big) + \frac {k_d C} {r_f} = C\\
C = \frac {m} {r} \big(1 + \frac {k_d} {k_a(l_0)} \big)
\end{eqnarray*}

Two new constants are useful:
\begin{eqnarray*} 
A_1 = \frac {k_a(l_1) C} {r_f} =  \frac {m} {r} \frac {1+(k_d/k_a(l_0))} 
{1+(k_d/k_a(l_1))} \\
M_1 = \frac {k_d C} {r_f} = \frac {m} {r} \frac {k_d} {k_a(l_1)} \frac 
{1+(k_d/k_a(l_0))}{1+(k_d/k_a(l_1))}
\end{eqnarray*}

Substituting the expressions for $A_1$ and $M_1$ into eqn. \ref{a_fast} and 
\ref{m_fast} we have obtain $A$ and $M$ on the fast time scale:
\begin{eqnarray}
A(\tau) = (A_s - A_1) e^{-r_f \tau} + A_1 \label{af_sol} \\
M(\tau) = (M_0 - M_1)e^{-r_f \tau} + M_1 \label{mf_sol}
\end{eqnarray}

We turn to the slow time scale, and we note that now $A = \frac {k_a(l_1)} 
{k_d} M$.  Substituting this expression into the original equations gives us:
\begin{eqnarray}
\frac {dM} {dt} = m   \label{ms} \\
\frac {k_a(l_1)} {k_d} \frac {dM} {dt} = -r \frac {k_a(l_1)} {k_d} M \label{ff}
\end{eqnarray}

By adding equations \ref{ms} and \ref{ff} we get:
\begin{eqnarray} 
\frac {dM} {dt} = \frac { m k_d} {k_a(l_1)+ k_d} - \frac {r k_a(l_1)}
{k_d + k_a(l_1)} M \nonumber \\
M(t) = \bar{c} e^{(-r_s t)} + \frac { m k_d} {(k_a(l_1)+ k_d) r_s} 
\nonumber
\\
M(t) = \bar{c} e^{(-r_s t)} + M_2 \label{ms_sol}
\end{eqnarray}

We have defined $$r_s = r \frac {k_a(l_1)} {k_d + k_a(l_1)} $$ and $$M_2 =
\frac {m k_d} {r k_a(l_1)}.$$

Similarly, we have an equation for $A$:
\begin{eqnarray}
A(t) =\frac {k_a(l_1)} {k_d} \bar{c} e^{(-r_s t)} + \frac{m} {r} 
\label{as_sol}
\end{eqnarray}

The solutions on the fast and slow times scale must match, so the limit of 
eqn \ref{af_sol} as $\tau \rightarrow \infty$ must be equal to eqn 
\ref{as_sol} at zero.  Similarly, $\lim_{\tau\to\infty}M_f = M_s(t=0)$. 
Matching the solutions allows us to find the expression for $\bar{c}.$
\begin{eqnarray}
\frac {k_a(l_1) \bar{c}} {k_d} + A_s = A_1 \\
\bar{c} = \frac {m} {r} \Big[\frac{k_d} {r_f} \big(1+\frac {k_d} {k_a(l_0)} 
\big) - \frac{k_d} {k_a(l_0)} \Big] = M_1 - M_0     
\end{eqnarray}

The full solution is the sum of the fast and slow terms with the common limit 
($A_1$ and $M_1$)  subtracted.  
\begin{eqnarray} 
A(t) = A_s + (A_s - A_1) e^{-r_f \lambda t} +(A_1 - A_s) e^{-r_s t}
\label{a_approx} \\
M(t) = M_2 + (M_0 - M_1) e^{-r_f \lambda t}+ (M_1 - M_2)e^{-r_s t} 
\label{m_approx}  
\end{eqnarray}

\section{Steady state solution}
\label{ss_calc}

We start with the four equations describing the two compartment model of
the signal transduction system again.
\begin{eqnarray}
\frac {dM_1} {dt} = m + \lambda (-k_a(l_1) M_1 + k_d A_1) + k_1(M_2-M_1)
\label{m1} \\
\frac {dA_1} {dt} = -r A_1 + \lambda ( k_a(l_1) M_1 - k_d A_1) +
k_2(A_2-A_1) \label{a1} \\
\frac {dM_2} {dt} = m + \lambda (-k_a(l_2) M_2 + k_d A_2) - k_1(M_2-M_1)
\label{m2} \\
\frac {dA_2} {dt} = -r A_2 + \lambda ( k_a(l_2) M_2 - k_d A_2) -
k_2(A_2-A_1) \label{a2} \\
M_1(0) = M_2(0) = \frac {m} {r} \frac {k_d} {k(l_0)} \ \
A_1(0) = A_2(0) = \frac {m} {r} \label{ic}
\end{eqnarray}

At the steady state the left hand side of these equations is zero.  We
start finding the expressions for the steady state by adding equations
\ref{m1} and \ref{m2}.  In the following calculations we write $A_1$,
$A_2$, $M_1$ and $M_2$ instead of $A_{1s}$, $A_{2s}$, $M_{1s}$ and
$M_{2s}.$
\begin{eqnarray}
2m - \lambda k_a(l_1)M_1 - \lambda k_a(l_2) M_2 + \lambda k_d (A_1+ A_2) 
=
0 \\
k_a(l_1) M_1 + k_a (l_2) M_2 = \frac {2 m} {\lambda} + k_d (A_1 + A_2)
\label{m1+m2}
\end{eqnarray}

Next we add equations \ref{a1} and \ref{a2}.
\begin{eqnarray}
-r(A_1 + A_2) + \lambda (k_a(l_1) M_1 + k_a(l_2)M_2) - \lambda k_d (A_1 
+
A_2) = 0  \\  
(A_1 + A_2) (\lambda k_d + r) = \lambda (k_a(l_1) M_1 + k_a (l_2) M_2)
\end{eqnarray}

We substitute the expression for $k_a(l_1) M_1 + k_a (l_2) M_2$ from 
\ref{m1+m2}:
\begin{eqnarray}
A_1 + A_2 = \frac {\lambda } {r + \lambda k_d} \Big( \frac{2m} {\lambda} 
+
k_d (A_1 + A_2) \Big) \nonumber \\
A_1 + A_2 = \frac {2m } {r+ \lambda k_d} + \frac {\lambda k_d} {r +
\lambda k_d} (A_1 + A_2)
\end{eqnarray}

Solving the equation for $A_1 + A_2$ we obtain
\begin{eqnarray}
A_1 + A_2 = \frac {2m} {r}
\label{a1s+a2s}
\end{eqnarray}

This allows us to express $A_2(A_1)=  \frac {2m } {r} - A_1.$

Now we return to \ref{m1+m2} to find an expression for $k_a(l_1) M_1 +
k_a (l_2) M_2$ explicitly.
\begin{eqnarray} k_a(l_1) M_1 + k_a (l_2) M_2 = \frac {2m} {\lambda} + 
k_d
\frac {2m} {r} \\
M_2 = \frac {2m} {k_a(l_2)} \Big( \frac {r+\lambda k_d} {\lambda r} 
\Big)
- \frac {k_a (l_1)} {k_a (l_2)} M_1 \label{m2(m1)}
\end{eqnarray}

Now we add equations \ref{m1} and \ref{a1}, and similarly, add \ref{m2}
and \ref{a2}.
\begin{eqnarray}
m-r A_1 + k_1 (M_2 - M_1) + k_2 (A_2 - A_1 ) = 0 \label{m1+a1} \\
m- rA_2 - k_1 (M_2 - M_1) - k_2 (A_2 - A_1) = 0 \label{m2+a2}
\end{eqnarray}

By subtracting \ref{m2+a2} from \ref{m1+a1} we arrive at
\begin{eqnarray}
A_1 - A_2 = \frac {2 k_1} {r + 2 k_2} (M_2 - M_1)
\label{a1-a2}
\end{eqnarray}

We want to express $A_1$ as function $M_1$, and this will allow us to 
find
$A_2(M_1)$.  In order to do this, we add \ref{a1s+a2s} and \ref{a1-a2}.
\begin{eqnarray*}
A_1 = \frac {m} {r} + \frac {k_1} {r+ 2 k_2} (M_2 - M_1)
\end{eqnarray*}

Now we use \ref{m2(m1)} to find both $A_1$  and $A_2$ as a function of
$M_1$, so we have $A_1(M_1)$, $A_2(M_1)$ and $M_2(M_1)$.
\begin{eqnarray}
A_1 = \frac {m} {r} + \frac {k_1} {r+2k_2} \Big( \frac {2m (r + \lambda
k_d)} {k_a(l_2) \lambda r} - \frac {k_a(l_1)} {k_a(l_2)} M_1 - M_1 \Big)
\label{a1(m1)}\\
A_2 = \frac {m} {r}- \frac {k_1} {r+2k_2} \Big( \frac {2m (r + \lambda
k_d)} {k_a(l_2) \lambda r} - \frac {k_a(l_1)} {k_a(l_2) }M_1 - M_1 \Big)
\label{a2(m1)}
\end{eqnarray}

We return to the equation \ref{m1} with its right hand side set to zero,  
and substitute $A_1(M_1)$, $A_2(M_1)$ and $M_2(M_1)$ from the equations
\ref{a1(m1)}, \ref{a2(m1)} and \ref{m2(m1)}, respectively, and solve the
equation for $M_1$.
$$ m+\lambda \Big[ -k_a(l_1)M_1 + k_d \Big(\frac {m} {r} + \frac {k_1}
{r+ 2 k_2}\big[\frac {2m (r+\lambda k_d)} {k_a(l_2)\lambda r} $$
$$- \frac {k_a(l_1) + k_a(l_2)} {k_a(l_2)} M_1 \big] \Big) \Big]
+ k_1 \Big( \frac{2m (r+\lambda k_d)} {k_a(l_2)\lambda r} - \frac
{k_a(l_1) + k_a(l_2)} {k_a(l_2)} M_1 \Big)  =   0 $$
$$ m+\frac { m \lambda k_d} {r} + k_1 \frac {2m (r+\lambda k_d)} 
{k_a(l_2)
\lambda r} \cdot \Big( \frac {\lambda k_d} {r+ 2k_2} + 1\Big) $$
$$ = \Big[\lambda k_a(l_1) + k_1 \cdot \frac {k_a(l_1)+ k_a(l_2)}
{k_a(l_2)}
\big( \frac{\lambda k_d} {r + 2k_2} + 1 \big) \Big] M_1 $$

If we simplify this expression, and substitute it back into 
\ref{a1(m1)},
\ref{a2(m1)} and \ref{m2(m1)} we arrive at the steady state solution:
\begin{eqnarray}
A_1 = \frac {m} {r} \cdot \Big[ 1 + \frac {r_1  k_1 k} {\lambda r_2
k_p + k_1 k_s( r_2 + \lambda k_d)}
\Big] \\
A_2 = \frac {m} {r} \cdot \Big[ 1 + \frac {- r_1 k_1 k } {\lambda r_2
k_p + k_1 k_s(r_2 + \lambda k_d)} \Big] \\
M_1 = \frac {m r_1} {\lambda r} \cdot \Big[\frac {r_2 ( \lambda 
k_a(l_2) + 2 k_1) + 2 \lambda k_d k_1} {r_2[\lambda k_p + k_1 k_s] +
\lambda k_d k_1 k_s} \Big]  \\
M_2 = \frac {m r_1} {\lambda r} \cdot \Big[\frac {r_2( \lambda
k_a(l_1) + 2 k_1) + 2 \lambda k_d k_1} {r_2 [\lambda k_p + k_1 k_s] +
\lambda k_d k_1 k_s} \Big]
\end{eqnarray}

where we have defined:
\begin{eqnarray*}
r_1= r+ \lambda k_d, \ \ r_2 = r + 2 k_2 \\
k = k_a(l_1)-k_a(l_2), \ \ k_s = k_a(l_1) + k_a(l_2) \\
k_p = k_a(l_1)k_a(l_2)
\end{eqnarray*} 

We want to verify that for $k_a(l_1) = k_a(l_2) =k_a$ we obtain the
original steady state.  This is clear for $A_1$ and $A_2$ by inspection,
but we need to simplify $M_1$ and $M_2$.  We show the calculations for
$M_1$.
\begin{eqnarray*}
M_1 = \frac {m r_1} {k_a \lambda r} \cdot \Big[ \frac {r_2 (\lambda
k_a + 2 k_1) + 2 \lambda k_d k_1} {r_2 [\lambda k_a + 2 k_1] + 2
\lambda k_d k_1} \Big] \\
M_1=\frac {m(r+\lambda k_d)} {\lambda r k_a} \\
\lim_{\lambda\to\infty} M_1 = \frac {m k_d} {r ka}
\end{eqnarray*}

Previously we have shown that if $k_1 =0$ then the two compartments
respond to stimulus as if they were not connected.   In the main text we
also mention that no important qualitative changes occur when $k_2=0.$  
In
this case our system of equations becomes:
\begin{eqnarray*}
A_1 = \frac {m} {r} \Big[ 1 + \frac {r_1 k_1 k} {\lambda r k_p +k_1
k_s r_1 } \Big] \\
A_2 = \frac {m} {r} \Big[ 1 - \frac {r_1 k_1 k} {\lambda r k_p +k_1 k_s
r_1 } \Big] \\
M_1 = \frac {m r_1} {\lambda r} \cdot \Big[\frac {r( \lambda k_a(l_2) +
2 k_1) + 2 \lambda k_d k_1} {r [\lambda k_p + k_1 k_s] + \lambda k_d k_1
k_s} \Big] \\
M_2 = \frac {m r_1} {\lambda r} \cdot \Big[\frac {r( \lambda k_a(l_1) +
2 k_1) + 2 \lambda k_d k_1} {r [\lambda k_p + k_1 k_s] + \lambda k_d k_1
k_s} \Big]
\end{eqnarray*}

It is clear from the above equations that $A_1$ and $A_2$ depend on the
ligand difference, so the system will respond to ligand gradients.  The
main text shows that the assumption $k_2 \gg 1$, on the other hand,
results in cells where $A_1 = A_2$, so the cell cannot maintain an
internal ligand gradient.

We also want to find the steady state of the system for $\lambda \gg 1$,
and show that the qualitative behavior remains the same as in the 
$\lambda
\approx O(1)$ case.  We rearrange $M_1$ and $A_1$ to show decreasing
powers of $\lambda$, and note that similar rearrangements can be made 
for
$M_2$ and $A_2$.
\begin{eqnarray*}
M_1 = \frac {m} {r} \frac {\lambda^2 a + \lambda b+ c}
{\lambda^2 d + \delta e } \\
a=k_a(l_2) k_d (r +2 k_2) + 2 k_1 k_d^2 \\
d=(r+2 k_2) k_a(l_1) k_a(l_2) + k_1 k_d(k_a(l_1) + k_a(l_2)) \\
A_1 = \frac {m} {r} + \frac {m} {r} \cdot \frac {\lambda \alpha} 
{\lambda
\beta + \gamma} \\
\alpha= k_1 k_d (k_a (l_2)- k_a(l_1))\\
\beta = k_d k_1 (k_a (l_2)+ k_a(l_1)) + (r + 2 k_2)k_a(l_1)k_a(l_2)
\end{eqnarray*}

We take the limit of these expressions as $\lambda \rightarrow \infty$:

\begin{eqnarray*}
M_1 = \frac {m} {r} \cdot \frac {a} {d} \\
A_1 = \frac {m} {r} + \frac {m} {r} \cdot \frac {\alpha} {\beta}
\end{eqnarray*}

We arrive at the new steady state:
\begin{eqnarray}
M_1 = \frac {m k_d} {r} \frac {2 k_1 k_d + k_a (l_2) r_2} {k_1
k_d k_s +  k_p r_2} \\
M_2 =  \frac {m k_d} {r} \frac {k_a(l_1) r_2 + 2 k_d k_1} {k_d 
k_1 k_s + k_p r_2} \\
A_1 = \frac {m} {r} \Big( 1 + \frac {k_1 k_d k}
{k_1 k_d k_s + k_p r_2} \Big) \\
A_2 = \frac {m} {r} \Big( 1 - \frac {k_1 k_d k}
{k_1 k_d k_s + k_p r_2} \Big)
\end{eqnarray}

It is simple to verify that similarly to the original system where
$\lambda \approx O(1)$, the qualitative behavior is the same with 
respect
to $k_1$ and $k_2$.  If $k_1 =0$, then the compartments reach the same
steady state as when they were not connected.
\begin{eqnarray*}
M_1 = \frac {m k_d} {r} \cdot \frac {k_a(l_2)r_2}
{k_p r_2} = \frac {mk_d} {r k_a(l_1)} \\
M_2 = \frac {mk_d} {r k_a(l_2)} \\
A_1 = A_2 = \frac {m} {r}
\end{eqnarray*}

The case $k_2 \gg 1$ and $\lambda \gg 1$ is discussed in the main text.
We conclude that regardless of the assumption we have make about
$\lambda$, we must have $k_1 \gg k_2$ in order to have the desired
dynamics.

Finally, we examine the absolute difference between $|A_1-A_2|$.  As
previously, we can assume without loss of generality that $k_2=0$.  Then
$|A_1-A_2|$ is bounded below by zero and above by $\frac {2m} {r}$.
\begin{eqnarray*}
A_1-A_2 = \frac {2m} {r} \frac {k_1 k_d k} {k_1 k_d k_s + r k_p}
\end{eqnarray*}

We consider $A_1-A_2$ as a function of $k$ and $k_s$, as before.  The 
same
analysis as in the $\lambda \approx O(1)$ case shows that for a fixed
concentration difference, $k_a(l_1)-k_a(l_2)$, the optimal concentration
range will be where $k_1 k_d (k_a(l_1)+k_a(l_2)) =r$.  We have shown 
that
for the particular cases we have considered, our equations have the same
qualitative behavior for $\lambda \approx O(1)$ and $\lambda \gg 1$.
This is sufficient for our purposes, but we note that in order to make 
this statement rigorous, we would have to compare the approximate 
solution
based on the separation of time scales (i.e. on the assumption that
$\lambda $ is large) with the exact analytical solution.

\section{Analytical solution and approximation}
\label{anal_sol}

We consider the system of equations:
\begin{eqnarray*}
\frac {dM_1} {dt} = m + \lambda (-k_a(l_1) M_1 + k_d A_1) + k_1(M_2-M_1)
\\
\frac {dA_1} {dt} = -r A_1 + \lambda ( k_a(l_1) M_1 - k_d A_1)  \\
\frac {dM_2} {dt} = m + \lambda (-k_a(l_2) M_2 + k_d A_2) - k_1(M_2-M_1)
\\ \frac {dA_2} {dt} = -r A_2 + \lambda ( k_a(l_2) M_2 - k_d A_2) \\
M_1(0) = M_2(0) = \frac {m} {r}  \frac {k_d} {k(l_0)} \ \
A_1(0) = A_2(0) = \frac {m} {r}
\end{eqnarray*}

Based on our previous analysis we assumed that the flux of $A$ was much   
smaller than the flux of $M$, $k_1$, so we set the rate of flux of $A$ 
to
be zero.  We rewrite the our equations in matrix form.  
\begin{eqnarray*}
\frac {d \vec{y}(t)} {dt} = D \vec{y} + h \\
\vec{y}(0) = \vec{y}_0
\end{eqnarray*}

with
\begin{eqnarray*}
\vec{y} = \left[ \begin{array}{c}
M_1(t)  \\ A_1(t) \\  M_2(t) \\ A_2(t)
\end{array} \right]
\end{eqnarray*}

\begin{eqnarray*}
D = \left[ \begin{array}{c c c c}
-(\lambda k_a(l_1)+k_1) & \lambda k_d & k_1 & 0 \\
\lambda k_a(l_1) & -(r+k_d) & 0 & 0\\
k_1 & 0 & -(\lambda k_a(l_2)+ k_1) & \lambda k_d \\     
0 & 0 & \lambda k_a(l_2) & -(r + k_d)
\end{array} \right]
\end{eqnarray*}

\begin{eqnarray*}
h = \left[ \begin{array}{c}
m  \\ 0 \\  m \\ 0
\end{array} \right]
\end{eqnarray*}  

and
\begin{eqnarray*}
\vec{y}_0 = \left[ \begin{array}{c}
\frac {m} {r} \frac {k_d} {k_a(l_0)} \\ \frac {m} {r}  \\ \frac {m} {r}
\frac {k_d} {k_a(l_0)}   \\ \frac {m} {r}
\end{array} \right]
\end{eqnarray*}

We want to find $X$, $\Lambda$ such that $D=X\Lambda X^{-1}.$  Now our
system becomes
\begin{eqnarray}
\frac {d \vec{y}(t)} {dt} = =X\Lambda X^{-1} \vec{y} + h \\
X^{-1} \frac {d \vec{y}(t)} {dt} = \Lambda X^{-1} y + X^{-1} h
\label{matrix}
\end{eqnarray}

Define $v=X^{-1} \vec{y}$ and $ \bar{h} = X^{-1} h$, so
\begin{eqnarray*}
\frac {dv} {dt} = \Lambda w + \bar{h} \\
v(0) = X^{-1} \vec{y}_0
\end{eqnarray*}

By making one more substitution, and letting $w = v + \Lambda^{-1}
\bar{h}$ we obtain
\begin{eqnarray*}
\frac {dw} {dt} = \Lambda w \\
w(0) = v(0) + \Lambda^{-1} \bar{h}
\end{eqnarray*} 

The solution to this equation is $w = w(0) e^{\Lambda t}$, and by making
the appropriate substitutions again, this gives the solution to the
equation \ref{matrix} to be
\begin{eqnarray*}
\vec{y} = (\vec{y}_0 + D^{-1} h) e^{\Lambda t} - D^{-1} h
\end{eqnarray*}  

In order to find the explicit formula for $\vec{y}$, we must find
$\Lambda$, the diagonal matrix of eigenvalues of $D$ and $D^{-1}$ in terms
of our parameters.  In spite the fact that this is a problem with some
symmetry, finding the eigenvalues and the inverse of the four-by-four
matrix, D is difficult even with Matlab's Symbolic Math Toolbox.  We leave
the exact solution in this form.

It is possible to approximate the exact solution to equations
\ref{m1}-\ref{ic} in case $\lambda \gg 1$.  The solution is similar to 
the solution of the equations for perfect adaptation in Appendix
\ref{approx_sol}.  

Depending on the size of the flux $k_1$ between the two compartments we
can consider two cases.  First, we assume that the flux between
the two compartments is very fast, and in fact, $k_1=O(\lambda)$.  Based
on our intuition developed by the steady state solution and the
approximate solutions, we expect in this case the greatest change to be
that the fast time scale, $r_f$ has to depend on both $k_a(l_1)$ and
$k_a(l_2)$.  If we write down the equations that apply on the fast time
scale, \begin{eqnarray*} \frac {dM_1} {d\tau} = -k_a(l_1)  M_1 + k_d A_1 +
k_1(M_2-M_1) \\ \frac {dA_1} {d\tau} = k_a(l_1) M_1 - k_d A_1 \\ \frac
{dM_2} {d\tau} = -k_a(l_2) M_2 + k_d A_2 - k_1(M_2-M_1)\\ \frac {dA_2}
{d\tau} = k_a(l_2) M_2 - k_d A_2 \end{eqnarray*} we see that now we must
sum all four components to get a constant, i.e., $M_1 + M_2 + A_1 + A_2
=C$, so the four equations are coupled.  In fact, solving this system of
equations is not simpler than providing the exact solution, therefore we
do not pursue this line of investigation.

Now we examine the case when $k_1 \approx O(1)$, so on the fast time scale
the it is still true that $A_1 + M_1 = C_1$ for a constant $C_1$, and
similarly, $A_2 + M_2 = C_2$.  The same calculations we used in Appendix
\ref{approx_sol} apply, and we can obtain the solution to on the fast time
scale: \begin{eqnarray} A_1(\tau) = (A_{s1} - A_{11}) e^{-r_{f1} \tau} +
A_{11} \label{a1f_sol} \\ M_1(\tau) = (M_{01} - M_{11})e^{-r_{f1} \tau} +
M_{11} \label{m1f_sol} \\ A_2(\tau) = (A_{s2} - A_{12}) e^{-r_{2f} \tau} +
A_{12} \label{a2f_sol} \\ M_2(\tau) = (M_{02} - M_{12})e^{-r_{2f} \tau} +
M_{12} \label{m2f_sol} \end{eqnarray}

where we have defined the constants
\begin{eqnarray*}
r_{f1} = k_a(l_1)+ k_d, \ \ r_{f2}=k_a(l_2)+ k_d \\
M_{01} = M_{02}= \frac {m k_d} {r k_a(l_0)} \\ 
M_{11}=\frac {m} {r} \frac {k_d} {k_a(l_1)} \frac
{1+(k_d/k_a(l_0))}{1+(k_d/k_a(l_1))} , \ \ M_{12}= \frac {m} {r} \frac {k_d}
{k_a(l_1)} \frac {1+(k_d/k_a(l_0))}{1+(k_d/k_a(l_2))}\\
A_{s1}=A_{s2}= \frac {m} {r} \\ 
A_{11}=\frac {m} {r} \frac {1+(k_d/k_a(l_0))} {1+(k_d/k_a(l_1))}, \ \ A_{12}
= \frac {m} {r} \frac {1+(k_d/k_a(l_0))} {1+(k_d/k_a(l_2))}
\end{eqnarray*}

On the slow time scale it remains true that $A_1 = \frac {k_a(l_1)} {k_d}
M_1$, and $A_2 = \frac {k_a(l_2)} {k_d} M_2$.

As before, substituting these expressions into
equations \ref{m1}-\ref{ic} we arrive at a new system of equations:
\begin{eqnarray} 
\frac {dM_1} {dt} = m + k_1 (M_2 - M_1) \label{m1a} \\
\frac {k_a(l_1)} {k_d} \frac {dM_1} {dt} = -r \frac {k_a(l_1)} {k_d} M_1
\label{f1a} \\
\frac {dM_2} {dt} = m - k_1 (M_2 - M_1) \label{m2a} \\
\frac {k_a(l_1)} {k_d} \frac {dM_1} {dt}= -r \frac {k_a(l_2)} {k_d} M_2
\label{f2a}
\end{eqnarray}

As before, we add equations \ref{m1a} and \ref{f1a}, and equations \ref{m2a} 
and \ref{f2a}.  We arrive at 
\begin{eqnarray}   
\frac {dM_1} {dt} = \frac {m k_d} {k_d + k_a(l_1)} + \frac {k_1 k_d}
{k_d + k_a(l_1)}(M_2 - M_1) - r \frac {k_a(l_1)} {k_d + k_a(l_1)} M_1
 \label{m1as} \\
\frac {dM_2} {dt} = \frac {m k_d} {k_d + k_a(l_2)} + \frac {k_1 k_d}
{k_d + k_a(l_2)}(M_2 - M_1) - r \frac {k_a(l_2)} {k_d + k_a(l_2)} 
M_2 \label{m2as}
\end{eqnarray}

The solution to eqns. \ref{m1as}, \ref{m2as} is given by:
\begin{eqnarray*}
\bar{M}(t) = (\bar{M}(0) + D^{-1} h) e^{\Lambda t} - D^{-1} h
\end{eqnarray*}

where
\begin{eqnarray*}  
D = \left[ \begin{array}{c c}
\frac {-r k_a(l_1) - k_1 k_d} {k_d + k_a (l_1)}  & \frac{k_1 k_d} {k_d +
k_a (l_1)} \\
\frac {k_1 k_d} {k_d + k_a(l_1)} &  \frac {-r k_a (l_2)- k_1 k_d} {k_d +
k_a (l_2)}
\end{array} \right]
\end{eqnarray*}

\begin{eqnarray*}
h = \left[ \begin{array}{c}
\frac {m k_d} {k_d + k_a(l_1)}  \\  \frac {m k_d} {k_d + k_a(l_2)}
\end{array} \right]
\end{eqnarray*}  

\begin{eqnarray*}
\bar{M}(0) = \left[ \begin{array}{c}
c_1  \\ c_2
\end{array} \right]
\end{eqnarray*}    

We define $\beta = Tr (D)$ and $\gamma = \det(D)$.  Then the
eigenvalues $\lambda_1$, $\lambda_2$ of $D$ can be found as follows:
\begin{eqnarray*}  
\lambda_{1,2} = \frac {-\beta \pm \sqrt{\beta^2 - 4 \gamma}} {2} 
\end{eqnarray*}
By making the appropriate substitutions and carrying out the calculations, 
the discriminant $\sqrt{\beta^2-4\gamma}$ can be reduced to:

$$\Big(\frac {(k_1 k_d +r k_a(l_1))} {(k_d + k_a(l_1))} - \frac {(k_1k_d + 
rk_a(l_2))} {(k_d + k_a(l_2))} \Big)^2 - 2 
\frac {(k_1 k_d)^2} {(k_d + k_a(l_1))(k_d + k_a(l_2))}  $$

If we assume that the second term is much smaller than the first one, then
the above expression greatly simplifies.  This is true if 
\begin{eqnarray} (r+k_1)k_d (k_a(l_1)-k_a(l_2)) \gg k_1 k_d \label{condition} 
\end{eqnarray}
Equation \ref{condition} implies that our approximation is appropriate when 
the ligand concentrations in the two compartment are very different, i.e. 
when $k_a(l_1)-k_a(l_2) \gg 1$.  Now we are able to find the 
two eigenvalues: 
\begin{eqnarray*}
\lambda_{1,2} \simeq - \frac {r k_a(l_1) + k_1 k_d} {k_d + k_a (l_1)} , - 
\frac {r k_a(l_2) + k_1 k_d} {k_d + k_a(l_2)}
\end{eqnarray*}

The two eigenvalues define the two slow time scales, 
$$ r_{s1} = - \frac {r k_a(l_1) + k_1 k_d} {k_d + k_a (l_1)}$$
$$ r_{s2} = - \frac {r k_a(l_2) + k_1 k_d} {k_d + k_a (l_2)}$$

We notice that if the flux $k_1$ is small, then we have recovered the slow 
time scale for each compartment independently of each other.  This result 
confirms conclusions we have drawn from our steady state analysis.  

Next, we find $D^{-1}h$:
\begin{eqnarray*}
d_1 = D^{-1}h_1 = -\frac{m} {r} \frac {k_d (r k_a(l_2)+k_1 k_d)} {r
k_a(l_1)k_a(l_2) + k_1 k_d (k_a(l_1)+ k_a(l_2))} \\
d_2 = D^{-1}h_2 = - \frac {m} {r} \frac {k_d k_1 (k_d + k_a(l_1))} { r
k_a(l_1) k_a(l_2) + k_1 k_d(k_a(l_1)+ k_a(l_2))}
\end{eqnarray*}    

The solution to eqns. \ref{m1as} - \ref{m2as} is:
\begin{eqnarray}
M_{1}= (c_1 + d_1) e^{(-r_{s1} t)} - d_1 \label{m1s_sol} \\
M_{2}= (c_2 + d_2) e^{(-r_{s2} t)} - d_2 \label{m2s_sol}
\end{eqnarray}

The equations for $M_1$ and $M_2$ also determine the expressions for $A_1$ 
and $A_2$:
\begin{eqnarray}
A_{1}= \frac {k_a(l_1)} {k_d}(c_1 + d_1) e^{(-r_{s1} t)} - \frac 
{k_a(l_1) d_1} {k_d}  \label{a1s_sol} \\
A_{2}= \frac {k_a(l_2)} {k_d}(c_2 + d_2) e^{(-r_{s2} t)} - \frac 
{k_a(l_2) d_2} {k_d}  \label{a2s_sol}
\end{eqnarray}

Now we can match the fast and slow solutions as before to determine the 
constants $c_1$ and $c_2$, and to give arrive at the full approximation. 
Taking the limit as $\tau \rightarrow \infty$ of the equations on the fast 
time scale, and setting this equal to the initial condition of the equations 
of the slow time scale gives us: 
\begin{eqnarray*}
c_1 + d_1 - d_1 = M_{11} \\
c_2 + d_2 - d_2 = M_{12} \\
\frac {k_a(l_1)} {k_d}(c_1 + d_1) - \frac {k_a(l_1) d_1} {k_d} = A_{11} \\
\frac {k_a(l_2)} {k_d}(c_2 + d_2) - \frac {k_a(l_2) d_2} {k_d} = A_{12}
\end{eqnarray*}

Thus the approximate solution to our system of equations 
\begin{eqnarray}
A_1 = (A_{s1} - A_{11}) e^{-r_{f1} \lambda t} + \frac {k_a(l_1)} {k_d}(A_{11} 
+ d_1) e^{(-r_{s1} t)} - \frac {k_a(l_1) d_1} {k_d}  \\
M_1 = (M_{01} - M_{11})e^{-r_{f1} \lambda t} + (M_{11} + d_1) e^{(-r_{s1} t)} 
- d_1 \\
A_2 = (A_{s2} - A_{12}) e^{-r_{f2} \lambda t} + \frac {k_a(l_2)} {k_d}(A_{12}
+ d_2) e^{(-r_{s2} t)} - \frac {k_a(l_2) d_2} {k_d}  \\
M_2 = (M_{02} - M_{12})e^{-r_{f2} \lambda t} + (M_{12} + d_2) e^{(-r_{s2} t)}  
- d_2 
\end{eqnarray}

\section{Calcium switch}
\label{switch}

We return to the experimental observation that the cytosolic calcium
concentration can change the turning behavior of a growth cone.  Let us  
assume that a netrin-1 gradient is set up outside the cell.  Recall that
in a cell with normal cytosolic calcium concentration a gradient 
develops
with the high calcium concentrations facing the source of netrin-1, and
the growth cone responds with attractive turning.  However, in the same
netrin-1 gradient a cell whose cytosolic calcium has been depleted 
before
the trial responds with repulsive turning (eventhough the high calcium
concentrations still face the source of netrin-1).  Such a behavior is 
possible by making some assumptions regarding $k_a$, the rate at which 
$A$
is produced in our model.

\begin{figure}[h!] \centering
\centerline{\includegraphics[width=0.8\textwidth]{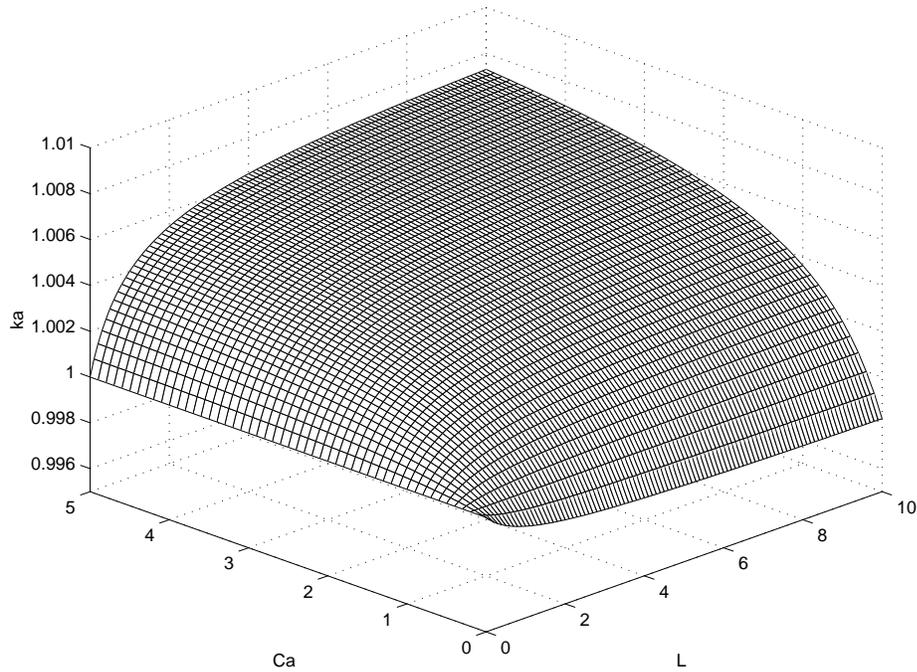}}
\caption[Calcium switch]{$k_a(l,Ca)$, the production rate of the 
modified
substance, $A$.
l: ligand, netrin-1, Ca: cytosolic calcium concentration.  $a=0.01$,   
$b=1$, $c=1$, $Ca_b=0.2$} \label{ka_1.ps}
\end{figure}

We assume that the production rate of $A$ depends both on the ligand      
concentration and the cytosolic calcium concentration, $k_a (l, Ca)$.  
Further, we want $k_a$ to be such that for large values of calcium 
$\frac
{\partial k_a} {\partial l} > 0 $, so $k_a$ is an increasing function of
the ligand, and for small values of calcium $\frac {\partial k_a}
{\partial l} < 0 .$ Without the constraints of experimental data, many
such functions can be found.  We chose $k_a = \exp(\frac {a l (Ca - 
Ca_b)}
{(l+ b)(Ca + c)} ) $ where $a$, $b$ and $c$ are new constants, and  
$Ca_b$ is the normal cytosolic calcium concentration.  $k_a$ is
illustrated in Figure \ref{ka_1.ps}.

\begin{figure}[h!]
\includegraphics[width=0.45\textwidth]{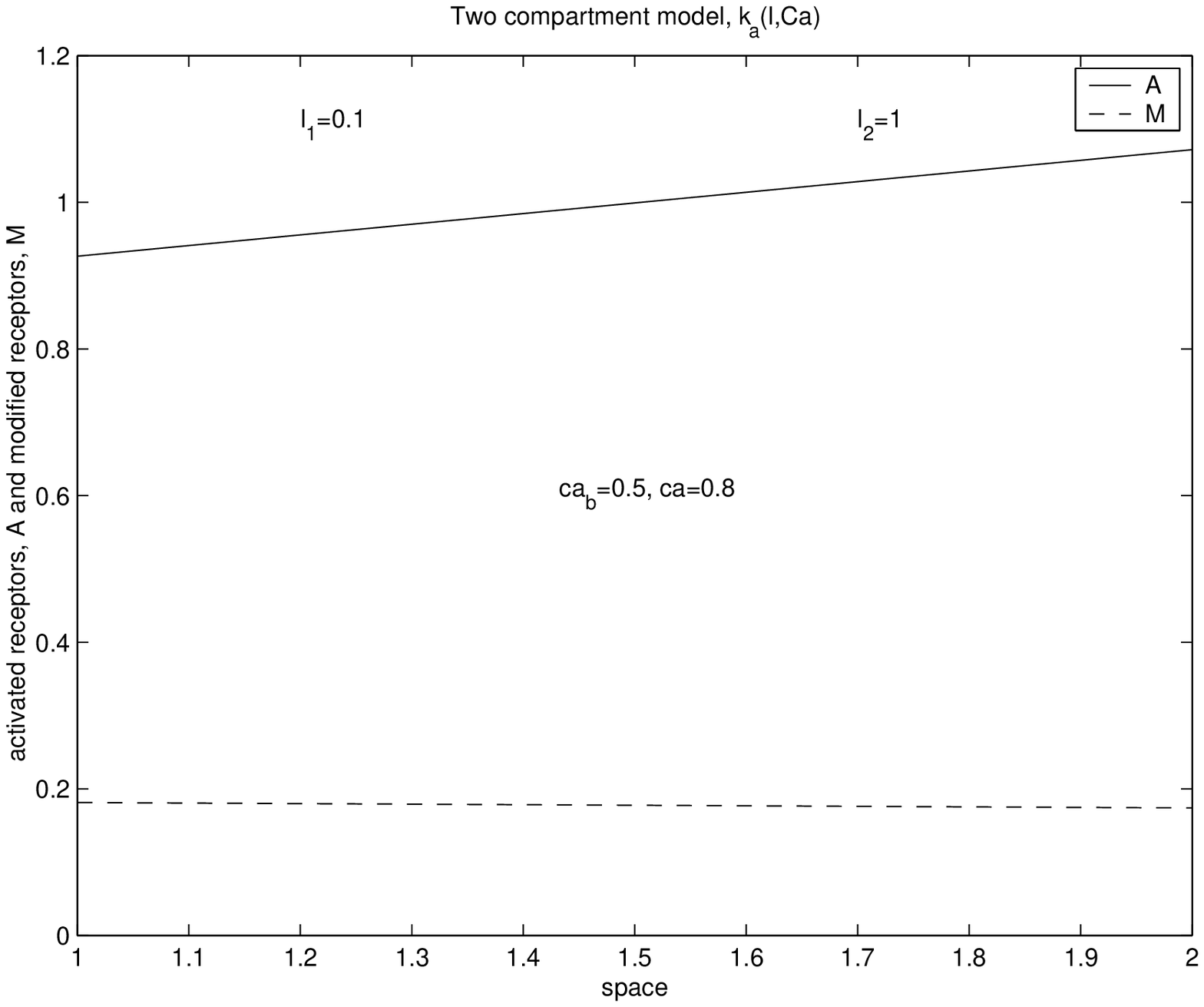}\qquad
\includegraphics[width=0.45\textwidth]{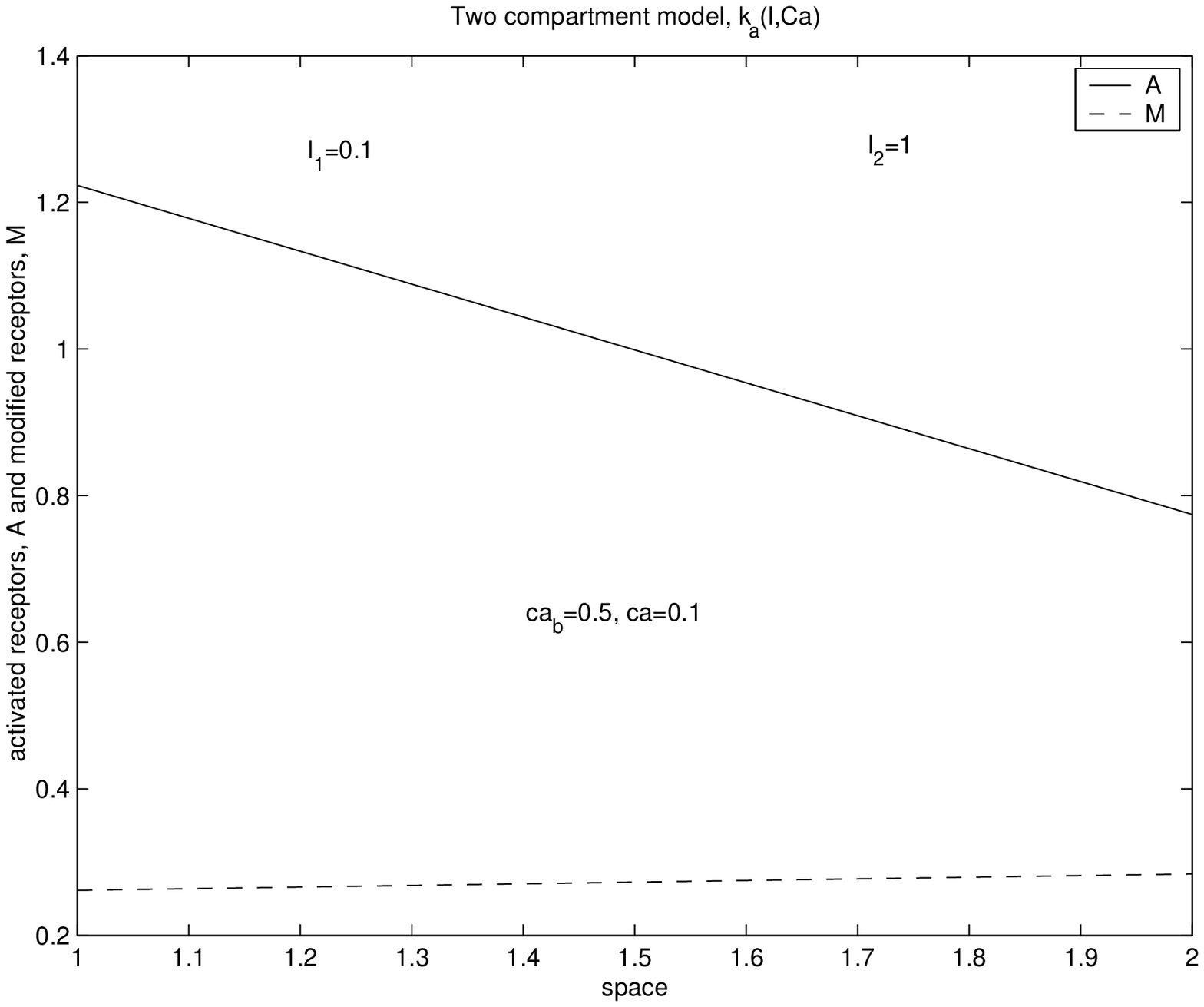}
\caption[Numerical results with the calcium switch.]{Numerical 
simulations
of the two compartment model using $k_a =
\exp(\frac {a l (Ca - Ca_b)}{(l+ b)(Ca + c)} )$. The first figure shows 
a
cell where the cytosolic calcium level is above the baseline.  The 
second 
figure shows a cell in the same ligand gradient where the calcium level 
is
below the baseline. }
\label{adapt_switch} \end{figure}

Using this rate we can simulate what happens with our system of equations
in the same ligand gradient when the calcium level is above and when the
calcium level is below the baseline, $Ca_b$.  Changing the calcium level
in the simulations corresponds to changing the value of the parameter $Ca$
in the expression for $k_a(l)$.  This implicitly implies that we take the
calcium level to be spatially uniform inside the cell.  Figure
\ref{adapt_switch} shows that in the normal cytosolic concentration we get
a gradient of the adapted substance, $A$ when a netrin-1 gradient is
presented, with the higher concentration of $A$ corresponding to the
higher concentration of $l$.  (As before, the "gradient" of $A$ is based
on two values only, $A_1=A(1)$ in the left hand compartment and $A_2=A(2)$
in the right hand compartment.) Lowering the cytosolic concentration
level, and presenting the cell the same ligand gradient results in a
gradient of $A$ where now the level of $A$ is lower in the compartment
corresponding to the higher ligand concentration.

Both panels of the figure are run to the same time, t=50 seconds.  The
simulations are based on the same Matlab code as the previous two
compartment models, except that now $k_a (l)$ is changed.

This model presents a very simple explanation of how calcium levels 
influence turning behavior.  Only the overall calcium concentration is   
considered by the simulations, the spatial and temporal gradients of
calcium are not.  Further investigations are necessary to create a more
realistic description of the behavior.


\newpage
\thispagestyle{myheadings}
\markright{  \rm \normalsize APPENDIX. \hspace{0.5cm}
  MATHEMATICAL MODELS IN BIOLOGY}
\chapter{Sample Matlab code} \label{num_sol}
\thispagestyle{myheadings}

We are interested in investigating the deformations of two coupled Kelvin
bodies in (i) steady flow, $F = F_0$ and (ii) oscillatory flow, $F = F_0
\cos (\omega t)$.  Equations \ref{maineqs1} and \ref{maineqs2} were solved
with a Matlab code given below.  The solution presented here is for
oscillatory flow.

First we can see the program \verb+twobodies_aF.m+ that asks for the input
from the user and displays the results of the simulations.  The user
defines the coefficients for the two Kelvin bodies, and the
program stores this information in the matrix \verb+B+.  In this code a
function, \verb+parallel2_aF.m+ is called which actually computes the
solution with a four-stage Runge-Kutta method.  Then, as the output of
\verb+parallel2_aF.m+, the solution of the differential equation is
returned to the matrix \verb+u_1+ in \verb+twobodies_aF.m+, the solution
(in this case $u(t)$) is plotted.

\begin{verbatim} 
r=input('Two-body system with the bodies connected in parallel');
r=input('Coefficients for ith body are:
[\mu_{0i} \mu_{1i} \eta_{1i}] ');

B=zeros(2,3);
B(1,:)=input('Coefficients for Body 1:- ');
B(2,:)=input('Coefficients for Body 2:- ');

h = 0.1;                        % size of the time step
N0=5;
N = N0/h;

x=[1:1:N];                      % length of a, u, F

u_1=parallel2_aF(B,h,N);        % function call

u1=u_1(1,1:1:N);
F=u_1(3,1:1:N);
a=u_1(2,1:1:N);
a=a./F;

plot(x,u,'r-.')
ylabel('u(t), deformation')
xlabel('time, t')
title('u vs t for different values of \eta_{12},
F=F_0 cos(\omega t)')

\end{verbatim}

Now we can look at the code for the actual ODE solver,
\verb+parallel2_aF.m+.  Here only the main loop of the program is
included which contains the fourth-order four-stage Runge-Kutta method.
In the actual program the matrices \verb+M1+, \verb+M2+, \verb+M3+ and
\verb+M+ are defined to be $A$, $D$, $\vec{c}$ and $A^{-1}D$,
respectively.  The full code contains the initialization of all the
appropriate variables.  This program also calls a function,
\verb+ve2_aF.m+, given below.  The input of \verb+parallel2_aF.m+
function is the matrix of coefficients of the Kelvin bodies, denoted by
\verb+B+ , the step size, \verb+h+, and the length of the solution
vector, \verb+N+.  The output of this function is also a matrix, called
\verb+u1_plot+ whose first row is $u(t)$, second row is $a(t)F(t)$, and
third row is $F(t)$.  $u(t)$ and $a(t)F(t)$ are obtained as the solutions
to the differential equations.  $F(t)$, the third row of the matrix, is
simply $F_0 \cos (\omega t)$ evaluated for each time step.

\begin{verbatim}
function f=parallel2_aF(B,h,N);

for k=1:N

        F(k+1)=F0*cos(w*(t+h));         % oscillatory

        k1=u;
        k2=u+(1/2)*h*ve2_aF(M,inM1,M3,k1,t+(h/2));
        k3=u+(1/2)*h*ve2_aF(M,inM1,M3,k2,t+(h/2));
        k4=u+h*ve2_aF(M,inM1,M3,k3,t+h);

        u_new=u+h*((1/6)*ve2_aF(M,inM1,M3,k1,t)...
        +(1/3)*ve2_aF(M,inM1,M3,k2,t+(h/2))...
        +(1/3)*ve2_aF(M,inM1,M3,k3,t+(h/2))...
        +(1/6)*ve2_aF(M,inM1,M3,k4,t+h));

        u = u_new;
        t = t+h;
        u1_plot(1:2,k+1)=u;
        u1_plot(3,k+1)=F(k+1);

end

f=u1_plot;

\end{verbatim}

Finally, we can look at \verb+ve2_aF.m+.  The input of this function is
$A^{-1} D$, denoted here by \verb+M+; $A^{-1}$, denoted by \verb+inM1+;
$\vec{c}$, denoted by \verb+M3+; a vector which consists of $u(t)$ and
$a(t)F(t)$ at the previous time step, denoted by \verb+u+, and finally,
$t$, the current time step.  The output of this function is the right hand
side of equation \ref{maineqs1} for the appropriate time step with the
given type of flow.  Here oscillatory flow is shown.
        
\begin{verbatim}
        
function f = ve2_aF(M,inM1,M3,u,t)

F0=1;
w=2*pi;
F=F0*cos(w*t);
dF=-F0*w*sin(w*t);

C=zeros(2,1);
C = M3*dF;
C(2,1) = F+C(2,1);

f = M*u+inM1*C;

\end{verbatim}


\end{document}